\documentclass[11pt,twoside]{article}

\usepackage[left=2cm, bottom=2cm, right=2cm, top=2cm]{geometry}
\usepackage[utf8]{inputenc}
\usepackage[T1]{fontenc}
\usepackage[french]{babel}
\usepackage[dvipsnames]{xcolor}
\usepackage[hidelinks]{hyperref}
\usepackage{setspace}
\usepackage{fancyhdr}
\usepackage[strict]{changepage}
\usepackage{mdframed}
\usepackage{ifthen}
\usepackage{amsmath}
\usepackage{amsthm}
\usepackage{amssymb}
\usepackage{tikz-cd}

\newlength\aliaslength
\newlength\defparindent
\newlength\defparskip
\theoremstyle{definition}

\newtheorem{mydefinition}{Definition}
\newenvironment{definition}[1][noref]
{\begin{mydefinition}\label{#1}}
    {\end{mydefinition}}

\newtheorem{myspedefinition}{Definition}

\newtheorem{mysituation}[mydefinition]{Situation}
\newenvironment{situation}[1][noref]
{\begin{mysituation}\label{#1}}
    {\end{mysituation}}

\theoremstyle{plain}

\newtheorem{mytheorem}[mydefinition]{Theorem}
\newenvironment{theorem}[1][noref]
{\begin{mytheorem}\label{#1}}
    {\end{mytheorem}}

\newtheorem{mylemma}[mydefinition]{Lemma}
\newenvironment{lemma}[1][noref]
{\begin{mylemma}\label{#1}}
    {\end{mylemma}}

\newtheorem{myproposition}[mydefinition]{Proposition}
\newenvironment{proposition}[1][noref]
{\begin{myproposition}\label{#1}}
    {\end{myproposition}}

\newtheorem{mycorollary}[mydefinition]{Corollary}
\newenvironment{corollary}[1][noref]
{\begin{mycorollary}\label{#1}}
    {\end{mycorollary}}

\newtheorem{myspetheorem}[myspedefinition]{Theorem}
\newenvironment{spetheorem}[1][noref]
{\begin{myspetheorem}\label{#1}}
    {\end{myspetheorem}}

\newtheorem{myspecorollary}[myspedefinition]{Corollary}
\newenvironment{specorollary}[1][noref]
{\begin{myspecorollary}\label{#1}}
    {\end{myspecorollary}}

\theoremstyle{remark}

\newtheorem{myremark}[mydefinition]{Remark}
\newenvironment{remark}[1][noref]
{\begin{myremark}\label{#1}}
    {\end{myremark}}

\newtheorem{myexample}[mydefinition]{Example}
\newenvironment{example}[1][noref]
{\begin{myexample}\label{#1}}
    {\end{myexample}}

\newenvironment{myproof}[1][\hspace{-\aliaslength}]
{\setlength\parskip{0.5\defparskip}%
    \vspace\parskip%
    \noindent\underline{\textit{%
            Proof #1.}}}
{\hfill\qedsymbol%
    \setlength\parskip\defparskip}

\newcommand\Definition[1]{%
    Definition~\ref{#1}}
\newcommand\Lemma[1]{%
    Lemma~\ref{#1}}
\newcommand\Proposition[1]{%
    Proposition~\ref{#1}}
\newcommand\Corollary[1]{%
    Corollary~\ref{#1}}
\newcommand\Theorem[1]{%
    Theorem~\ref{#1}}
\newcommand\Remark[1]{%
    Remark~\ref{#1}}
\newcommand\Example[1]{%
    Example~\ref{#1}}
\newcommand\Situation[1]{%
    Situation~\ref{#1}}

\newcommand{\CC}{\mathbb C}
\newcommand{\FF}{\mathbb F}
\newcommand{\GG}{\mathbb G}
\newcommand{\NN}{\mathbb N}
\newcommand{\QQ}{\mathbb Q}
\newcommand{\ZZ}{\mathbb Z}

\newcommand{\A}{\mathcal A}
\newcommand{\B}{\mathcal B}
\newcommand{\C}{\mathcal C}
\newcommand{\J}{\mathcal J}

\newcommand{\F}{\mathcal F}

\renewcommand{\P}{\mathcal P}
\renewcommand{\O}{\mathcal O}
\newcommand{\T}{\mathcal T}
\newcommand{\U}{\mathcal U}
\newcommand{\W}{\mathcal W}

\newcommand{\eps}{\varepsilon}
\renewcommand{\phi}{\varphi}
\newcommand{\Ohm}{\Omega}
\newcommand{\ohm}{\omega}

\newcommand{\Adh}[1]{\overline{#1}}
\newcommand{\aff}{\mathrm{aff}}
\DeclareMathOperator{\Alg}{\mathrm{Alg}}
\newcommand{\Ab}{\mathrm{Ab}}
\DeclareMathOperator\Br{\mathrm{Br}}
\DeclareMathOperator\cd{\mathrm{cd}}
\DeclareMathOperator{\characteristic}{\mathrm{char}}
\DeclareMathOperator{\coker}{\mathrm{coker}}
\newcommand\cond{\mathrm{cond}}
\newcommand\cores{\mathrm{cores}}
\newcommand\et{\mathrm{et}}
\newcommand\Et{\mathrm{Et}}
\DeclareMathOperator\Ext{\mathrm{Ext}}
\DeclareMathOperator\EXT{\underline{\mathrm{Ext}}}
\newcommand\fin{\mathrm{fin}}

\newcommand\fppf{\mathrm{fppf}}
\DeclareMathOperator\Gal{\mathrm{Gal}}
\DeclareMathOperator\Hom{\mathrm{Hom}}
\DeclareMathOperator\HOM{\underline{{\mathrm{Hom}}}}
\newcommand{\id}{\mathrm{id}}
\newcommand{\Ind}{\mathrm{I}}
\newcommand{\IndAb}{\mathrm{IAb}}
\newcommand{\IndAlg}{\mathrm{IAlg}}
\newcommand{\IndPro}{\mathrm{IP}}
\newcommand{\IndProAb}{\mathrm{IPAb}}
\newcommand{\IndProAlg}{\mathrm{IPAlg}}
\newcommand\indrat{\mathrm{indrat}}
\DeclareMathOperator{\im}{\mathrm{im}}

\newcommand{\LCAb}{\mathrm{LCAb}}

\newcommand\opp{\mathrm{op}}
\newcommand\mm{\mathfrak{m}}
\newcommand{\myhat}[1]{\widehat #1}
\newcommand\perar{\mathrm{perar}}
\newcommand\perf{\mathrm{perf}}
\newcommand\Pro{\mathrm{P}}
\newcommand\ProAb{\mathrm{PAb}}
\newcommand\ProAlg{\mathrm{PAlg}}
\newcommand\proet{\mathrm{proet}}
\newcommand\profppf{\mathrm{profppf}}
\newcommand\PRPS{\mathrm{PRPS}}
\newcommand\PSh{\mathrm{PSh}}
\newcommand\qcqs{\mathrm{qcqs}}
\newcommand\res{\mathrm{res}}
\newcommand\Res{\mathrm{Res}}
\newcommand\RPS{\mathrm{RPS}}
\newcommand\RPSch{\mathrm{RPSch}}
\newcommand\RPSSch{\mathrm{RPSSch}}
\newcommand\RP{\mathrm{RP}}
\newcommand\Sch{\mathrm{Sch}}
\newcommand\set{\mathrm{set}}

\DeclareMathOperator\Sh{\mathrm{Sh}}
\DeclareMathOperator\Spec{\mathrm{Spec}}

\newcommand{\Tr}{\mathrm{Tr}}
\newcommand{\tr}{\mathrm{tr}}
\renewcommand{\top}{\mathrm{top}}
\newcommand{\TopAb}{\mathrm{TopAb}}
\newcommand{\tur}{\mathrm{tur}}
\newcommand{\ur}{\mathrm{ur}}

\title{Arithmetic duality for finite Galois modules over two-dimensional local fields of mixed characteristic}
\author{Antoine GALET\footnote{Sorbonne Université, Université Paris Cité, CNRS, IMJ-PRG, F-75005 Paris, France}, 2026}
\date{}

\begin{document}

\maketitle

\renewcommand\abstractname{Abstract}
\begin{abstract}
    \noindent A field $K$ is $d$-local if there exist fields $K=k_d,\dots,k_0$ where each $k_{i+1}$ is a complete discrete valuation field with residue field $k_i$, and $k_0$ is a finite field of characteristic $p$. By work of Deninger and Wingberg, the Galois cohomology of such fields with coefficients in finite Galois modules satisfies a duality generalizing Tate duality when either $d=0$, $\characteristic k_1=0$ or the coefficients have no $p$-torsion. Based on recent progress by Kato and Suzuki, we obtain duality statements for arbitrary finite Galois modules, under the weaker assumption that either $d\leq 1$ or $\characteristic k_2=0$. We also get duality for the étale cohomology of $K$-varieties with coefficients in finite étale groups, in the style of Artin-Verdier duality. These dualities are stated in terms of condensed structures (in fact, locally compact Hausdorff topologies) on the cohomology groups. More generally we obtain results for any perfect $k_0$, endowing the totally unramified cohomology groups of $K$ with the structure of ind-pro-quasi-algebraic $k_0$-groups.
\end{abstract}

\scriptsize
\renewcommand\contentsname{Contents}
\setlength\parskip{0pt}
\tableofcontents
\setlength\parskip\defparskip
\normalsize


\section*{Introduction}
\label{Introduction}
\addcontentsline{toc}{section}{Introduction}

\emph{Higher local fields} are a higher-dimensional generalization of nonarchimedian local fields, classically defined as follows : a field $K$ is \emph{$d$-local} if there exists a sequence of fields $K=k_d,\dots,k_0$ such that $k_0$ is finite and $k_i$ is a complete discrete valuation field with residue field $k_{i-1}$, for $1\leq i\leq d$. It is often not a great sacrifice to allow $k_0$ to only be \emph{quasi-finite}, \emph{i.e.} perfect with absolute Galois group $\myhat\ZZ$, or even just perfect, and $k_{i+1}$ to only be Henselian when $\characteristic k_{i+1}=0$.

The cohomology of a $p$-adic field $K$ with coefficients in finite Galois modules satisfies a duality, by Tate's local duality theorem : namely the cup-product $H^q(K,M)\otimes H^{2-q}(K,M^\vee)\to \Br(K)\cong\QQ/\ZZ$ is a perfect pairing of finite groups for all $q\in\ZZ$ and finite Galois module $M$, where $M^\vee=\Hom(M,\GG_m)$ is the Cartier dual. If instead $K$ is local of characteristic $p>0$, the same statement holds if $M$ has no $p$-torsion. This has a higher local analogue by \cite[Prop. 1.2]{DenWin} : if $K=k_d,\dots,k_1$ is $d$-local with $k_0$ quasi-finite of exponent characteristic $p$, and $M$ is a finite Galois module of order $m$ prime to $p$, then we have perfect pairings of finite groups :
$$H^q(K,M)\otimes H^{d+1-q}(K,M^\vee)\to H^{d+1}(K,\mu_m^{\otimes d})\cong\ZZ/m$$
where $M^\vee=\Hom(M,\mu_m^{\otimes d})$ is the Tate-twisted Cartier dual of $M$. If $\characteristic k_1=0$ and $k_0$ is finite then this holds even if $M$ has $p$-torsion, but in general the cohomology of $K$ with $p$-torsion coefficients is not finite. In this paper, we consider such pairings with finite coefficients of arbitrary order.

One approach is by induction on $d$ : considering a complete discrete valuation field $L$ with residue field $l$, we have Hochschild-Serre spectral sequences of Galois cohomology :
$$E^{ij}_2=H^i(l,H^j(L^\ur,-))\Rightarrow H^{i+j}(L,-)$$
where $L^\ur$ is the maximal unramified extension of $L$ with respect to its discrete valuation, so duality for a $d$-local field $K=k_d,\dots,k_0$ reduces to duality for $k_0$ and for the unramified cohomology functors $R\Gamma(k_i^\ur,-)$ for $i\geq 1$. Three cases are known :
\begin{enumerate}
    \item[a.] When $m\geq 2$ is invertible in $l$, $H^q(L^\ur,-)$ sends perfect pairings of finite $m$-torsion Galois modules $M\otimes N\to\mu_m^{\otimes i}$ over $L$ to perfect pairings $H^r(L^\ur,M)\otimes H^{1-r}(L^\ur,N)\to H^1(L^\ur,\mu_m^{\otimes i})=\mu_m^{\otimes(i-1)}$ of finite Galois modules over $l$, for any $i\in\ZZ$. This recovers the statement above for higher local fields.
    \item[b.] When $L$ has mixed characteristic $(0,p)$ and $[l:l^p]<+\infty$, Bertapelle and Suzuki \cite{BertapelleSuzuki} showed that $H^q(L^\ur,\mu_{p^n}^{\otimes r})$ has the structure of the \emph{relative perfection of an algebraic unipotent (RPAU)} $l$-group $R^q\Psi(\mu_{p^n}^{\otimes r})$. Such groups, when viewed as sheaves on the \emph{relatively perfect} or \emph{relatively perfectly smooth} site of $l$, admit a notion of duality first developed by Kato \cite{Kato86, Kato87} and expanded upon in \cite{KatoSuzuki2019, BertapelleSuzuki}. Kato and Suzuki \cite{KatoSuzuki2019} obtained perfect pairings $R\Psi(\mu_{p^n}^{\otimes r})\otimes^LR\Psi(\mu_{p^n}^{\otimes (i-r)})\to\nu_n(i-1)[-i]$ where $p^{i-1}=[l:l^p]$ and $\nu_n(r)=W_n\Ohm^r_{\log}$ is the sheaf of logarithmic de Rham-Witt differentials \cite{Kato86}.
    \item[c.] When $l$ is perfect, $L$ has equal characteristic $p>0$, and $G$ is an RPAU $L$-group, Suzuki \cite[Sec. 4]{SuzukiCFT} showed that $H^q(L^\ur,G)$ has the structure of an ind-pro-RPAU $l$-group $R^q\pi_{L,\RP,*}(G)$. Furthermore it belongs to a class of good ind-pro-RPAU groups $\W_l$ to which duality for RPAU $l$-groups extends, and we have a perfect pairing $R\pi_{L,\RP,*}(G)\otimes^LR\pi_{L,\RP,*}(G^\vee)\to\QQ_p/\ZZ_p[-1]$. This pairing is defined by viewing objects of $\W_l$ over the \emph{perfect Artinian étale} or \emph{indrational proétale} sites of $l$ defined in \cite{SuzukiGN, SuzukImp}. If $l$ is finite, composing with $R\Gamma(l,-)$ yields a perfect pairing of locally compact Hausdorff groups (viewed as condensed groups), recovering approaches of \cite[Sec. III.7 and III.11]{Milne2006} and \cite{Shatz64}.
\end{enumerate}

In \textbf{(b.)}, if $L$ has a $p$-th root of unity then the RPAU structure on $R^q\Psi(\mu_p^{\otimes r})$ is given by Kato's unit filtration on Milnor $K$-theory $K_q^M(L)/p\cong H^q(L,\mu_p^{\otimes q})\cong H^q(L,\mu_p^{\otimes r})$, which is bounded with graded objects isomorphic to groups of differentials over $l$. In \textbf{(c.)}, the ind-pro-RPAU structure of $R\pi_{L,\RP,*}\GG_{a,L}^\RP$ is essentially given by the Tate vector space isomorphism $L\cong l(\!(\pi)\!)\cong\bigoplus_{n<0}lt^n\oplus\prod_{n\geq 0}lt^n$, where $\GG_a^\RP$ is the relative perfection of the additive group.

In this paper we generalize, detail and assemble such results into a statement on a higher local field $K$ and varieties over it, and interpret them concretely in terms of étale cohomology. This improves on \cite{DenWin} by allowing $K$ to have $\characteristic(k_2)=0$ instead of $\characteristic(k_1)=0$. This improves on \cite{SuzukiCFT} by allowing arbitrary coefficients and slightly more general $K$.


\subsection*{Main results}
\label{Main results}
\addcontentsline{toc}{subsection}{Main results}

We need some notations and definitions to state the results.

For $F$ a perfect field, write $\Spec F_\proet^\indrat$ the indrational proétale site of $F$ (\Definition{sites of perfect field}). Write $D^b_{\W_F}(F_\proet^\indrat)$ the full subcategory of bounded objects of $D(F_\proet^\indrat)$ with cohomology representable by "good" ind-pro-RPAU $F$-groups (\Definition{good indpro} and \Definition{Wk as sheaves}).

For $X$ a scheme and $S_X$ one of the sites $X_\et$ (in general) or $X_\proet^\indrat$ (if $X=\Spec F$), write $D_\fin(S_X)$ (resp. $D_\ell(S_X)$ with $\ell$ a prime) the full subcategory of bounded objects of $D(S_X)$ with cohomology representable by finite (resp. $\ell$-primary) étale $X$-groups.

In any symetric monoidal closed category $\C$, given objects $X,Y,Z\in\C$ we say a morphism $X\otimes Y\to Z$ is a perfect pairing if both adjunct maps $X\to\HOM(Y,Z)$ and $Y\to\HOM(X,Z)$ are isomorphisms.

Let $K=k_d,\dots,k_0$ be a sequence of fields such that $k_0$ is perfect of characteristic $p>0$, and for $i>0$ the field $k_i$ is a Henselian discrete valuation field with residue field $k_{i-1}$ and either $\characteristic k_i=0$ or $k_i$ is complete. Let $X$ be a proper, smooth, geometrically integral $K$-scheme of relative dimension $e$ and $U\subseteq X$ a nonempty open. We consider the object $\QQ/\ZZ(r)=\varinjlim_m\mu_m^{\otimes r}$ of $\Sh(U_\et)$.

\begin{spetheorem}[intro general]
    (\Theorem{higher local duality}) Assume $d\geq 2$ and $\characteristic k_2=0$.
    
    There exist functors $R\Psi_{U/k_0},R\Psi_{U/k_0,c}:D(U_\et)\to D(k_{0,\proet}^\indrat)$ with the following properties.
    \begin{enumerate}
        \item For prime $\ell\neq p$, $R\Psi_{K/k_0}$ and $R\Psi_{U/k_0,c}$ send $D_\ell(K_\et)$ to $D_\ell(k_{0,\proet}^\indrat)$, and $D_p(K_\et)$ to $D^b_{\W_{k_0}}(k_{0,\proet}^\indrat)$.
        \item There exists a canonical trace map $\tr:R\Psi_{U/k_0,c}(\QQ/\ZZ(d+e))\to\QQ/\ZZ[-d-2e]$ in $D(k_{0,\proet}^\indrat)$.
        \item If $G\in D_\fin(U_\et)$ and $G^\vee=R\HOM_{U_\et}(G,\QQ/\ZZ(d))$, then we have a perfect pairing in $D(k_{0,\proet}^\indrat)$ :
        $$R\Psi_{U/k_0}(G)\otimes^LR\Psi_{U/k_0,c}(G^\vee)\to R\Psi_{U/k_0,c}(\QQ/\ZZ(d+e))\xrightarrow\tr\QQ/\ZZ[-d-2e].$$
        \item If $G\in D_\fin(U_\et)$ and $q\in\ZZ$ then $R^q\Psi_{U/k_0}(G)_{\Adh k_0}=H^q(U_{K^\tur},G)$ where $K^\tur$ is the maximal totally unramified extension of $K$ (that is, the union of all finite extensions of $K$ which induce unramified extensions of $k_i$ for all $i\geq 1$), and $(-)_{\Adh k_0}=\varinjlim_{l/k_0 \text{ finite}}\Gamma(l,-)$ denotes the étale stalk at $\Adh k_0$.
    \end{enumerate}
\end{spetheorem}

When $d=1$, still $\characteristic K=0$, and $U=X$, the analogous theorem is proven in \cite[Th. 2.8.1]{SuzukiDualLoc} using the étale rational site $\Spec k_\et^{\mathrm{rat}}$ instead. The essential case of this theorem, with $d=2$, $e=0$ and $G=\mu_{p^n}^{\otimes r}$, is treated in \cite[Prop. 5.10]{SuzukiCFT} and \cite[Prop. 6.2.2]{Suzuki2024}, albeit in the language of relative sites,.

The duality of \Theorem{intro general} is not a purely derived abstraction. If $G\in D_\ell(U_\et)$ with $\ell\neq p$, then for $q\in\ZZ$ we have perfect pairing of finite étale $k_0$-groups :
$$R^q\Psi_{U/k_0}(G)\otimes R^{d+2e-q}\Psi_{U/k_0,c}(G^\vee)\to R^{d+2e}\Psi_{U/k_0,c}(\QQ/\ZZ(d+e))\xrightarrow\tr\QQ/\ZZ.$$
If $G\in D_p(U_\et)$, then for $q\in\ZZ$ we have so-called Serre dualities of the neutral components and component groups of the cohomology, in the form :
\begin{align*}
    R^q\Psi_{U/k_0}(G)^0&=\EXT^1_{k_{0,\proet}^\indrat}(R^{d+2e+1-q}\Psi_{U/k_0,c}(G^\vee)^0,\QQ/\ZZ),\\
    \pi_0(R^q\Psi_{U/k_0}(G))&=\HOM_{k_{0,\proet}^\indrat}(\pi_0(R^{d+2e-q}\Psi_{U/k_0,c}(G^\vee)),\QQ/\ZZ).
\end{align*}

\Theorem{intro general} gives a concrete statement for the cohomology of $U$ when $k_0$ is algebraically closed. If $k_0$ is finite, $R\Gamma(k_0,-)$ sends ind-pro-RPAU $k_0$-groups to ind-pro-finite groups, so we can give concrete statements also in that case. This fails if $k_0$ is only quasi-finite, so in constrast with the case of prime-to-$p$ coefficients, a generalization to arbitrary quasi-finite $k_0$ seems out of reach.

Let $*_\proet$ be the proétale site of a point (\Definition{proetale point definition} and \cite{ScholzeClausen}). Recall that locally compact Hausdorff topological spaces can be fully faithfully embedded as sheaves of sets over $*_\proet$. Write $D_\fin(*_\proet$ (resp. $D_\ell(*_\proet)$ with $\ell$ a prime) the subcategory of bounded objects of $D(*_\proet)$ with cohomology representable by finite discrete (resp. and $\ell$-primary) abelian groups. Write $D^b_{\W_p}(*_\proet)$ the subcategory of bounded objectsof $D(*_\proet)$ with cohomology representable by "good" ind-pro-finite $p$-groups, which are locally compact Hausdorff topological abelian groups (\Definition{topological abelian groups} and \Definition{proetale point definition}).

\begin{spetheorem}[intro finite]
    (\Theorem{higher local duality finite}) Assume $d\geq 2$, $\characteristic k_2=0$ and $k_0$ is finite.
    
    There exist functors $R\Psi_{U/*},R\Psi_{U/*,c}:D(U_\et)\to D(*_\proet)$ with the following properties.
    \begin{enumerate}
        \item For any prime $\ell\neq p$, $R\Psi_{U/*}$ and $R\Psi_{U/*,c}$ send $D_\ell(U_\et)$ to $D_\ell(*_\proet)$, and $D_p(U_\et)$ to $D^b_{\W_p}(*_\proet)$. In particular, if $G\in D_\fin(U_\et)$ then $R^q\Psi_{U/*}(G)$ and $R^q\Psi_{U/*,c}(G)$ are representable by locally compact Hausdorff topological abelian groups for all $q\in\ZZ$.
        \item There exists a canonical trace map $\tr:R\Psi_{U/*,c}(\QQ/\ZZ(d+e))\to\QQ/\ZZ[-d-2e-1]$ in $D(*_\proet)$.
        \item If $G\in D_\fin(U_\et)$ and $G^\vee=R\HOM_{U_\et}(G,\QQ/\ZZ(d+e))$ then we have a perfect pairing in $D(k_{0,\proet}^\indrat)$ :
        $$R\Psi_{U/*}(G)\otimes^LR\Psi_{U/*,c}(G^\vee)\to R\Psi_{U/*,c}(\QQ/\ZZ(d+e))\xrightarrow\tr\QQ/\ZZ[-d-2e-1].$$
        For $q\in\ZZ$, this induces perfect Pontryagin dualities of locally compact Hausdorff abelian groups :
        $$R^q\Psi_{U/*}(G)\otimes R^{d+2e+1-q}\Psi_{U/*,c}(G^\vee)\to R^{d+2e+1}\Psi_{U/*,c}(\QQ/\ZZ(d+e))\xrightarrow\tr\QQ/\ZZ.$$
        \item If $G\in D_\fin(U_\et)$ and $q\in\ZZ$ then $R^q\Psi_{U/*}(G)(*)=H^q(U,G)$ and $R^q\Psi_{U/*,c}(G)(*)=H^q_c(U,G)$. In other words, $R^q\Psi_{U/*}(G)$ and $R^q\Psi_{U/*,c}(G)$ can be viewed as the groups $H^q(U,G)$ and $H^q_c(U,G)$ equipped with certain locally compact Hausdorff group topologies.
    \end{enumerate}
\end{spetheorem}

\begin{specorollary}[intro cor] (\Corollary{higher local interpretation bis})
    For $G\in D_\fin(U_\et)$ with $G^\vee=R\HOM_{U_\et}(G,\QQ/\ZZ(d+e))$ and $q\in\ZZ$ the cup-product pairing $H^q(U,G)\otimes H^{d+2e+1-q}_c(U,G^\vee)\to H^{d+2e+1}(U,\QQ/\ZZ(d+e))\to \QQ/\ZZ$ is nondegenerate.
\end{specorollary}


\subsection*{Structure of the paper}
\label{structure of the paper}
\addcontentsline{toc}{subsection}{Structure of the paper}

The first three parts of the paper focus on recalling, detailing and expanding upon the various formalisms used to obtain the duality results of \cite{KatoSuzuki2019, SuzukiCFT} mentioned in the introduction. The fourth focuses on generalizing these duality results, assembling them into points \textbf{(1.-3.)} of \Theorem{intro general} and \Theorem{intro finite}, and relating them to actual étale cohomology to obtain points \textbf{(4.)} of these theorems and \Corollary{intro cor}.

In part \ref{Relatively perfect group schemes} we give a review on commutative relatively perfect group schemes after \cite[Sec. 2 and 8]{BertapelleSuzuki} and \cite[Sec. 3]{SuzukiCFT}. Unipotent such groups, abbreviated RPAU groups, satisfy a duality first outlined in \cite{Kato86} and expanded upon in \cite{KatoSuzuki2019, SuzukiCFT}. We especially detail the dévissage of RPAU groups, the equivalent descriptions of their derived category, and their behaviour with regards to finite field extensions.

In part \ref{Ind-pro-RPAU groups} we recall the formalism of the sites $\Spec k_\proet^\indrat$ and $\Spec k_\et^\perar$ for a perfect field $k$ due to Suzuki \cite{SuzukiGN, SuzukImp, SuzukiDualLoc}. We focus on their main features of interest : the class of good ind-pro-RPAU $k$-groups $\W_k$ can be viewed as sheaves over them in a derived sense, and satisfies a duality as such. We prove $\W_k$ is stable under extensions, a property used in \cite{SuzukiCFT} to prove duality in the equal characteristic $p$ step.

In part \ref{Ind-pro-finite groups} we detail some (possibly folklore) comparisons of ind-pro-finite groups, locally compact groups, and condensed groups : the latter fully faithfully contain the former two in a derived sense. The class of good ind-pro-finite groups $\W_\fin$ is in the intersection ; it is stable under Pontryagin duality and extensions of condensed groups, and such extensions respect topology. If $k$ is a finite field of characteristic $p>0$, $R\Gamma(k,-)$ sends $\W_k$ to $\W_\fin$ and satisfies a duality - this is the step separating \Theorem{intro general} and \Theorem{intro finite}.

Finally in part \ref{Duality for higher local fields} we prove the main results. We write down the steps corresponding to \textbf{(a.)}, \textbf{(b.)} and \textbf{(c.)} in the introduction, largely citing work of Kato and Suzuki. Extra work is provided in extending the mixed characteristic step to general $p$-primary coefficients, adapting techniques used in \cite{BertapelleSuzuki}. We deal with prime-to-$p$ coefficients with the same languages. The higher local duality theorem follows nicely but the interpretation for Galois cohomology still requires some work, to which we dedicate a section.

In an appendix, we recall Kato's canonical lifting and study its structure.

\emph{Acknowledgements.} The author would like to thank Cyril Demarche for his suggestion and continuous support of this project. He would also like to thank Takashi Suzuki for gracefully sharing his insights and comments : several results in this paper are directly based on his advice.


\subsection*{Notations and conventions}
\label{Notations and conventions}
\addcontentsline{toc}{subsection}{Notations and conventions}

All sites considered are defined by pretopologies. Namely for our purposes, a site is a category $S$ equipped with a class of families of maps with common target called \emph{covering families}, such that : $S$ has fiber products $X\times_ZY$ for any $X\to Z$ appearing in a covering family and any $Y\to Z$ ; any isomorphism in $S$ is a one-term covering family ; and covering families are stable under composition and arbitrary base change. See \cite[Def. 00VH]{stacks-project}.

For $S$ such a site and $\Lambda$ a ring, we write $\PSh(S,\Lambda)$ and $\Sh(S,\Lambda)$ the categories of presheaves and sheaves of $\Lambda$-modules on $S$. We write $D(S,\Lambda)=D(\Sh(S,\Lambda))$ its derived category (and $D^b(S,\Lambda)$, $D^+(S,\Lambda)$, $D^-(S,\Lambda)$ the bounded subcategories). We write $\otimes_\Lambda$ and $\HOM_{S,\Lambda}$ the tensor product and inner Hom of $\Sh(S,\Lambda)$. If $\Lambda=\ZZ$, the ring $\Lambda$ is omitted from notations and we write simply $\PSh(S)$, $\Sh(S)$, $D(S)$, $\otimes$ and $\HOM_S$. On occasion we write $\PSh^\set(S)$ and $\Sh^\set(S)$ to denote (pre)sheaves of sets on $S$. Note that $\PSh^\set(S)$ and $\PSh(S,\Lambda)$ only depend on the underlying category of $S$, and we will use those notations also for categories not equipped with any site structure.

A \emph{premorphism} between two sites defined by pretopologies $v:S\to S'$ is a functor of underlying categories $v^{-1}:S'\to S$ which sends covering families to covering families, and $v^{-1}(X)\times_{v^{-1}(Z)}v^{-1}(Y)$ for any $X\to Z$ appearing in a covering family and any $Y\to Z$ in $S'$. Precomposing by $v^{-1}$ defines a functor $v_*:\Sh^\set(S)\to \Sh^\set(S')$ which restricts to $v_*:\Sh(S,\Lambda)\to\Sh(S',\Lambda)$ for any ring $\Lambda$, and both admit left adjoints $v^\set:\Sh^\set(S')\to\Sh^\set(S)$ and $v^*:\Sh(S',\Lambda)\to\Sh(S,\Lambda)$. The functor $v^*$ admits a left derived functor $Lv^*:D(S',\Lambda)\to D(S,\Lambda)$, which is right adjoint to $Rv_*:D(S)\to D(S')$, and there exists a natural morphism $Rv_*(-)\otimes^LRv_*(-)\to Rv_*(-\otimes^L-)$ in $D(S')$ called the cup-product. When $v^\set$ is exact we say $v$ is a \emph{morphism} of site : in that case $v^*$ is an exact monoidal functor. See \cite[Sec. 2]{SuzukImp}.

For $X$ a scheme, we write $X_\Et$ the \emph{big étale site} of $X$ \cite[Def. 021A and 021B]{stacks-project}. For $x=\Spec k$ a point of $X$, we write $\Adh x=\Spec \Adh k$ a fixed choice geometric point (separable closure) over $x$.

Here a \emph{filtration} of an object $G$ of an abelian category $\A$ is a sequence of subobjects $\dots\leq G_{i+1}\leq G_i\leq\dots\leq G_1\leq G_0=G$ in $\A$. Its $i$-th quotient is $G_i/G_{i+1}$. It is \emph{finite separated} if $G_i=0$ for some $i\geq 0$.


\subsection*{Reminders on ind-pro-categories}
\label{Reminders on ind-pro-categories}

For $\A$ a category, we write $\Ind\A$ its ind-category and $\Pro\A=(\Ind\A^\opp)^\opp$ its pro-category. We refer to \cite[Sec. 6, 8.6 and 15]{KashiwaraSchapira} and \cite[Sec. 2.2]{SuzukiGN} for all details pertaining to them ; we recall the key points here.

\textbf{Definition.} Let $\U$ be a fixed universe containing $\NN$. Objects of $\Ind\A$ are filtered $\U$-small diagrams $\{X_i\}_{i\in I}$ of objects of $\A$, viewed as formal filtered colimits $\varinjlim_{i\in I}X_i$. Given such $X=\varinjlim_iX_i$ and $Y=\varinjlim_jY_j$ in $\Ind\A$, with $X_i,Y_j\in\A$, morphisms $X\to Y$ are defined by :
$$\Hom_{\Ind\A}\left(\varinjlim_iX_i,\varinjlim_jY_j\right) =\varprojlim_i\varinjlim_j\Hom_\A(X_i,Y_j).$$

The category $\A$ fully faithfully embeds in $\Ind\A$. The category $\Ind\A$ has all $\U$-small filtered colimits, and objects of $\A$ are exactly the compacts of $\Ind\A$. If $X=\{X_i\}_i$ is a filtered diagram of $\A$ then we genuinely have a colimit $X=\varinjlim_iX_i$ in $\Ind\A$. The formation of $\Ind\A$ is functorial in $\A$, and if $\A$ has $\U$-small filtered colimits then their computation in $\A$ defines a functor $\varinjlim:\Ind\A\to\A$ left adjoint to the inclusion $\A\to\Ind\A$.

\textbf{Abelian case.} If $\A$ is abelian then so is $\Ind\A$ and any short exact sequence in $\Ind\A$ arises as a filtered colimit of short exact sequences of $\A$ \cite[Th. 8.6.5 and Prop. 8.6.6]{KashiwaraSchapira}. Filtered colimits are exact in $\Ind\A$, and the formation of $\Ind\A$ takes left (resp. right) exact functors to left (resp. right) exact functors of ind-categories \cite[Cor. 8.6.8]{KashiwaraSchapira}. If $\A$ is abelian and small, then $\Ind\A$ is a Grothendieck category hence admits enough injectives \cite[Th. 8.6.5]{KashiwaraSchapira}. In general $\Ind\A$ does not have enough injectives \cite[Cor. 15.1.3]{KashiwaraSchapira}. However if $F:\A\to\B$ is a left (resp. right) exact functor of abelian categories and $\A$ has enough $F$-injectives (resp. $F$-projectives) then $\Ind\A$ has enough $\Ind F$-injectives (resp. $\Ind F$-projectives) and $R^q(\Ind F)=\Ind(R^qF)$ (resp. $L_q(\Ind F)=\Ind(L_qF)$) for $q\in\ZZ$ \cite[Prop. 15.3.2 and 15.3.7]{KashiwaraSchapira}.

\textbf{Derived Hom-functor.} We will primarily consider categories of the form $\IndPro\A$ where $\A$ is an essentially small abelian category (either the category of RPAU groups over a perfect field, or a category of finite abelian groups). Then \cite[Prop. 2.2.2]{SuzukiGN} shows the derived functor $R\Hom_{\IndPro\A}:D^-(\IndPro\A)^\opp\times D^+(\IndPro\A)\to D^+(\Ab)$ exists and provides a method to compute it. Namely we have spectral sequences and equalities :
\begin{align*}
    E^{pq}_2= R^p\varprojlim_i\Ext^q_{\IndPro\A}(X_i,Y)\Rightarrow\Ext^{p+q}_{\IndPro\A}(\varinjlim_iX_i,Y)&&\Ext^n_{\IndPro\A}(X_0,\varinjlim_kY_k)=\varinjlim_k\Ext^n_{\Pro\A}(X_0,Y_k)\\
    'E^{pq}_2= R^p\varprojlim_l\Ext^q_{\Pro\A}(X_0,Y_l)\Rightarrow\Ext^{p+q}_{\Pro\A}(X_0,\varprojlim_lY_l)&&\Ext^n_{\Pro\A}(\varprojlim_jX_j,Y_0)=\varinjlim_j\Ext^n_\A(X_j,Y_0)
\end{align*}
whenever $Y\in\IndPro\A$, $X_i,X_0,Y_k\in\Pro\A$, and $X_j,Y_0\in\A$.

This will often be summarized by writing, for $X_{ij},Y_{kl}\in\A$ :
$$R\Hom_{\IndPro\A}(\varinjlim_i\varprojlim_jX_{ij},\varinjlim_k\varprojlim_lY_{kl}) = R\varprojlim_i\varinjlim_kR\varprojlim_l\varinjlim_jR\Hom_\A(X_{ij},Y_{kl}).$$
Note that the inverse limits inside $R\Hom_{\IndPro\A}$ are underived. This equality leads to a correct understanding when considering Ext's, as illustrated by the previous spectral sequences, but may be misleading.

First, $\Hom_\A$ may not genuinely have a derived functor $R\Hom_\A:D^-(\A)^\opp\times D^+(\A)\to D^+(\Ab)$. Instead for $X_0,Y_0\in\A$ we take $R\Hom_\A(X_0,Y_0)=R\Hom_{\Pro\A}(X_0,Y_0)\in D^+(\Ab)$. With this we still have $H^nR\Hom_\A(X_0,Y_0)=\Ext^n_\A(X_0,Y_0)$ (Ext-group in the sense of Yoneda), so this recovers the actual derived $R\Hom_\A$ if it exists, and is no problem when computing cohomology.

Second, $R\varprojlim$ and $\varinjlim$ do not denote genuine (derived) (co)limits in the derived category $D^+(\Ab)$. Varying $i,j,k,l$ we can obtain diagrams of objects $R\Hom_\A(X_{ij},Y_{kl})\in D^+(\Ab)$, and the right-hand term of the above equality misleadingly suggests taking successive (derived) (co)limits of these diagrams, as a composite :
$$(\IndPro\A)^\opp\times\IndPro\A\xrightarrow{R\Hom_\A}\Pro\Ind\Pro\Ind D^+(\Ab)\xrightarrow{R\varprojlim (\Pro\varinjlim) (\Pro\Ind R\varprojlim) (\Pro\Ind\Pro\varinjlim)}D^+(\Ab).$$
where $\varinjlim$ and $R\varprojlim$ denote (derived) (co)limit-computing functors. We do not claim that such a computation makes sense or recovers $R\Hom_{\IndPro\A}$. What is true is we have a factorization of $R\Hom_{\IndPro\A}$ as :
\begin{align*}
    D^-(\IndPro\A)^\opp\times D^+(\IndPro\A)
    &\xrightarrow{R\Pro\Ind\Pro\Ind\Hom_\A}D^+(\Pro\Ind\Pro\Ind\Ab)
    \xrightarrow{\Pro\Ind\Pro\varinjlim} D^+(\Pro\Ind\Pro\Ab)\\
    &\xrightarrow{R\Pro\Ind\varprojlim}D^+(\Pro\Ind\Ab)
    \xrightarrow{\Pro\varinjlim}D^+(\Pro\Ab)
    \xrightarrow{R\varprojlim}D^+(\Ab)
\end{align*}
where $\varinjlim:\Ind\B\to\B$ and $\varprojlim:\Pro\B\to\B$ are the (co)limit-computing functors for $\B=\Ab,\Pro\Ab,\Ind\Pro\Ab,\dots$ Then $\varinjlim,\Pro\varinjlim,\Ind\Pro\varinjlim,\dots$ are exact and $\varprojlim,\Ind\varprojlim,\Pro\Ind\varprojlim,\dots$ are left-exact, so they admit derived functors as in the above composition. If $F$ is a left-exact functor of abelian categories with enough $F$-injectives, we have $R^n\Ind F=\Ind R^nF$ and $R^n\Pro F=\Pro R^nF$, so again subtleties vanish when only considering cohomology.


\section{Relatively perfect group schemes}
\label{Relatively perfect group schemes}


\subsection{Relatively perfect schemes}
\label{Relatively perfect schemes}

Fix $S$ an $\FF_p$-scheme. For $X/S$, write $X^{(p/S)}$ the base change along the absolute Frobenius $F_S:S\to S$. In particular $S^{(p)}=S^{(p/S)}$ is $S$ with the $S$-scheme structure given by $F_S$. Then $F_X$ factors uniquely through a map $F_{X/S}:X\to X^{(p/S)}$, called the \emph{relative Frobenius of $X$ over $S$}. When $S=\Spec k$ and $X=\Spec R$ are affine, we use similar notations $k^{(p)}$, $R^{(p/k)}$, $F_k:k\to k$ and $F_{R/k}:R^{(p/k)}\to R$ for the corresponding $k$-algebras and ring maps.

\begin{definition}[relatively perfect]
    We say $X/S$ is \emph{relatively perfect over $S$} if its relative Frobenius $F_{X/S}$ is an isomorphism. The category of relatively perfect $S$-schemes with all $S$-scheme morphisms is written $\RPSch/S$. The \emph{relatively perfect site of} $S$ is category $\RPSch/S$ equipped with the étale topology, which we write $S_\RP$.
\end{definition}

More precisely, $S_\RP$ is the full subcategory of the big étale site $S_\Et$ of \cite[Def. 021B]{stacks-project} of objects relatively perfect over $S$, with the étale topology.

If $S=\Spec k$ is affine, we also say an $S$-scheme is relatively perfect over $k$ if it is so over $\Spec k$, and write $\RPSch/k=\RPSch/S$. If $R$ is a $k$-algebra, we say $R$ is relatively perfect over $k$ if $\Spec R$ is so over $\Spec k$, or equivalently $F_{R/k}$ is an isomorphism. For $X$ a $k$-scheme, the relative Frobenius $F_{X/S}$ correspond to the map of ring sheaves $F_{\O_X/k}:\O_X\otimes_kk^{(p)}\to\O_X$, defined by $F_{\O_X/k}(U)=F_{\O_X(U)/k}$ for $U\subseteq X$ open. In particular, if $X\in\RPSch/k$ then the global sections $\O(X)$ form a relatively perfect $k$-algebra.

Relatively perfect morphisms are stable under composition and base change \cite[Lem. 1.2]{Kato86} ; in particular, the formation of the site $S_\RP$ is functorial with respect to morphisms of $\FF_p$-schemes. A relatively perfect morphism is formally étale, and the relatively perfect morphisms locally of finite presentation are exactly étale morphisms \cite[Lem. 1.3]{Kato86}. In particular for $X\in S_\RP$, $X_\et$ is contained in $S_\RP$.

Many examples of relatively perfect morphisms can be described using $p$-bases.

\begin{definition}[p-basis def]
    Let $k$ be an $\FF_p$-algebra and $B\subseteq k$ a family. We denote by $\tilde B$ the family $\left(\prod_{x\in X} x^{m(x)}\right)$ indexed by functions $m:B\to\{0,\dots,p-1\}$ such that $m(x)=0$ for all but finitely many $x\in B$.
    \begin{itemize}
        \item We say $B$ is a \emph{weak $p$-basis} of $k$ if $\tilde B$ is a free basis of $k$ over the subalgebra of $p$-powers $k^p\subseteq k$.
        \item We say $B$ is a \emph{strong $p$-basis} of $k$ if $\tilde B$ is a free basis of $k^{(p)}$ as a $k$-module.\footnote{This coincides with the notion of $p$-basis used in \cite[Def. 1.3]{Kato91}, by \Lemma{p-basis comparison}.2.}
    \end{itemize}
    We say the structure sheaf of an $\FF_p$-scheme $S$ admits a (resp. finite) strong $p$-basis Zariski-locally if there is a covering by affine opens $S=\bigcup_i\Spec R_i$ such that each $R_i$ admits a (resp. finite) strong $p$-basis.
\end{definition}

In the literature the term "$p$-basis" more commonly means the weak notion, but for our purposes the strong one is more useful. The prototypical example of a strong $p$-basis is $(t_1,\dots,t_n)$ for $\FF_p[t_1,\dots,t_n]$. Any field of characteristic $p>0$ admits a strong $p$-basis with $\log_p[k:k^p]\in\NN\cup\{+\infty\}$ elements.

\begin{lemma}[p-basis comparison]
    Let $k$ be a ring of prime characteristic $p>0$ and $R$ a $k$-algebra.
    \begin{enumerate}
        \item A family $B$ of $k$ is a strong $p$-basis if and only if it is a weak $p$-basis and $k$ is reduced.
        \item If $k$ has a strong $p$-basis $B$ then $R$ is relatively perfect if and only if $B$ is a strong $p$-basis of $R$.
        \item If $k$ is reduced and $R$ is flat and reduced over $k$ then $R$ is reduced.
    \end{enumerate}
    In particular if $k$ admits a strong $p$-basis, then any relatively perfect $k$-algebra is reduced.
\end{lemma}

In general a weak $p$-basis may not be strong : $t$ is a one-term weak $p$-basis of $\FF_p[t]/(t^p)$ but cannot be a strong $p$-basis as the latter is not reduced. Note that if $k$ is regular Noetherian then any relatively perfect $k$-algebra is $k$-flat \cite[Prop. 5.2]{Kato86}.

\begin{myproof}[of \Lemma{p-basis comparison}]
    \textbf{1.} If $B$ is a strong $p$-basis and $a\in\ker(F_k)$ then $a^p\cdot 1$ and $0^p\cdot 1$ are two decompositions of $a^p$ in the $k$-free basis $\tilde B$ of $k^{(p)}$, so by uniqueness $a=0$ and $k$ is reduced. If $k$ is reduced, then $F_k$ induces an isomorphism $k\xrightarrow\sim k^p$, thus identifies free bases of $k$ over $k^p$ and of $k^{(p)}$ over $k$. This concludes.
    
    \textbf{2.} The strong $p$-basis assumption gives a $k$-module isomorphism $k^{(p)}\cong k^{\oplus\tilde B}$ hence an $R$-module isomorphism $R\otimes_kk^{(p)}\cong R^{\oplus\tilde B}$. The relative Frobenius is an $R$-module morphism when viewed as a map $R\otimes_kk^{(p)}\to R^{(p)}$, given as $(x,y)\mapsto y^px$, thus is an isomorphism if and only if the free $R$-basis $\tilde B$ of $R^{\oplus\tilde B}$ is mapped to a free $R$-basis of $R^{(p)}$. Hence $R$ is relatively perfect over $k$ if and only if the image of $B$ in $R$ under the structure map $k\to R$ is a strong $p$-basis of $R$.
    
    \textbf{3.} An $\FF_p$-algebra is reduced if and only if its absolute Frobenius is injective. The absolute Frobenius of $R$ decomposes as $F_{R/k}\circ(1\otimes F_k):R=R\otimes_kk\to R\otimes_kk^{(p)}\to R$, where $1\otimes F_k$ is injective because $k$ is reduced and $R$ is $k$-flat, and $F_{R/k}$ is an isomorphism. Thus $F_R$ is injective and $R$ is reduced.
\end{myproof}

\begin{proposition}[relative perfection]
    \cite[Prop. 1.4 and Lem. 1.5]{Kato86} Assume $F_S$ is finite locally free.
    \begin{enumerate}
        \item The functor $X\mapsto X^{(p/S)}$ on $\Sch/S$ has a right adjoint, \emph{i.e.} the Weil restriction $\Res_{F_S}(-)$ exists.
        \item The inclusion $\RPSch/S\hookrightarrow\Sch/S$ admits a right adjoint, called the \emph{relative perfection over $S$}, written $(-)^{\RP/S}$ or simply $(-)^\RP$, given by $X^{\RP/S}=\varprojlim_{n\geq 0}\Res_{F_S}^n(X)$.
    \end{enumerate}
\end{proposition}

The transitions in $\Res_{F_S}^{n+1}(X)\to\Res_{F_S}(X)$ are given by the natural map $\Res_{F_S}\to\id$ corresponding, under the adjunction between $(-)^{(p/S)}$ and $\Res_{F_S}$, to the relative Frobenius $\id\to (-)^{(p/S)}$. In other words, $\Res_{F_S}^{n+1}(X)\to\Res_{F_S}^n(X)$ is the map $\Res_{F_S}^n(\varepsilon_X\circ F_{X/S})$ where $\varepsilon_X:\Res_{F_S}(X^{(p/S)})\to X$ is the counit.

The assumption on $F_S$ is true, for example, if $S$ is smooth over a field $k$ of characteristic $p>0$ such that $[k:k^p]<+\infty$ ; or more generally, if the structure sheaf of $S$ Zariski-locally admits a finite strong $p$-basis. Relative perfection preserves any property preserved by projective limits with affine transitions and Weil restrictions. For instance, the relative perfection of an affine (resp. quasi-compact, quasi-separated) $S$-scheme is again affine (resp. quasi-compact, quasi-separated).

\begin{example}[affine RP]
    Let $k$ be an $\FF_p$-algebra which admits a finite strong $p$-basis $B=(x_1,\dots,x_r)$ and consider $\tilde B=\{x^m : m\in I\}$ be the corresponding free $k$-pasis of $k^{(p)}$, where $I=\{(m_1,\dots,m_r), 0\leq m_i<p\}$ and we write $x^m=x_1^{m_1}\cdots x_r^{m_r}$ for $m=(m_1,\dots,m_r)$. More generally, the $k$-algebra $k^{(p^n)}$ given by the $n$-fold Frobenius $F_k^n:k\to k$ admits a free $k$-linear basis given by :
    $$\tilde B_n=\{x^M : M\in I_n\}, \qquad I_n=\{(M_1,\dots,M_r) : 0\leq M_i<p^n\}.$$
    By \Lemma{p-basis comparison}.2, the $k$-algebra $k=k[t^{1/p^\infty}]=k[t^{1/p^n}, n\geq 0]$ is relatively perfect since clearly $B$ is a strong $p$-basis of it. However unless $k$ is perfect, $k[t^{1/p^\infty}]$ is not the relative perfection of $k[t]$ : any element of $k\setminus k^p$ defines a map $k[t]\to k$ which does not factor through $k[t^{1/p^\infty}]$, showing failure of the universal property. Rather, $k[t]^\RP$ is the smallest $k[t]$-algebra "closed under taking coordinates in $\tilde B$" : one can show that $\Res_{F_k}^n(k[t])$ is the polynomial $k[t]$-algebra $k[t_M, M\in I_n]$ whose structure morphism "decomposes $t$ in $\tilde B_n$" :
    \begin{align*}
        k[t] &\to k[t_M:M\in I_n]\\
        t & \mapsto \sum_{M\in I_n}t_M^{p^n}x^M
    \end{align*}
    and the transition morphism $\Res_{F_k}^n(k[t])\to\Res_{F_k}^{n+1}(k[t])$ is defined by sending $t_M$ to $\sum_{m\in I}t_{pM+m}^px^m$, making the variables in $\Res_{F_k}^{n+1}(k[t])$ the coordinates of those of $\Res_{F_k}^n(k[t])$ in the $p$-basis $B$. Then :
    $$k[t]^\RP = \left.k\left[t_{(M,n)}, (M,n)\in\bigsqcup_{n\geq 1} I_n\right]\middle/\left\langle t_{(M,n)}-\sum_{m\in I}t_{(pM+m,n+1)}^px^m, (M,n)\in\bigsqcup_{n\geq 1} I_n\right\rangle\right..$$
    In particular $k[t]^\RP$ is not Noetherian : the ideal $\langle t_{(M,n)}, (M,n)\in \bigsqcup_nI_n\rangle$ is not finitely generated.
\end{example}

\begin{definition}[RPS scheme]
    Assume $F_S$ is finite locally free. We say a scheme $X/S$ is \emph{relatively perfectly smooth (resp. relatively perfectly locally of finite presentation) over $S$} if, Zariski-locally on $X$ and $S$, $X$ is isomorphic to the relative perfection of a smooth (resp. finite presentation) $S$-scheme. We write $\RPSSch/S$ the category of relatively perfectly smooth $S$-schemes. The \emph{relatively perfectly smooth site} of $S$ is the category $\RPSSch/S$ equipped with the étale topology. We write it $S_\RPS$.
\end{definition}

\begin{lemma}[base change RP]
    Let $S$ be an $\FF_p$-scheme with and $T$ an $S$-scheme
    \begin{enumerate}
        \item Assume $F_S$ is finite locally free, $T$ is relatively perfect over $S$ and $X$ is an $S$-scheme. Then we have a canonical isomorphism $(X\times_ST)^{\RP/T} = X^{\RP/S}\times_ST$.
        \item Assume there exists a map $g:S\to T$ which forms a commutative diagram :
        $$\begin{tikzcd}
            T\arrow[r,"F_T^n"]\arrow[d] & T \arrow[d] \\
            S\arrow[ru,"g"] \arrow[r,"F_S^n"] & S
        \end{tikzcd}$$
        for some $n\geq 1$. Then the base change functor $(-)\times_ST : \RPSch/S\to \RPSch/T$ is an equivalence of categories with essential inverse given by base change by $g$.
        \item Consider $S$, $T$ and $g:S\to T$ as in the previous point. Assume furthermore that $F_S$, $F_T$ and $g$ are finite locally free. Then for $X\in\Sch/S$ there is an isomorphism of $T$-schemes :
        $$X^{\RP/S}\times_ST = (\Res_gX)^{\RP/T}=\Res_g(X^{\RP/S})$$
        and the base change functor restricts to an equivalence $(-)\times_ST : \RPSSch/S\to \RPSSch/T$ with essential inverse given by base change by $g$.
        \item Assume $S$ is the spectrum of a field with finite Frobenius, and $T$ is the spectrum of a finite field extension. Then $(-)\times_ST$ sends $\RPSSch/S$ to $\RPSSch/T$.
    \end{enumerate}
\end{lemma}

As seen in the proof, it is equivalent in \textbf{(3.)} to assume $F_S$, $F_T$ and $T\to S$ (rather than $g$) are finite locally free. This lemma implies that the site $S_\RPS$ is functorial with respect to maps of $\FF_p$-schemes with finite locally free Frobenius, which either are étale, satisfy the assumptions of \textbf{(3.)}, or are finite extensions of fields with finite Frobenius. Points \textbf{(2.-3.)} hold for instance if $S$ is reduced and quasi-compact, $T\to S$ is a finite universal homeomorphism, and $F_S$ and $F_T$ are finite locally free \cite[Cor. 20]{Kelly}.

\begin{myproof}[of \Lemma{base change RP}]
    \textbf{1.} Since $X\times_ST$ is relatively perfect over $X$, by \cite[Cor. 1.9]{Kato86} we have identifications  $(X\times_ST)^{\RP/S}=X^{\RP/S}\times_X(X\times_ST)=X^{\RP/S}\times_ST$.
    
    \textbf{2.} If $U$ is any $\FF_p$-scheme, by definition the $n$-fold relative Frobenius $Y\to Y^{(p^n/U)}$ is a natural isomorphism of $U$-schemes for $Y$ relatively perfect over $U$, where $Y^{(p^n/U)}$ is the base change of $Y$ along $F_U^n$. Thus $(-)^{(p^n/U)}:\RPSch/U\to\RPSch/U$ is naturally isomorphic to identity. We have a commutative diagram of functors :
    $$\begin{tikzcd}
        \RPSch/S\arrow[r,"(-)^{(p^n/S)}"]\arrow[d,"(-)\times_ST"'] & \RPSch/S \arrow[d,"(-)\times_ST"] \\
        \RPSch/T \arrow[r,"(-)^{(p^n/T)}"']\arrow[ru,"(-)\times_TS"] & \RPSch/T
    \end{tikzcd}$$
    where $(-)\times_TS$ is base change by $g$, and by the above both horizontal arrows are naturally isomorphic to identity. This shows $(-)\times_ST$ and  $(-)\times_TS$ form an equivalence $\RPSch/S\cong\RPSch/T$.
    
    \textbf{3.} Remark that under these assumptions, $f:T\to S$ is also finite locally free. First $g$ is quasi-compact and flat, and it is surjective because $g\circ f = F_T^n$ is surjective, so $f$ is flat because $F_S$ is \cite[Lem. 02L2]{stacks-project}. By \cite[Lem. 0560 and 00F4 (4)]{stacks-project}, $f$ is also finite and of finite presentation because $F_T^n$ and $g$ are, so $f$ is finite locally free. In particular $f$ and $g$ play symetric roles.
    
    The functor $(-)\times_TS:\Sch/T\to\Sch/S$ admits a right adjoint $\Res_g$ by assumption. For $X$ an $S$-scheme and $Y\in\RPSch/T$, we have :
    \begin{align*}
        \Hom_{\RPSch/T}(Y,X^{\RP/S}\times_ST)
        &= \Hom_{\RPSch/S}(Y\times_TS,X^{\RP/S}) &\text{by equivalence,}\\
        &= \Hom_{\Sch/S}(Y\times_TS,X) &\text{by definition of }(-)^{\RP/S},\\
        &= \Hom_{\Sch/T}(Y,\Res_gX) &\text{by definition of }\Res_g\\
        &= \Hom_{\RPSch/T}(Y,(\Res_gX)^{\RP/T}) &\text{by definition of }(-)^{\RP/T}.
    \end{align*}
    thus $X^{\RP/S}\times_ST=(\Res_{S/T}X)^{\RP/T}$. We also have $Z\times_ST = \Res_{S/T}(Z)$ for any $Z\in\Sch/S$, by uniqueness of the right adjoint of $(-)\times_TS$.
    
    By \textbf{(1.)} and symetry of $f$ and $g$, it remains only to show $(-)\times_ST$ sends $\RPSSch/S$ into $\RPSSch/T$. By compatibility of base change with Zariski coverings, we need to see $X^{\RP/S}\times_ST$ is relatively perfectly smooth over $T$ if $X$ is $S$-smooth. This follows easily from the formula $X^{\RP/S}\times_ST=(\Res_{S/T}X)^{\RP/T}$ and the fact that the Weil restriction $\Res_g$ sends $S$-smooth schemes to $T$-smooth ones.
    
    \textbf{4.} Write $S=\Spec k$ and $T=\Spec l$. Consider $X\in\RPSSch/S$. After taking a suitable Zariski covering of $X$, we can assume $X=X_0^\RP$ for some $S$-smooth $X_0$. By \cite[Lem. 030K]{stacks-project} and the associativity of fiber products, we reduce to the cases where $l/k$ is either finite separable or finite purely inseparable. The separable case is by \textbf{(1.)}, since $X_0\times_kl$ is $l$-smooth. In the purely inseparable case, there exists some $n\geq 1$ such that $l^{p^n}$ is contained in $k$ by \cite[Lem. 09HI]{stacks-project} (for instance, we can take $p^n=[l:k]$). This gives a commutative diagram :
    $$\begin{tikzcd}
        k\arrow[r,"F_k^n"]\arrow[d] & k \arrow[d] \\
        l  \arrow[r,"F_l^n"]\arrow[ur] & l
    \end{tikzcd}$$
    and we conclude that $(-)\times_kl$ sends $\RPSSch/k$ to $\RPSSch/l$ by \textbf{(3.)}.
\end{myproof}

The following definition will be used in section \ref{Mixed characteristic step}.

\begin{definition}[RPS site of O]
    Let $K$ be a Henselian discrete valuation field of mixed characteristic $(0,p)$, with residue field $k$ such that $[k:k^p]<+\infty$. Define $\RPSSch/\O_K$ as the category of schemes over the ring of integers $\O_K$ which are $\O_K$-flat and whose special fiber belongs to $\RPSSch/k$. The \emph{relatively perfectly smooth site} of $\Spec\O_K$ is the category $\RPSSch/\O_K$ equipped with the étale topology. We write it $\Spec \O_{K,\RPS}$.
\end{definition}

In the rest of the paper, we will mostly use $S_\RP$ and $S_\RPS$ for $S$ the spectrum of a field. For more properties of relatively perfect schemes over a Dedekind scheme, see \cite{OverkampSuzuki}.

\begin{remark}[affine or qcqs restriction]
    For $\tau=\RP,\RPS$, define $S_{\tau,\qcqs}$ (resp. $S_{\tau,\aff}$) as the category of objects of $S_\tau$ quasi-compact quasi-separated (resp. affine) over $S$, with the étale topology. Most of this part \ref{Relatively perfect group schemes} holds verbatim if we replace $S_\tau$ with $S_{\tau,\qcqs}$ (resp. $S_{\tau,\aff}$), because relative perfection preserves affine, quasi-compact or quasi-separated morphisms, and the RPA (resp. RPAU) group schemes considered thereafter are quasi-compact (resp. affine). One checks that \cite[Sec. 8]{BertapelleSuzuki} and \cite[Sec. 3]{SuzukiCFT} still hold for these sites when restricting the discussion to only quasi-compact quasi-separated (resp. affine) group schemes.
\end{remark}


\subsection{RPAU group schemes}
\label{RPAU group schemes}

From now on $S=\Spec k$ is the spectrum of a field $k$ of characteristic $p>0$ such that $[k:k^p]<+\infty$. Recall that a commutative group scheme $G/k$ is \emph{unipotent} if it is affine and every nontrivial closed subgroup $H\leq G$ admits a nontrivial $k$-group morphism to the additive group $H\to\GG_{a,k}$. See \cite[Sec. IV.2]{DemazureGabriel}.

\begin{lemma}[RPAU equivalence]
    Let $G/k$ be a commutative group scheme. The following are equivalent :
        \begin{enumerate}
            \item $G$ is (resp. unipotent) quasi-compact and relatively perfectly locally of finite presentation over $k$ ;
            \item $G$ is (resp. unipotent) quasi-compact and relatively perfectly smooth over $k$ ;
            \item $G$ is the relative perfection of a (resp. unipotent) quasi-compact, smooth $k$-group.
        \end{enumerate}
    We then say $G$ is \emph{relatively perfectly algebraic (resp. unipotent) over $k$} (abbreviated RPA, resp. RPAU).
\end{lemma}

Over perfect $k$, the perfection of an algebraic group is sometimes called \emph{quasi-algebraic}. RPA groups, where perfection is replaced with relative perfection, generalize quasi-algebraicity to imperfect $k$.

\begin{myproof}[of \Lemma{RPAU equivalence}]
    See \cite[Prop. 8.7]{BertapelleSuzuki} and the discussion preceding it, for the non-unipotent version. A commutative unipotent group satisfying \textbf{(2.)} is the same as the relative perfection of a unipotent, commutative, quasi-compact, smooth $k$-group by \cite[Prop. 3.1]{SuzukiCFT}.
\end{myproof}

\begin{definition}[RPAU group scheme]
    We write $\Alg(k)$ (resp. $\Alg_u(k)$) the category of (resp. unipotent) commutative, quasi-compact smooth $k$-group schemes with $k$-group scheme morphisms. We write $\Alg^\RP(k)$ (resp. $\Alg^\RP_u(k)$) the category of commutative RPA (resp. RPAU) $k$-group schemes with $k$-group scheme morphisms. For $\tau=\RP,\RPS$, we write $\Sh_a(k_\tau)$ (resp. $\Sh_0(k_\tau)$) the full subcategory of $\Sh(k_\tau)$ of sheaves representable by objects of $\Alg^\RP(k)$ (resp. $\Alg_u^\RP(k)$).
\end{definition}

If instead we consider all finite presentation commutative $k$-groups, the resulting category would be an abelian subcategory of $\Sh(k_\fppf)$ \cite[Th. 3.2 of $\mathrm{VI}_{\mathrm A}$]{SGA3}. By considering only smooth $k$-groups, we ensure the exact structure is étale-local\footnote{In other words, the Yoneda functor $\Alg(k)\to\Sh(k_\Et)$ is exact. It is clearly fully faithful and left-exact. An admissible epimorphism $G\to H$ in $\Alg(k)$ is a smooth surjection, so by \cite[Cor. 17.16.3 (ii)]{EGA4} there exists an étale cover $H'\to H$ and a section of $H$-schemes $H'\to G$. Then in $\Sh^\set(k_\Et)$, $H'\to H$ is surjective and decomposes as $H'\to G\to H$, so $G\to H$ is surjective.}, but $\Alg(k)$ is not abelian. Thus it is remarkable that $\Alg^\RP(k)$ and $\Alg_u^\RP(k)$ are abelian categories of étale sheaves, in the sense of \Proposition{RPAU characterization}.

\begin{lemma}[yoneda and pullback]
    Let $u:\C_2\to\C_1$ be a premorphism of sites. Assume $\C_1$ and $\C_2$ admit finite products and $u^{-1}$ preserves those products. Let $\A_i$ be the category of commutative group objects of $\C_i$, $Y_i:\A_i\to\Sh(\C_i)$ be the Yoneda functor followed by sheafification. Then $u^*\circ Y_1 = Y_2\circ u^{-1}$ as functors $\A_1\to\Sh(\C_2)$.
\end{lemma}

\begin{myproof}[of \Lemma{yoneda and pullback}]
    The same statement about the sheafified Yoneda functors $\C_i\to\Sh^\set(\C_i)$ and $u^\set$ is true with no particular assumption on $u$, however in that case $u^{-1}$ does not preserve group objects and does not restrict to a functor $\A_1\to\A_2$. Write $Y_i^p:\A_i\to\PSh(\C_i)$ the Yoneda functor and $a_i:\PSh(\C_i)\to\Sh(\C_i)$ the sheafification functor. For $G\in\A_1$ and $F\in\Sh(\C_2)$ we have :
    \begin{align*}
        \Hom_{\Sh(\C_2)}(u^*Y_1(G),F)
        &=\Hom_{\Sh(\C_1)}(Y_1(G),u_*F)
        =\Hom_{\PSh(\C_1)}(Y_1^p(G),u_*F)
        =(u_*F)(G)\\
        &=F(u^{-1}G)
        =\Hom_{\PSh(\C_2)}(Y_2^p(u^{-1}G),F)
        =\Hom_{\Sh(\C_2)}(Y_2(u^{-1}G),F)
    \end{align*}
    by various adjunctions. Hence $u^*Y_1(G)=Y_2(u^{-1}G)$ canonically.
\end{myproof}

\begin{proposition}[RPAU characterization]
    For $k$ a field of characteristic $p>0$ such that $[k:k^p]<+\infty$, the following holds.
    \begin{enumerate}
        \item $\Alg^\RP(k)$ is abelian, and $\Alg_u^\RP(k)$ is an abelian subcategory stable under extensions.
        \item For $\tau=\RP,\RPS$, $\Sh_0(k_\tau)$ and $\Sh_a(k_\tau)$ are full subcategories of $\Sh(k_\tau)$ stable under extensions.
        \item The Yoneda functor, and pushforward and pullback of the premorphism $k_\RP\to k_\RPS$ given by identity, give equivalences $\Alg^\RP(k)\cong\Sh_a(k_\RP)\cong\Sh_a(k_\RPS)$ and $\Alg^\RP_u(k)\cong\Sh_0(k_\RP)\cong\Sh_0(k_\RPS)$.
    \end{enumerate}
\end{proposition}

\begin{myproof}[of \Proposition{RPAU characterization}]
    Write $v:k_\RP\to k_\RPS$ the premorphism given by identity. The Yoneda functor $\Alg^\RP(k)\to\Sh(k_\tau)$ is fully faithful by Yoneda's lemma and \Lemma{RPAU equivalence}, and essentially surjective by definition. The rest of the proposition for RPA groups, except Noetherianity and properties of $v^*$, is \cite[Prop. 8.12 and 8.17]{BertapelleSuzuki}. The category $k_\RPS$ has finite products and $v^{-1}$ preserves them, so $v^*$ sends $\Sh_a(k_\RPS)$ to $\Sh_a(k_\RP)$ by \Lemma{yoneda and pullback}. The restriction $v^*:\Sh_a(k_\RPS)\to\Sh_a(k_\RP)$ is left adjoint to $v_*$, hence its essential inverse. Clearly $\Alg_u^\RP(k)$ is closed under kernels, cokernels and extensions in $\Alg^\RP(k)$, and $v^*$ sends $\Sh_0(k_\RPS)$ to $\Sh_0(k_\RP)$, specializing the results to RPAU groups.
\end{myproof}

The following says that RPA $k$-groups are even characterized by only their sections on fields.

\begin{proposition}[yoneda fields]
    Let $k$ be a field of characteristic $p>0$ such that $[k:k^p]<+\infty$. Let $\F_k$ be the category of finite disjoint unions of spectra of relatively perfect field extensions of $k$. Then the Yoneda functor $Y:\Alg^\RP(k)\to\PSh(\F_k)$ is fully faithful.
\end{proposition}

\begin{myproof}[of \Proposition{yoneda fields}]
    We adapt a proof of Suzuki. For $X$ a reduced $k$-scheme with finitely many irreducible components (such as $X=H,H\times_kH$ for $H$ in $\Alg^\RP(k)$ or $\Alg(k)$, by \cite[Rem. 2.2.(2)]{BertapelleSuzuki}), we write $\xi_X$ the disjoint union of its generic points. Consider $H,G\in\Alg^\RP(k)$. 
    
    If $f:H\to G$ is a $k$-group morphism which becomes $0$ in $\PSh(\F_k)$, then for each generic point $\xi$ of $H$ the restriction $f|_\xi:\xi\to G$, as the image of the inclusion $\xi\to H$ under $H(\xi)\to G(\xi)$, is zero. In particular, $f$ is zero on the dense subset $\xi_H$, hence zero everywhere because $\{0\}$ is closed in $G$. Thus $Y$ is faithful.
    
    Consider $\phi:Y(H)\to Y(G)$ in $\PSh(\F_k)$. In particular it induces an element $\phi'$ of the group :
    $$N(H,G)=\ker\left[G(\xi_H)\xrightarrow{p_1^*+p_2^*-m^*}G(\xi_{H\times_kH})\right]$$
    where $p_1,p_2,m:\xi_{H\times_kH}\to\xi_H$ are induced by the projections and group operations. Consider groups $H_0,G_0\in\Alg(k)$ such that $G=G_0^\RP$ and $H=H_0^\RP$. Then $H=\varprojlim_nH_n$, where $H_n=\Res_{F_k}^n(H_0)$ again belongs to $\Alg(k)$ for all $n\geq 0$, and also $H=H_n^\RP$. The transitions $H_{n+1}\to H_n$ are faithfully flat (because $H\to H_{n+1}$ and $H\to H_n$ are, by \cite[Prop. 8.13]{BertapelleSuzuki}), in particular they restrict to the generic points and we have $\xi_H=\varprojlim_n\xi_{H_n}$ and $\xi_{H\times_kH}=\varprojlim_n\xi_{H_n\times_kH_n}$. We have :
    \begin{align*}
        N(H,G)
        &=\ker\left[G_0(\xi_H)\to G_0(\xi_{H\times_kH})\right] &\text{because }\xi_H,\xi_{H\times_kH}\in\RPSch/k,\\
        &=\varinjlim_n\ker\left[G_0(\xi_{H_n})\to G_0(\xi_{H_n\times_kH_n})\right]&\text{because }G_0\text{ is finite presentation},\\
        &=\varinjlim_n\ker\left[(G_0)_{\Adh k}(\xi_{(H_n)_{\Adh k}})\to (G_0)_{\Adh k}(\xi_{(H_n)_{\Adh k}\times_{\Adh k}(H_n)_{\Adh k}})\right]^{\Gal(\Adh k/k)}&\text{by Galois descent}
    \end{align*}
    where $(-)_{\Adh k}$ denotes base change to a separable closure $\Adh k$, and $(-)^{\Gal(\Adh k/k)}$ the subset of Galois-fixed elements. Hence $\phi'\in N(H,G)$ is induced by a Galois-fixed rational $\Adh k$-map $\phi'':(H_n)_{\Adh k}\to (G_0)_{\Adh k}$ for some $n\geq 0$, which is a homomorphism on generic points. By \cite[Lem. 6 of Sec. V.1.5]{Serre88} (though stated for connected groups over algebraically closed $k$, the proof works identically in our case), $\phi''$ is everywhere regular, and descends to a $k$-group morphism $\phi''':H_n\to G_0$. Then $\phi=Y((\phi''')^\RP)$, showing $Y$ is full.
\end{myproof}


\subsection{Dévissage for RPAU groups}
\label{Dévissage for RPAU groups}

Let $k$ be again a field of characteristic $p>0$ such that $[k:k^p]<+\infty$.

\begin{proposition}[devissage RPAU]
    Let $G\in\Alg_u^\RP(k)$.
    \begin{enumerate}
        \item The group $G$ is $p^r$-torsion for some $r\geq 0$, and there exists a finite separated filtration of $G$ in $\Alg_u^\RP(k)$ with $p$-torsion quotients.
        \item If $G$ is $p$-torsion connected, there is an exact sequence $0\to G\to(\GG_{a,k}^\RP)^r\to\GG_{a,k}^\RP\to 0$ with $r\geq 1$.
        \item If $G$ is $p$-torsion, there exists a monomorphism $G\to(\GG_{a,k}^\RP)^r$ with $r\geq 1$.
    \end{enumerate}
\end{proposition}

\begin{myproof}[of \Proposition{devissage RPAU}]
    \textbf{1.} By \Lemma{RPAU equivalence}, $G=\tilde G_0^\RP$ for some unipotent commutative finite presentation $k$-group $\tilde G_0$. Because $\tilde G_0$ is unipotent, we can construct a decreasing sequence of subgroups $\tilde G_i$ of $\tilde G_0$ as follows : for $i\geq 0$, either $\tilde G_i$ is trivial, or it admits a nontrivial map $\phi_i:\tilde G_i\to\GG_{a,k}$ and we set $\tilde G_{i+1}=\ker(\phi_i)$. In particular this sequence decreases strictly, either indefinitely or until it stabilizes at $0$. Because $\tilde G_0$ is Noetherian, we have $\tilde G_r=0$ for some $r\geq 0$. Thus $\tilde G_0$ is $p^r$-torsion as a successive extension of $r$ subgroups of $\GG_{a,k}$, and $G$ is $p^r$-torsion by additivity of $(-)^\RP$.
    
    Then the subgroups $G_i=\im(p^i:G\to G)$ define a decreasing filtration of $G$ in the abelian category $\Alg_u^\RP(k)$, with $p$-torsion quotients, and $G_i=0$ for large enough $i$ because $G$ is finite $p$-primary exponent.
    
    \textbf{2.} If $G=G_0^\RP$ with $G_0\in\Alg_u(k)$, then by \cite[Cor. B.1.13]{CGP} (which holds with no assumption on $k$ because $G_0$ is connected) there is a nonzero morphism $f:\GG_{a,k}^r\to\GG_{a,k}$ with kernel $G_0$. It is surjective because $\GG_{a,k}$ is its only own smooth connected nontrivial split subgroup. We conclude for $G$ by \Proposition{RP is exact}.1.
    
    \textbf{3.} If $G=G_0^\RP$ with $G_0\in\Alg_u(k)$, then $G_0$ is affine so it embeds in some $\GG_{a,k}^r$ by \cite[Lem. B.1.10]{CGP}.
\end{myproof}

\begin{proposition}[RP is exact]
    Consider the relative perfection functor $(-)^\RP$.
    \begin{enumerate}
        \item $(-)^\RP:\Alg(k)\to\Alg^\RP(k)$ is an exact functor, and $(-)^\RP:\Sch/k\to\RPSch/k$ preserves limits.
        \item $(-)^\RP:\Alg(k)\to\Alg^\RP(k)$ is essentially surjective. Any morphism in $\Alg^\RP(k)$ is the image of a morphism in $\Alg(k)$ under $(-)^\RP$.
        \item If $A_0\in\Alg(k)$ and $0\to A_0^\RP\to B\to C\to 0$ is an exact sequence in $\Alg^\RP(k)$, then it is the relative perfection of an exact sequence $0\to A_0\to B_0\to C_0\to 0$ in $\Alg^\RP(k)$.
        \item For $G\in\Alg(k)$ we have canonical isomorphisms $\pi_0(G^\RP)\cong\pi_0(G)$ and $(G^0)^\RP\cong(G^\RP)^0$.
        \item All of the above holds again if we replace $\Alg(k)$ and $\Alg^\RP(k)$ by $\Alg_u(k)$ and $\Alg_u^\RP(k)$.
    \end{enumerate}
\end{proposition}

Exactness of $(-)^\RP$ on the nonabelian category $\Alg(k)$ means if $0\to N\to G\to H\to 0$ is an exact sequence in $\Sh(k_\Et)$, with $N,G,H\in\Alg(k)$, then $0\to N^\RP\to G^\RP\to H^\RP\to 0$ is exact in $\Sh(k_\tau)$ for $\tau=\RP,\RPS$. If we consider all finite presentation commutative groups, then $(-)^\RP$ is still left exact as a right adjoint, but the relative perfection of a non-smooth fppf map may not be surjective \cite[Ex. 8.11]{BertapelleSuzuki}.

\begin{myproof}[of \Proposition{RP is exact}]
    \textbf{1.} This functor preserves limits (and in particular is left exact) as right adjoint to the inclusion. It is right exact because it takes smooth surjections to surjections by \cite[Prop. 8.8]{BertapelleSuzuki}.
    
    \textbf{2.} Essential surjectivity is by definition of RPA $k$-groups, \textit{i.e.} \Proposition{RPAU equivalence}. For the second part, we can more specifically prescribe the target of the morphism in $\Alg(k)$, as well as its source up to Weil restriction. Consider $G,H\in\Alg^\RP(k)$ and $G_0,H_0\in\Alg(k)$ such that $G=G_0^\RP$, $H=H_0^\RP$. Then $G=\varprojlim_nG_n$ where $G_n=\Res_{F_k}^n(G_0)$. By definition of the relative perfection and \cite[Prop. 01ZC]{stacks-project} :
    $$\Hom_{\RPSch/k}(G,H)=\Hom_{\Sch/k}(G,H_0)=\varinjlim_n\Hom_{\Sch/k}(G_n,H_0)$$
    and these identifications restrict to the subsets of morphisms of group objects. By \Lemma{base change RP}.3 (with $S=T=\Spec k$ and $g=F_k^n$), we have $G_n^\RP = G_0^\RP\times_{k,F_k^n}\Spec k$, and applying the relative Frobenius we get $G_n^\RP=G_0^\RP=G$. Thus any morphism $G\to H$ is the relative perfection of a morphism $G_n\to H_0$.
    
    It only remains to show $G_n\in\Alg(k)$, which reduces to showing $\Res_{F_k}(G_0)\in\Alg(k)$ for $G_0\in\Alg(k)$, after an induction on $n$, which holds because Weil restriction preserves smoothness and quasi-compactness.
    
    \textbf{3.} First we show any short exact sequence $(E):0\to A\to B\to C\to 0$ in $\Alg^\RP(k)$ is the image of a short exact sequence $0\to A_1\to B_1\to C_1\to 0$ in $\Alg(k)$ without prescribing $A_1$.
    
    By \textbf{(2.)}, the morphism $f:B\to C$ is the relative perfection of a morphism $f_1: B_1\to \tilde C_1$ in $\Alg(k)$. Its kernel $\tilde A_1$ is a finite type $k$-group but not necessarily smooth. Its maximal smooth subgroup $A_1=(\tilde A_1)^\sharp$ \cite[Def. 8.3]{BertapelleSuzuki} belongs to $\Alg(k)$, and by \cite[Prop. 8.4]{BertapelleSuzuki} and left exactness of relative perfection, we have $A_1^\RP = (\tilde A_1)^\RP=\ker(f)$. Set $C_1=B_1/A_1$. Then $C_1$ belongs to $\Alg(k)$ as quotient of a smooth group \cite[Ex. A.1.12]{CGP}. By \textbf{(1.)} we have an exact sequence $0\to A_1^\RP\to B_1^\RP\to C_1^\RP\to 0$. It follows that $C_1^\RP=\coker(A\to B)=\tilde C_1^\RP$, so $(E)$ is the relative perfection of the exact sequence $0\to A_1\to B_1\to C_1\to 0$.
    
    Now consider some fixed $A_0\in\Alg(k)$ such that $A=A_0^\RP$. Using the previous construction, we have $A=\varprojlim_nA_{1,n}$ where $A_{1,n}=\Res_{F_k}^n(A_1)$, so by \cite[Prop. 01ZC]{stacks-project} the natural map $A\to A_0$ factors as $A\to A_{1,n}\to A_0$ for some $n$. Similarly write $B_{A,n}=\Res_{F_k}^n(B_1)$ and $C_{1,n}=\Res_{F_k}^n(C_1)$ and consider the commutative diagram :
    $$\begin{tikzcd}
       0 \arrow[r]& A_{1,n} \arrow[r]\arrow[d]& B_{1,n} \arrow[r]\arrow[d]& C_{1,n} \arrow[r]\arrow[d]& 0\\
       0 \arrow[r]& A_0 \arrow[r]\arrow[r]& B_0 \arrow[r]& C_0\arrow[r] & 0
    \end{tikzcd}$$
    where $B_0$ is defined as the pushout of $B_{1,n}$ and $A_0$, and $C_0=C_{1,n}$. The Weil restriction of a smooth quasi-compact group is again smooth quasi-compact hence $C_0$ and (as quotient of the smooth quasi-compact group $B_{A,n}\oplus A_0$) $B_0$ belong to $\Alg(k)$. By the same proof as \cite[Prop. 8.8]{BertapelleSuzuki} (applied to the functor $\Res_{F_k}^n:\Alg(k)\to\Alg(k)$ instead of $(-)^\RP:\Alg(k)\to\Alg^\RP(k)$), the first row in the above diagram is exact. By pushout the second row is also exact. Since $A=A_{1,n}^\RP=A_0^\RP$ and $C=C_{1,n}^\RP=C_0^\RP$, by \textbf{(1.)} the map $B=B_{1,n}^\RP\to B_0^\RP$ is an isomorphism. Thus $0\to A_0\to B_0\to C_0\to 0$ has relative perfection $(E)$.
    
    \textbf{4.} By \cite[Prop. 8.13]{BertapelleSuzuki}, the canonical morphism $G^\RP\to G$ induces a bijection between irreducible components, meaning $\pi_0(G^\RP)\cong\pi_0(G)$, thus by left exactness $(G^0)^\RP\cong(G^\RP)^0$.
    
    \textbf{5.} Point \textbf{(1.)} and \textbf{(4.)} clearly specialize to the exact subcategory $\Alg_u^\RP(k)\subseteq\Alg^\RP(k)$.
    
    For \textbf{(2.)} and \textbf{(3.)} the proofs are identical with the additional fact that $\Res_{F_k}(G_0)$ is unipotent for $G_0\in\Alg_u(k)$. If $G_0$ is $p$-torsion, by \Proposition{devissage RPAU}.3 and left exactness of Weil restriction we have a closed immersion of groups $\Res_{F_k}(G_0)\subseteq \Res_{F_k}(\GG_{a,k}^r)$ for some $r\geq 0$, where $\Res_k(\GG_{a,k}^r)\cong\GG_{a,k}^{[k:k^p]r}$ is unipotent by \Example{affine RP}, hence $\Res_{F_k}(G_0)$ is unipotent. In general we have exact sequences for $i\geq 0$ :
    $$0\to \Res_{F_k}(p^{i+1}G_0)\to \Res_{F_k}(p^iG_0)\to\Res_{F_k}(p^iG_0/p^{i+1}G_0).$$
    For large $i\geq 0$, $p^iG_0=0$ is unipotent by \Proposition{devissage RPAU}.1. Also for $i\geq 0$, $\Res_{F_k}(p^iG_0/p^{i+1}G_0)$ is unipotent by the $p$-torsion case, so the subgroup $\im(\Res_{F_k}(p^iG_0)\to\Res_{F_k}(p^iG_0/p^{i+1}G_0))$ is unipotent. Hence $\Res_{F_k}(p^iG_0)$ is unipotent by induction on $i$, as an extension of unipotent groups, so $\Res_{F_k}(G_0)\in\Alg_u(k)$.
\end{myproof}

\begin{definition}[wound and split]
    A group $G\in\Alg_u(k)$ (resp. $G\in\Alg_u^\RP(k)$) is called \emph{(resp. relatively perfect) split} if it admits a finite separated filtration with quotients $0$ or $\GG_{a,k}$ (resp. $0$ or $\GG_{a,k}^\RP$). It is called \emph{(resp. relatively perfect) wound} if all morphisms $\GG_{a,k}\to G$ (resp. $\GG_{a,k}^\RP\to G$) are trivial.
\end{definition}

We also call RPAU $k$-groups \emph{split} or \emph{wound} when it is clear from context that they belong to $\Alg_u^\RP(k)$.

\begin{lemma}[wound and split equivalence]
    Consider $G,H\in\Alg_u^\RP(k)$.
    \begin{enumerate}
    \item If $G$ and $H$ are respectively split and wound, then any morphism $G\to H$ is trivial.
    \item The group $G$ is split (resp. wound) if and only if $G=G_0^\RP$ with $G_0\in\Alg_u(k)$ split (resp. wound).
    \item There exists a unique exact sequence $0\to G_s\to G\to G_w\to 0$ where $G_s$ (resp. $G_w$) is relatively perfect split (resp. wound) unipotent. Any morphism $H\to G$ with $H$ split (resp. $G\to H$ with $H$ wound) factors uniquely through $G_s$ (resp. $G_w$).
    \item If $G$ is split and $G\to H$ is a surjection in $\Alg_u^\RP(k)$ then $H$ is split. If $G$ is wound and $H\to G$ is an injection in $\Alg_u^\RP(k)$ then $H$ is wound.
    \end{enumerate}
\end{lemma}

If $G\in\Alg_u^\RP(k)$ is split then it is connected, and if it is finite étale then it is wound. If $k$ is perfect then the converses are true, since they are true for split and wound objects of $\Alg_u(k)$.

\begin{myproof}[of \Lemma{wound and split equivalence}]
    \textbf{1.} Consider a finite separated filtration $\dots\leq G_1\leq G_0=G$ with quotients $0$ or $\GG_{a,k}^\RP$. The restriction map $\Hom_{\Alg_u^\RP(k)}(G_i,H)\to\Hom_{\Alg_u^\RP(k)}(G_{i+1},H)$ has kernel $\Hom_{\Alg_u^\RP(k)}(G_i/G_{i+1},H)$, trivial because $H$ is wound. Thus $\Hom_{\Alg_u^\RP(k)}(G,H)\to\Hom_{\Alg_u^\RP(k)}(G_i,H)$ is injective, but $G_i=0$ for large $i$, so $\Hom_{\Alg_u^\RP(k)}(G,H)=0$, proving the claim.
    
    \textbf{2.} By \Proposition{RP is exact}.1, the relative perfection of a split group is split. If $G_0\in\Alg_u(k)$ is wound unipotent then  $\Hom(\GG_{a,k}^\RP,G_0^\RP)=\Hom(\GG_{a,k}^\RP,G_0)=\varinjlim_n\Hom(\Res^n_{F_k}(\GG_{a,k}),G_0)$ where all Hom-sets are of $k$-group schemes, by \cite[Prop. 8.10]{BertapelleSuzuki}. As a $k$-group, $\Res^n_{F_k}(\GG_{a,k})$ is isomorphic to $\GG_{a,k}^{[k:k^p]^n}$ by \Example{affine RP}, we conclude that $G_0^\RP$ is relatively perfect wound unipotent.
    
    Now take any $G_0\in\Alg_u(k)$ such that $G=G_0^\RP$, and consider $0\to G_{0,s}\to G_0\to G_{0,w}\to 0$ an exact sequence in $\Alg_u(k)$ with split $G_{0,s}$ and wound $G_{w,0}$ ; such a sequence exists by \cite[Th. B.3.4]{CGP}. By \Proposition{RP is exact}.1, we have an exact sequence in $\Alg_u^\RP(k)$ :
    $$0\to G_{0,s}^\RP\to G\to G_{0,w}^\RP\to 0$$
    where we have just seen that $G_{0,s}^\RP$ (resp. $G_{0,w}^\RP$) was relatively perfect split (resp. wound) unipotent. If $G$ is relatively perfect wound unipotent, then $G_{0,s}^\RP\to G$ is trivial and $G=G_{0,w}^\RP$. If $G$ is relatively perfect split unipotent, then $G\to G_{0,w}^\RP$ is trivial and $G=G_{0,s}^\RP$. This concludes.
    
    \textbf{3.} The existence was seen in the proof of \textbf{(2.)}. For uniqueness it suffices to prove the universal properties, which follow from \textbf{(1.)}.
    
    \textbf{4.} The wound statement is clear from the definition. If $G$ is split then the composite $G\to H\to H_w$ is $0$ by \textbf{(1.)}, so $G\to H$ is surjective and factors as $G\to H_s\to H$. This implies $H_s=H$.
\end{myproof}

\begin{proposition}[RPAU generation]
    Let $\tau=\RP,\RPS$. Then $\Sh_0(k_\tau)$ is the smallest full abelian subcategory of $\Sh(k_\tau)$ closed under extensions which contains $\GG_{a,k}^\RP$.
\end{proposition}

\begin{myproof}[of \Proposition{RPAU generation}]
    By \Proposition{RPAU characterization}, $\Sh_0(k_\tau)$ is a full abelian subcategory of $\Sh(k_\tau)$ closed under extensions, and it clearly contains $\GG_{a,k}^\RP$. Consider $\A\subseteq\Sh(k_\tau)$ a full abelian subcategory closed under extensions containing $\GG_{a,k}^\RP$. By \Proposition{devissage RPAU}.2, it contains all $p$-torsion connected objects as well. If $G\in\Sh_0(k_\tau)$ is $p$-torsion, it fits in an exact sequence $0\to G\to(\GG_{a,k}^\RP)^r\to H\to 0$ for some $H\in\Alg_u^\RP(k)$ by \Proposition{devissage RPAU}.3, and $H$ is automatically $p$-torsion and connected, hence in $\A$, as the image of the $p$-torsion connected group $(\GG_{a,k}^\RP)^r$. By extension, $\A$ contains all $p$-torsion objects of $\Sh_0(k_\tau)$. In turn it contains all of $\Sh_0(k_\tau)$ by \Proposition{devissage RPAU}.1, which concludes.
\end{myproof}


\subsection{Derived category of RPAU groups}
\label{Derived category of RPAU groups}

The goal of this section is to prove a derived version of \Proposition{RPAU characterization}. Its content leads up to \Proposition{fully faithful Yon RPAU} and finally \Proposition{definition of D0}. To prove \Proposition{fully faithful Yon RPAU}, we reproduce a proof communicated by Suzuki. This proof takes inspiration from \cite[Prop. 8 of Sec 8.8]{Serre1961} but uses $\Ind\Alg_u^\RP(k)$ instead of $\Pro\Alg_u^\RP(k)$.

In this section let $\tau=\RP,\RPS$ be fixed. We write $\Alg_u^\RP(k,\ZZ/p)$ the full subcategory of $p$-torsion objects of $\Alg_u^\RP(k)$. From \Proposition{RPAU characterization}, it follows that $\Alg_u^\RP(k,\ZZ/p)$ is an abelian subcategory of $\Alg_u^\RP(k)$, equivalent to the full subcategory $\Sh_0(k_\tau,\ZZ/p)$ of $p$-torsion objects of $\Sh_0(k_\tau)$.

\begin{lemma}[injective Ga]
    The additive group $\GG_{a,k}^\RP$ is injective in $\Alg_u^\RP(k,\ZZ/p)$.
\end{lemma}

\begin{myproof}[of \Lemma{injective Ga}]
    Consider $f:\GG_{a,k}^\RP\to G$ a monomorphism in $\Alg_u^\RP(k,\ZZ/p)$. By \Proposition{RP is exact}.3, this is the image of a closed immersion $f_0:\GG_{a,k}\to G_0$ with $G_0\in\Alg_u(k)$. By \cite[Cor. B.1.12]{CGP} this $f_0$ splits, so $f$ is a split monomorphism.
\end{myproof}

\begin{proposition}[fully faithful p]
    The Yoneda functor $Y:\Alg_u^\RP(k,\ZZ/p)\to\Sh(k_\RP,\ZZ/p)$ is exact and its derived functor $RY:D^b(\Alg_u^\RP(k,\ZZ/p))\to D^b(\Sh(k_\RP,\ZZ/p))$ is fully faithful.
\end{proposition}

\begin{myproof}[of \Proposition{fully faithful p}]
    The exactness of $Y$ is by \Proposition{RPAU characterization}. We need to show equalities : $$\Ext^q_{\Alg_u^\RP(k,\ZZ/p)}(G,H)=\Ext^q_{k_\tau,\ZZ/p}(Y(G),Y(H))$$
    for $q\geq 0$ and $G,H\in\Alg_u^\RP(k,\ZZ/p)$. By \Proposition{RPAU generation}, we can assume $G=H=\GG_{a,k}^\RP$. For $q=0$ this is \Proposition{RPAU characterization}. Then for $q\geq 1$, the right-hand side vanishes by \cite[Prop. 2.1]{Kato86}, and the left-hand side vanishes by \Lemma{injective Ga}, which concludes.
\end{myproof}

\begin{proposition}[p-divisible injective]
    Any injective object of $\IndAlg_u^\RP(k)$ is $p$-divisible.
\end{proposition}

\begin{myproof}[of \Proposition{p-divisible injective}]
    For $n\geq 0$ and $R$ a $k$-algebra, we consider the Cohen ring $C_n(R)=\J_n^\B(R)$ of \cite[Def. 7.4]{BertapelleSuzuki} and \cite[Sec. 2.1]{Schoeller}, which depends on the choice of a fixed $p$-basis $\B$ of $k$. Here $C_n(R)$ is trivial for $n=0$ by convention. As a presheaf on affine $k$-schemes, $C_n$ is representable by an affine $k$-scheme, in particular it has the structure of a commutative $k$-group scheme \cite[Prop. 7.6]{BertapelleSuzuki}. By \cite[Prop. 2.9]{Schoeller} we have exact sequences (of presheaves of groups on affine $k$-schemes) :
    $$0\to C_n\to C_{n+1}\to\Res_{F_k^n}(\GG_{a,k})\to 0$$
    so $C_n$ is unipotent and finite presentation by induction on $n\geq 0$. Hence $C_n^\RP\in\Alg_u^\RP(k)$.
    
    For $0\leq m\leq n$ we have a closed immersion $j_{m,n}:C_m\to C_n$ and for $R$ a $k$-algebra, an identification $C_m(R)=p^{n-m}C_n(R)$ of subgroups of $C_n(R)$ \cite[Prop. 7.8]{BertapelleSuzuki}. In particular, as $R$ runs through relatively perfect $k$-algebras, it follows that $C_m^\RP=p^{n-m}C_n^\RP$. Thus $C_\infty^\RP=\varinjlim_nC_n^\RP$ is $p$-divisible in $\IndAlg_u^\RP(k)$.
    
    With these preliminaries in mind, let $I\in\IndAlg_u^\RP(k)$ be injective. Consider $G\in\Alg_u^\RP(k)$.
    
    Write $G=G_0^\RP$ with $G_0\in\Alg_u(k)$. By \cite[Prop. 4.7]{Schoeller}, there exists a closed immersion $G_0\to C_n^r$ for some $n,r\geq 0$. By left exactness of relative perfection, we have a monomorphism $G\to (C_n^\RP)^r$ in $\Alg_u^\RP(k)$, hence a monomorphism $G\to(C_\infty^\RP)^r$ in $\IndAlg_u^\RP(k)$. Then any morphism $G\to I$ factors as $G\to (C_\infty^\RP)^r\to I$ by injectivity. In particular, $G\to I$ factors as $G\to pI\to I$ by $p$-divisibility of $(C_\infty^\RP)^r$. It follows that the natural map :
    $$\Hom_{\IndAlg_u^\RP(k)}(G,pI)\to\Hom_{\Ind\Alg_u^\RP(k)}(G,I)$$
    is surjective for $G\in\Alg_u^\RP(k)$. Since injectivity is clear, and Hom functors preserve limits, it follows that the inclusion $pI\to I$ induces an isomorphism for any $H\in\IndAlg_u^\RP(k)$
    $$\Hom_{\IndAlg_u^\RP(k)}(H,pI)\to\Hom_{\Ind\Alg_u^\RP(k)}(H,I).$$
    Thus $pI=I$ by Yoneda's lemma, which concludes.
\end{myproof}

To work with $\ProAlg_u^\RP(k)$ instead we would need every projective objects in $\ProAlg_u^\RP(k)$ to be $p$-torsion free (as in as in \cite[Prop. 8 of Sec 8.8]{Serre1961} for the case $k$ perfect). With a dual approach, we would want every $G\in\Alg_u^\RP(k)$ to be a quotient of some $Q_n^r$, where $Q_n\in\Alg_u^\RP(k)$ fits in a codirected sequence $\{Q_n\}_{n\geq 0}$ with surjective transitions such that $\ker(Q_n\to Q_m) = Q_n[p^{n-m}]$ for $0\leq m\leq n$.

\begin{lemma}[p-torsion properties]
    Let $\T:\IndAlg_u^\RP(k)\to\IndAlg_u^\RP(k,\ZZ/p)$ (resp. $\Sh(k_\tau)\to\Sh(k_\tau,\ZZ/p)$) be the $p$-torsion functor.
    \begin{enumerate}
        \item Both functors $\T$ are left exact and we have commutative diagrams of Yoneda functors :
        $$\begin{tikzcd}
            \IndAlg_u^\RP(k)\arrow[r]\arrow[d,"\T"]& \Sh(k_\tau)\arrow[d,"\T"]\\
            \IndAlg_u^\RP(k,\ZZ/p)\arrow[r]& \Sh(k_\tau,\ZZ/p)
        \end{tikzcd}
        \begin{tikzcd}
            D^+(\IndAlg_u^\RP(k))\arrow[r]\arrow[d,"R\T"]& D^+(k_\tau)\arrow[d,"R\T"]\\
            D^+(\IndAlg_u^\RP(k,\ZZ/p))\arrow[r]& D^+(k_\tau,\ZZ/p)
        \end{tikzcd}$$
        \item For $G\in\IndAlg_u(k)$ (resp. $G\in\Sh(k_\tau)$), we have $R^1\T(G)=G/pG$ and $R^q\T(G)=0$ for $q\geq 2$.
        \item For $G\in\IndAlg_u(k,\ZZ/p)$ (resp. $G\in\Sh(k_\tau,\ZZ/p)$) we have :
        \begin{align*}
            R\Hom_{\IndAlg_u(k,\ZZ/p)}(G,R\T(-)) &= R\Hom_{\IndAlg_u(k)}(G,-)\\ \text{(resp. }R\Hom_{k_\tau,\ZZ/p}(G,R\T(-)) &= R\Hom_{k_\tau}(G,-)\text{)}.
            \end{align*}
        \item The functor $R\T:D^+(\IndAlg_u(k))\to D^+(\IndAlg_u(k,\ZZ/p))$ restricts to a functor between full subcategories $D^b(\Alg_u(k))\to D^b(\Alg_u(k,\ZZ/p))$.
    \end{enumerate}
\end{lemma}

In \textbf{(4.)} above, recall that by \cite[Th. 15.3.1.(i)]{KashiwaraSchapira} the inclusion $\A\to\Ind\A$, for $\A$ an abelian category, identifies $D^b(\A)$ with the full subcategory of $D^b(\Ind\A)$ of objects with cohomology in $\A$.

\begin{myproof}[of \Lemma{p-torsion properties}]
    \textbf{1. Underived case.} Left exactness is obvious. The functor $\Alg_u^\RP(k)\to\Sh(k_\tau)$ is exact by \Proposition{RPAU characterization}, thus so is its restriction to abelian subcategories $\Alg_u^\RP(k,\ZZ/p)\to\Sh(k_\tau,\ZZ/p)$. By ind-completion and the exactness of filtered colimits in $\Sh(k_\tau)$, the Yoneda functors $\IndAlg_u^\RP(k)\to\Sh(k_\tau)$ and $\IndAlg_u^\RP(k,\ZZ/p)\to\Sh(k_\tau,\ZZ/p)$ are exact, giving the desired commutation.
    
    \textbf{2. Sheaf case.} Consider the inclusion $J:\Sh(k_\tau,\ZZ/p)\to\Sh(k_\tau)$. Since $\HOM_{k_\tau}(\ZZ,-)=\id$, the exact sequence $0\to\ZZ\to\ZZ\to\ZZ/p\to 0$ gives identifications for $G\in\Sh(k_\tau)$ and $n\geq 2$ :
    $$\HOM_{k_\tau}(\ZZ/p,G)=J\circ\T(G),\qquad \EXT^1_{k_\tau}(\ZZ/p,G)=G/pG,\qquad \EXT^n_{k_\tau}(\ZZ/p,G)=0$$
    Because $J$ is exact, we have $J\circ R^n\T(G)=\EXT^n_{k_\tau}(\ZZ/p,G)$ for $n\geq 1$. We conclude by full faithfulness of $J$.
    
    \textbf{1. Derived case.} The Yoneda functor is exact, so it suffices to show it sends injectives of $\IndAlg_u^\RP(k)$ to $\T$-acyclics in $\Sh(k_\tau)$. This is due to \Proposition{p-divisible injective} and \textbf{(2.)} for $\T:\Sh(k_\tau)\to\Sh(k_\tau,\ZZ/p)$.
    
    \textbf{2. Group case.} For $\T:\IndAlg_u^\RP(k)\to\IndAlg_u^\RP(k,\ZZ/p)$ this now follows from the diagram of derived categories in \textbf{(1.)} and the case of $\T:\Sh(k_\tau)\to\Sh(k_\tau,\ZZ/p)$.
    
    \textbf{3.} The functor $\T:\IndAlg_u^\RP(k)\to\IndAlg_u^\RP(k,\ZZ/p)$ is right adjoint to the inclusion $J:\IndAlg_u^\RP(k,\ZZ/p)\to\IndAlg_u^\RP(k)$, which is exact (as the ind-completion of an exact functor). In particular $\T$ preserves injectives. The equality now follows by derivation of composite functors and the identifications :
    $$\Hom_{\IndAlg_u^\RP(k)}(G,-)=\Hom_{\IndAlg_u^\RP(k)}(G,\T(-))=\Hom_{\IndAlg_u^\RP(k,\ZZ/p)}(G,\T(-))$$
    for $G\in\IndAlg_u^\RP(k,\ZZ/p)$. The argument in $\Sh(k_\tau)$ is identical.
    
    \textbf{4.} For $G\in D^b(\Alg_u^\RP(k))$, by \textbf{(2.)}, $R\T$ sends $D^b(\Alg_u^\RP(k))$ to $D^b(\IndAlg_u^\RP(k,\ZZ/p))$ and we have a convergent spectral sequence :
    $$E^{ij}_2=R^i\T(H^j(G))\Rightarrow R^{i+j}\T(G)=H^{i+j}(R\T(G)).$$
    By \textbf{(2.)} again $R^i\T(H^j(G))\in\Alg_u^\RP(k,\ZZ/p)$ for $i,j\in\ZZ$. Since $\Alg_u^\RP(k,\ZZ/p)$ is an abelian full subcategory of $\IndAlg_u^\RP(k,\ZZ/p)$, it follows that $R\T(G)$ belongs to $D^b(\Alg_u^\RP(k,\ZZ/p))$.
\end{myproof}

\begin{proposition}[fully faithful Yon RPAU]
    The exact, fully faithful Yoneda functor $Y:\Alg_u^\RP(k)\to\Sh(k_\RP)$ induces an exact, fully faithful functor of triangulated categories $RY:D^b(\Alg_u^\RP(k))\to D(k_\RP)$.
\end{proposition}

\begin{myproof}[of \Proposition{fully faithful Yon RPAU}]
    The functor $D^b(\Alg_u^\RP(k))\to D^b(\IndAlg_u^\RP(k))$ is fully faithful and the composite $\varinjlim\Ind Y:\IndAlg_u^\RP(k)\to\Ind\Sh(k_\RP)\to\Sh(k_\RP)$ (still written $Y$) is exact by the proof of \Lemma{p-torsion properties}.1. Hence we want to show : $$R\Hom_{\IndAlg_u^\RP(k)}(A,B)=R\Hom_{k_\RP}(Y(A),Y(B))$$
    for $A,B\in\Alg_u^\RP(k)$. By \Lemma{devissage RPAU}.1 we can assume $A\in\Alg_u^\RP(k,\ZZ/p)$. Then we have :
    \begin{align*}
        R\Hom_{\IndAlg_u^\RP(k)}(A,B)
        &=R\Hom_{\IndAlg_u^\RP(k,\ZZ/p)}(A,R\T(B)) &\text{(\Lemma{p-torsion properties}.3)}\\
        &=R\Hom_{\Sh(k_\RP,\ZZ/p)}(Y(A),RY(R\T(B))) &\text{(\Proposition{fully faithful p})}\\
        &= R\Hom_{\Sh(k_\RP,\ZZ/p)}(Y(A),R\T(Y(B))) &\text{(\Lemma{p-torsion properties}.1 and exactness of }Y\text{)}\\
        &= R\Hom_{\Sh(k_\RP)}(Y(A),Y(B)) &\text{(\Lemma{p-torsion properties}.3)}
    \end{align*}
    which shows $RY:D^b(\Alg_u^\RP(k))\to D(k_\RP)$ is fully faithful.
\end{myproof}

\begin{proposition}[definition of D0]
    Let $\tau=\RP,\RPS$. Then the following full subcategories of $D(k_\tau)$ coincide :
    \begin{enumerate}
        \item the essential image of $D^b(\Alg_u^\RP(k))$ by the derived Yoneda functor ;
        \item the smallest full triangulated subcategory containing $\GG_{a,k}^\RP$ (as a complex concentrated in degree $0$) ;
        \item the full subcategory of bounded objects with cohomology in $\Sh_0(k_\tau)$.
    \end{enumerate}
    This subcategory is written $D_0(k_\tau)$. The Yoneda functor $D^b(\Alg_u^\RP(k))\to D(k_\tau)$, and the pushforward of the premorphism of sites $v:k_\RP\to k_\RPS$ given by identity induce exact equivalences of triangulated categories $D^b(\Alg_u^\RP)\cong D_0(k_\RP)\cong D_0(k_\RPS)$.
\end{proposition}

\begin{myproof}[of \Proposition{definition of D0}]
    Because $\Sh_0(k_\tau)$ is closed under extensions in $\Sh(k_\tau)$, \textbf{(3.)} is a full triangulated subcategory of $D(k_\tau)$. It contains $\GG_{a,k}^\RP$, so \textbf{(3.)} contains \textbf{(2.)}. To prove conversely that \textbf{(2.)} contains \textbf{(3.)}, it suffices to show every object of $\Sh_0(k_\tau)$ (in degree $0$) is contained in \textbf{(2.)}, which follows from (the same proof as) \Proposition{RPAU generation}. Write $D_0(k_\tau)$ for the equivalent definitions \textbf{(2.)} and \textbf{(3.)}.
    
    The functor $v_*$ is exact and defines an equivalence $v_*:D_0(k_\RP)\xrightarrow\sim D_0(k_\RPS)$ by \cite[Prop. 2.2]{KatoSuzuki2019} : though this is stated for the $p^n$-torsion versions $D_0(-,\Lambda_n)$ with $\Lambda_n=\ZZ/p^n$, the proof for $D_0(-)$ is identical upon replacing $\Lambda_n$ with $\ZZ$ and MacLane's resolution with the Deligne-Scholze resolution of \cite[Th. 4.5]{ScholzeClausen}.
    
    Clearly the composite $\Alg_u^\RP(k)\to \Sh(k_\RP)\to\Sh(k_\RPS)$ is the Yoneda functor, and since both factors are exact, the composite $D^b(\Alg_u^\RP(k))\to D_0(k_\RP)\to D_0(k_\RPS)$ is also the derived Yoneda functor. Hence it only remains to show that $D^b(\Alg_u^\RP(k))\to D(k_\RP)$ defines an equivalence $D^b(\Alg_u^\RP(k))\cong D_0(k_\RP)$.
    
    The Yoneda functor $Y:\Alg_u^\RP(k)\to\Sh(k_\RP)$ is exact and fully faithful by \Proposition{RPAU characterization}.2. By \cite[Eq. 06UR and adjacent discussion]{stacks-project} its derived functor $RY:D^b(\Alg_u^\RP(k))\to D(k_\RP)$ is an exact functor of triangulated categories, and is fully faithful by \Proposition{fully faithful Yon RPAU}. For $G\in D^b(\Alg_u^\RP(k))$, $H^q(RY(G))=Y(H^q(G))\in\Sh_0(k_\RP)$ for all $q\in\ZZ$. Hence the essential image of $RY$ is contained in $D_0(k_\RP)$.
    
    Conversely the essential image of $RY$ is a full triangulated subcategory of $D(k_\RP)$. Indeed given an exact triangle $Y(A)\to B\to Y(C)\xrightarrow{\delta} Y(A)[1]$ in $D(k_\RP)$ with $A,C\in D^b(\Alg_u^\RP(k))$, then $B$ is determined by the morphism $\delta\in\Hom_{D(k_\RP)}(Y(C),Y(A)[1])$ ; the fullness of $RY$ ensures $\delta=Y(\delta')$ for some $\delta':C\to A[1]$, and the exactness of $RY$ implies that $B=RY(B')$ where $B'$ is the fiber of $\delta'$ in $D^b(\Alg_u^\RP(k))$. The essential image of $RY$ is a full triangulated subcategory containing $\GG_{a,k}^\RP$, so it contains $D_0(k_\RP)$. This shows $RY:D^b(\Alg_u^\RP(k))\to D_0(k_\RP)$ is an equivalence.
\end{myproof}


\subsection{Duality for RPAU groups}
\label{Duality for RPAU groups}

Let $k$ be again a field of characteristic $p>0$ such that $[k:k^p]<+\infty$ and let $d=\log_p[k:k^p]$. Recall that $d$ is an integer, the cardinal of any $p$-basis of $k$.

Recall the degree $r\geq 0$, level $n\geq 1$, logarithmic part of the de Rham-Witt complex of \cite[Sec. 4]{Kato86}, $\nu_n(r)(R)$ for $R$ an $\FF_p$-algebra which admits a finite strong $p$-basis (\textit{e.g.} $\O(X)$ for $X\in\RPSch/k$, by \Lemma{p-basis comparison}.2). For $r<0$, we define $\nu_n(r)(R)=0$ by convention. The formation of $\nu_n(r)(R)$ is functorial in $R$ and satisfies étale descent, thus we can view it as an object $\nu_n(r)_k$ (or simply $\nu_n(r)$) of $\Sh(k_\tau)$ for $\tau=\et,\Et,\RP,\RPS$. We write $\ZZ/p^n(r)=\nu_n(r)[-r]\in D(k_\tau)$ and :
$$\nu_\infty(r)=\varinjlim_n\nu_n(r),\qquad \QQ_p/\ZZ_p(r)=\varinjlim_n\ZZ/p^n(r)=\nu_\infty(r)[-r].$$
This is intended to unify notations with the Tate twists $\ZZ/m(r)=\mu_m^{\otimes r}$ for $m$ invertible on the base. We actually define $\ZZ/m(r)=\mu_{m'}^{\otimes r}\oplus\nu_n(r)[-r]$ when $m=p^nm'$ with $m'$ prime to $p$. This notation is also partially motivated by \cite{Geisser2005}. For our purposes, where only objects of $\RPSch/k$ are considered, these conventions are sufficient. Better candidates for $\ZZ/p^n(r)_X$ are studied in \textit{e.g.} \cite{Kato87, Ren} for $X$ separated of finite type over a perfect field possibly not regular.

\begin{lemma}[log dRW is RPAU]
    For $r\in\ZZ$ and $\tau=\RP,\RPS$, $\nu_n(r)$ defines a relatively perfect wound object of $\Sh_0(k_\tau)$.
\end{lemma}

\begin{myproof}[of \Lemma{log dRW is RPAU}]
    Any $p$-basis of $k$ of order $d$ gives a strong $p$-basis of $\O(X)$ for any $X\in\RPSch/\FF_p$, by \Lemma{p-basis comparison}.2, which gives a differential basis of that ring \cite[Sec. (38.A)]{Matsumura80}, that is $\Ohm^1\cong(\GG_{a,k}^\RP)^d$. Hence $\Ohm^r\cong(\GG_{a,k}^\RP)^{\binom dr}$ belongs to $\Alg_u^\RP(k)$, and so does $\Ohm^r/d\Ohm^{r-1}$ as the cokernel of a map in $\Alg_u^\RP(k)$. Using \Proposition{RPAU characterization}.2 and the following exact sequences \cite[Eq. (3.1.5) and (4.1.8)]{Kato86} :
    $$0\to\nu_{n-1}(r)\to\nu_n(r)\to\nu_1(r)\to 0,\qquad 0\to \nu_1(r)\to \Ohm^r\xrightarrow{C^{-1}-1}\Ohm^r/d\Ohm^{r-1}\to 0$$
    we can inductively reduce to $n=1$ by \Lemma{RPAU characterization}.1 and conclude that $\nu_n(r)\in\Alg_u^\RP(k)$. That it is wound, \textit{i.e.} $\Hom_{\Alg_u^\RP(k)}(\GG_{a,k}^\RP,\nu_n(r))=0$, follows from \cite[Th. 3.2.(ii)]{Kato86}.
\end{myproof}

\begin{lemma}[change n if you like]
    Consider $\tau=\RP,\RPS$, $r\in\ZZ$ and $G\in D^b(k_\tau,\ZZ/p^n)$ for some $n\geq 1$. Consider $J_n:\Sh(k_\tau,\ZZ/p^n)\to\Sh(k_\tau)$ the inclusion. Then we have a canonical isomorphism in $D(k_\tau)$ :
    $$R\HOM_{k_\tau}(J_n(G),\nu_\infty(r)) = J_n(R\HOM_{k_\tau,\ZZ/p^n}(G,\nu_n(r))).$$
\end{lemma}

\begin{myproof}[of \Lemma{change n if you like}]
    We can assume $G\in\Sh(k_\tau,\ZZ/p^n)$. Consider $\T_n:\Sh(k_\tau)\to\Sh(k_\tau,\ZZ/p^n)$ the $p^n$-torsion functor. It is right adjoint to $J_n$, which is exact, hence $\T_n$ is left exact and preserves injectives. By composition of derived functors it follows that :
    $$R\HOM_{k_\tau}(J_n(G),-)=J_n(R\HOM_{k_\tau,\ZZ/p^n}(G,R\T_n(-)).$$
    We now want to show $R\T_n(\nu_\infty(r))=\nu_n(r)$. As in the proof of \Lemma{p-torsion properties}.2, using the full faithfulness and exactness of $J_n$ and the exact sequence $0\to\ZZ\to\ZZ\to\ZZ/p^n\to 0$, we get :
    $$R^1\T_n(H)=H/p^nH,\qquad R^q\T_n(H)=0$$
    for $q\geq 2$ and $H\in\Sh(k_\tau)$. For $s,t\geq 0$ we have exact sequences :
    $$0\to\nu_s(r)\xrightarrow{j_{s,t}}\nu_{s+t}(r)\xrightarrow{q_{s,t}}\nu_t(r)\to 0$$
    where $j_{t,s}\circ q_{s,t}=p^n\id_{\nu_{s+t}(r)}$. In particular we have identifications $\nu_s(r) = \nu_{s+t}(r)[p^s] = p^t\nu_{s+t}(r)$. Taking the colimit with respect to $t$ we get $\nu_s(r)=\nu_\infty(r)[p^s]$ hence $\T_n(\nu_\infty(r))=\nu_n(r)$. Taking the colimit with respect to $t$ we get $\nu_\infty(r)=p^t\nu_\infty(r)$ hence $R^1\T_n(\nu_\infty(r))=0$. Thus we have a quasi-isomorphism
    $\nu_n(r)=RT_n(\nu_\infty(r))$ in $D(k_\tau,\ZZ/p^n)$, as desired.
\end{myproof}

\begin{theorem}[duality for RPAU groups]
    Let $k$ be a field of characteristic $p>0$ such that $[k:k^p]=p^d<+\infty$. Let $\tau=\RP,\RPS$. For $G\in D(k_\tau)$, define $G^\vee = R\HOM_{k_\tau}(G,\QQ_p/\ZZ_p(d))$.
    \begin{enumerate}
        \item For $G\in D_0(k_\tau)$, we have $G^\vee\in D_0(k_\tau)$ and the canonical morphism $G\to(G^\vee)^\vee$ is an isomorphism.
        \item For $G\in\Sh_0(k_\tau)$ and $q\in\ZZ$, we have $H^q(G^\vee)=0$ for $q\neq d,d+1$, $H^d(G^\vee)=H^d((G_w)^\vee)$ is relatively perfect wound, and $H^{d+1}(G^\vee)=H^{d+1}((G_s)^\vee)$ is relatively perfect split.
        \item If $G\in D_0(k_\tau)$ then for $q\in\ZZ$ we have isomorphisms in $\Sh(k_\tau)$ :
        $$H^q(G)_s = \EXT^1(H^{d+1-q}(G^\vee)_s,\nu_\infty(d)),\qquad H^q(G)_w=\HOM(H^{d-q}(G^\vee)_w,\nu_\infty(d)).$$
        We say that $H^q(G)_s$ and $H^{d+1-q}(G^\vee)_s$ (resp. $H^q(G)_w$ and $H^{d-q}(G^\vee)_w$) are perfect Serre duals as split (resp. wound) RPAU $k$-groups.
        \item For $G_0$ a vector group, $(G_0^\RP)^\vee\cong\HOM_{k_\tau,\O_k}(G_0^\RP,\Ohm^d)[-d-1]$ (where $\HOM_{k_\tau,\O_k}$ is the inner-Hom of modules over the structure sheaf $X\mapsto \O(X)$ of $\Spec k_\tau$). For $r,n\geq 0$ we have $\ZZ/p^n(r)^\vee\cong\ZZ/p^n(d-r)$.
    \end{enumerate}
\end{theorem}

If $(t_1,\dots,t_d)$ is a $p$-basis of $k$ then by \Lemma{p-basis comparison}.2 we have $\Ohm^d\cong\GG_{a,k}^\RP\cdot dt_1\wedge\dots\wedge dt_d$ in $\Sh(k_\tau)$, hence Serre duality for split groups is a generalization of linear duality for vector groups. If $k$ is perfect then $\nu_\infty(d)=\QQ_p/\ZZ_p$, relatively perfect wound groups are the same as finite étale $p$-primary groups, and Serre duality for them coincides with Pontryagin duality.

\begin{myproof}[of \Theorem{duality for RPAU groups}]
    The case $\tau=\RPS$ is equivalent to the case $\tau=\RP$ via \Proposition{RPAU characterization} and \Proposition{definition of D0}. By \Lemma{change n if you like} and \Proposition{devissage RPAU}.1, it suffices to prove the theorem with $D(k_\RP)$ replaced by $D(k_\RP,\ZZ/p^n)$, $D_0(k_\RP)$ replaced by the subcategory of $p^n$-torsion objects $D_0(k_\RP,\ZZ/p^n)$, $(-)^\vee$ replaced by $R\HOM_{\ZZ/p^n}(-,\nu_n(r))$ and correspondingly $\EXT^q(-,\nu_\infty(r))$ by $\EXT^q_{\ZZ/p^n}(-,\nu_\infty(r))$.
    
    Then the isomorphism $G\to (G^\vee)^\vee$ for $G\in D_0(k_\tau,\ZZ/p^n)$ and \textbf{(4.)} come from \cite[Th. 3.2 and 4.3]{Kato86}, with our convention $\ZZ/p^n(r)=\nu_n(r)[-r]$. Points \textbf{(2.)} and \textbf{(3.)} come from \cite[Prop. 3.3 and 3.4]{SuzukiCFT}, and it follows that $G^\vee\in D_0(k_\tau,\ZZ/^n)$ for $G\in D_0(k_\tau,\ZZ/p^n)$.
\end{myproof}

\begin{remark}[versions of RPAU duality]
    We can instead take $(-)^\vee = R\HOM(-,\nu_\infty(d))$ (or $R\HOM_{\ZZ/p^n}(-,\nu_n(d))$ for the $p^n$-torsion version), as in \cite{Kato86,SuzukiCFT}. Then \Theorem{duality for RPAU groups}.1 holds as-is, while the other points are true after obvious shifts. For instance, \Theorem{duality for RPAU groups}.4 becomes $\nu_n(r)^\vee\cong\nu_n(d-r)$.
\end{remark}


\subsection{Change of field and trace maps}
\label{Change of field and trace maps}

As before let $k$ be a field of characteristic $p>0$ such that $[k:k^p]=p^d<+\infty$. If $l$ is a finite field extension of $k$, we again have $[l:l^p]=[k:k^p]=p^d$ (indeed $[l:k^p]=[l:l^p][l^p:k^p]=[l:k][k:k^p]$ and $[l^p:k^p]=[k:k^p]$ by injectivity of the Frobenius).

For $R$ an $\FF_p$-algebra with Zariski-local strong $p$-bases, we write $\Ohm^q_R$, $W_n\Ohm^q_R$, $\nu_n(q)_R$ the sheaves of Kähler differentials, de Rham-Witt differentials, logarithmic de Rham-Witt differentials, on $\Spec R_\RP$ (or $\Spec R_\RPS$ if defined). We write $\Ohm^q(R)$, $W_n\Ohm^q(R)$, $\nu_n(q)(R)$ their respective groups of global sections.

\begin{proposition}[pushforward exactness]
    Consider $\tau=\RP,\RPS$ and $\pi:\Spec l\to\Spec k$ a finite field extension.
    \begin{enumerate}
        \item Base change by $\pi$ defines a morphism of sites $\Spec l_\tau\to \Spec k_\tau$, again written $\pi$.
        \item The functor $\pi_*:\Sh(l_\tau)\to\Sh(k_\tau)$ is both left and right adjoint to $\pi^*$. In particular, $\pi_*$ is exact.
        \item If $l/k$ is separable then $\pi^*$ is given by restricting sheaves to the subcategory $\Spec l_\tau\subseteq\Spec  k_\tau$.
        \item If $l/k$ is purely inseparable then $\pi_*$ and $\pi^*$ define an equivalence $\Sh(k_\tau)\cong\Sh(l_\tau)$. Furthermore there exists a purely inseparable finite field extension $\rho:\Spec k\to\Spec l$ such that $\pi^*=\rho_*$.
    \end{enumerate}
\end{proposition}

\begin{myproof}[of \Proposition{pushforward exactness}]
    The well-definition of $\pi^{-1}:\Spec k_\tau\to\Spec l_\tau$ is by \Lemma{base change RP}.4, and this defines a premorphism of sites by \Lemma{base change RP}.1. By \cite[Lem. 030K]{stacks-project}, it suffices to prove the other statements in the case $l/k$ is separable or purely inseparable.
    
    \textbf{Separable case.} Because $\pi$ is an étale map, $\pi^\set$ is just given by restriction to the subcategory $\Spec l_\tau\subseteq\Spec k_\tau$. In particular it is exact, hence $\pi:\Spec l_\tau\to\Spec k_\tau$ is a morphism of sites and $\pi^*$ is described similarly. The proof that $\pi_*$ is left adjoint to $\pi^*$ is identical to \cite[Lem. 1.12 of Ch. V]{Milne80}.
    
    \textbf{Purely inseparable case.} For some $n\geq 0$ (say, $n=\log_p[l:k]$) we have $l^{p^n}\subseteq k$ hence the $n$-fold Frobenius on $\Spec l$ factors as a morphism $\rho:\Spec k\to \Spec l$. By \Lemma{base change RP}.2-3, $\pi^{-1}:\Spec k_\tau\to\Spec l_\tau$ is an equivalence with inverse $\rho^{-1}:\Spec l_\tau\to\Spec k_\tau$. Thus $\pi_*$ and $\rho_*$ define an equivalence $\Sh^\set(l_\tau)\cong\Sh^\set(k_\tau)$, so $\pi^\set=\rho_*$ is exact as an equivalence functor, and $\pi$ is a morphism of sites. Similarly, $\pi_*$ and $\rho_*$ define an equivalence $\Sh(l_\tau)\cong\Sh(k_\tau)$, so $\pi^*$ is an essential inverse to $\pi_*$ thus right adjoint to it.
\end{myproof}

\begin{proposition}[weil restriction RPAU]
    Consider $\tau=\RP,\RPS$ and $\pi:\Spec l\to\Spec k$ a finite field extension. Then $\pi_*:\Sh(l_\tau)\to\Sh(k_\tau)$ and $\pi^*:\Sh(k_\tau)\to\Sh(l_\tau)$ restrict to functors :
    $$\pi_*:\Sh_0(l_\tau)\to\Sh_0(k_\tau),\qquad \pi^*:\Sh_0(k_\tau)\to\Sh_0(l_\tau).$$
    If $l/k$ is purely inseparable, these define an equivalence $\Sh_0(k_\tau)\cong\Sh_0(l_\tau)$.
\end{proposition}

\begin{myproof}[of \Corollary{weil restriction RPAU}]
    By \cite[Lem. 030K]{stacks-project}, we can assume $l/k$ is either separable or purely inseparable. To show $\pi_*$ and $\pi^*$ send $\Sh_0(k_\tau)$ and $\Sh_0(l_\tau)$ to each other, by \Proposition{RPAU generation} and exactness of $\pi_*$ and $\pi^*$ it suffices to show $\pi_*(\GG_{a,l}^\RP)\in\Sh_0(k_\tau)$ and $\pi^*(\GG_{a,k}^\RP)\in\Sh_0(l_\tau)$. For any $X\in\Spec k_\tau$, we have :
    $$\pi_*(\GG_{a,l}^\RP)(X)=\GG_{a,l}^\RP(X\times_k\Spec l)=\GG_{a,l}(X\times_k\Spec l)=(\Res_{l/k}\GG_{a,l})(X) = (\Res_{l/k}\GG_{a,l})^\RP(X)$$
    by definitions of $\pi_*$, $(-)^\RP$ and $\Res_{l/k}$, and by permanence of relative perfectness under base change, where $X\times_k\Spec l\in\Spec l_\tau$ by \Lemma{base change RP}.4. The group $\Res_{l/k}\GG_{a,l}$ is isomorphic to $\GG_{a,k}^{[l:k]}$ by \Example{affine RP}, thus $\pi_*(\GG_{a,l}^\RP)=(\Res_{l/k}\GG_{a,l})^\RP\in\Sh_0(k_\tau)$ as desired. 
    
    If $l/k$ is separable then $\pi^*(\GG_{a,k}^\RP) = \GG_{a,k}^\RP\times_k\Spec l = (\GG_{a,k}\times_k\Spec l)^\RP = \GG_{a,l}^\RP$ by \Lemma{yoneda and pullback} and \Lemma{base change RP}.1, thus $\pi^*$ maps $\Sh_0(k_\tau)$ to $\Sh_0(l_\tau)$. If $l/k$ is purely inseparable, by \Proposition{pushforward exactness}.4 we have $\pi^*=\rho_*$ for some purely inseparable finite field extension $\rho:\Spec k\to\Spec l$ and $\rho_*$ sends $\Sh_0(k_\tau)$ to $\Sh_0(l_\tau)$ by the previous paragraph. Since $\pi_*$ and $\pi^*$ define an equivalence $\Sh(k_\tau)\cong\Sh(l_\tau)$, their restrictions to subcategories are again an equivalence $\Sh_0(k_\tau)\cong\Sh_0(l_\tau)$.
\end{myproof}

The rest of the section is devoted to \Theorem{weil duality RP}, which compares Serre duality over $k$ and over $l$. It is straightforward when $l/k$ is separable, but requires some preliminaries for the purely inseparable case.

For $R\to S$ a finite syntomic ring map, recall from \cite[Lem. 0FLB]{stacks-project} the map $\Theta_{S/R}^q:\Ohm^q(S)\to\Ohm^q(R)$. It is compatible with base change, $\Theta^0_{S/R}:S\to R$ is the trace morphism, and $\Theta_{S/R}^\bullet$ is a map of differential graded $\Ohm^\bullet(R)$-modules. For $\tau=\RP,\RPS$ and $\pi:\Spec l\to\Spec k$ a finite field extension, as $\pi$ is finite syntomic we have a morphism $\Theta_{l/k}^\bullet:\pi_*\Ohm^\bullet_l\to\Ohm^\bullet_k$ of differential graded $\Ohm^\bullet_k$-modules in $\Sh(k_\tau)$.

\begin{proposition}[trace log dRW]
    Consider $\tau=\RP,\RPS$, $\pi:\Spec l\to\Spec k$ a finite field extension, and $q\geq 0$. There exist morphisms $\Tr^q_n=\Tr^q_{n,l/k}:\pi_*\nu_n(q)_l\to\nu_n(q)_k$ in $\Sh(k_\tau)$, 
    for $n\geq 1$, which satisfy the following.
    \begin{enumerate}
        \item For $n=1$ we have a commutative square with the inclusions $\nu_1(q)\subseteq\Ohm^q$ :
        $$\begin{tikzcd}
            \pi_*\nu_1(q)_l\arrow[r,"\Tr_1^q"]\arrow[d] & \nu_1(q)_k \arrow[d] \\
            \pi_*\Ohm^q_l\arrow[r,"\Theta_{l/k}^q"] & \Ohm^q_k
        \end{tikzcd}$$
        \item For $n\geq 1$ we have commutative squares with the surjections $\nu_{n+1}(q)\to\nu_n(q)$ :
        $$\begin{tikzcd}
            \pi_*\nu_{n+1}(q)_l\arrow[r,"\Tr_{n+1}^q"]\arrow[d] & \nu_{n+1}(q)_k \arrow[d] \\
            \pi_*\nu_n(q)_l\arrow[r,"\Tr_n^q"] & \nu_n(q)_k
        \end{tikzcd}$$
        \item For $n\geq 1$ we have commutative squares with the injections $"p":\nu_n(q)\to\nu_{n+1}(q)$ :
        $$\begin{tikzcd}
            \pi_*\nu_n(q)_l\arrow[r,"\Tr_n^q"]\arrow[d,"{"p"}"] & \nu_n(q)_k \arrow[d,"{"p"}"] \\
            \pi_*\nu_{n+1}(q)_l\arrow[r,"\Tr_{n+1}^q"] & \nu_{n+1}(q)_k
        \end{tikzcd}$$
        \item If $\pi':\Spec m\to\Spec l$ is a further finite field extension, then $\Tr^q_{n,m/k}=\Tr^q_{n,l/k}\circ\pi_*(\Tr^q_{n,m/l})$.
    \end{enumerate}
\end{proposition}

\begin{myproof}[of \Proposition{trace log dRW}]
    We have $\nu_n(q)_k,\pi_*\nu_n(q)_l,\Ohm^q_k,\pi_*\Ohm^q_l\in\Alg^\RP(k)$ by \Proposition{log dRW is RPAU} and \Proposition{weil restriction RPAU}. Hence by \Proposition{yoneda fields} it is enough to define the maps $\Tr^q_{n,l/k}$ and check \textbf{(1.-3.)} in the category $\PSh(\F_k)$, where $\F_k$ is the category of schemes that are finite disjoint unions of spectra of relatively perfect field extensions of $k$.
    
    \textbf{Construction.} Recall the norm maps in Milnor $K$-theory $N^q_{E/F}:K^M_q(E)\to K^M_q(F)$
    for any finite field extension $E/F$ (\cite{BassTate} ; \cite[Prop. 5 of Sec. 1.7]{Kato80} ; \cite[Sec. 7.3]{GillesSzamuely}). We have differential symbols $c_n^q(F):K^M_q(F)\to\nu_n(q)(F)$ for any field $F$ of characteristic $p>0$ and any $n\geq 0$, defined on symbols by :
    $$c_n^q(F)(\{x_1,\dots,x_q\})=\frac{d[x_1]_n}{[x_1]_n}\wedge\dots\wedge\frac{d[x_q]_n}{[x_q]_n}$$
    where $[-]_n:F^\times\to (W_nF)^\times$ is the Teichmüller character. By \cite[Cor. 2.8]{BlochKato86}, it induces isomorphisms $K^M_q(F)/p^n\cong\nu_n(q)(F)$. Thus for any finite extension $E/F$ of fields of characteristic $p>0$ and $n\geq 1$ we have a transfer map $\Tr^q_n(E/F):\nu_n(q)(E)\to\nu_n(q)(F)$ defined by the commutative square :
    $$\begin{tikzcd}
        K_q^M(E)/p^n\arrow[r,"N^q_{E/F}"]\arrow[d,"c_n^q(E)"', "\sim"] & K_q^M(F)/p^n \arrow[d,"c_n^q(F)","\sim"'] \\
        \nu_n(q)(E)\arrow[r,"\Tr_n^q(E/F)"] & \nu_n(q)(F)
    \end{tikzcd}$$
    
    For $k'/k$ a relatively perfect field extension, $l\otimes_kk'$ is a finite product of finite field extensions of $k'$ because $k'/k$ is geometrically reduced \cite[Rem. 2.2.(2)]{BertapelleSuzuki}. Therefore for any $S\in\F_k$, we have a transfer map $\Tr^q_{n,l/k}(S):\pi_*\nu_n(q)_l(S)=\nu_n(q)(S\times_k\Spec l)\to\nu_n(q)(S)$. This morphism is natural in $S\in\F_k$ by \cite[Lem. 7.3.6]{GillesSzamuely} (again using that $k'/k$ is geometrically reduced), so $\Tr^q_{n,l/k}$ defines a morphism in $\PSh(\F_k)$.
    
    \textbf{1.} We show we have, for $E/F$ a finite extension of fields of characteristic $p>0$, a commutative square :
    $$\begin{tikzcd}
        K^M_q(E)\arrow[r,"N^q_{E/F}"]\arrow[d,"c_1^q(E)"'] & K^M_q(F) \arrow[d,"c_1^q(F)"] \\
        \Ohm^q(E)\arrow[r,"\Theta^q_{E/F}"] & \Ohm^q(F)
    \end{tikzcd}$$
    Let $F'/F$ be the extension associated to any $p$-Sylow of $\Gal(\Adh F/F)$. Then it is separable over $F$ so $E'=F'\otimes_FE$ is a finite product of finite field extensions $E_i'/F'$, $1\leq i\leq r$. We have commutative squares :
    $$\begin{tikzcd}
        K^M_q(E)\arrow[r,"N^q_{E/F}"]\arrow[d] & K^M_q(F) \arrow[d] \\
        \bigoplus_{i=1}^rK^M_q(E_i') \arrow[r,"\sum_iN^q_{E_i'/F'}"] & K^M_q(F')
    \end{tikzcd}
    \qquad
    \begin{tikzcd}
        \Ohm^q(E)\arrow[r,"\Theta^q_{E/F}"]\arrow[d] & \Ohm^q(F) \arrow[d] \\
        \bigoplus_{i=1}^r\Ohm^q(E_i') \arrow[r,"\sum_i\Theta^q_{E_i'/F'}"] & \Ohm^q(F')
    \end{tikzcd}$$
    by \cite[Lem. 7.3.6]{GillesSzamuely} and by compatibility of $\Theta^q$ with base change. Since $\Ohm^q_F\to\Ohm^q_{F'}$ is injective by separability of $F'/F$, we are reduced to the case $F=F'$. Then $F$ has no prime-to-$p$ extension so by \cite[Cor. 7.2.10]{GillesSzamuely}, $K^M_\bullet(E)$ is generated by $K^M_1(E)$ as a $K^M_\bullet(F)$-module. By $K^M_\bullet(F)$-linearity of $N^\bullet_{E/F}$ and $\Ohm^\bullet(F)$-linearity of $\Theta^\bullet_{E/F}$, we are reduced to the case $q=1$, which follows from \cite[Cor. in Sec. 7]{Garel}.
        
    \textbf{2.} This property follows from the commutative squares, for $n\geq 1$ and $F$ a field of characteristic $p>0$ :
    $$\begin{tikzcd}
        K^M_q(F)\arrow[r,equals]\arrow[d,"c^q_{n+1}(F)"'] & K^M_q(F) \arrow[d,"c_n^q(F)"] \\
        \nu_{n+1}(q)(F)\arrow[r] & \nu_n(q)(F)
    \end{tikzcd}$$
    which are clear from the definitions of the differential symbols.
    
    \textbf{3.} This follows formally from \textbf{(2.)}. Indeed the injections $"p"$ are defined by the commutative squares :
    $$\begin{tikzcd}
        \nu_{n+1}(q)\arrow[r,"p"]\arrow[d,"s_n"] &\nu_{n+1}(q)\arrow[d,equals]\\
        \nu_n(q)\arrow[r,"{"p"}"] & \nu_{n+1}(q)
    \end{tikzcd}$$
    where $s_n$ denotes the canonical surjection, so by \textbf{(2.)} we have :
    $$\Tr^q_{n+1}\circ\pi_*("p")\circ\pi_*(s_n)=\Tr^q_{n+1}\circ p = p\circ \Tr^q_{n+1} = "p"\circ s_n\circ\Tr^q_{n+1}="p"\circ\Tr_n^q\circ \pi_*(s_n).$$
    By \Proposition{pushforward exactness}.2, $\pi_*(s_n)$ is an epimorphism in $\Sh(k_\tau)$, hence $\Tr_{n+1}^q\circ\pi_*("p")="p"\circ\Tr_n^q$, giving \textbf{(3.)}.

    \textbf{4.} This follows from the functoriality of norms in Milnor $K$-theory, that is the relation $N^q_{E'/F}=N^q_{E/F}\circ N^q_{E'/E}$ for $E'/E/F$ a tower of finite field extensions.
\end{myproof}

\begin{lemma}[trace isomorphism]
    Consider $\tau=\RP,\RPS$ and $\pi:\Spec l\to\Spec k$ a purely inseparable finite field extension. Then the map $\Tr=\varinjlim_n\Tr^d_n:\pi_*\nu_\infty(d)_l\to\nu_\infty(d)_k$ of \Proposition{trace log dRW} (where $d=\log_p[k:k^p]$) is an isomorphism in $\Sh(k_\tau)$.
\end{lemma}

Note that the same lemma generally does not hold for $l/k$ separable or for $\Tr^q_n$ with $0\leq q<d$.

\begin{myproof}[of \Lemma{trace isomorphism}]
    By \Proposition{trace log dRW}.4 and \cite[Lem. 09HI]{stacks-project}, we can assume $l=k[\beta]$ with $\beta^p=\alpha\in k\setminus k^p$. By \Proposition{trace log dRW}.2-3 and \cite[Lem. 4.1.5]{Kato86}, we have commutative diagrams with exact rows  :
    $$\begin{tikzcd}
        0 \arrow[r]& \pi_*\nu_{n-1}(d)_l\arrow[d,"\Tr^d_{n-1}"] \arrow[r]& \pi_*\nu_n(d)_l\arrow[r]\arrow[d,"\Tr^d_n"] & \pi_*\nu_1(d)_l\arrow[d,"\Tr^d_1"]\arrow[r] & 0\\
        0\arrow[r] & \nu_{n-1}(d)_k\arrow[r] & \nu_n(d)_k\arrow[r] & \nu_1(d)_k\arrow[r] & 0
    \end{tikzcd}$$
    which reduces inductively to the case $n=1$. This can be checked at stalks, \textit{i.e.} for $\Tr^d_1(R):\nu_1(d)(S)\to\nu_1(d)(R)$ where $R$ is a strictly Henselian relatively perfect $k$-algebra and $S=R\otimes_kl=R[\beta]$. Then by \Proposition{trace log dRW}.2, the map $\Tr^d_1(R)$ coincides with a morphism $\Theta=\Theta^d_{S/R}:\Ohm^d(S)\to \Ohm^d(R)$ such that :
    \begin{align*}
        \Theta(u\beta^id\log(\beta)\wedge dx_2\wedge\dots\wedge dx_d)
        &= 0 &\text{if }i\neq 0\\
        &= ud\log(\alpha)\wedge dx_2\wedge\dots\wedge dx_d&\text{if }i=0
    \end{align*}
    for all $0\leq i<p$ and $u,x_2,\dots,x_d\in R$ (see \cite[Rem. 0FLF]{stacks-project} with the facts that $\beta\in S^\times$, $\alpha\in R^\times$).
    
    Write $B^d=\im(d:\Ohm^{d-1}\to\Ohm^d)$ the cocycles, $q:\Ohm^d\to\Ohm^d/B^d$ the quotient, and $C^{-1}:\Ohm^d\to\Ohm^d/B^d$ the inverse Cartier operator. Because $R$ and $S$ are strictly Henselian local rings, we have exact sequences :
    $$0\to\nu_1(d)(A)\to\Ohm^q(A)\xrightarrow{C^{-1}_A-q_A}\Ohm^d(A)/B^d(A)\to 0$$
    for $A=R,S$, by \cite[Lem. 4.1.5]{Kato86}. Consider the following claims :
    \begin{enumerate}
        \item[(i)] The morphism $\Theta:\Ohm^d(S)\to \Ohm^d(R)$ is surjective.
        \item[(ii)] There is a map $\tilde\Theta:\Ohm^d(S)/B^d(S)\to\Ohm^d(R)/B^d(R)$ such that $q_R\circ\Theta = \tilde\Theta\circ q_S$ and $C^{-1}_R\circ\Theta=\tilde\Theta\circ C^{-1}_S$.
        \item[(iii)] We have $\ker(\Theta)\subseteq B^d(S)\subseteq\Ohm^d(S)$.
        \item[(iv)] The morphisms $C^{-1}_R:\Ohm^d(R)\to\Ohm^d(R)/B^d(R)$ and $C^{-1}_S:\Ohm^d(S)\to\Ohm^d(S)/B^d(S)$ are isomorphisms.
    \end{enumerate}
    Assuming these, the bijectivity of $\Tr_1^d(R):\nu_1(d)(S)\to\nu_1(d)(R)$ follows. Indeed we have a commutative diagram with exact rows :
    $$\begin{tikzcd}
        0\arrow[r] & \nu_1(d)(S)\arrow[r]\arrow[d,"\Tr_1^d(R)"] & \Ohm^d(S)\arrow[r,"C^{-1}_S-q_S"]\arrow[d,"\Theta"] & \Ohm^d(S)/B^d(S)\arrow[r]\arrow[d,"\tilde\Theta"] & 0\\
        0\arrow[r] & \nu_1(d)(R)\arrow[r] & \Ohm^d(R)\arrow[r,"C^{-1}_R-q_R"] & \Ohm^d(R)/B^d(R)\arrow[r] & 0
    \end{tikzcd}$$
    by (ii), which by (i) and the snake lemma, induces an exact sequence :
    $$0\to\ker(\Tr_1^d(R))\to\ker(\Theta)\xrightarrow{C^{-1}_S-q_S}\ker(\tilde\Theta)\to\coker(\Tr_1^d(R))\to 0.$$
    By (iii) the map $\ker(\Theta)\xrightarrow{C^{-1}_S-q_S}\ker(\tilde\Theta)$ is simply the inverse Cartier operator $\ker(\Theta)\xrightarrow{C^{-1}_S}\ker(\tilde\Theta)$, which is an isomorphism by (ii) and (iv). Hence $\ker(\Tr_1^d(R))=\coker(\Tr_1^d(R))=0$.
    
    We now prove (i-iv). We start with preliminary computations of the modules of differentials. By \cite[Lem. 07P2]{stacks-project}, the $p$-independant one-term family $\{\alpha\}$ of $k$ can be extended to a $p$-basis $(\alpha,x_2,\dots,x_d)$ of $k$. One easily checks that $(\beta,x_2,\dots,x_d)$ is a $p$-basis of $l$. Hence $(\alpha,x_2,\dots,x_d)$ is a strong $p$-basis of $R$ and $(\beta,x_2,\dots,x_d)$ is a strong $p$-basis of $S$, by \Lemma{p-basis comparison}.2. We shall write :
    $$x^m=x_2^{m_2}\cdots x_d^{m_d},\qquad dx = dx_2\wedge\cdots\wedge dx_d,\qquad d\log(x) = x_2^{-1}\cdots x_d^{-1}dx$$
    for $m=(m_2,\dots,m_d)\in I$ where $I=\{(m_2,\dots,m_d) ; 0\leq m_i<p\}$.
    
    \textbf{Modules $\Ohm^d$ and $\Ohm^d/B^d$ and map $q$.} Fix $A=R,S$ and $(y_1,\dots,y_d)=(\alpha,x_2,\dots,x_d),(\beta,x_2,\dots,x_d)$ correspondingly. Then $y_i\in A$ is a unit, since it comes from either field $k$ or $l$. We write :
    $$y^m=y_1^{m_1}\dots y_d^{m_d},\qquad d\log(y_J)=d\log(y_{j_1})\wedge\cdots\wedge d\log(y_{j_r}),\qquad d\log(y)=d\log(y_1)\wedge\cdots\wedge d\log(y_d)$$
    for $m=(m_1,\dots,m_d)\in I'=\{(m_1,\dots,m_d);0\leq m_i<p\}$ and $J=\{1\leq j_1<\dots<j_r\leq d\}$ a set of distinct indices. By \cite[Sec. (38.A)]{Matsumura80}, $(dy_1,dy_2,\dots,dy_d)$ is a free basis of the $A$-module $\Ohm^1(A)$ ; this is unchanged by multiplying each term by a unit, so $(d\log(y_1),\dots,d\log(y_d))$ is again a free $A$-basis of $\Ohm^1(A)$. It follows that we have identifications for $0\leq r\leq d$ :
    \begin{align*}
        \Ohm^r(A)
        &= \bigoplus_{|J|=r}A\cdot d\log(y_J) &\text{as a module over }A,\\
        &= \bigoplus_{|J|=r}\bigoplus_{m\in I'}A^p\cdot y^md\log(y_J) & \text{as a module over }A^p.
    \end{align*}
    In particular we have identifications :
    \begin{align*}
        \Ohm^d(S)
        &= S\cdot d\log(\beta)d\log(x) &\text{as a module over }S,\\
        &= \bigoplus_{0\leq n<p}R\cdot \beta^nd\log(\beta)d\log(x) &\text{as a module over }R,\\
        &= \bigoplus_{0\leq n<p}\bigoplus_{m\in I}S^p\cdot \beta^nx^md\log(\beta)d\log(x) &\text{as a module over }S^p=R^p[\alpha],
    \end{align*}
    and similarly :
    \begin{align*}
        \Ohm^d(R)
        &= R\cdot d\log(\alpha)d\log(x) &\text{as a module over }R,\\
        &= \bigoplus_{0\leq n<p}\bigoplus_{m\in I}R^p\cdot \alpha^nx^md\log(\alpha)d\log(x) &\text{as a module over }R^p.\\
    \end{align*}
    Consider an $A^p$-linear basis element $\ohm=y^md\log(y_{\hat i})$ in $\Ohm^{d-1}(A)$, where $\hat i=\{1\leq j\leq d ; j\neq i\}$. Then : $$d\ohm = (-1)^{i-1}m_iy^md\log(y)$$
    where $(-1)^{i-1}m_i$ is either $0$ or a unit in $A$. Hence $B^d(A)$ is exactly the $A^p$-linear submodule of $\Ohm^d(A)$ spanned by those basis elements $y^md\log(y)$ with $m\neq(0,\dots,0)$. It follows that we have identifications :
    $$\Ohm^d(S)/B^d(S) = S^pd\log(\beta)d\log(x),\qquad \Ohm^d(R)/B^d(R) = R^pd\log(\alpha)d\log(x)$$
    such that the projection $q:\Ohm^d\to\Ohm^d/B^d$ is the projection onto the $(n,m)=(0,(0,\dots,0))$-th coordinate.
    
    \textbf{Inverse Cartier operators.} Take as before $A=R,S$ and $(y_1,\dots,y_d)$ the fixed strong $p$-basis of $A$. The inverse Cartier operator $C^{-1}(A):\Ohm^\bullet(A)\to\Ohm^\bullet(A)/B^\bullet(A)$ is characterized as a ring map satisfying :
    $$C^{-1}(u)=u^p,\qquad C^{-1}(du)=u^{p-1}du$$
    for any $u\in A$ ; if $u\in A^\times$ we thus have $C^{-1}(d\log(u))=d\log(u)$. Hence in degree $d$ it is the map :
    $$C^{-1}_A:\Ohm^d(A)=A\cdot d\log(y)\to\Ohm^d(A)/B^d(A)=A^p\cdot d\log(y),\qquad ud\log(y)\mapsto u^pd\log(y).$$
    In particular it is clearly bijective, proving (iv).
    
    \textbf{Trace morphisms.} By \cite[Lem. 0AX5]{stacks-project}, the map $\Theta:\Ohm^d(S)\to\Ohm^d(R)$ is given by projecting onto the $R$-linear coordinate $n=0$ :
    \begin{align*}
        \Theta:\Ohm^d(S)=\bigoplus_{0\leq n<p}R\cdot \beta^nd\log(\beta)d\log(x)&\to\Ohm^d(R)=R \cdot d\log(\alpha)d\log(x)\\
        \sum_{n=0}^{p-1}u_n\beta^nd\log(\beta)d\log(x)&\mapsto u_0d\log(\alpha)d\log(x).
    \end{align*}
    From this description, $\Theta$ is clearly surjective, proving (i). The kernel of $\Theta$ is spanned by basis elements $\beta^nx^md\log(\beta)d\log(x)$ with $n\neq 0$, thus contained in the kernel of $q_S$ (spanned by those $\beta^nx^md\log(\beta)d\log(x)$ with $n\neq 0$ or $m\neq(0,\dots,0)$), proving (iii). One checks that, under the identification of free $R^p$-modules :
    $$\Ohm^d(S)/B^d(S)=S^p\cdot d\log(\beta)d\log(x)=R^p[\alpha]\cdot d\log(\beta)d\log(x)=\bigoplus_{0\leq n < p}R^p\cdot \alpha^nd\log(\beta)d\log(x)$$
    one can choose $\tilde\Theta$ to be the projection onto the coordinate $n=0$ :
    \begin{align*}
        \tilde\Theta:\Ohm^d(S)/B^d(S)=\bigoplus_{0\leq n<p}R^p\cdot\alpha^nd\log(\beta)d\log(x)&\to\Ohm^d(R)/B^d(R)=R^p\cdot d\log(\alpha)d\log(x),\\
        \sum_{n=0}^{p-1}u_n^p\alpha^nd\log(\beta)d\log(x)&\mapsto u_0^pd\log(\alpha)d\log(x).
    \end{align*}
    from which (ii) follows easily. This concludes.
\end{myproof}

\begin{theorem}[weil duality RP]
    Consider $\tau=\RP,\RPS$ and $\pi:\Spec l\to\Spec k$ a finite field extension. There exist maps :
    $$\pi^*\nu_\infty(d)_k\xrightarrow\eta\nu_\infty(d)_l,\quad \pi_*\nu_\infty(d)_l\xrightarrow\Tr\nu_\infty(d)_k,$$
    in $\Sh(l_\tau)$ and $\Sh(k_\tau)$, and corresponding maps $\pi^*\QQ_p/\ZZ_p(d)_k\to\QQ_p/\ZZ_p(d)_l)$ in $D(l_\tau)$ and $\pi_*\QQ_p/\ZZ_p(d)_l\to\QQ_p/\ZZ_p(d)_k$ in $D(k_\tau)$, with the following properties.
    \begin{enumerate}
        \item The map $\eta:\pi^*\nu_\infty(d)_k\to\nu_\infty(d)_l$ is an isomorphism, and the
        composite :
        $$\pi_*\pi^*\nu_\infty(d)_k\xrightarrow{\pi_*\eta}\pi_*\nu_\infty(d)_l\xrightarrow\Tr\nu_\infty(d)_k$$
        coincides with the counit morphism $\pi_*\pi^*\to\id$.
        \item If $G\otimes^L H\to\QQ_p/\ZZ_p(d)_k$ is a perfect pairing in $D(k_\tau)$ then we have a perfect pairing in $D(k_\tau)$ :
        $$(\pi^*G)\otimes ^L(\pi^*H) = \pi^*(G\otimes^LH)\to\pi^*\QQ_p/\ZZ_p(d)_k\xrightarrow\eta\QQ_p/\ZZ_p(d)_l.$$
        \item If $G\otimes^L H\to\QQ_p/\ZZ_p(d)_l$ is a perfect pairing in $D(l_\tau)$ then we have a perfect pairing in $D(k_\tau)$ :
        $$(\pi_*G)\otimes ^L(\pi_*H) \to \pi_*(G\otimes^LH)\to\pi_*\QQ_p/\ZZ_p(d)_l\xrightarrow\Tr\QQ_p/\ZZ_p(d)_k.$$
        \item Given $G\otimes^LH\to\QQ_p/\ZZ_p(d)_k$ a pairing in $D(k_\tau)$, we have a commutative diagram :
        $$\begin{tikzcd}
            G\arrow[d] &[-3em]\otimes^L &[-3em] H\arrow[rrr] &&& \QQ_p/\ZZ_p(d)_k\arrow[d,equals]\\
            (\pi_*\pi^*G) &[-3em]\otimes^L &[-3em] (\pi_*\pi^*H)\arrow[u]\arrow[r] & \pi_*\pi^*\QQ_p/\ZZ_p(d)_k \arrow[r,"\pi_*\eta"]& \pi_*\QQ_p/\ZZ_p(d)_l \arrow[r,"\Tr"]& \QQ_p/\ZZ_p(d)_k
        \end{tikzcd}$$
        where the vertical maps $\pi_*\pi^*\to \id$ and $\id\to\pi_*\pi^*$  are the counit and unit.
    \end{enumerate}
\end{theorem}

\begin{myproof}[of \Theorem{weil duality RP}]
    By \cite[Lem. 030K]{stacks-project}, we can assume $l/k$ is either separable or purely inseparable.
    
    \textbf{Separable case.} Recall that $\nu_n(q)_k$ for $q\in\ZZ$ is the subsheaf of the sheaf of de Rham-Witt differentials :
    $$W_n\Ohm^q_k:X\in\Spec k_\tau\mapsto W_n\Ohm^q(\Gamma(X,\O_X))$$
    whose sections are elements which are étale-locally sums of elements form $d\log(x_1)\wedge\dots\wedge d\log(x_q)$ with $x_i\in \Gamma(X,\GG_m)$. By \Proposition{pushforward exactness}.3, $\pi^*W_n\Ohm^q_k$ is the restriction of $W_n\Ohm^q_k$ to the subcategory $\Spec l_\tau\subseteq\Spec k_\tau$. Evidently, this coincides with $W_n\Ohm^q_l$ and the subsheaves $\pi_*\nu_n(q)_k\subseteq\pi^*W_n\Ohm^q_k$ and $\nu_n(q)_l\subseteq W_n\Ohm^q_l$ coincide also. This defines an isomorphism $\eta^q_n:\pi^*\nu_n(q)_k\xrightarrow\sim\nu_n(q)_l$ for all $q\in\ZZ$ and $n\geq 1$, which further induces an isomorphism $\eta=\varinjlim_n\eta^d_n:\pi^*\nu_\infty(d)_k\to\nu_\infty(d)_l$.
    
    Let $\Tr$ be the composite $\pi_*\nu_\infty(d)_l\xrightarrow{\pi_*\eta^{-1}}\pi_*\pi^*\nu_\infty(d)_k\to\nu_\infty(d)_k$. Then \textbf{(1.)} is clear.
    
    For \textbf{(2.)} we have :
    \begin{align*}
        \Hom_{D(l_\tau)}(T,\pi^*R\HOM_{k_\tau}(A,B))
        &= \Hom_{D(k_\tau)}((\pi_*T)\otimes^L A,B)\\
        &= \Hom_{D(k_\tau)}(\pi_*(T\otimes^L\pi^*A),B)\\
        &= \Hom_{D(l_\tau)}(T,R\HOM_{l_\tau}(\pi^*A,\pi^*B))
    \end{align*}
    for $A,B\in D(k_\tau)$ and $T\in D(l_\tau)$, by adjunctions (using \Proposition{pushforward exactness}.2), and the projection formula (which can be proven in this context by the same argument as \cite[Lem. 0DK5]{stacks-project}). It follows that
    $\pi^*R\HOM_{k_\tau}(A,B) = R\HOM_{l_\tau}(\pi^*A,\pi^*B)$. By taking $A=G$, $B=\QQ_p/\ZZ_p(d)_k$, and composing with the isomorphism $\eta:\pi^*\QQ_p/\ZZ_p(d)_k\to\QQ_p/\ZZ_p(d)_l$, we get \textbf{(2.)}
    
    For \textbf{(3.)}, the composite :
    $$\pi_*R\HOM_{l_\tau}(A,\pi^*B)\xrightarrow\sim R\HOM_{k_\tau}(\pi_*A,\pi_*\pi^*B)\to R\HOM_{k_\tau}(\pi_*A,B)$$
    is an isomorphism for $A\in D(l_\tau)$ and $B\in D(k_\tau)$, where the first map is obtained by deriving the composite functor $\pi_*\HOM_{l_\tau}(A,\pi^*(-))=\HOM_{k_\tau}(\pi_*A,\pi_*\pi^*(-))$, wherein $\pi_*$ and $\pi^*$ are exact and preserve injectives (by \Proposition{pushforward exactness}.2), and the second is the counit $\pi_*\pi^*\to\id$ ; the proof is identical to \cite[Prop. 1.13 of Ch. V]{Milne80}.
    
    Point \textbf{(4.)} follows from \textbf{(1.)} because we have a commutative diagram :
    $$\begin{tikzcd}
        G\arrow[d] &[-3em]\otimes^L &[-3em] H\arrow[r] & \QQ_p/\ZZ_p(d)_k\\
        \pi_*\pi^*G &[-3em]\otimes^L &[-3em] \pi_*\pi^*H\arrow[u]\arrow[r] & \pi_*\pi^*\QQ_p/\ZZ_p(d)_k \arrow[u]
    \end{tikzcd}$$
    where the vertical maps are (co)units and the lower pairing is induced by the upper pairing. Indeed it suffices to check the commutativity of :
    $$\begin{tikzcd}
        A_0(X)\arrow[d] &[-3em]\otimes &[-3em] B_0(X)\arrow[r] & C_0(X)\\
        \pi_*\pi^*A_0(X) &[-3em]\otimes &[-3em] \pi_*\pi^*B_0(X)\arrow[u]\arrow[r] & \pi_*\pi^*C_0(X) \arrow[u]
    \end{tikzcd}$$
    for $A_0\otimes B_0\to C_0$ a pairing in $\Sh(k_\tau)$ and $X\in\Spec k_\tau$. This can be checked explicitly by picking a Galois extension $l'/k$ containing $l$ and using the explicit description of $\Tr$ as in the proof of \cite[Prop. 1.12 of Ch. V]{Milne80}.
    
    \textbf{Purely inseparable case.} The morphism $\Tr:\pi_*\nu_\infty(d)_l\to\nu_\infty(d)_k$ is the isomorphism of \Lemma{trace isomorphism}. Recall that $\pi_*$ and $\pi^*$ form an equivalence $\Sh(k_\tau)\cong\Sh(l_\tau)$ and define :
    $$\eta:\pi^*\nu_\infty(d)_k\xrightarrow{\pi^*\Tr^{-1}}\pi^*\pi_*\nu_\infty(d)_l\xrightarrow\sim\nu_\infty(d)_l.$$
    Clearly $\eta$ is an isomorphism and we have a commutative diagram :
    $$\begin{tikzcd}
        \pi_*\pi^*\nu_\infty(d)_k\arrow[d,"\sim"]\arrow[r,"\pi_*\pi^*\Tr^{-1}"]\arrow[rr,bend left=20,"\pi_*\eta"] & \pi_*\pi^*\pi_*\nu_\infty(d)_l\arrow[d,"\sim"]\arrow[r,"\sim"] & \pi_*\nu_\infty(d)_l\arrow[d,equals]\\
        \nu_\infty(d)_k\arrow[r,"\Tr^{-1}"] & \pi_*\nu_\infty(d)_l\arrow[r,equals] & \pi_*nu_\infty(d)_l
    \end{tikzcd}$$
    so that $\Tr\circ\pi_*\eta$ is the counit $\pi_*\pi^*\nu_\infty(d)_k\to\nu_\infty(d)_k$ followed by $\Tr\circ\Tr^{-1}=\id_{\nu_\infty(d)_k}$, proving \textbf{(1.)}. Because $\pi:\Spec l_\tau\to\Spec k_\tau$ is a morphism of sites, $\pi^*:D(k_\tau)\to D(l_\tau)$ is monoidal hence its essential inverse $\pi_*:D(l_\tau)\to D(k_\tau)$ is monoidal as well. It follows that we have natural isomorphisms :
    \begin{align*}
        \pi^*R\HOM_{k_\tau}(-,-)&=R\HOM_{l_\tau}(\pi^*(-),\pi^*(-))\\
        \pi_*R\HOM_{k_\tau}(-,\pi^*(-))&=R\HOM_{l_\tau}(\pi_*(-),\pi_*\pi^*(-))=R\HOM_{l_\tau}(\pi_*(-),-)
    \end{align*}
    and the proofs of \textbf{(2.)} and \textbf{(3.)} follow just like in the separable case. Point \textbf{(4.)} follows from the monoidality of $\pi_*$ and $\pi^*$ and the fact that the counit and unit morphisms $\pi_*\pi^*\to\id$, $\id\to\pi_*\pi^*$ are inverse maps.
\end{myproof}


\section{Ind-pro-RPAU groups}
\label{Ind-pro-RPAU groups}


\subsection{Ind-pro-RPAU groups and $\W_k$}
\label{Ind-pro-RPAU groups and Wk}

\begin{definition}[good indpro]
    Let $k$ be a perfect field of characteristic $p>0$. We say a group $G\in\IndProAlg_u^\RP(k)$ is :
        \begin{itemize}
            \item \emph{Type 1} if $G\in\ProAlg_u^\RP(k)\subseteq\IndProAlg_u^\RP(k)$ is an $\NN$-indexed limit $\varprojlim_nG_n$ where each $G_n\in\Alg^\RP_u(k)$ is connected and the transition morphisms $G_{n+1}\to G_n$ are surjective with connected kernels ;
            \item \emph{Type 2} if $G\in\IndAlg_u^\RP(k)\subseteq\IndProAlg_u^\RP(k)$ is an $\NN$-indexed colimit $\varinjlim_nG_n$ where each $G_n\in\Alg^\RP_u(k)$ is connected and the transition morphisms $G_n\to G_{n+1}$ are injective.
        \end{itemize}
    Such an inverse (resp. direct) $\NN$-indexed system will be called a \emph{type 1 inverse system} (resp. \emph{type 2 direct system}). We define $\W_k$ as the full subcategory of $\IndProAlg_u^\RP(k)$ of objects $G$ which admit a filtration $G''\leq G'\leq G$ such that $G''$ is type 1, $G'/G''$ is type 2, and $G/G'$ is finite étale $p$-primary. We say such an object of $\W_k$ is \emph{connected} if we can in addition choose $G=G'$.
\end{definition}

The category $\W_k$ is not abelian. It is however stable under extensions in $\IndProAlg_u^\RP(k)$, which is \Lemma{stability of Wk} below. The proof breaks down into many similar cases, making it long, but it is not very technical.\footnote{We thank Takashi Suzuki for explaining this fact in private communication, which was omitted from his papers.}

\begin{lemma}[stability of Wk]
    Let $k$ be a perfect field of characteristic $p>0$. Then the full subcategory $\W_k\subseteq\IndProAlg_u^\RP(k)$ is stable under extensions. An extension in $\IndProAlg_u^\RP(k)$ of type 1, type 2, or connected objects of $\W_k$ is again type 1, type 2 or connected in $\W_k$, respectively.
\end{lemma}

\begin{myproof}[of \Lemma{stability of Wk}]
    Let $0\to H\to E\to G\to 0$ be an exact sequence in $\IndProAlg_u^\RP(k)$ with $G,H\in\W_k$. We need to show $E\in\W_k$. We work case by case.
    
    \textbf{1. Case $H$ and $G$ are connected.} Consider $\C\subseteq\Alg_u^\RP(k)$ the full subcategory of connected objects. It is closed under extensions and as such is an exact category, where all monics are admissible and admissible epics are surjections with connected kernel. Type 1 (resp. type 2) objects of $\W_k$ are the objects of the full subcategory $\mathrm{Pro}^a_{\aleph_0}(\C)\subseteq\Pro\C$ (resp. $\mathrm{Ind}^a_{\aleph_0}(\C)\subseteq\Ind\C$) of admissible objects of at most countable size, in the sense of \cite[Def. 3.2, Def. 4.1]{BOGTate}, and the full subcategory $\W_k^0\subseteq \W_k$ of connected objects is exactly the category of elementary Tate objects of $\C$ of at most countable size in the sense of \cite[Th. 5.4]{BOGTate} (the proof of which shows any elementary Tate object admits an admissible pro-subobject with admissible ind-quotient). By the same result, $\W_k^0$ is closed under extensions in the category $\mathrm{Ind}^a_{\aleph_0}(\mathrm{Pro}^a_{\aleph_0}(\C))$.

    By \cite[Th. 3.7]{BOGTate}, for $\A$ an exact category, $\mathrm{Ind}^a_{\aleph_0}(\A)$ is closed under extensions in the category of left-exact presheaves on $\A$, which contains $\Ind\A$ as an exact subcategory, therefore $\mathrm{Ind}^a_{\aleph_0}(\A)$ is closed under extensions in $\Ind\A$, and similarly for the pro-categories. Thus $\W_k^0$ is closed under extensions in $\IndPro\C$. By \cite[Prop. 8.6.12]{KashiwaraSchapira}, $\IndPro\C$ is closed under extensions in $\IndProAlg_u^\RP(k)$, hence so is $\W_k^0$.
    
    \textbf{2. Case $G$ is type 1 and $H$ finite étale.} Take a type 1 inverse system $G=\varprojlim_nG_n$. By \cite[Th. 2.2.2]{SuzukiGN}, for some $n_0\in\NN$ we have a commutative diagram with exact rows and cartesian right square :
    $$\begin{tikzcd}[row sep=10]
        0 \arrow[r] & H \arrow[r]\arrow[d,equals] & E \arrow[r]\arrow[d] & G \arrow[r]\arrow[d] & 0 \\
        0 \arrow[r] & H \arrow[r] & E_{n_0} \arrow[r] & G_{n_0} \arrow[r] & 0 
    \end{tikzcd}$$
    Define $E_n$ for $n\geq n_0$ as the pullback of $E_{n_0}$ and $G_n$ over $G_{n_0}$. This way $E=\varprojlim_n E_n$, with $E_n\in\Alg_u^\RP(k)$ and $E_{n+1}\to E_n$ surjective. There is an exact sequence :
    $$0\to H/H'_n\to E_n/E_n^0\to G_n/G_n'\to 0$$
    where $H_n'$ is the inverse image of $E_n^0$ in $H$, and $G_n'$ the image of $E_n^0$ in $G_n$. Then $G_n/G_n'$ is finite étale as a quotient of the finite étale groups $H/H_n'$ and $E_n/E_n^0$, and connected as a quotient of the connected group $G_n$, hence $G_n/G_n'=0$. Thus $G_n=G_n'$ meaning we also have an exact sequence :
    $$0\to H'_n\to E_n^0\to G_n\to 0$$
    Now $(H'_n)_n$ is a decreasing sequence of closed subgroups of the finite étale (thus Noetherian) group $H$, so for large $n\geq n_0$, it is constant at some $H'\leq H$. We have commutative diagrams with exact rows :
    $$\begin{tikzcd}[row sep=10]
        0 \arrow[r] & H' \arrow[r]\arrow[d,equals] & E_{n+1}^0 \arrow[r]\arrow[d] & G_{n+1} \arrow[r]\arrow[d] & 0 \\
        0 \arrow[r] & H' \arrow[r] & E_n^0 \arrow[r] & G_n \arrow[r] & 0 
    \end{tikzcd}$$
    so $E_{n+1}^0\to E_n^0$ has the same kernel and cokernel as $G_{n+1}\to G_n$. It follows that $E'=\varprojlim_{n\geq n_1}E_n^0$ is type 1. We have also seen that $E_n/E_n^0=H/H_n'$, so taking the limit in $n$ we have an exact sequence $0\to E'\to E\to H/H'\to 0$, and $E$ is in $\W_k$.
    
    \textbf{3. Case $G$ is type 2 and $H$ finite étale.} Take a type 2 direct system $G=\varinjlim_nG_n$. Let $E_n\in\Alg_u^\RP(k)$ be the pullback of $G_n$ in $E$. Again define $H'_n\leq H$ as the inverse image of $E_n^0$ in $H$. Then $(H_n')$ is an increasing sequence of reduced closed subgroups of the finite étale group $H$, so for large $n$ is is constant at some $H'\leq H$. As in the previous case we have exact sequences :
    $$0\to H'\to E_n^0\to G_n\to 0,\qquad 0\to E_n^0\to E_n\to H/H'\to 0$$
    because $\coker(E_n^0\to G_n)=(E_n/E_n^0)/(H/H_n')$ is simultaneously connected and finite étale, thus trivial. The snake lemma applied to the first sequence, for varying $n\geq n_0$, shows $E_n^0\to E_{n+1}^0$ is injective thus $E'=\varinjlim_{n\geq n_0}E_n^0$ is type 2. Passing to the colimit in the second sequence shows $E'$ is a subgroup of $E=\varinjlim_{n\geq n_0}E_n$, with finite étale quotient $H/H'$, thus $E$ is in $\W_k$.
    
    \textbf{4. Case $G$ is connected and $H$ is finite étale.} Let $G'\leq G$ be a type 1 subgroup such that $G/G'$ is type 2. Let $E'$ be the inverse image of $G'$ in $E$, then we have a short exact sequence $0\to H\to E'\to G'\to 0$ so by case \textbf{(2.)}, $E'$ admits a type 1 subgroup $E_1$ such that $E'/E_1$ is finite étale. We have a short exact sequence :
    $$0\to E'/E_1\to E/E_1\to E/E'\to 0$$
    where $E/E'=G/G'$ is type 2, so by case \textbf{(3.)}, $E/E_1$ admits a type 2 subgroup such that the quotient is finite étale. Taking $E_2\leq E$ to be the inverse image of such a subgroup in $E$, we have $E_1\leq E_2\leq E$ where $E_1$ is type 1, $E_2/E_1$ is type 2, and $E/E_2$ is finite étale. Thus $E$ is in $\W_k$.
    
    \textbf{5. General case.} Let $G'\leq G$ and $H'\leq H$ be connected objects of $\W_k$ with $G/G'$ and $H/H'$ finite étale. Let $E'$ be the inverse image of $G'$ in $E$. The exact sequence $0\to H/H'\to E'/H'\to G'\to 0$ and case \textbf{(2.)} show $E'/H'$ belongs to $\W_k$. Hence $E'$ admits a subobject $E''$ containing $H'$ such that $E''/H'$ is a connected object of $\W_k$ and $(E'/H')/(E''/H')$ is finite étale. Then case \textbf{(1.)} shows $E''$ is a connected object of $\W_k$. The exact sequence $E'/E''\to E/E''\to E/E'\to 0$ and identifications $E/E'=G/G'$ and $E'/E''=(E'/H')/(E''/H')$ show $E/E''$ is finite étale, proving $E$ belongs to $\W_k$.
\end{myproof}


\subsection{Ind-pro-RPAU groups as sheaves}
\label{Ind-pro-RPAU groups as sheaves}

\begin{definition}[sites of perfect field]
    Let $k$ be a perfect field.
    \begin{itemize}
        \item A $k$-algebra is \emph{rational} if it is a finite product of perfections of finite type field extensions of $k$, and \emph{indrational} if it is a filtered colimit of rational algebras. The \emph{indrational étale (resp. proétale) site} $\Spec k_\et^\indrat$ (resp. $\Spec k_\proet^\indrat$) is the opposite of the category of spectra of indrational $k$-algebras, with all $k$-scheme morphisms, with the étale (resp. proétale) topology.
        \item A $k$-algebra is \emph{perfect Artinian} if it is a finite product of (arbitrary) perfect field extensions of $k$. The \emph{perfect Artinian étale site} $\Spec k_\et^\perar$ is the category of spectra of perfect Artinian $k$-algebras with all $k$-scheme morphisms, with the étale topology.
        \item A family of maps $\{X_i\to X\}_i$ of perfect affine $k$-schemes $X=\Spec R,X_i=\Spec R_i$ is a \emph{pro-fppf covering} if each $R_i$ is a filtered colimit of perfections of flat finitely presented $R$-algebras, and $R\to\prod_iR_i$ is faithfully flat. The \emph{perfect pro-fppf site} $\Spec k_\profppf^\perf$ is the category of perfect affine $k$-schemes with all $k$-scheme morphisms, with the topology defined by pro-fppf coverings.
    \end{itemize}
    We write $f_k:\Spec k_\profppf^\perf\to\Spec k_\proet^\indrat$, $g_k:\Spec k_\proet^\indrat\to\Spec k_\et^\perar$ and $h_k=g\circ f$ the premorphisms of defined by the identity. We define the \emph{change of site functor} $\alpha_k = Rf_{k,*}Lh^*_k:D(k_\et^\perar)\to D(k_\proet^\indrat)$.
\end{definition}

See \cite[Sec. 2.1]{SuzukiGN}, \cite[Sec. 3]{SuzukImp} and \cite[Sec. 3.1]{SuzukiDualLoc} for references on these sites. A sheaf $F\in\Sh(k_\proet^\indrat)$ is a fortiori an object of $\Sh(k_\et^\indrat)$ ; we write $R\Gamma((-)_\proet,F)$ and $R\Gamma((-)_\et,F)$ to differentiate between the resulting cohomologies. Note that the other sites in this paper ($S_\et$, $S_\Et$, $S_\RP$, $S_\RPS$, $\Spec k_\et^\perar$, $*_\proet$) are categories only ever viewed with one topology and we write $R\Gamma(-,-)$ unambiguously for their cohomology.

We will usually omit the $\Spec$ from notations concerning these sites, writing \textit{e.g.} $\Sh(k_\et^\perar)$, $D(k_\et^\perar)$, $\HOM_{k_\et^\perar}(-,-)$ instead of $\Sh(\Spec k_\et^\perar)$, $D(\Spec k_\et^\perar)$, $\HOM_{\Spec k_\et^\perar}(-,-)$.

\begin{definition}[Wk as sheaves]
    Let $k$ be a perfect field. Let $\C=\Spec k_\et^\perar,\Spec k_\proet^\indrat,\IndProAlg^\RP_u(k)$.
    \begin{itemize}
        \item We write $D^b_{\W_k}(\C)$ the full subcategory of objects of $D^b(\C)$ with cohomology in (or representable by objects of) $\W_k$.
        \item We write $\langle\W_k\rangle_\C$ the smallest full triangulated subcategory of $D(\C)$ closed under taking split subobjects, and which contains (objects representable by) objects of $\W_k$ in degree $0$.
    \end{itemize}
\end{definition}

Again we will usually write \textit{e.g.} $\langle\W_k\rangle_{k_\et^\perar}$ instead of $\langle\W_k\rangle_{\Spec k_\et^\perar}$.

The reason for defining these sites is that they properly capture ind-pro-RPAU groups over $k$, a statement made precise in \Proposition{Yoneda for Wk}. Our eventual interest lies in $\Spec k_\proet^\indrat$, which understands all of $\IndProAlg_u^\RP(k)$. The site $\Spec k_\et^\perar$ captures at least $\W_k$ and we can pass from it to $\Spec k_\proet^\indrat$ via the change of site functor. Many of the results we use rely on the study of $\Spec k_\profppf^\perf$, but we will not work with it directly ; we mention it only to define $\alpha_k$.

\begin{proposition}[Yoneda for Wk]
    Let $k$ be a perfect field of characteristic $p>0$.
    \begin{enumerate}
        \item The Yoneda functor induces fully faithful, exact functors of abelian and triangulated categories :
        $$\IndProAlg_u^\RP(k)\to\Sh(k_\proet^\indrat),\qquad D^b(\IndProAlg_u^\RP(k))\to D(k_\proet^\indrat),$$
        identifying $\IndProAlg_u^\RP(k)$ as a full abelian subcategory of $\Sh(k_\proet^\indrat)$ closed under extensions and $D^b(\IndProAlg_u^\RP(k))$ as a full triangulated subcategory of $D(k_\proet^\indrat)$.
        \item The Yoneda functor induces fully faithful, exact functors of exact and triangulated categories :
        $$\W_k\to\Sh(k_\et^\perar),\qquad \langle\W_k\rangle_{\IndProAlg_u^\RP(k)}\to D(k_\et^\perar)$$
        identifying $\W_k$ as a full subcategory of $\Sh(k_\et^\perar)$ closed under extensions, and induces an exact equivalence of triangulated categories $\langle\W_k\rangle_{\IndProAlg_u^\RP(k)}\cong\langle\W_k\rangle_{k_\et^\perar}$. We have commutative squares :
        $$\begin{tikzcd}
            \Alg_u^\RP(k)\arrow[r,"\subseteq"]\arrow[d] & \W_k \arrow[d] \\
            \Sh(k_\RP) \arrow[r,"w_*"] & \Sh(k_\et^\perar)
        \end{tikzcd}
        \qquad
        \begin{tikzcd}
            D^b(\Alg_u^\RP(k))\arrow[r,"\subseteq"]\arrow[d] & \langle\W_k\rangle_{\IndProAlg_u^\RP(k)} \arrow[d] \\
            D(k_\RP) \arrow[r,"w_*"] & D(k_\et^\perar)
        \end{tikzcd}$$
        where $w:\Spec k_\RP\to\Spec k_\et^\perar$ is the premorphism of sites induced by identity.
        \item The change of site $\alpha_k:D(k_\et^\perar)\to D(k_\proet^\indrat)$ induces a commutative square of exact equivalences :
        $$\begin{tikzcd}
            \langle\W_k\rangle_{\IndProAlg_u^\RP(k)}\arrow[r,equals]\arrow[d] & \langle\W_k\rangle_{\IndProAlg_u^\RP(k)} \arrow[d] \\
            \langle\W_k\rangle_{k_\et^\perar}  \arrow[r,"\alpha_k"] & \langle\W_k\rangle_{k_\proet^\indrat}
        \end{tikzcd}$$
        \item For all $k'\in\Spec k^\perar_\et$ and $G\in\langle\W_k\rangle_{k_\et^\perar}$, we have $R\Gamma(k'_\proet,\alpha_kG)=R\Gamma(k'_\et,G)$. If additionally $H^q(G)\in\W_k$ for some $q\in\ZZ$ then $H^q(\alpha_kG)\in\Sh(k_\proet^\indrat)$ and $H^q(G)\in\Sh(k_\et^\perar)$ are representable by the same object of $\W_k$. In particular, $\alpha_k$ also induces a commutative square of equivalences :
        $$\begin{tikzcd}
            D^b_{\W_k}(\IndProAlg_u^\RP(k))\arrow[r,equals]\arrow[d] & D^b_{\W_k}(\IndProAlg_u^\RP(k)) \arrow[d] \\
            D^b_{\W_k}(k_\et^\perar) \arrow[r,"\alpha_k"] & D^b_{\W_k}(k_\proet^\indrat)
        \end{tikzcd}$$
    \end{enumerate}
\end{proposition}

\begin{myproof}[of \Proposition{Yoneda for Wk}]
    \textbf{1.} See \cite[Prop. 2.3.4]{SuzukiGN} for the second functor's full faithfulness. The first's follows because $\Hom_{D(k^\indrat_\proet)}(G,H)=\Hom_{\Sh(k_\proet^\indrat)}(G,H)$ for $G,H\in\Sh(k^\indrat_\proet)$, viewed as complexes concentrated in degree $0$. The stability of $\IndProAlg_u^\RP(k)$ under extensions in $\Sh(k_\proet^\indrat)$ follows from the isomorphism, for $G,H\in\IndProAlg_u^\RP(k)$ :
    $$\Ext^1_{\IndProAlg_u^\RP(k)}(G,H)=\Hom_{D^b(\IndProAlg_u^\RP(k))}(G,H[1])=\Hom_{D(k^\indrat_\proet)}(G,H[1])=\Ext^1_{\Sh(k_\proet^\indrat)}(G,H).$$
    
    \textbf{3.} See \cite[Prop. 3.1.5]{Suzuki2024}. 
    
    \textbf{4.} See \cite[Prop. 3.1.4]{Suzuki2024} and \cite[Prop. 2.1]{SuzukiCFT}.
    
    \textbf{2.} For the first part, the full faithfulness $\langle\W_k\rangle_{\IndProAlg_u^\RP(k)}\to D(k_\et^\perar)$ follows from \textbf{(1.)} and \textbf{(3.)}, that of $\W_k\to\Sh(k_\et^\perar)$ follows as in \textbf{(1.)}, and the stability of $\W_k$ under extensions in $\Sh(k_\et^\perar)$ follows as in \textbf{(1.)} from its stability under extensions in $\IndProAlg_u^\RP(k)$ (\Lemma{stability of Wk}). For the commutative squares, note that $w_*$ takes a sheaf to its restriction (as a presheaf) to the full subcategory $\Spec k_\et^\perar\subseteq\Spec k_\RP$. As such it commutes with the Yoneda functor on the underived level. Since both $w_*$ and the Yoneda functor are exact, they commute on the derived level as well. 
\end{myproof}


\subsection{Serre duality and duality for $\W_k$}
\label{Serre duality and duality for Wk}

There is a suitable notion of duality for many ind-pro-RPAU groups and especially objects of $\W_k$, which generalizes in an intuitive way the duality for RPAU groups (\Theorem{duality for RPAU groups}), called Serre duality in \cite{SuzukiGN}.

Consider the functors  $(-)^0:\Alg_u^\RP(k)\to\Alg_u^\RP(k)$ and $\pi_0:\Alg_u^\RP(k)\to\Alg_u^\RP(k)$ of neutral component and component group. Their extension between ind-pro-categories are again written $(-)^0,\pi_0:\IndProAlg_u^\RP(k)\to\IndProAlg_u^\RP(k)$. Namely we have, by definition :
$$G^0=\varinjlim_i\varprojlim_jG_{ij}^0,\qquad \pi_0(G)=\varinjlim_i\varprojlim_j\pi_0(G_{ij}).$$
whenever $G_{ij}\in\Alg_u^\RP(k)$ and $G=\varinjlim_i\varprojlim_jG_{ij}\in\IndProAlg_u^\RP(k)$. The natural maps $(-)^0\to\id\to\pi_0$ extend to the ind-pro-functors similarly. We will say $G$ is \emph{connected} if $\pi_0(G)=0$ and \emph{ind-pro-finite étale} if $G^0=0$.

\begin{lemma}[connected Wk]
    Consider $G\in\IndProAlg_u^\RP(k)$.
    \begin{enumerate}
        \item We have an exact sequence $0\to G^0\to G\to \pi_0(G)\to 0$ formed by the natural maps, functorial in $G$. In particular $G$ is connected (resp. ind-pro-finite étale) if and only if $G=G^0$ (resp. $G=\pi_0(G)$).
        \item If $G\to H$ is a surjection in $\IndProAlg_u^\RP(k)$ and $G$ is connected then $H$ is connected. If $H\to G$ is an injection in $\IndProAlg_u^\RP(k)$ and $G$ is ind-pro-finite étale, then $H$ is ind-pro-finite-étale.
        \item If $G$ is connected and $H$ is ind-pro-finite étale, then any morphism $G\to H$ is trivial.
        \item The object $G^0$ is the unique subobject of $G$ which is connected with ind-pro-finite étale quotient.
        \item If $G\in\W_k$ then $G^0$ is connected in $\W_k$ (in the sense of \Definition{good indpro}) and $\pi_0(G)$ is finite étale.
    \end{enumerate}
\end{lemma}

\begin{myproof}[of \Lemma{connected Wk}]
    \textbf{1.} This follows from the exactness of cofiltered limits in $\ProAlg_u^\RP(k)$ and of filtered colimits in $\IndProAlg_u^\RP(k)$. Functoriality follows from the naturality of the maps $(-)^0\to \id\to \pi_0$.
    
    \textbf{2.} From \textbf{(1.)}, any surjection $G\to H$ (resp. injection $H\to G$) induces a surjection $\pi_0(G)\to\pi_0(H)$ (resp. injection $H^0\to G^0$). The claim  then follows from the definitions.
    
    \textbf{3.} As $G$ is connected, the map factors as $G=G^0\to H^0\to H$ by functoriality of $(-)^0$, and $H^0=0$.
    
    \textbf{4.} For $H\leq G$ a subobject in $\IndProAlg_u^\RP(k)$ with $H$ connected and $G/H$ ind-pro-finite étale, as in \textbf{(3.)} we have a map $H=H^0\to G^0$ which induces $G/H\to G/G^0=\pi_0(G)$, and we have a map $\pi_0(G)\to\pi_0(G/H)=G/H$. These maps exhibit an isomorphism $G/H=G/G^0$, hence $H=G^0$ as subobjects of $G$.
    
    \textbf{5.} There exists $G^{(0)}\subseteq G$ connected in $\W_k$ such that $G/G^{(0)}$ is finite étale. Then $G/G^{(0)}$ is in particular ind-pro-finite étale, so by \textbf{(4.)} it suffices to show $G^{(0)}$ is connected in $\IndProAlg_u^\RP(k)$. Consider an exact sequence $0\to G^{(1)}\to G^{(0)}\to G^{(2)}\to 0$ with $G^{(1)}$ type 1 and $G^{(2)}$ type 2. Given type 1 and 2 systems $G^{(1)}=\varprojlim_m G^{(1)}_m$ and $G^{(2)}=\varinjlim_nG^{(2)}_n$, consider a following commutative diagram with exact rows :
    $$\begin{tikzcd}
       0 \arrow[r]& G^{(1)}   \arrow[r]\arrow[d,equals]& G^{(0)}       \arrow[r]& G^{(2)} \arrow[r]& 0\\
       0 \arrow[r]& G^{(1)}   \arrow[r]\arrow[d]& G^{(0)}_n\arrow[u]\arrow[d]     \arrow[r]& G^{(2)}_n\arrow[u] \arrow[r]& 0\\
       0 \arrow[r]& G^{(1)}_m \arrow[r]& G^{(0)}_{n,m} \arrow[r]& G^{(2)}_n\arrow[u,equals] \arrow[r]& 0
    \end{tikzcd}$$
    with cartesian upper-right square and cocartesian lower-left square. Then $G^{(0)}_{n,m}\in\Alg_u^\RP(k)$ is connected as an extension of connected $k$-groups. Hence $G^{(0)}=\varinjlim_n\varprojlim_mG^{(0)}_{n,m}$ is connected by definition.
\end{myproof}

\begin{theorem}[Serre duality]
    Let $k$ be a perfect field of characteristic $p>0$. Let $S_k$ be one of the sites $\Spec k_\proet^\indrat$ or $\Spec k_\et^\perar$. We define the \emph{Serre duality} functor as $(-)^\vee_{S_k}=R\HOM_{S_k}(-,\QQ_p/\ZZ_p):D(S_k)\to D(S_k)$.
    \begin{enumerate}
        \item The functor $(-)^\vee_{S_k}$ sends $\langle\W_k\rangle_{S_k}$ to itself and we have a commutative square :
        $$\begin{tikzcd}
           \langle\W_k\rangle_{k_\et^\perar}^\opp \arrow[r,"\alpha_k"]\arrow[d,"(-)^\vee_{k_\et^\perar}"'] & \langle\W_k\rangle_{k_\proet^\indrat}^\opp \arrow[d,"(-)^\vee_{k_\proet^\indrat}"] \\
           \langle\W_k\rangle_{k_\et^\perar} \arrow[r,"\alpha_k"] & \langle\W_k\rangle_{k_\proet^\indrat}
        \end{tikzcd}$$
        where $\alpha_k:D(k_\et^\perar)\to D(k_\proet^\indrat)$ is the change of site functor. Thus $(-)^\vee_{k_\et^\perar}$ and $(-)^\vee_{k_\proet^\indrat}$ restrict via the Yoneda embeddings to the same functor on $\langle\W_k\rangle_{\IndProAlg_u^\RP(k)}$, simply denoted $(-)^\vee$.
        \item For $G\in\langle\W_k\rangle_{\IndProAlg_u^\RP(k)}$, the natural map $G\to (G^\vee)^\vee$ is an isomorphism. In other words, the pairing  $G\otimes^L G^\vee\to\QQ_p/\ZZ_p$ is perfect when viewed in either $D(k_\et^\perar)$ or $D(k_\proet^\indrat)$.
        \item For $G\in D^b_{\W_k}(\IndProAlg_u^\RP(k))$, the dual $G^\vee$ is computed as follows.
        \begin{enumerate}
            \item If $G\in\Alg_u^\RP(k)$ then $G^\vee=w_*R\HOM_{k_\RP}(G,\QQ_p/\ZZ_p)$ where $w:\Spec k_\RP\to \Spec k_\et^\perar$ is the premophism of sites induced by identity.
            \item If $G\in\W_k$ is connected then $G^\vee$ is concentrated in degree $1$ and $\EXT^1(G,\QQ_p/\ZZ_p)$ is a connected object of $\W_k$. The functor $\EXT^1(-,\QQ_p/\ZZ_p)$ sends type 1 systems to type 2 systems and \textit{vice versa}.
            \item If $G\in\W_k$, $H^0(G^\vee)=\HOM(\pi_0(G),\QQ_p/\ZZ_p)$ is finite étale, $H^1(G^\vee)=\EXT^1(G^0,\QQ_p/\ZZ_p)$ is connected in $\W_k$, and $H^q(G^\vee)=0$ for $q\neq 0,1$. In particular, $G^\vee\in D^b_{\W_k}(\IndProAlg_u^\RP(k))$.
            \item If $G\in D^b_{\W_k}(\IndProAlg_u^\RP(k))$ and $H=G^\vee$, then for all $q\in\ZZ$ we have isomorphisms :
            $$H^q(G)^0=\EXT^1(H^{1-q}(H)^0,\QQ_p/\ZZ_p),\qquad\pi_0(H^q(G))=\HOM(\pi_0(H^{-q}(H)),\QQ_p/\ZZ_p).$$
            We say $H^q(G)^0$ and $H^{1-q}(H)^0$ (resp. $\pi_0(H^q(G))$ and $\pi_0(H^{-q}(H))$) satisfy a perfect Serre duality of connected (resp. finite étale) objects of $\W_k$.
        \end{enumerate}
    \end{enumerate}
\end{theorem}

\begin{myproof}[of \Theorem{Serre duality}]
    \textbf{1.} That $\alpha_k(G^\vee_{k_\et^\perar})=(\alpha_kG)^\vee_{k_\proet^\indrat}$ for $G\in\W_k$ is \cite[Prop. 3.1.8]{Suzuki2024}. The fact that $(-)^\vee_{S_k}$ sends $\W_k$ to $\langle\W_k\rangle_{S_k}$, hence $\langle\W_k\rangle_{S_k}$ to itself, is part of \textbf{(3.c)} below.
    
    \textbf{2.} If $G=\varprojlim_nG_n$ is a type 1 inverse system then by \textbf{(3.b)} below, $G^\vee$ is concentrated in degree $1$ with $H^1(G^\vee)=\varinjlim_n H^1(G_n^\vee)$ a type 2 direct system, and $(G^\vee)^\vee$ is concentrated in degree $0$ with $H^0((G^\vee)^\vee)=\varprojlim_n H^1(H^1(G_n^\vee)^\vee)$. By \textbf{(3.a)} and \Theorem{duality for RPAU groups}.2, since each $G_n$ is connected hence split (because $k$ is perfect), the maps $G_n\to H^1(H^1(G_n^\vee)^\vee)$ are isomorphisms. Thus $G\to (G^\vee)^\vee$ is an isomorphism in $D^b(\IndProAlg_u^\RP(k))$. If $G$ is type 2, by an identical argument $G\to (G^\vee)^\vee$ is an isomorphism. If $G$ is finite étale $p$-primary, by \textbf{(3.a)} and \Theorem{duality for RPAU groups}.1 the morphism $G\to (G^\vee)^\vee$ is an isomorphism.
    
    Objects of $\W_k$ are extensions of type 1, type 2 and finite étale $p$-primary groups, so by exactness of $(-)^\vee$, the natural map $G\to (G^\vee)^\vee$ is an isomorphism for $G\in\W_k$, and in turn for $G\in\langle\W_k\rangle_{\IndProAlg_u^\RP(k)}$.
    
    \textbf{3.a.} This is part of \cite[Prop. 2.4.1.b]{SuzukiGN}.
    
    \textbf{3.b.} Let $G=\varprojlim_nG_n$ be a type 1 inverse system. By \cite[Prop. 2.3.3.(a)]{SuzukiGN}, for $q\in\ZZ$ we have :
    \begin{align*}
        H^q(G^\vee)
        &= \EXT^q_{k_\proet^\indrat}(G,\QQ_p/\ZZ_p)\\
        &=\varinjlim_m\varinjlim_n\EXT^q_{k_\proet^\indrat}(G_n,\ZZ/p^m)\\
        &=\varinjlim_n\varinjlim_m\EXT^q_{k_\proet^\indrat}(G_n,\ZZ/p^m)\\
        &=\varinjlim_n\EXT^q_{k_\proet^\indrat}(G_n,\QQ_p/\ZZ_p)\\
        &=\varinjlim_n H^q(G_n^\vee)
    \end{align*}
    Hence $H^q(G^\vee)=\varinjlim_nH^q(G_n^\vee)$. By \textbf{(3.a)} and \Theorem{duality for RPAU groups}.2, $H^q(G_n^\vee)$ is in $\Alg_u^\RP(k)$, trivial if $q\neq 1$, and connected if $q=1$. Thus $H^q(G^\vee)=0$ for $q\neq 1$. For $n\in\NN$ we have an exact sequence :
    $$0\to N_n\to G_{n+1}\to G_n\to 0$$
    where $N_n\in\Alg_u^\RP(k)$ is some connected group by the type 1 assumption, so taking the long exact sequence of Ext's and applying the above yields an exact sequence :
    $$0\to H^1(G^\vee_n)\to H^1(G^\vee_{n+1})\to H^1(N_n^\vee)\to 0$$
    showing that $H^1(G^\vee_n)\to H^1(G^\vee_{n+1})$ is injective, so $H^1(G^\vee)=\varinjlim_nH^q(G_n^\vee)$ is a type 2 direct system.
    
    Let $G=\varinjlim_nG_n$ be a type 2 direct system. Again by \cite[Th. 2.3.3.(a)]{SuzukiGN} we have a spectral sequence :
    $$E^{ij}_2=(R^i\varprojlim_n)H^j(G_n^\vee) \Rightarrow \EXT^{i+j}_{k_\proet^\indrat}(G,\QQ/\ZZ)=H^{i+j}(G^\vee)$$
    where as before, $H^j(G_n^\vee)$ is trivial if $j\neq 1$, and connected in $\Alg_u^\RP(k)$ if $j=1$. Also by \cite[Prop. 2.1.2.f]{SuzukiGN}, $E^{ij}_2=0$ if $i\neq 0$. This shows $H^q(G^\vee)=0$ if $q\neq 1$, and $H^1(G^\vee)=E^{0,1}_2=\varprojlim_nH^1(G_n^\vee)$. For $n\in\NN$ we have an exact sequence :
    $$0\to G_n\to G_{n+1}\to Q_n\to 0$$
    where $Q_n\in\Alg_u^\RP(k)$ is some connected group as a quotient of a connected group, so taking the long exact sequence of Ext's yields an exact sequence :
    $$0\to H^1(Q^\vee_n)\to H^1(G^\vee_{n+1})\to H^1(G^\vee_n)\to 0$$
    where $H^1(Q^\vee_n)$ is connected, so $H^1(G^\vee)=\varinjlim_nH^1(G_n)$ is a type 1 inverse system.
    
    If $G\in\W_k$ is connected, writing $G$ as an extension of type 1 and type 2 groups and taking the long exact sequence of Ext's shows $H^q(G^\vee)=0$ if $q\neq 1$, and $H^1(G^\vee)$ is a connected object of $\W_k$.
    
    \textbf{3.c.} This is obtained by taking the long exact sequence of Ext's for $0\to G^0\to G\to \pi_0(G)\to 0$ and noting that, by \textbf{(3.a)} and \Theorem{duality for RPAU groups}.2, $\pi_0(G)^\vee$ is concentrated in degree $0$.
    
    \textbf{3.d.} This follows from \textbf{(3.c)} as the spectral sequence :
    $$E^{ij}_2 = \EXT^i(H^{-j}(H),\QQ_p/\ZZ_p)\Rightarrow H^{i+j}(G)$$
    degenerates into a short exact sequence for each $q\in\ZZ$ :
    $$0\to \EXT^1(H^{1-q}(H),\QQ_p/\ZZ_p)\to H^q(G)\to \HOM(H^{-q}(H),\QQ_p/\ZZ_p)\to 0$$
    whose first term is connected in $\W_k$ and whose third term is finite étale, by \textbf{(3.c)}.
\end{myproof}


\section{Ind-pro-finite groups}
\label{Ind-pro-finite groups}


\subsection{Ind-pro-finite groups as condensed groups}
\label{Ind-pro-finite groups as condensed groups}

\begin{definition}[topological abelian groups]
    Let $S$ be a fixed set of primes.
    \begin{itemize}
        \item We write $\Ab_S$ the category of finite abelian groups $G$ such that the $p$-torsion of $G$ is trivial for all primes $p\notin S$. When $S$ is the set of all primes, write $\Ab_S=\Ab_\fin$ the category of all finite abelian groups. When $S=\{p\}$ has one element, write $\Ab_S=\Ab_p$ the category of finite abelian $p$-groups.
        \item Define $\W_S\subseteq\Ind\ProAb_S$ as the full subcategory of objects $G\in\IndProAb_S$ such that there exists an exact sequence $0\to G'\to G\to G''\to 0$ where $G'\in\ProAb_S$ is an $\NN$-indexed limit of objects of $\Ab_S$, and $G''\in\IndAb_S$ is an $\NN$-indexed colimit of objects of $\Ab_S$. Similarly we write $\W_S=\W_\fin$ if $S$ is the set of all primes and $\W_S=\W_p$ if $S=\{p\}$.
    \end{itemize}
\end{definition}

Note that $\Ab_S=\bigoplus_{p\in S}\Ab_p$ and $D^b(\Ab_p)=\bigoplus_{p\in S}D^b(\Ab_p)$, as $\Ext^i_{\Ab}(A,B)=0$ for $i\geq 0$ if $A\in\Ab_p$ and $B\in\Ab_q$ with primes $p\neq q$. The same decomposition holds for $\Ind\Ab_S$ and $\Pro\Ab_S$ if $S$ is finite, but in general $\IndPro\Ab_S$ is larger than $\bigoplus_{p\in S}\IndProAb_p$. For instance, $\Ind\Ab_\fin\cong\prod_p\Ind\Ab_p$ is the category of torsion abelian groups, but $\bigoplus_p\Ind\Ab_p$ is the category of finite exponent abelian groups.\footnote{Here for additive categories $\{\A_i\}_{i\in I}$, the direct sum $\bigoplus_i\A_i$ is the full subcategory of the product $\prod_i\A_i$ whose objects are those families $(A_i)_{i\in I}$ such that $A_i=0$ for all but finitely many $i$'s.}

\begin{definition}[proetale point definition]
    Let $S$ be a fixed set of primes.
    \begin{itemize} 
        \item The \emph{proétale site of the point} $*_\proet$ is the category of pro-finite spaces with coverings given by jointly surjective finite families of continuous maps. A \emph{condensed group} is an object of $\Sh(*_\proet)$.
        \item We write $\TopAb$ the category of topological abelian groups with continuous morphisms, and $\LCAb$ the full subcategory of locally compact Hausdorff topological abelian groups.
        \item Consider three functors :
        $$Y:\IndProAb_S\to\Sh(*_\proet),\qquad (-)_\top:\IndProAb_S\to\TopAb,\qquad (-)_\cond:\TopAb\to\Sh(*_\proet)$$
        respectively given as the ind-completion of the Yoneda functor  $\ProAb_S\to\Sh(*_\proet)$ ; the ind-pro-completion of the functor $\Ab_S\to\TopAb$ giving the discrete topology to a finite group ; by $G_\cond(X)$ being the set of continuous maps $X\to G$ for $X$ pro-finite and $G\in\TopAb$.
        \item We say a sequence $0\to A\xrightarrow v B\xrightarrow u C\to 0$ in $\TopAb$ is \emph{topologically exact} if $u\circ v=0$, $v$ is an embedding (\textit{i.e.} a homeomorphism between $A$ and $v(A)\subseteq B$ with the subspace topology), and the map $B/A\to C$ induced by $u$ is a homeomorphism.
    \end{itemize}
\end{definition}

See \cite{ScholzeClausen} for more on $*_\proet$ and condensed groups ; our $(-)_\cond$ is written $G\mapsto\underline G$ there. Equivalently, $0\to A\xrightarrow v B\xrightarrow u C\to 0$ is topologically exact if it is exact, $v$ is an embedding, and $u$ is open. In the literature, such a sequence is also called \emph{strictly exact}. This endows $\TopAb$ and $\LCAb$ with exact structures.

In this section, we compare objects of $\IndProAb_S$, $\LCAb$ and $\Sh(*_\proet)$ : the former two categories fully faithfully exactly embed in the latter, and their intersection contains the category of "good" ind-pro-finite groups $\W_S$ (\Proposition{topology of W0}). This allows us to work in the nice category $\Sh(*_\proet)$ without difficult topological considerations, and at the end translate condensed statements in terms of locally compact groups. The category $\TopAb$ is used mostly to define the functors of \Definition{proetale point definition} in the greatest possible generality, though only objects of $\LCAb$ will eventually be of interest to us.

\begin{remark}[topological subtleties]
    The category $\TopAb$ has all limits, given by computing the limit in $\Ab$ and equipping it with the initial topology with respect to the projections. The category $\TopAb$ also has all colimits, and the object of $\Ab$ underlying a colimit $G=\mathrm{colim}_iG_i$ in $\TopAb$ is the same colimit computed in $\Ab$. The topology on $G$ is given as the initial topology with respect to (an essentially small, appropriate class of) maps $G\to H$ with $H\in\TopAb$ such that $G_i\to H$ is continuous for all $i$'s. This topology is hard to pin down, and usually different from the the final topology with respect to the maps $G_i\to G$ : the formation of the final topology on a colimit of sets does not commute with products in general, so $G$ equipped with the final topology is not generally a topological group. However, if the colimit is indexed by a countable directed set, and the $G_i$'s are locally compact Hausdorff, then $\mathrm{colim}_i(G_i\times G_i)=(\mathrm{colim}_iG_i)\times(\mathrm{colim}_iG_i)$ in the category of topological spaces, so the colimits in $\TopAb$ and in topological spaces coincide in this case \cite[Th. 4.1]{HHN}. Thus for  $G=\varinjlim_iG_i\in\IndProAb_\fin$ with $G_i\in\Pro\Ab_\fin$, if the index category for $i$ is a countable directed set, then the topology on the colimit $G_\top=\varinjlim_i(G_i)_\top$ is the same as the final topology with respect to the maps $(G_i)_\top\to G_\top$. This in particular applies to $G\in\W_\fin$. 
\end{remark}

\begin{lemma}[surjection in Pro]
    For $S$ a set of primes and $G\to H$ a surjection in $\ProAb_S$, the map $G_\top\to H_\top$ is surjective.
\end{lemma}

\begin{myproof}[of \Lemma{surjection in Pro}]
    Write $N=\ker(G\to H)\in\ProAb_S$, and choose cofiltered systems $N=\varprojlim_iN_i$ and $H=\varprojlim_jH_j$ with $N_i,H_j\in\Ab_S$.
    Let $G_i$ be the pushout of $G$ along $N\to N_i$, so we have an exact sequence $0\to N_i\to G_i\to H\to 0$ and an isomorphism $G=\varprojlim_iG_i$.
    
    By \cite[Prop. 2.2.2]{SuzukiGN}, we have an isomorphism $\Ext^1_{\ProAb_S}(H,N_i) = \varinjlim_j\Ext^1_{\ProAb_S}(H_j,N_i)$ so for some $j_0$ there exists an exact sequence $0\to N_i\to G_{i,j_0}\to H_{j_0}\to 0$ which pulls back to $0\to N_i\to G_i\to H_{j_0}\to 0$. Taking $G_{i,j}$ to be the pullback of $G_{i,j_0}$ along $H_j\to H_{j_0}$ for $j\geq j_0$, we get an exact sequence in $\Ab_S$, $0\to N_i\to G_{i,j}\to H_j\to 0$ where $G_i=\varprojlim_{j\geq j_0}G_{i,j}$.
    
    Clearly the sequence $0\to (N_i)_\top\to (G_{i,j})_\top\to (H_j)_\top\to 0$ is exact. Taking the cofiltered limit in $j\geq j_0$, because the system of kernels $\{(N_i)_\top\}_{j\geq j_0}$ is constant, one can adapt the proof of \cite[Lem. 0598]{stacks-project} to see the map $(G_i)_\top\to H_\top$ is surjective for all $i$.
    
    To deduce $G_\top\to H_\top$ is surjective, consider $x\in H_\top$ and $Z\subseteq G_\top$ its inverse image. We want to see $Z$ is nonempty. Consider also $Z_i$ the inverse image of $x$ in $(G_i)_\top$, for each $i$ : these are nonempty by surjectivity of $(G_i)_\top\to H_\top$, and we have $Z=\varprojlim_iZ_i$ because limits commute with fiber products. The spaces $H_\top$ and $(G_i)_\top$ are compact Hausdorff and the map $(G_i)_\top\to H_\top$ is continuous, so $Z_i$ is compact Hausdorff as a closed subset of $(G_i)_\top$. By Tychonoff's theorem, as a cofiltered limit of nonempty compact Hausdorff spaces, $Z$ is (compact Hausdorff and) nonempty. This concludes.
\end{myproof}

\begin{proposition}[W0 closed under ext]
    For $S$ a set of primes, the full subcategory $\W_S\subseteq\IndProAb_S$ is closed under extensions.
\end{proposition}

\begin{myproof}[of \Lemma{W0 closed under ext}]
    Consider an exact sequence $0\to H\to E\to G\to 0$ in $\IndProAb_S$ such that $H,G\in\W_S$.
    
    \textbf{Case $H$ and $G$ ind-finite.} Write $H=\varinjlim_mH_m$ and $G=\varinjlim_nG_n$ with $H_m,G_n\in\Ab_S$. Let $E_n$ be the pullback of $E$ along $G_n\to G$, so we have a short exact sequence $0\to H\to E_n\to G_n\to 0$ with $E=\varinjlim_nE_n$. By \cite[Prop. 2.2.2]{SuzukiGN} we have :
    $$\Ext^1_{\IndProAb_S}(G_n,H)=\Ext^1_{\IndAb_S}(G_n,H)=\varinjlim_m\Ext^1_{\Ab_S}(G_n,H_m)$$
    so for some $m_0$ there exists an exact sequence $0\to H_{m_0}\to E_{n,m_0}\to G_n\to 0$ which pushes out to $0\to H\to E_n\to G_n\to 0$. For $m\geq m_0$ let $E_{n,m}$ be the pushout of $E_{n,m_0}$ along $H_{m_0}\to H_m$ ; then we have exact sequences $0\to H_m\to E_{n,m}\to G_n\to 0$ so $E_{n,m}\in\Ab_S$, and $E=\varinjlim_n\varinjlim_{m\geq m_0}E_{n,m}$. Choosing $\{H_m\}_m$ and $\{G_n\}_n$ to be $\NN$-indexed systems, we have $E=\varinjlim_{m\geq m_0}E_{m,m}$ hence $E\in\W_S$.
    
    \textbf{Case $H$ and $G$ pro-finite.} The proof is dual to the previous case.
    
    \textbf{General case.} Consider exact sequences $0\to H'\to H\to H''\to 0$ and $0\to G'\to G\to G''\to 0$ where $H',G'$ are $\NN$-indexed objects of $\ProAb_S$ and $H'',G''$ are $\NN$-indexed objects of $\IndAb_S$. Write $H''=\varinjlim_mH''_m$ with $H_m''\in\Ab_S$ and let $H_m\subseteq H$ be the pullback of $H''_m$. The exact sequence $0\to H'\to H_m\to H''_m\to 0$ and the pro-finite case show $H_m\in\W_S$ is pro-finite, and we have $H=\varinjlim_mH_m$.
    
    By \cite[Prop. 2.2.2]{SuzukiGN} , the sequence $0\to H\to E'\to G'\to 0$ is the pushout of some exact sequence $0\to H_m\to E_m'\to G'\to 0$. By the pro-finite case, $E_m'\in\W_S$ is pro-finite. By definition of $H_m$ we have :
    $$H/H_m = H''/H''_m = \varinjlim_{n\geq m}H_n''/H_m''$$
    so $H/H_m\in\W_S$ is ind-finite. The exact sequence $0\to H/H_m\to E/E'_m\to G''\to 0$ shows $E/E_m'\in\W_S$ is ind-finite by the ind-finite case. Hence $E$ belongs to $\W_S$.
\end{myproof}

The bounded derived category $D^b(\LCAb)$ of the nonabelian category $\LCAb$ is defined in \cite{HoffmannSpitzweck} ; for $G,H\in\LCAb$ and $q\in\geq 0$, the group $\Hom_{D^b(\LCAb)}(G,H[q])$ coincides with the group of Yoneda extensions $\Ext^q_{\LCAb}(G,H)$ (where one replaces the notion of exact sequences with that of strict exact sequences).

\begin{proposition}[topology of W0]
    Let $S$ be a fixed set of primes.
    \begin{enumerate}
        \item The functor $(-)_\cond$ induces a fully faithful exact functor of triangulated categories $D^b(\LCAb)\to D(*_\proet)$. In particular, it makes $\LCAb$ into a full subcategory of $\Sh(*_\proet)$ closed under extensions.
        \item The Yoneda functor induces fully faithful exact functors of abelian and triangulated categories :
        $$Y:\IndProAb_S\to\Sh(*_\proet),\qquad RY:D^b(\IndProAb_S)\to D(*_\proet)$$
        identifying $\IndProAb_S$ as a full abelian subcategory of $\Sh(*_\proet)$ closed under extensions and $D^b(\IndProAb_S)$ as a full triangulated subcategory of $D(*_\proet)$.
        \item For $G\in\W_S$, $G_\top$ is locally compact Hausdorff and $(G_\top)_\cond=Y(G)$ in $\Sh(*_\proet)$.
        \item The functor $(-)_\top:\W_S\to\LCAb$ is a fully faithful exact functor of exact categories which identifies $\W_S$ as a full subcategory closed under extensions.
    \end{enumerate}
\end{proposition}

The topological realization of a general $G\in\IndProAb_\fin$ may not be locally compact Hausdorff. For example, the topological realization of $(\prod_iF_i)/(\bigoplus_iF_i)$ for any infinite family of nontrivial finite groups $F_i$, coincides with the quotient topology of $\prod_iF_i$ by the nonclosed subgroup $\bigoplus_iF_i$ hence is not Hausdorff. This example also shows $\W_S$ is not stable under quotients.

\begin{myproof}[of \Proposition{topology of W0}]
    \textbf{1.} This is \cite[Cor. 4.9]{ScholzeClausen}.
    
    \textbf{2.} \textbf{Exactness of $Y$.} Because filtered colimits are exact, it suffices to prove that the restricted functor $Y:\ProAb_S\to\Sh(*_\proet)$ is exact. Left-exactness is generally true of the Yoneda functor. For right-exactness, if $G\in\ProAb_S$, then $Y(G)\in\Sh(*_\proet)$ is precisely the sheaf representable by the pro-finite group $G_\top$, so by \Lemma{surjection in Pro}, a surjection $G\to H$ in $\ProAb_S$ gives a continuous surjection $Y(G)\to Y(H)$.
    
    \textbf{Full faithfulness of $RY$.} By \cite[Prop. 15.3.2]{KashiwaraSchapira}, the functor $RY:D^b(\IndProAb_S)\to\Sh(*_\proet)$ is well defined. By \cite[Prop. 2.2.3]{SuzukiGN} applied to the functor $\Hom_{\Sh(*_\proet)}(Y(-),Y(-))$, for $A=\varinjlim_i\varprojlim_jA_{ij}$ and $B=\varinjlim_k\varprojlim_lB_{kl}$ with $A_{ij},B_{kl}\in\Ab_S$, we get a sequence of morphisms and isomorphisms :
    \begin{align*}
        &R\varprojlim_i\varinjlim_kR\varprojlim_l\varinjlim_jR\Hom_{*_\proet}(Y(A_{ij}),Y(B_{kl}))\\
        \to& R\varprojlim_i\varinjlim_kR\varprojlim_lR\Hom_{*_\proet}(\varprojlim_j Y(A_{ij}),Y(B_{kl}))\\
        =&R\varprojlim_i\varinjlim_kR\Hom_{*_\proet}(\varprojlim_j Y(A_{ij}),R\varprojlim_l Y(B_{kl}))\\
        \to& R\varprojlim_iR\Hom_{*_\proet}(\varprojlim_j Y(A_{ij}),\varinjlim_kR\varprojlim_l Y(B_{kl}))\\
        =&R\Hom_{*_\proet}(\varinjlim_i\varprojlim_j Y(A_{ij}),\varinjlim_kR\varprojlim_l Y(B_{kl}))
    \end{align*}
    As seen in the previous point, $R\varprojlim_lY(B_{kl})=\varprojlim_lY(B_{kl})$, so the last term is simply $R\Hom_{*_\proet}(Y(A),Y(B))$. On the other hand, we have : $$R\Hom_{*_\proet}(Y(A_{ij}),Y(B_{kl}))=RY(R\Hom_{\Ab_S}(A_{ij},B_{kl}))$$ for $A_{ij},B_{kl}\in\Ab_S$ as a very particular case of \textbf{(1.)} since $R\Hom_\LCAb(A_{ij},B_{kl})$ is the complex of discrete groups $R\Hom_\Ab(A_{ij},B_{kl})$ for the discrete abelian groups $A_{ij}$ and $B_{kl}$, by \cite[Rem. 4.17]{HoffmannSpitzweck}. Hence the first term in the sequence is $R\Hom_{\IndProAb_S}(A,B)$.
    
    Hence it remains to show that we have isomorphisms for $A_{ij},B_{kl}\in\Ab_S$ :
    \begin{align*}
        \varinjlim_kR\Hom_{*_\proet}(\varprojlim_j Y(A_{ij}),R\varprojlim_l Y(B_{kl}))&\to R\Hom_{*_\proet}(\varprojlim_j Y(A_{ij}),\varinjlim_kR\varprojlim_l Y(B_{kl})),\\
        \varinjlim_jR\Hom_{*_\proet}(Y(A_{ij}),Y(B_{kl}))&\to R\Hom_{*_\proet}(\varprojlim_j Y(A_{ij}),Y(B_{kl})).
    \end{align*}
    Since $R\varprojlim_l Y(B_{kl})=\varprojlim_l Y(B_{kl})$ by \textbf{(1.)}, it is enough to show for $X=\varprojlim_jX_j$ any cofiltered limit in $\ProAb_S$ and $Z=\varinjlim_iZ_i$ any filtered colimit in $\Sh(*_\proet)$, the equality :
    $$\varinjlim_i\varinjlim_jR\HOM_{*_\proet}(Y(X_j),Z_i)=R\HOM_{*_\proet}(Y(X),Z)$$
    For $G$ an abelian group, recall $M_\bullet(G)\to G$ the Deligne-Scholze resolution \cite[Th. 4.5]{ScholzeClausen} : this is a resolution concentrated in nonnegative (homological) degree, such that $M_n(G)=\bigoplus_{i=1}^{r_n}\ZZ[G^{m_{i,n}}]$ for some integers $r_n,m_{i,n}\geq 0$, where $\ZZ[G^m]$ means the free abelian group generated by the set $G^m$. The resolution $M_\bullet(G)\to G$ and isomorphisms $M_n(G)=\bigoplus_{i=1}^{r_n}\ZZ[G^{m_{i,n}}]$ are functorial in $G$ : thus for any sheaf $F\in\Sh(*_\proet)$, we similarly have a resolution $M_\bullet(F)\to F$ with $M_n(F)\to\bigoplus_{i=1}^{r_n}\ZZ[F^{m_{i,n}}]$, where $M_n(F)$ and $\ZZ[F^m]$ mean the sheafifications of the presheaves $U\mapsto M_n(F(U))$ and $U\mapsto\ZZ[F(U)]$. We have a functorial spectral sequence and natural identifications on $\ProAb_S\times\Sh(*_\proet)$ :
    \begin{align*}
        E^{pq}_1(-,-)=&\Ext^p_{*_\proet}(M_q(Y(-)),-)\Rightarrow \Ext^{p+q}_{*_\proet}(M_\bullet(Y(-)),-)=\Ext^{p+q}_{*_\proet}(Y(-),-),\\
        E^{pq}_1(-,-)=&\bigoplus_{i=1}^{r_q}\Ext^p_{*_\proet}(\ZZ[Y(-)^{m_{i,q}}],-)=\bigoplus_{i=1}^{r_q}H^q((-)^{m_{i,q}},-)
    \end{align*}
    and the claim reduces to the equalities $\varinjlim_i\varinjlim_j H^q(X_j^m,Z_i)=H^q(X^m,Z)$ of \cite[Lem. 0739]{stacks-project}, for $m\geq 0$.
    
    \textbf{Full faithfulness of $Y$ and stability under extensions.} These follow from the full faithfulness of $RY$ and exactness of $Y$, which imply $\Hom_{*_\proet}(Y(A),Y(B))=\Hom_{\IndProAb_S}(A,B)$ and $\Ext^1_{*_\proet}(Y(A),Y(B))=\Ext^1_{\IndProAb_S}(A,B)$ for $A,B\in\IndProAb_S$.
    
    \textbf{3.} Consider an exact sequence $0\to G'\to G\to G''\to 0$ such that $G'=\varprojlim_mG_m'$ and $G''=\varinjlim_nG''_n$ with $G'_m,G''_n\in\Ab_S$ indexed by $\NN$. We will show that the sequence $0\to G'_\top\to G_\top\to G''_\top\to 0$ is topologically exact. Then $G'_\top$ is a pro-finite open subgroup of $G_\top$ because $G''_\top$ is discrete, so $G_\top$ is locally pro-finite (hence locally compact Hausdorff). 
    
    Define $G_n$ and $G_{n,m}$ by the following commutative diagram with exact rows, cartesian top-right square and cocartesian bottom-left square :
    $$\begin{tikzcd}[row sep=10]
        0 \arrow[r] & G' \arrow[r]\arrow[d,equals] & G \arrow[r] & G'' \arrow[r]& 0 \\
        0 \arrow[r] & G' \arrow[r]\arrow[d] & G_n \arrow[r]\arrow[d]\arrow[u] & G''_n \arrow[r]\arrow[d,equals]\arrow[u]& 0 \\
        0 \arrow[r] & G'_m \arrow[r] & G_{n,m} \arrow[r] & G''_n \arrow[r]& 0
    \end{tikzcd}$$
    where $G_n=\varprojlim_mG_{n,m}$ and $G=\varinjlim_nG_n$ in $\IndProAb_S$. By definition $(G'_m)_\top$, $(G_{n,m})_\top$ and $(G''_n)_\top$ are finite discrete, so the sequences $0 \to (G'_m)_\top \to (G_{n,m})_\top \to (G''_n)_\top \to 0$ are trivially topologically exact.
    
    Let $N_m'=\ker(G'\to G_m')$ and $N_{n,m}=\ker(G_n\to G_{n,m})$. The snake lemma gives an isomorphism $N_m'=N_{n,m}$. By left exactness of limits, the induced map $u_n:G'_\top\to (G_n)_\top$ is injective and we have an identification of subgroups (without topology) :
    $$\ker(G'_\top\to(G'_m)_\top)=\varprojlim_{m\leq m'}\ker((G'_{m'})_\top\to(G'_m)_\top)=\ker(G'\to G'_m)_\top=(N_m')_\top$$
    and similarly $\ker((G_n)_\to(G_{n,m})_\top)=(N_{n,m})_\top$. Thus by definition of the pro-finite topologies on $(G_n)_\top$ and $G'_\top$, the groups $(N_m')_\top=(N_{n,m})_\top$ for varying $m$ form a basis of open neighborhoods of $0$ in $G'_\top$ and in $(G_n)_\top$ simultaneously. Since $u_n$ is a group homomorphism, this implies $u_n$ is open. As an open injection, it is an open embedding, thus $(G_n)_\top/G'_\top$ is discrete, hence coincides with $(G''_n)_\top$ as a topological group. Thus the sequences $0\to G'_\top\to(G_n)_\top\to (G_n'')_\top\to 0$ are topologically exact.
    
    By exactness of filtered colimits of abelian groups, the sequence $0\to G'_\top\to G_\top\to G''_\top\to 0$ is exact in $\Ab$ (without topologies). Let $u:G'_\top\to G_\top$ be the corresponding injection. By definition of the colimit topology (using here that $G=\varinjlim_nG_n$ is a countable filtered colimit, see also \Remark{topological subtleties}), $u$ is an open embedding if and only if $u_n:G'_\top\to (G_n)_\top$ is an open embedding for each $n$, which we have just seen. In particular $G_\top/G'_\top$ is discrete.  The group $G''_\top$ is discrete, thus the isomorphism of abelian groups $G_\top/G'_\top\cong G''_\top$ is a homeomorphism, and this sequence is topologically exact. This concludes.
    
    Now we show $(G_\top)_\cond = Y(G)$. Consider $G=\varinjlim_nG_n$ as above. For $X$ a pro-finite set, we have by definition :
    $$(G_\top)_\cond(X) = \Hom(X,G_\top),\qquad Y(G)(X) = \varinjlim_n\Hom(X,(G_n)_\top)$$
    where the Hom's here denote continuous maps. We can assume the transitions in $G''_\top=\varinjlim_n(G_n'')_\top$ (hence in $G_\top=\varinjlim_n(G_n)_\top$) are injective, because the $G_n''$'s are finite. Then the canonical map :
    $$\varinjlim_n\Hom(X,(G_n)_\top)\to \Hom(X,G_\top)$$
    is injective. If $f:X\to G_\top$ is a continuous map then the composite $X\to G_\top\to G''_\top$ has compact and discrete, hence finite image, hence factors through one of the subgroups $(G''_n)_\top\leq G''_\top$. Thus $f$ factors through some $(G_n)_\top$. Therefore the canonical map above is bijective, and $(G_\top)_\cond(X) = Y(G)(X)$.
    
    \textbf{4.} The fully faithful exactness follows from \textbf{(1.)}, \textbf{(2.)} and \textbf{(3.)}. Stability under extensions similarly reduces to $\W_S$ being closed under extensions in $\IndProAb_S$, which is \Lemma{W0 closed under ext}.
\end{myproof}


\subsection{Duality for ind-pro-finite groups}
\label{Duality for ind-pro-finite groups}

\begin{definition}[derived good indpro]
    Let $S$ be a fixed set of primes. Let $\C$ be one of the categories $*_\proet$ or $\IndProAb_S$.
    \begin{itemize} 
        \item We write $D^b_{\W_S}(\C)\subseteq D^b(\C)$ be the full sucategory of bounded objects with cohomology in (or representable by objects of) $\W_S$.
        \item We write $\langle\W_S\rangle_{\C}\subseteq D(\C)$ the smallest full triangulated sucategories which contains (objects representably by) objects of $\W_S$ in degree $0$.
    \end{itemize}
\end{definition}

\begin{proposition}[indproab over proet]
    Let $S$ be a fixed set of primes.
    \begin{enumerate}
        \item The functor $\Gamma(*_\proet,-):\Sh(*_\proet)\to\Ab$ commutes with all limits and colimits. For $G\in\TopAb$, $\Gamma(*,G_\cond)$ is the abelian group underlying $G$ (without topology).
        \item The derived Yoneda functor $RY:D^b(\IndProAb_S)\to D(*_\proet)$ induces exact equivalences :
        $$D^b_{\W_S}(\IndProAb_S)\xrightarrow\sim D^b_{\W_S}(*_\proet),\qquad \langle\W_S\rangle_{\IndProAb_S}\xrightarrow\sim\langle\W_S\rangle_{*_\proet}.$$
        \item For objects $A=\varinjlim_i\varprojlim_jA_{ij}$ and $B=\varinjlim_k\varprojlim_lB_{kl}$ of $\IndProAb_S$ with $A_{ij},B_{kl}\in\Ab_S$, we have :
        $$R\HOM_{*_\proet}(A,B)=R\varprojlim_i\varinjlim_k R\varprojlim_j\varinjlim_l R\HOM_{*_\proet}(A_{ij},B_{kl})$$
        where $R\HOM_{*_\proet}(A_{ij},B_{kl})$ is representable by $R\Hom_\Ab(A_{ij},B_{kl})\in D^b(\Ab_S)\subseteq D^b(\IndProAb_S)$.
        \item Let $G\in D^b_{\W_S}(*_\proet)$ and $H=R\HOM_{*_\proet}(G,\QQ_p/\ZZ_p)$. Then $H\in D^b_{\W_S}(*_\proet)$ and for $q\in\ZZ$ :
        $$H^q(H) = \HOM_{*_\proet}(H^{-q}(G),\QQ/\ZZ)$$
        making $H^q(H)_\top$ and $H^{-q}(G)_\top$ perfect Pontryagin dual locally compact Hausdorff groups.
    \end{enumerate}
\end{proposition}

\begin{myproof}[of \Proposition{indproab over proet}]
    \textbf{1.} This is \cite[proof of Th. 2.2]{ScholzeClausen} ; in fact $\Gamma(S,-)$ commutes with all limits and colimits for any extremely disconnected profinite set $S$.
    
    \textbf{2.} These follow from \Proposition{topology of W0}.2.
    
    \textbf{3.} This reduces to the analogous formula for $R\Hom_{*_\proet/U}((-)|_U,(-)|_U)$, for each $U\in *_\proet$, for which the proof is identical to the proof of full faithfulness of $RY$ in \Proposition{topology of W0}.2, and the identification of $R\HOM_{\LCAb}(A,B)$ with the complex of discrete groups $R\Hom_{\Ab}(A,B)$ for $A,B\in\Ab_S$.
    
    \textbf{4.} If $G=\varprojlim_iG_i\in\ProAb_S$ with $G_i\in\Ab_S$, then by \textbf{(3.)} and exactness of $\Hom_\Ab(-,\QQ/\ZZ)$ we have :
    $$R\HOM_{*_\proet}(Y(G),\QQ/\ZZ)=\varinjlim_iR\HOM_{*_\proet}(Y(G_i),\QQ/\ZZ)=\varinjlim_i Y(\HOM_\Ab(G_i,\QQ/\ZZ))= Y(H)$$
    where $H=\varinjlim_i\Hom_{\Ab_S}(G_i,\QQ/\ZZ)$ belongs to $\IndAb_S$, and clearly $H_\top$ is the Pontryagin dual of $G_\top$. We conclude similarly in the case $G\in\IndAb_S$ (using the fact that $R\varprojlim_iY(H_i)=\varprojlim_iY(H_i)$ in $D(*_\proet)$ for $H_i\in\Ab_S$).  The general case $G\in\W_S$ follows immediately.
    
    Now for $G\in D^b_{\W_S}(*_\proet)$, letting $H=R\HOM_{*_\proet}(G,\QQ/\ZZ)$, the spectral sequence :
    $$E^{ij}_2=\EXT^i_{*_\proet}(H^{-j}(G),\QQ/\ZZ	)\Rightarrow H^{i+j}(H)$$
    degenerates at the first page ($E^{ij}_2=0$ for $i\neq 0$ by the case $G\in\W_S$), giving the desired isomorphism. This proves $H\in D^b_{\W_S}(*_\proet)$, and the Pontryagin duality statement follows.
\end{myproof}


\subsection{Cohomology of a finite field}
\label{Cohomology of a finite field}

\begin{situation}[finite field]
    Let $k$ be a finite field of characteristic $p$. Consider $\psi_{k/*}:\Spec k^\indrat_\proet\to *_\proet$ the premorphism given on the underlying categories by :
    $$\psi_{k/*}^{-1}(S)=\varprojlim_i\left(\bigsqcup_{s\in S_i}\Spec k\right)\in \Spec k^\indrat_\proet$$
    for $S=\varprojlim S_i$ a pro-finite set with finite $S_i$. We write $\Psi_{k/*}:\Sh(k_\proet^\indrat)\to\Sh(*_\proet)$ its pushforward.
\end{situation}

In \cite{SuzukiCFT}, the functor $R\Psi_{k/*}$ is simply written $R\Gamma(k,-)$. This notation is justified by the equality $R^q\Psi_{k/*}(G)(*)=H^q(k_\proet,G)$
for $G\in D^b(k_\proet^\indrat)$ and $q\in\ZZ$, due to the exactness of $R\Gamma(*_\proet,-)$. Hence $R^q\Psi_{k/*}(G)$ can truly be thought of as the group $H^q(k_\proet,G)$ enhanced with the structure of a condensed group. Nevertheless we prefer a notation that does not conflict with the (non-condensed) cohomology groups, and reminds one of the functors $R\Psi_{L/l}$ defined in the next part, though in the case of \Situation{finite field} the functor $R\Psi_{k/*}$ is unrelated to nearby cycles.

\begin{lemma}[pushforward of sites]
    Let $f:X\to Y$ be a premorphism of sites. For $G\in\Sh(X)$ and $q\in\ZZ$, the sheaf $R^qf_*(G)$ is the sheafification of the presheaf $U\in Y\mapsto H^q(f^{-1}(U),G)$.
\end{lemma}

\begin{myproof}[of \Lemma{pushforward of sites}]
    Let $a:\PSh(Y)\to\Sh(Y)$ be the sheafification functor and $F:\Sh(X)\to\PSh(Y)$ the functor given by $F(G)(U)=G(f^{-1}(U))$ for $U\in Y$ and $G\in\Sh(X)$. By exactness of $\Gamma(U,-)$ on $\PSh(Y)$ :
    $$\Gamma(U,-)\circ R^qF = R^q[\Gamma(U,-)\circ F]=H^q(f^{-1}(U),-).$$
    Thus $aR^qF(G)$ is sheafification of $U\mapsto H^q(f^{-1}(U),G)$, while on the other hand $aF(G)$ is $f_*(G)$ by definition of the pushforward. The equality $R^q(aF)=aR^qF$ for $q\geq 0$ follows from the case $q=0$ by exactness of $a$, proving the claim.
\end{myproof}

\begin{proposition}[topological cohomology]
    Consider $k$ and $p$ as in \Situation{finite field}.
    \begin{enumerate}
        \item The functor $R\Psi_{k/*}$ sends $D^b_{\W_k}(k_\proet^\indrat)$ to $D^b_{\W_p}(*_\proet)$. In particular for $G\in D^b_{\W_k}(k_\proet^\indrat)$ and $q\in\ZZ$, $R^q\Psi_{k/*}(G)$ is representable by a locally compact Haudorff topological abelian group.
        \item For $G\in D^b_{\W_k}(k_\proet^\indrat)$, the object $R\Psi_{k/*}(G)\in\W_p$ is computed as follows.
        \begin{enumerate}
            \item If $G\in\Alg_u^\RP(k)$ then for all $q\in\ZZ$ the object $R^q\Psi_{k/*}(G)$ is representable by the discrete finite group $H^q(k_\et,G)$. In particular $R^q\Psi_{k/*}(G)=0$ for $q\neq 0,1$.
            \item If $G\in\W_k$ is connected then $R\Psi_{k/*}(G)$ is concentrated in degree $0$ and $R^0\Psi_{k/*}(G)\in\W_p$. The functor $R^0\Psi_{k/*}$ sends type 1 and type 2 systems to limits and colimits in $\Sh(*_\proet)$, respectively.
            \item If $G\in\W_k$ then $R^q\Psi_{k/*}(G)=0$ for $q\neq 0,1$, we have $R^1\Psi_{k/*}(G)=R^1\Psi_{k/*}(\pi_0(G))$ and this object belongs to $\Ab_p$, and $R^0\Psi_{k/*}(G)$ belongs to $\W_p$.
            \item For general $G\in D^b_{\W_k}(k_\proet^\indrat)$ and $q\in\ZZ$, there is an exact sequence :
            $$0\to R^1\Psi_{k/*}(\pi_0(H^{q-1}(G)))\to R^q\Psi_{k/*}(G)\to R^0\Psi_{k/*}(H^q(G))\to 0.$$
        \end{enumerate}
        \item There is a canonical trace map $\tr:R\Psi_{k/*}(\QQ_p/\ZZ_p)\to \QQ_p/\ZZ_p[-1]$.
        \item For $G\in D^b_{\W_k}(k_\proet^\indrat)$ and $G^\vee=R\HOM_{k_\proet^\indrat}(G,\QQ_p/\ZZ_p)$, we have a perfect pairing in $D(*_\proet)$ :
        $$R\Psi_{k/*}(G)\otimes^LR\Psi_{k/*}(G^\vee)\to R\Psi_{k/*}(\QQ_p/\ZZ_p)\xrightarrow{\tr}\QQ_p/\ZZ_p[-1]$$
        which induces perfect Pontryagin duality of locally compact Hausdorff abelian groups for $q\in\ZZ$ :
        $$R^q\Psi_{k/*}(G)\otimes R^{1-q}\Psi_{k/*}(G^\vee)\to R^1\Psi_{k/*}(\QQ_p/\ZZ_p)\to \QQ_p/\ZZ_p.$$
    \end{enumerate}
\end{proposition}

\begin{myproof}[of \Proposition{topological cohomology}]
    \textbf{1.} This will follow from the computations of \textbf{(2.)} and \Proposition{topology of W0}.3.
    
    \textbf{2.a.} Consider $G_0\in\Alg_u(k)$ and $G=G_0^\RP$. Clearly $H^0(k_\et,G)=G_0(k)^S$ is finite. We have $H^1(k_\et,G)=H^1(k_\et,\pi_0(G_0))$ by \cite[Prop. 8.9]{BertapelleSuzuki} and \cite[Cor. 7.10 of Sec. III.5]{DemazureGabriel}, and this group is finite as a subquotient of the finite group $\pi_0(G_0)(\Adh k)$ by \cite[Prop 1 of Ch. XIII]{SerreLC}. For $q\geq 2$, $H^q(k_\et,G)=H^q(k_\et,G)$ is trivial because $G$ is $p$-primary torsion. Hence $H^q(k_\et,G)$ is finite for $q\in\ZZ$. By \Proposition{Yoneda for Wk}.4 the group $H^q(k_\proet,G)=H^q(k_\et,G)$ is finite. For $S=\varprojlim_iS_i$ a profinite set with $S_i$ finite discrete sets, we have : $$H^q(\psi_{k/*}^{-1}(S),G)=\varinjlim_i H^q(k_\proet,G)^{S_i}=\varinjlim_i \Hom(S_i,H^q(k_\proet,G))=\Hom(S,H^q(k_\proet,G))$$
    where the Hom's here denotes continuous maps of topological spaces and the finite set $H^q(k_\proet,G)$ has the discrete topology. By \Lemma{pushforward of sites}, $R^q\Psi_{k/*}(G)$ is representable by the finite discrete group $H^q(k_\proet,G)$.
        
    \textbf{2.b.} By \cite[Prop. 2.2.4]{SuzukiGN} for $G=\varinjlim_i\varprojlim_jG_{ij}\in\IndProAlg_u^\RP(k)$ with $G_{ij}\in\Alg^\RP_u(k)$ we have :
    $$R\Psi_{k/*}(G)=\varinjlim_iR\varprojlim_j R\Psi_{k/*}(G_{ij}).$$
    If $G=\varinjlim_nG_n$ is type 2 then for all $q\in\ZZ$ we have $R^q\Psi_{k/*}(G)=\varinjlim_nR^q\Psi_{k/*}(G_n)$, and the cancellation for $q\geq 1$ follows from \textbf{(2.a)}. If $G=\varprojlim_nG_n$ is type 1 then we have a spectral sequence :
    $$E^{ij}_2=R^i\varprojlim_iR^j\Psi_{k/*}(G_i)\Rightarrow R^{i+j}\Psi_{k/*}(G).$$
    By \textbf{(2.a)} we have $R^j\Psi_{k/*}(G_i)=0$ if $j\geq 1$. In particular the transitions $R^0\Psi_{k/*}(G_{i+1})\to R^0\Psi_{k/*}(G_i)$ are surjective, so by \cite[Lem. 0CQA and 07KW (1)]{stacks-project}, $E^{i0}_2=0$ for $i\geq 1$. This implies $R^q\Psi_{k/*}(G)=0$ for $q\geq 1$, and $R^0\Psi_{k/*}(G)=\varprojlim_iR^0\Psi_{k/*}(G_i)$, as desired.
    
    \textbf{2.c} The long exact sequence for $R\Psi_{k/*}$ associated to $0\to G^0\to G\to \pi_0(G)\to 0$, combined with the cancellations $R^q\Psi_{k/*}(\pi_0(G))=0$ for $q\geq 2$ and $R^q\Psi_{k/*}(G^0)=0$ for $q\geq 1$ given by \textbf{(2.a-b)}, give an isomorphism $R^1\Psi_{k/*}(G)=R^1\Psi_{k/*}(\pi_0(G))$ and an exact sequence :
    $$0\to R^0\Psi_{k/*}(G^0)\to R^0\Psi_{k/*}(G)\to R^0\Psi_{k/*}(\pi_0(G))\to 0.$$
    We conclude by stability of $\W_p$ under extensions in $\Sh(*_\proet)$ (\Lemma{W0 closed under ext} and \Proposition{topology of W0}.2).
    
    \textbf{2.d} We have a spectral sequence :
    $$E^{ij}_2=R^i\Psi_{k/*}(H^j(G))\Rightarrow R^{i+j}\Psi_{k/*}(G)$$
    where $H^j(G)\in\W_k$ for all $j\in\ZZ$. By \textbf{(2.c)}, $E^{ij}_2=0$ for all $j\geq 2$ so this spectral sequence degenerates to short exact sequences $0\to R^1\Psi_{k/*}(H^{q-1}(G))\to R^q\Psi_{k/*}(G)\to R^0\Psi_{k/*}(H^q(G))\to 0$ where $R^1\Psi_{k/*}(H^{q-1}(G))=R^1\Psi_{k/*}(\pi_0(H^{q-1}(G)))$ by \textbf{(2.c)}.
    
    \textbf{3.} The isomorphism $H^1(k_\et,\ZZ/p^n)=\ZZ/p^n$ for $k$ finite is well known, given by identifying the Galois group of $k$ with $\hat\ZZ$ and $H^1(k_\et,\ZZ/p^n)$ with the set of continuous, $p^n$-torsion characters of this Galois group. The result for $R^1\Psi_{k/*}(\QQ_p/\ZZ_p)=\varinjlim_nR^1\Psi_{k/*}(\ZZ/p^n)=\varinjlim H^1(k_\et,\ZZ/p^n)$ follows. The morphism $R\Psi_{k/*}(\QQ_p/\ZZ_p)\to\QQ_p/\ZZ_p[-1]$ then comes from the fact that $R^q\Psi_{k/*}(\QQ_p/\ZZ_p)=0$ for $q>1$.
    
    \textbf{4.} See \cite[Prop. 12.2]{SuzukiCFT} and \Proposition{indproab over proet}.4.
\end{myproof}


\subsection{Finite coefficients}
\label{Finite coefficients}

\begin{definition}[finite coeffs]
    Consider $X$ a scheme and $S_X$ any of the following sites :
    \begin{enumerate}
        \item $S_X=X_\et$ with no assumption ;
        \item $S_X=X_\Et$ if $X=\Spec k$ with $k$ a field ;
        \item $S_X=X_\RP,X_\RPS$ if $X=\Spec k$ with $k$ a field of characteristic $p>0$ such that $[k:k^p]<\infty$ ;
        \item $S_X=X_\et^\perar,X_\proet^\indrat$ if $X=\Spec k$ with $k$ a perfect field of characteristic $p>0$.
    \end{enumerate}
    Define $\Sh_\fin(S_X)$ (resp. $\Sh_\ell(S_X)$ for $\ell$ a prime) as the full subcategory of $\Sh(S_X)$ of sheaves representable by finite étale (resp. and $\ell$-primary) $X$-group schemes. Define $D_\fin(S_X)$ (resp. $D_\ell(S_X)$) as the full subcategory of $D(S_X)$ of bounded objects with cohomology in $\Sh_\fin(S_X)$ (resp. in $\Sh_\ell(S_X)$).
\end{definition}

It is clear that $\Sh_\fin(S_X)=\bigoplus_\ell\Sh_\ell(S_X)$ and $D_\fin(S_X)=\bigoplus_\ell D_\ell(S_X)$. The following lemma indicates that these categories contain the same data regardless of the choice of site.

\begin{lemma}[finite coeffs equivalence]
    Consider $X$ and $S_X$ as in any case of \Definition{finite coeffs} and $\ell$ a prime.
    \begin{enumerate}
        \item $\Sh_\fin(S_X)$ (resp. $\Sh_\ell(S_X)$) is an abelian subcategory of $\Sh(S_X)$ closed under extension and $D_\fin(S_X)$ (resp. $D_\ell(S_X)$) is the full triangulated subcategory of $D(S_X)$ generated by objects of $\Sh_\fin(S_X)$ (resp. $\Sh_\ell(S_X)$) viewed as complexes concentrated in degree $0$.
        \item Consider $u:S_X\to X_\et$ the morphism of sites defined by identity. Then its pullback $u^*$ defines exact equivalences of abelian and triangulated categories :
        $$\Sh_\fin(X_\et)\xrightarrow\sim\Sh_\fin(S_X),\quad \Sh_\ell(X_\et)\xrightarrow\sim\Sh_\ell(S_X),\quad D_\fin(X_\et)\xrightarrow\sim D_\fin(S_X),\quad D_\ell(X_\et)\xrightarrow\sim D_\ell(S_X)$$
        which commute with the Yoneda functor from the category of finite ($\ell$-primary) étale $X$-groups.
    \end{enumerate}
\end{lemma}

\begin{myproof}[of \Lemma{finite coeffs equivalence}]
    \textbf{1.} First observe that a sheaf $G\in\Sh(S_X)$ is representable by a finite étale $k$-group if and only if it is finite locally constant, \textit{i.e.} there exists a covering $\{X_i\}$ of $X$ such that each $G|_{S_X/X_i}$ is a constant sheaf associated to a finite group. For $S_X=X_\et,X_\Et,X_\RP,X_\RPS,X_\et^\perar,X^\indrat_\et$ the proof is identical to \cite[Lem. 03RV]{stacks-project}. The case $S_X=X^\indrat_\proet$ reduces to the case $S_X=X_\et^\indrat$ by an argument identical to \cite[Lem. 099Y]{stacks-project}.
    
    Clearly finite (resp. $\ell$-primary) locally constant sheaves form an abelian subcategory of $\Sh(S_X)$ closed under extensions. It follows that $D_\fin(S_X)$ (resp. $D_\ell(S_X)$) is a triangulated subcategory of $D(S_X)$, which is clearly the smallest such subcategory containing objects of $\Sh_\fin(S_X)$ (resp. $\Sh_\ell(S_X)$) in degree $0$.
    
    \textbf{2.} Note that $u$ is a morphism of sites because $X_\et$ admits all finite products and $u^{-1}$ preserves them. By \Lemma{yoneda and pullback}, $u^*:\Sh(X_\et)\to\Sh(S_X)$ commutes with the Yoneda functors from the category of finite étale $X$-groups (since $u^{-1}:X_\et\to S_X$ is the inclusion). Thus $u^*$ restricts to $\Sh_\fin(X_\et)\to\Sh_\fin(S_X)$.
    
    The unit map $F\to u_*u^*(F)$ is an isomorphism for all $F\in\Sh(X_\et)$, because $u$ is defined by the inclusion $X_\et\to S_X$. It follows that $u^*:\Sh(X_\et)\to\Sh(S_X)$ is fully faithful. Its restriction $\Sh_\fin(X_\et)\to\Sh_\fin(S_X)$ is essentially surjective by compatibility with the Yoneda functors, thus an equivalence. The pushforward $u_*:\Sh(S_X)\to\Sh(X_\et)$ admits $u^*$ as left adjoint, thus it preserves injectives, and it is clearly exact. Thus by derivation of the composite functors $u_*u^*=\id$ and $\Hom_{X_\et}(F,u_*(-))=\Hom_{S_X}(u^*F,-)$, we have :
    $$R\Hom_{X_\et}(F,G)=R\Hom_{X_\et}(F,Ru_*u^*G)=R\Hom_{S_X}(u^*F,u^*G)$$
    for $F,G\in\Sh(X_\et)$. Hence $u^*:D_\fin(X_\et)\to D(S_X)$ is fully faithful. We conclude that it is an equivalence $D_\fin(X_\et)\cong D_\fin(S_X)$ by an argument identical to that in the proof of \Proposition{definition of D0}.
\end{myproof}

Recall that for $r\in\ZZ$ and $X$ a $\ZZ[1/m]$-scheme, with $m>0$ an integer, we define the Tate twist $\ZZ/m(r)\in D(S_X)$ as the sheaf $\mu_m^{\otimes r}$ concentrated in degree $0$ ; in particular $\ZZ/m(0)=\ZZ/m$ and $\ZZ/m(r)$ is an invertible $\ZZ/m$-module with inverse $\ZZ/m(-r)=\HOM_{S_X}(\ZZ/m,\ZZ/m(r))$. For $\ell$ a prime invertible on $X$, define $\QQ_\ell/\ZZ_\ell(r)=\varinjlim_n\ZZ/\ell^n(r)$. If $X$ is a $\QQ$-scheme, define $\QQ/\ZZ(r)=\bigoplus_\ell\QQ_\ell/\ZZ_\ell(r)=\varinjlim_m\ZZ/m(r)$.

\begin{proposition}[cartier duality]
    Consider $X$ and $S_X$ as in any case of \Definition{finite coeffs}. Let $\ell$ be a prime invertible on $X$ and $d$ a fixed integer. Consider the functor $(-)^\vee_{S_X}=R\HOM_{S_X}(-,\QQ_\ell/\ZZ_\ell(d))$ on $D(S_X)$.
    \begin{enumerate}
        \item For $G\in D_\ell(S_X)$, we have $H^q(G^\vee_{S_X})=\HOM_{S_X}(H^{-q}(G),\QQ_\ell/\ZZ_\ell(d))$ for $q\in\ZZ$, and $G^\vee_{S_X}\in D_\ell(S_X)$.
        \item For any $G\in D_\ell(S_X)$, the canonical map $G\to(G^\vee_{S_X})^\vee_{S_X}$ is an isomorphism in $D(S_X)$.
        \item Consider $u:S_X\to X_\et$ the premorphism of sites defined by identity. We have a commutative square :
        $$\begin{tikzcd}
           D_\ell(X_\et) \arrow[r,"(-)^\vee_{X_\et}"]\arrow[d,"u^*"] & D_\ell(X_\et) \arrow[d,"u^*"] \\
           D_\ell(S_X)^\opp \arrow[r,"(-)^\vee_{S_X}"] & D_\ell(S_X)^\opp
        \end{tikzcd}$$
        \item If $X$ is a $\QQ$-scheme, all of the above holds with $\Sh_\ell(S_X)$, $D_\ell(S_X)$, $\QQ_\ell/\ZZ_\ell(d)$ respectively replaced by $\Sh_\fin(S_X)$, $D_\fin(S_X)$, $\QQ/\ZZ(d)$.
    \end{enumerate}
    We say that $G\in\Sh_\ell(S_X)$ and $\HOM_{S_X}(G,\QQ_\ell/\ZZ_\ell(d))$ satisfy a perfect Cartier duality of finite étale groups.
\end{proposition}

\begin{myproof}[of \Proposition{cartier duality}]
    \textbf{1.} Because $\ell$ is invertible on $X$, the sheaf has injective stalks $\QQ_\ell/\ZZ_\ell$ at geometric points. Thus $\HOM_{S_X}(-,\QQ_\ell/\ZZ_\ell(d))$ is exact, giving the desired equality. Consider $G\in\Sh_\ell(S_X)$ of $\ell^n$-torsion, and an étale covering $\{U_i\}_i$ which trivializes $G$ and $\mu_{\ell^n}$. The sheaf $\HOM_{S_X}(G,\QQ_\ell/\ZZ_\ell(d))=\HOM_{S_X}(G,\ZZ/\ell^n(d))$ is trivialized by this étale covering and thus is locally constant. By the same argument as \Lemma{finite coeffs equivalence}.1, this sheaf belongs to $\Sh_\ell(S_X)$. Hence $(-)^\vee_{S_X}$ sends objects of $D_\ell(S_X)$ to $D_\ell(S_X)$.
    
    \textbf{2.} For any $n>0$ we have :
    $$\HOM_{S_X}(\HOM_{S_X}(\ZZ/\ell^n,\QQ_\ell/\ZZ_\ell(d)),\QQ_\ell/\ZZ_\ell(d))=\HOM_{S_X}(\ZZ/\ell^n(d),\QQ_\ell/\ZZ_\ell^n(d))=\ZZ/\ell^n$$
    meaning $G\to (G^\vee_{S_X})^\vee_{S_X}$ is an isomorphism when $G=\ZZ/\ell^n$. If $G\in\Sh_\ell(S_X)$ is general, then étale-locally we have $G\cong\bigoplus_{i=1}^r\ZZ/\ell^{n_r}$ for some integers $n_1,\dots,n_r\geq 1$ thus $G\to (G^\vee_{S_X})^\vee_{S_X}$ is again an isomorphism. The same follows for $G\in D_\ell(S_X)$ by \Lemma{finite coeffs equivalence}.1.
    
    \textbf{3.} For any $n>0$ we have :
    $$u^*\HOM_{X_\et}(\ZZ/\ell^n,\QQ_\ell/\ZZ_\ell(d))=u^*\ZZ/\ell^n(d) = \ZZ/\ell^n(d)=(\ZZ/\ell^n)^\vee_{S_X}=(u^*\ZZ/\ell^n)^\vee_{S_X}$$
    by \Lemma{finite coeffs equivalence}.2. Thus $u^*(G^\vee_{X_\et})=u^*(G)^\vee_{S_X}$ when $G=\ZZ/\ell^n$, and as above this generalizes to $G\in D_\fin(X_\et)$.
    
    \textbf{4.} This is clear from $\Sh_\fin(S_X)=\bigoplus_\ell\Sh_\ell(S_X)$, $D_\fin(S_X)=\bigoplus_\ell D_\ell(S_X)$ and $\QQ/\ZZ(d)=\bigoplus_\ell\QQ_\ell/\ZZ_\ell(d)$.
\end{myproof}

The following lemma will not be used in this section, but will appear in the next part.

\begin{lemma}[devissage]
    Let $K$ be a field and $S_K=S_{\Spec K}$ be one of the sites of \Definition{finite coeffs}. Let $\ell\neq\characteristic K$ be a prime. For any $G\in\Sh_\ell(S_K)$ there exists a finite extension $L/K$ such that the following holds, where $S_L$ is the analogous site and $\pi_{L/K}:S_L\to S_K$ the morphism induced by base change.
    \begin{enumerate}
        \item The degree $[L:K]$ is coprime to $\ell$.
        \item There exists a finite separated filtration of $G|_L$ in $\Sh_\ell(S_L)$ with quotients $0$ or $\ZZ/\ell$.
        \item The unit $G\to\pi_{L/K,*}(G|_L)$ is split monic and the counit $\pi_{L/K,*}(G|_L)\to G$ is split epic in $\Sh_\ell(S_L)$.
    \end{enumerate}
\end{lemma}

\begin{myproof}[of \Lemma{devissage}]
    We reproduce the method used in \cite[Prop. 9.6]{BertapelleSuzuki}. By \Lemma{finite coeffs equivalence}.2, $G$ can be viewed as a Galois module over $K$. Let $L'/K$ be a finite Galois extension which contains a primitive $\ell$-th root of unity, over which $G$ is trivial, and $L/K$ be a subextension such that $\Gal(L'/L)$ is a $\ell$-Sylow of $\Gal(L'/K)$.
    
    Then $\Gal(L'/L)$ is a finite $\ell$-group and $G|_L$ is a finite $\ell$-primary $\Gal(L'/L)$-module, so $G|_L$ either is zero or has a nonzero $\Gal(L'/L)$-fixed submodule, in which case it admits a submodule isomorphic to $\ZZ/\ell$ in $\Sh_\ell(S_L)$. Inductively (reasoning the same way on $G|_L/(\ZZ/\ell)$), because $G|_L$ is finite, it admits a finite separated $\Gal(L'/L)$-module filtration with nontrivial quotients isomorphic to $\ZZ/\ell$, proving \textbf{(1.)}.
    
    The composite of the unit and counit maps $G\xrightarrow\res\pi_{L/K,*}(G|_L)\xrightarrow\cores G$ is multiplication by $[L:K]$. By construction $[L:K]$ is coprime to $\ell$, thus $\res$ is a split monomorphism with section $[L:K]^{-1}\cores$, and $\cores$ is a split epimorphism with section $[L:K]^{-1}\res$, proving \textbf{(2.)}.
\end{myproof}

\begin{remark}[restriction and separability]
    The above \Lemma{devissage} appears in the proof of \cite[Prop. 9.6]{BertapelleSuzuki} for $p$-primary $G$, in the case where $K$ is a Henselian discrete valuation field with residue field of charcateristic $p$ (there our $L'$ and $L$ are written $L$ and $M$). In particular, the residue field extension $l/k$ has degree prime to $p$ hence is separable. Later in the same proof, the authors use the Weil restriction $\Res_{l/k}:\Sh(l_\RPS)\to\Sh(k_\RPS)$, \textit{i.e.} the pushforward of the premorphism defined by base change, and the fact that it maps $\Sh_0(l_\RPS)$ to $\Sh_0(k_\RPS)$. The well-definition of the premorphism $\Spec l_\RPS\to\Spec k_\RPS$, and stability of RPAU groups under $\Res_{l/k}$, are easily checked in the case where $l/k$ is separable, which covers the use case of \cite{BertapelleSuzuki} (see \cite[Cor. 1.9]{Kato86}, \cite[Prop. 2.5]{BertapelleSuzuki}, or \Lemma{base change RP}.1). Without observing that $l/k$ is separable, one needs to treat the case that $l/k$ is any finite extension ; this is \Lemma{base change RP}.4 and \Corollary{weil restriction RPAU}.
\end{remark}


\section{Duality for higher local fields}
\label{Duality for higher local fields}


\subsection{Positive equal characteristic step}
\label{Positive equal characteristic step}

Recall the definition of the Kato lift $h^{K/k}(T)$ from \Proposition{Kato lift complete} and \Definition{Kato lift definition}, for $K$ the fraction field of an integral, complete Noetherian local ring with residue field $k$ of positive characteristic and $T$ an affine, relatively perfect $k$-scheme.

\begin{situation}[positive equal char setup]
    Let $K$ be a complete discrete valuation field of characteristic $p>0$ with perfect residue field $k$. Set $S_K^p=\Spec K_\RP$ and $S_k^p=\Spec k_\et^\perar$. Define a premorphism of sites $\psi_{K/k}^p:S_K^p\to S_k^p$ by :
    $$(\psi_{K/k}^p)^{-1}:T/k\mapsto h^{K/k}(T).$$
    We write $\Psi_{K/k}^p:\Sh(S_K^p)\to\Sh(S_k^p)$ the pushforward of $\psi_{K/k}^p$.
\end{situation}

The objects $K,k,\psi_{K/k}^p$ and $\Psi_{K/k}^p$ correspond to $k,F,\pi_{k,\RP}$ and $\pi_{k,\RP,*}$ in the notations of \cite[Sec. 4]{SuzukiCFT}. The superscript $(-)^p$ stands for "characteristic $p$", and will be omitted when no confusion is possible.

For $F=K,k$ and $G\in\Sh(S_F)$ we write $G_{\Adh F}=\varinjlim_{E/F}\Gamma(E,G)$ the étale stalk at the separable closure $\Adh F$ (where $E/F$ runs over finite separable extensions). We view it as a Galois module over $F$, that is an object of $\Sh(F_\et)$ ; this defines an exact functor $(-)_{\Adh F}:\Sh(S_F)\to \Sh(F_\et)$. For either $F=K,k$ which are characteristic $p>0$, we write $\ZZ/p^n(r)=\nu_n(r)[-r]\in D(S_F)$ and $\QQ_p/\ZZ_p(r)=\nu_\infty(r)[-r]$ as in section \ref{Duality for RPAU groups}.

The following results are from \cite[Sec 4]{SuzukiCFT}. Recall Serre duality per \Theorem{duality for RPAU groups} and \Theorem{Serre duality}.

\begin{theorem}[equal char p]
    Consider \Situation{positive equal char setup}.
    \begin{enumerate}
        \item $R\Psi_{K/k}$ sends $D_0(S_K)$ to $D^b_{\W_k}(S_k)$. If $G\in\Sh_0(S_K)$ then $R^q\Psi_{K/k}(G)=0$ for $q\neq 0,1$.
        \item There exists a canonical trace map $\tr:R\Psi_{K/k}(\QQ_p/\ZZ_p(1))\to\QQ_p/\ZZ_p[-1]$ in $D(S_k)$.
        \item For any $G\in D_0(S_K)$ with Serre dual $G^\vee=R\HOM_{S_K}(G,\QQ_p/\ZZ_p(1))$, we have a perfect pairing :
        $$R\Psi_{K/k}(G)\otimes^LR\Psi_{K/k}(G^\vee)\to R\Psi_{K/k}(\QQ_p/\ZZ_p(1))\xrightarrow\tr\QQ_p/\ZZ_p[-1]$$
        which gives perfect Serre dualities of connected and finite étale parts for $q\in\ZZ$ :
        \begin{align*}
            R^q\Psi_{K/k}(G)^0&=\EXT^1_{S_k}(R^{1-q}\Psi_{K/k}(G^\vee)^0,\QQ_p/\ZZ_p)
            ,\\
            \pi_0(R^q\Psi_{K/k}(G))&=\HOM_{S_k}(\pi_0(R^{-q}\Psi_{K/k}(G^\vee)),\QQ_p/\ZZ_p).
        \end{align*}
        \item For $G\in\Sh_0(S_K)$ and $q\in\ZZ$, $R^q\Psi_{K/k}(G)$ is the sheafification of the presheaf $T/k\mapsto H^q(h^{K/k}(T),G)$. In particular we have $R\Psi_{K/k}(G)_{\Adh k} = R\Gamma(K^\ur,G_{\Adh K})$ in $D(k_\et)$ for $G\in D_0(S_K)$.
    \end{enumerate}
\end{theorem}

\begin{myproof}[of \Theorem{equal char p}]
    \textbf{1.} If $G\in\Sh_0(S_K)$ then we have $R^q\Psi_{K/k}(G)=0$ if $q\neq 0,1$ and $R^q\Psi_{K/k}(G)\in\W_k$ for any $q\in\ZZ$ by \cite{SuzukiCFT} (this is Prop. 4.1 for $q=0$ ; Prop. 4.3 for $q\neq 0,1$ ; Prop. 4.6 for $q=1$ ; note that this reference uses our \Lemma{stability of Wk}). For $G\in D_0(S_K)$, the spectral sequence :
    $$E^{ij}_2=R^i\Psi_{K/k}(H^j(G))\Rightarrow R^{i+j}\Psi_{K/k}(G)$$
    gives short exact sequences $0\to R^1\Psi_{K/k}(H^{q-1}(G))\to R^q\Psi_{K/k}(G)\to R^0\Psi_{K/k}(H^q(G))\to 0$. Hence $R^q\Psi_{K/k}(G)$ is trivial whenever $H^q(G)=H^{q-1}(G)=0$, which holds for all but finitely many $q\in\ZZ$. Also by \Lemma{stability of Wk} and the case $G\in\Sh_0(S_K)$, $R^q\Psi_{K/k}(G)$ belongs to $\W_k$. Thus $R\Psi_{K/k}(G)\in D^b_{\W_k}(S_k)$.
    
    \textbf{2.} By \cite[Prop. 6.1.6]{Suzuki2024}, we have $R^q\Psi_{K/k}(\ZZ/p^n(1))=0$ if $q>1$, and identifications for $n\geq 1$ :
    $$R^1\Psi_{K/k}(\ZZ/p^n(1))=R^0\Psi_{K/k}(\GG_{m,K}^\RP)/p^n$$
    The latter is the sheafification of the presheaf
    $\Spec F\in S_k\mapsto h^{K/k}(F)^\times /p^n$. \Proposition{Kato lift structure}.2 gives a valuation map $R^0\Psi_{K/k}(\GG_{m,k}^\RP)\to\ZZ$ which defines the morphism $R^1\Psi_{K/k}(\ZZ/p^n(1))\to\ZZ/p^n$ for $n\geq 1$.
    
    We have $\varinjlim_nR^q\Psi_{K/k}(\ZZ/p^n(1))=R^q\Psi_{K/k}(\QQ_p/\ZZ_p(1))$ for all $q\in\ZZ$, by \cite[Prop. 2.2.4.(b)]{SuzukiGN}. Hence we have $R^q\Psi_{K/k}(\QQ_p/\ZZ_p(1))=0$ for $q>1$ and a morphism $R^1\Psi_{K/k}(\QQ_p/\ZZ_p(1))\to\QQ_p/\ZZ_p$, which defines the morphism $\tr:R\Psi_{K/k}(\QQ_p/\ZZ_p(1))\to\QQ_p/\ZZ_p[-1]$ in $D(S_k)$.
    
    \textbf{3.} The duality in $D(S_k)$ is \cite[Prop. 4.4]{SuzukiCFT}. The rest follows by \textbf{(1.)} and \Theorem{Serre duality}.3.(d).
    
    \textbf{4.} The first statement is by \Lemma{pushforward of sites}. By exactness of $(-)_{\Adh k}$ and \Proposition{definition of D0}, it suffices to show $R^q\Psi_{K/k}(G)_{\Adh k}=H^q(K^\ur,G_{\Adh K})$ for $G=\GG_{a,K}^\RP$ and $q\in\ZZ$. For $q\neq 0$, both sides of the equality are then trivial by \cite[Prop. 8.9]{BertapelleSuzuki} and \cite[Prop. 3.7 of Ch. III]{Milne80}. For $q=0$ we have :
    \begin{align*}
        \Psi_{K/k}(G)_{\Adh k}
        &=\varinjlim_{l/k}G(h^{K/k}(l)) &\text{because stalks before and after sheafification coincide,}\\
        &=\varinjlim_{l/k}\GG_{a,k}(h^{K/k}(l)) &\text{by \cite[Prop. 8.9]{BertapelleSuzuki},}\\
        &=\GG_{a,k}(\varinjlim_{l/k} h^{K/k}(l))=\GG_{a,k}(K^\ur) &\text{because }\GG_{a,k}\text{ is of finite presentation over }k,\\
        &=G(K^\ur)=H^0(K^\ur,G_{\Adh K})&\text{by \cite[Prop. 8.9]{BertapelleSuzuki} again,}
    \end{align*}
    where the colimits runs over finite separable extensions of $k$. This concludes.
\end{myproof}

\begin{remark}[for indrat proet]
    By \Proposition{Yoneda for Wk}.3-4 and \Proposition{Serre duality}.3, the change of site $\alpha_k:D(S_k)\to D(S_k')$, where $S_k'=\Spec k_\proet^\indrat$, preserves representability by objects of $\W_k$, Serre duality, and cohomology for such objects. Thus all of \Theorem{positive equal char setup}.1-3 remains true if we replace $S_k$ by $S_k'$ and $R\Psi_{K/k}$ by $\alpha_kR\Psi_{K/k}$. \Theorem{positive equal char setup}.4 is partially true in that $(\alpha_kR\Psi_{K/k}(G))_{\Adh k}=R\Gamma(K^\ur,G_{\Adh K})$ for $G\in D_0(S_K)$.
\end{remark}

\begin{remark}[adh k-sections]
    Since $\Spec\Adh k$ has no nontrivial covering in $S_k$, the functor $\Gamma(\Adh k,-)$ is exact on $\Sh(S_k)$ and given any presheaf $P\in\PSh(S_k)$ with sheafification $\tilde P$, we have $P(\Adh k)=\tilde P(\Adh k)$. The completion $\widehat{K^\ur}$ of $K^\ur$ with respect to its discrete valuation is identified with $h^{K/k}(\Adh k)$, so \Theorem{equal char p}.4 implies the equality $R\Psi_{K/k}(G)(\Adh k)=R\Gamma(\widehat{K^\ur},G)$ for $G\in D_0(S_K)$. 
\end{remark}


\subsection{Mixed characteristic step}
\label{Mixed characteristic step}

\begin{situation}[mixed char setup]
    Let $K$ be a Henselian discrete valuation field with residue field $k$. Assume $K$ has mixed characteristic $(0,p)$ and $[k:k^p]=p^d$ for some finite $d$. Set $S_K^m=\Spec K_\Et$, $S_{\O_K}^m=\Spec \O_{K,\RPS}$ and $S_k^m=\Spec k_\RPS$. Consider the immersions $j:\Spec K\to\Spec\O_K$ and $i:\Spec k\to\Spec\O_K$ and write $j:S_K^m\to S_{\O_K}^m$ and $i:S_{\O_K}^m\to S_k^m$ the premorphisms defined by base change. We write $R\Psi_{K/k}^m$ the composite functor $i^*Rj_*:D(S_K^m)\to D(S_k^m)$.
\end{situation}

The superscript $(-)^m$ stands for "mixed characteristic", and as before will often be omitted from notation. This definition of $R\Psi_{K/k}^m$ connects it to the theory of nearby cycles, whence the notation stems. We abusively retained similar notation in other cases (even ones that have no clear relation to nearby cycles), to highlight the similarities of each setup.

As previously, we write $(-)_F:D(S_F)\to D(F_\et)$ the étale stalk functor for $F=K,k$. Mind that over $S_K$ where $\characteristic K=0$, $\ZZ/p^n(r)$ denotes the sheaf $\mu_{p^n}^{\otimes r}$ (as in section \ref{Finite coefficients}), but over $S_k$ where $\characteristic k=p$ we mean $\ZZ/p^n(r)=\nu_n(r)[-r]$ (as in section \ref{Duality for RPAU groups}), and similarly for $\QQ_p/\ZZ_p(r)$.

\Situation{mixed char setup} can be unified with \Situation{positive equal char setup} as follows.

\begin{lemma}[interpretation char mixed]
    Consider \Situation{mixed char setup}. Let $\hat K$ be the completion of $K$. Define :
    \begin{itemize}
        \item $\tilde S_k^m$ the site $k_{\RPS,\aff}$ of \Remark{affine or qcqs restriction} ;
        \item $u_k:S_k^m\to\tilde S_k^m$ the premorphism of sites defined by identity on the underlying categories ;
        \item $\psi_{K/k}^m:S_K^m\to\tilde S_k^m$ the premorphism of sites defined by $(\psi_{K/k}^m)^{-1}(T)=h^{\hat K/k}(T)$ for $T\in\tilde S_k^m$ ;
        \item $R\tilde\Psi_{K/k}^m:D(S_K^m)\to D(\tilde S_k^m)$ the derived pushforward of $\psi_{K/k}^m$.
    \end{itemize}
    Then we have the following.
    \begin{enumerate}
        \item For $G\in\Sh(S_K^m)$ and $q\in\ZZ$, $R^q\tilde\Psi_{K/k}^m(G)$ is the sheafification of the presheaf $T/k\mapsto H^q(h^{\hat K/k}(T),G)$.
        \item We have $u_{k,*}R\Psi_{K/k}^m = R\tilde\Psi_{K/k}^m$. In particular for $T\in\tilde S_k^m$, $R\Gamma(T,R\Psi_{K/k}^m(-))=R\Gamma(T,R\tilde\Psi_{K/k}^m(-))$.
        \item For $G\in D_\fin(S_K)$ we have $R\Psi_{K/k}^m(G)_{\Adh k}=R\Gamma(K^\ur,G_{\Adh K})$.
    \end{enumerate}
\end{lemma}

The premorphism $\psi_{K/k}^m$ cannot be defined to target all of $S_k^m$ because the infinite-level Kato lift does not generalize to non-affine schemes (see \Remark{Kato lift generalization}.2). We could work with $\tilde S_k^m$ and $R\tilde\Psi_{K/k}^m$ from the beginning in \Situation{mixed char setup}, and for our purposes no issue would come up (see also \Remark{affine or qcqs restriction}).

\begin{myproof}[of \Lemma{interpretation char mixed}]
    Note that $u_{k,*}$ is exact because coverings in $S_k^m$ of objects of $\tilde S_k^m$ admit refinements by coverings in $\tilde S_k^m$. Point \textbf{(1.)} follows from \Lemma{pushforward of sites}. Point \textbf{(2.)} follows from \textbf{(1.)} by exactness of $u_{k,*}$ and by \cite[Cor. 3.3]{KatoSuzuki2019} (though stated for sheaves of $\ZZ/p$-modules, the same proof works for abelian sheaves). For $G\in\Sh_\fin(S_K)$, taking the colimit over finite separable extensions $T=\Spec l\to \Spec k$ in \textbf{(2.)} and applying \textbf{(1.)} gives :
    $$R\Psi_{K/k}(G)_{\Adh k}=\varinjlim_{l/k}R\Gamma(l,R\tilde\Psi_{K/k}(G))=\varinjlim_{l/k} R\Gamma(h^{\hat K/k}(l),G)=R\Gamma(K^\ur,G|_{\Adh K})$$
    where the last equality is by \cite[Lem. 1.16 of Ch. III]{Milne80} (using that $G$
    is the pullback of a sheaf on $\Spec K_\et$). Since $R\Psi_{K/k}(-)_{\Adh k}$ and $R\Gamma(K^\ur,(-)_{\Adh K})$ are exact, the same holds for $G\in D_\fin(S_K)$, proving \textbf{(3.)}.
\end{myproof}

\begin{remark}[also works for n]
    Consider $n\geq 1$ an integer.
    \begin{enumerate}
        \item \Lemma{interpretation char mixed}.1-2 still applies if we consider the categories of $n$-torsion sheaves, meaning we replace the categories $\Sh(S)$ and $D(S)$ with $\Sh(S,\ZZ/n)$ and $D(S,\ZZ/n)$ for $S=S_K^m,\tilde S_K^m,S_k^m$, and the functors $R\Psi_{K/k}^m$ and $R\tilde\Psi_{K/k}^m$ with the analogous functors $R\Psi_{K/k,n}^m$ and $R\tilde\Psi_{K/k,n}^m$ between those categories.
        \item For $S=S_K^m,\tilde S_K^m,S_k^m$, let $J_n:\Sh(S)\to\Sh(S,\ZZ/n)$ be the inclusion functor. For $G\in D^+(S_K^m,\ZZ/n)$, we have a map, induced by the isomorphism $J_n\circ\Psi_{K/k,n}^m=\Psi_{K/k}^m\circ J_n$ and exactness of $J_n$ :
        $$J_n\circ R\Psi_{K/k,n}^m=R(J_n\circ\Psi_{K/k,n}^m)\to R\Psi_{K/k}^m\circ J_n$$
        This is an isomorphism. To see this, for $G\in\Sh(S_K^m,\ZZ/n)$, for $q\in\ZZ$ and for any $T\in S_K^m$ with a geometric point $\Adh x\to T$, by \Lemma{interpretation char mixed}.1-2 we have :
        $$R^q\Psi_{K/k}^m(J_n(G))_{\Adh x}
        = \varinjlim_{T'/T} H^q(h^{\hat K/k}(T'),J_n(G))
        = \varinjlim_{T'/T} J_n(H^q(h^{\hat K/k}(T'),G))
        = (J_n\circ R^q\Psi_{K/k,n}^m(G))_{\Adh x}$$
        where $T'$ ranges over factorizations $\Adh x\to T'\to T$ with $T'\in S_k^m$ affine and étale over $T$. Hence we have $R\Psi_{K/k}^m\circ J_n(G)=J_n\circ R\Psi_{K/k,n}(G)$. For $G\in D^+(S_K^m,\ZZ/n)$ we have spectral sequences :
        $$E^{ij}_2=R^i\Psi_{K/k}^m(H^j(J_n(G)))\Rightarrow R^{i+j}\Psi_{K/k}^m(J_n(G)),\qquad 'E^{ij}_2=R^i\Psi^m_{K/k,n}(H^j(G))\Rightarrow R^{i+j}\Psi_{K/k}^m(G))$$
        from which the equality $R\Psi_{K/k}^m\circ J_n(G)=J_n\circ R\Psi_{K/k,n}(G)$ follows in that case.
    \end{enumerate}
\end{remark}

\begin{lemma}[duality restriction]
    Consider $k$ a field of characteristic $p\geq 0$, $S_k$ a site and $d\in\ZZ$ as in one of the following :
    \begin{itemize}
        \item $p=0$, $d\in\ZZ$ is arbitrary and $S_k=\Spec k_\et$ ;
        \item $p>0$, $k$ is perfect, $d=0$ and $S_k=\Spec k_\et^\perar$ ;
        \item $p>0$, $[k:k^p]=p^d$ and $S_k$ is either $\Spec k_\RP$ or $\Spec k_\RPS$.
    \end{itemize}
    Let $\pi:\Spec l\to\Spec k$ be a finite field extension and $S_l$ the analogous site.
    \begin{enumerate}
        \item Base change by $\pi$ defines a morphism of sites $S_l\to S_k$, also written $\pi$, where $\pi_*:\Sh(S_l)\to\Sh(S_k)$ and $\pi^*:\Sh(S_k)\to\Sh(S_l)$ are both left and right adjoint to each other.
        \item There is an isomorphism $\eta:\pi^*\QQ/\ZZ(d)_k\xrightarrow\sim\QQ/\ZZ(d)_l$ in $D(S_l)$ and a map $\Tr:\pi_*\QQ/\ZZ(d)_l\to\QQ/\ZZ(d)_k$ in $D(S_k)$ such that the composite $\Tr\circ\pi_*\eta$ is the counit $\pi_*\pi^*\QQ/\ZZ(d)_k\to\QQ/\ZZ(d)_k$.
        \item If $G\otimes^LH\to\QQ/\ZZ(d)_l$ is a perfect pairing in $D(S_l)$, then we have a perfect pairing in $D(S_k)$ :
        $$\pi_*G\otimes^L\pi_*H\to\pi_*\QQ/\ZZ(d)_l\xrightarrow\Tr\QQ/\ZZ(d)_k.$$
        \item If $G\otimes^LH\to\QQ/\ZZ(d)_k$ is a pairing in $D(S_k)$ then we have a commutative diagram :
        $$\begin{tikzcd}
            G\arrow[d] &[-3em]\otimes^L &[-3em] H\arrow[rrr] &&& \QQ/\ZZ(d)_k\arrow[d,equals]\\
            \pi_*\pi^*G &[-3em]\otimes^L &[-3em] \pi_*\pi^*H\arrow[u]\arrow[r] & \pi_*\pi^*\QQ/\ZZ(d)_k \arrow[r,"\pi_*\eta"]& \pi_*\QQ/\ZZ(d)_l \arrow[r,"\Tr"]& \QQ/\ZZ(d)_k
        \end{tikzcd}$$
        where the vertical maps $\pi_*\pi^*\to \id$ and $\id\to\pi_*\pi^*$ are the counit and unit.
    \end{enumerate}
\end{lemma}

\begin{myproof}[of \Lemma{duality restriction}]
    \textbf{1.} If $S_k=\Spec k_\RP,\Spec k_\RPS$, this is \Proposition{pushforward exactness}.
    
    If $S_k=\Spec k_\et$ then $\pi:S_l\to S_k$ is a well-defined morphism of sites because $S_k$ has finite fiber products which are preserved by $\pi^{-1}$.
    
    If $S_k=\Spec k_\et^\perar$ with $k$ perfect, then for $T=\Spec R\in\Spec k_\et^\perar$ the algebra $R\otimes_kl$ is perfect Artinian over $l$ : since a perfect field extension $k'/k$ is geometrically reduced, $k'\otimes_kl$ is a finite product of perfect field extensions of $l$. So $\pi:\Spec l_\et^\perar\to\Spec k_\et^\perar$ is a well-defined premorphism of sites. Because $\pi:\Spec l\to\Spec k$ is a covering in $\Spec k_\et^\perar$, $\pi^\set$ is simply given by restricting a sheaf of sets to the subcategory $S_l\subseteq S_k$. In particular $\pi^\set$ is exact. Thus $\pi$ defines a morphism of sites.
    
    In the cases $S_k=\Spec k_\et,\Spec k_\et^\perar$, because $k$ is perfect the functor $\pi_*$ is left adjoint to $\pi^*$ by the same proof as \cite[Lem. 1.12 of Ch. V]{Milne80}, because $\pi$ is a finite étale morphism.
    
    \textbf{2.} If $\ell\neq p$ is a prime, then $\ZZ/\ell^n(d)_k$ is representable by the finite étale $k$-group $\mu_n^{\otimes d}$ concentrated in degree $0$. By \Lemma{yoneda and pullback}, $\pi^*\ZZ/\ell^n(d)_k$ is representable by the base change of $\mu_n^{\otimes d}$ to $l$ concentrated in degree $0$, which is again representable by $\mu_n^{\otimes d}$. Hence we have an isomorphism  $\eta_\ell:\pi^*\QQ_\ell/\ZZ_\ell(d)\xrightarrow\sim\QQ_\ell/\ZZ_\ell(d)_k$ and we can arbitrarily define $\Tr_\ell$ as the composite $\pi_*\QQ_\ell/\ZZ_\ell(d)_l\xrightarrow{\sim}\pi_*\pi^*\QQ_\ell/\ZZ_\ell(d)_k\to\QQ_\ell/\ZZ_\ell(d)_k$.
    
    In case $S_k=\Spec k_\et^\perar$ where $k$ is perfect and $d=0$, the group $\ZZ/p^n(d)=\ZZ/p^n$ remains finite étale and we construct $\eta_p$ and $\Tr_p$ the same way. In cases $S_k=\Spec k_\RP,\Spec k_\RPS$, the $p$-primary parts $\eta_p:\pi^*\QQ_p/\ZZ_p(d)_k\xrightarrow\sim\QQ_p/\ZZ_p(d)_l$ and $\Tr_p:\pi_*\QQ_p/\ZZ_p(d)_l\to\QQ_p/\ZZ_p(d)_k$ can be taken as in \Theorem{weil duality RP}.
    
    In all cases, we can take $\eta=\bigoplus_\ell\eta_\ell$ and $\Tr=\bigoplus_\ell\Tr_\ell$, where $\ell$ ranges over all primes, which concludes.
    
    \textbf{3.-4.} By \cite[Lem. 030K]{stacks-project} and by construction of $\eta$ and $\Tr$, we can assume either $l/k$ is separable, or $S_k=\Spec k_\RP,\Spec k_\RPS$ and $l/k$ is purely inseparable which only occurs in cases $S_k=\Spec k_\RP,\Spec k_\RPS$. These can be treated by the same proofs as \Theorem{weil duality RP}, using \textbf{(1.-2.)}.
\end{myproof}

\begin{remark}[need for pure inseparability]
    We will not use the full generality of \Lemma{duality restriction}. In the proof of \Theorem{mixed char} we will consider an extension $l/k$ which is prime-to-$p$ hence separable. In \Theorem{equal char 0}, we will consider only $\QQ_\ell/\ZZ_\ell(d)$ and coefficients $G\in D_\ell(S_l)$, $G,H\in D(S_k)$ with $\ell\neq p$ ; by \Lemma{finite coeffs equivalence}, the lemma then reduces to claims over small étale sites, which follow from duality for finite étale maps (when $l/k$ is separable) and topological invariance (when $l/k$ is purely inseparable). The hard case of \Lemma{duality restriction}, \textit{i.e.} \Theorem{weil duality RP} for $l/k$ purely inseparable, goes unused, as well as \Proposition{trace log dRW}, \Lemma{trace isomorphism}, and the equivalences $\Spec l_\RP\cong\Spec k_\RP$ and $\Spec l_\RPS\cong\Spec k_\RPS$ of \Lemma{base change RP}.2-3.
\end{remark}

\begin{lemma}[mixed char restriction]
    Consider $K,k,S_K,S_k,R\Psi_{K/k}$ as in \Situation{mixed char setup}. Consider $L/K$ a finite extension with residue field $l$ (which is a finite, possibly not separable extension of $k$). Then one can analogously consider $S_L,S_l,R\Psi_{L/l}$ as in \Situation{mixed char setup}. There is a canonical isomorphism in $D(S_k)$ : $$R\Psi_{K/k}(\pi_{L/K,*}G)=\pi_{l/k,*}(R\Psi_{L/l}(G))$$
    natural in $G\in D^+(S_L)$, where $\pi_{L/K}:S_L\to S_K$ and $\pi_{l/k}:S_l\to S_k$ are induced by base change.
\end{lemma}

For this proof we reduce to showing the equality on affine objects of $S_k$, which we do using \Lemma{interpretation char mixed}. We will use this strategy again several times, such as in \Lemma{RPsi computation} and \Lemma{char 0pp-0p interpretation}. 

\begin{myproof}[of \Lemma{mixed char restriction}]
    The functor $\pi_{L/K,*}:\Sh(S_L)\to\Sh(S_K)$ is exact because $L/K$ is finite étale, and the functor $\pi_{l/k,*}:\Sh(S_l)\to\Sh(S_k)$ is exact by \Proposition{pushforward exactness}.2. The commutative diagram :
    $$\begin{tikzcd}
        \Spec L\arrow[r]\arrow[d] & \Spec\O_L\arrow[d] & \Spec l\arrow[l]\arrow[d]\\
        \Spec K\arrow[r] & \Spec\O_K & \Spec k\arrow[l]
    \end{tikzcd}$$
    induces a map $\pi_{l/k,*}R\Psi_{L/l}\to R\Psi_{K/k}\pi_{L/K,*}$. We need to show $R^q\Psi_{K/k}(\pi_{L/K,*}G)\to \pi_{l/k,*}(R^q\Psi_{L/l}(G))$ is an isomorphism for $q\in\ZZ$ and $G\in\Sh(S_L)$. Both sides are Zariski sheaves, so it suffices to show they coincide on affine objects of $S_k$. By \Lemma{interpretation char mixed}.2, it suffices to show $R^q\tilde\Psi_{K/k}(\pi_{L/K,*}G)= \pi_{l/k,*}(R^q\tilde\Psi_{L/l}(G))$.
    
    By \Lemma{interpretation char mixed}.1, for $(E,F)=(K,k),(L,l)$ we have $R^q\tilde\Psi_{E/F}=a_F\circ P^q_{E/F}$, where $a_F$ is sheafification over $\tilde S_F$ and $P^q_{E/F}(H)\in\PSh(\tilde S_F)$ is the presheaf $T\mapsto H^q(h^{E/F}(T),H)$, for $H\in\Sh(S_F)$. By \cite[Lem. 9.7.(2)]{BertapelleSuzuki}, $\pi_{E/F,*}\circ a_E=a_F\circ\pi_{E/F,p}$ where $\pi_{E/F,p}$ is the pushforward of presheaves. For $T=\Spec R\in\tilde S_k$ :
    $$\pi_{l/k,p}P^q_{L/l}(G)(T)=H^q(h^{L/l}(R\otimes_kl),G)=H^q(h^{K/k}(R)\otimes_K L,G)=P^q_{K/k}(\pi_{L/K,*}(G))(T)$$
    by \Lemma{Kato lift restriction}. By applying $a_k$ we get : $\pi_{l/k,*}R^q\tilde\Psi_{L/l}(G)=R^q\tilde\Psi_{K/k}(\pi_{L/K,*}(G))$
    which concludes.
\end{myproof}

The following theorem is by \cite{KatoSuzuki2019} for coefficients $\ZZ/p^n(r)$ ; we generalize to finite coefficients. Recall Cartier and Serre duality from \Proposition{cartier duality} and \Theorem{duality for RPAU groups}.

\begin{theorem}[mixed char]
    Consider \Situation{mixed char setup}. Recall that here, $d=\log_p[k:k^p]$ is finite.
    \begin{enumerate}
            \item $R\Psi_{K/k}$ sends $D_p(S_K)$ to $D_0(S_k)$. For $G\in\Sh_p(S_K)$, we have $R^q\Psi_{K/k}(G)=0$ if $q\neq 0,\dots,d+1$.
            \item There is a canonical trace map $R\Psi_{K/k}(\QQ_p/\ZZ_p(d+1))\to\QQ_p/\ZZ_p(d)[-1]$ in $D(S_k)$.
            \item For $G\in D_p(S_K)$ and $G^\vee=R\HOM_{S_K}(G,\QQ_p/\ZZ_p(d+1))$, there is a perfect pairing in $D(S_k)$ :
            $$R\Psi_{K/k}(G)\otimes^LR\Psi_{K/k}(G^\vee)\to R\Psi_{K/k}(\QQ_p/\ZZ_p(d+1))\xrightarrow\tr \QQ_p/\ZZ_p(d)[-1]$$
            which gives perfect Serre dualities of split and wound parts for $q\in\ZZ$ :
            \begin{align*}
                R^q\Psi_{K/k}(G)_s&=\EXT_{S_k}^1(R^{d+1-q}\Psi_{K/k}(G^\vee_K)_s,\nu_\infty(d)),\\ R^q\Psi_{K/k}(G)_w&=\HOM_{S_k}(R^{d-q}\Psi_{K/k}(G^\vee_K)_w,\nu_\infty(d)).
            \end{align*}
        \end{enumerate}
\end{theorem}

\begin{myproof}[of \Theorem{mixed char}]
    \textbf{1.} Consider $G\in\Sh_p(S_K)$, and let $L/K$ be as in \Lemma{devissage}, for $G$ and the prime $p\neq\characteristic K$. Let $l/k$ be the residue extension. Consider $S_L,S_l,R\Psi_{L/l}$ as in \Situation{mixed char setup} analogously.
    
    Let $\C\subseteq\Sh(S_L)$ be the subcategory of objects $H\in\Sh_p(S_L)$ such that $R^q\Psi_{L/l}(H)$ belongs to $\Sh_0(S_l)$ for all $q\in\ZZ$, and is trivial for $q\neq 0,\dots,d+1$. Then $\C$ is closed under extensions by \Proposition{RPAU characterization} and exactness of $R\Psi_{L/l}$. The constant sheaf $\ZZ/p$ belongs to $\C$ by \cite[Th. 3.4.(2)]{KatoSuzuki2019}. Thus $G|_L$, which admits a finite separated filtration with quotients $0$ or $\ZZ/p$, also belongs to $\C$.
    
    Write $P=\pi_{L/K,*}(G|_L)$. By \Lemma{devissage}.3 there is a decomposition $P = G\oplus H$ for some $H\in\Sh_p(S_K)$, so we have $R^q\Psi_{K/k}(P) = R^q\Psi_{K/k}(G)\oplus R^q\Psi_{K/k}(H)$. By \Lemma{mixed char restriction}, $R^q\Psi_{K/k}(P)=\pi_{l/k,*}R^q\Psi_{L/l}(G|_L)$, so $R^q\Psi_{K/k}(P)$ is trivial for $q\neq 0,\dots,d+1$, and belongs to $\Sh_0(S_k)$ by the above and \Corollary{weil restriction RPAU}, for all $q\in\ZZ$. As a direct factor of $R^q\Psi_{K/k}(P)$, $R^q\Psi_{K/k}(G)=0$ if $q\neq 0,\dots,d+1$ and for any $q\in\ZZ$ is the kernel of some endomorphism of $R^q\Psi_{K/k}(P)$, thus belongs to $\Sh_0(S_k)$ because $\Sh_0(S_k)$ is abelian.
    
    \textbf{2.} We use the notations of \Remark{also works for n}. For $n\geq 1$, by \cite[Th. 3.4.(1)]{KatoSuzuki2019}, we have a map :
    $$R\Psi_{K/k}(\ZZ/p^n(d+1))\cong J_n\circ R\Psi_{K/k,p^n}(\ZZ/p^n(d+1))\to \ZZ/p^n(d)[-1]$$
    using the equalities $\ZZ/p^n(d)=\nu_n(d)[-d]$ over $k$ and, by exactness of $J_{p^n}$, $\ZZ/p^n(r)=J_{p^n}(\ZZ/p^n(r))$.
    
    The functor $i^*:\Sh(S_{\O_K})\to\Sh(S_k)$ preserves colimits as a right adjoint functor, and for $q\in\ZZ$, $R^qj_*:\Sh(S_K)\to\Sh(S_{\O_K})$ preserves filtered colimits by \cite[Prop. 2.2.4.(b)]{SuzukiGN}. Hence $R^q\Psi_{K/k}$ preserves filtered colimits for all $q\in\ZZ$, which defines the morphism $R\Psi_{K/k}(\QQ_p/\ZZ_p(d+1))\to\QQ_p/\ZZ_p(d)[-1]$.
    
    \textbf{3.} We can take $G\in\Sh_p(S_K)$ by exactness of $R\Psi_{K/k}$ and $(-)^\vee$. Let $L/K$ and $l/k$ be as in \textbf{(1.)}, and write $(-)^\vee_F=R\HOM_{S_F}(-,\QQ_p/\ZZ_p(d+1))$ for $F=K,L$. By \Lemma{change n if you like}, \Remark{also works for n}.2 and \cite[Th. 3.4.(3)]{KatoSuzuki2019}, the pairing :
    $$R\Psi_{L/l}(H)\otimes^LR\Psi_{L/l}(H^\vee_L)\to R\Psi_{L/l}(\QQ_p/\ZZ_p(d+1))\xrightarrow\tr\QQ_p/\ZZ_p(d)[-1]$$
    is perfect when $H=\ZZ/p$. By exactness of $(-)^\vee_L=R\HOM_{S_L}(-,\QQ_p/\ZZ_p(d+1))$ and $R\Psi_{L/l}$, it is also perfect for $H=G|_L$ which has a finite separated filtration with quotients $0$ or $\ZZ/p$.
    
    Write $P=\pi_{L/K,*}(G|_L)$. The pairing :
    $$R\Psi_{K/k}(P)\otimes^LR\Psi_{K/k}(P^\vee_K)\to\QQ_p/\ZZ_p(d)[-1]$$
    is perfect in $D(S_l)$ by \Lemma{mixed char restriction} and \Lemma{duality restriction}.3 for the extension $l/k$. By \Lemma{devissage}.3, there is a decomposition $P = G\oplus H$ for some $H\in\Sh_p(S_K)$, so $R\Psi_{K/k}(P) = R\Psi_{K/k}(G)\oplus R\Psi_{K/k}(H)$ in $D(S_k)$. Similarly $P^\vee_K=G^\vee_K\oplus H^\vee_K$ and $R\Psi_{K/k}(P^\vee_K) = R\Psi_{K/k}(G^\vee_K)\oplus R\Psi_{K/k}(H^\vee_K)$, and the previous pairing decomposes as the orthogonal direct sum of the corresponding pairings :
    $$R\Psi_{K/k}(G)\otimes^LR\Psi_{K/k}(G^\vee_K)\to\QQ_p/\ZZ_p(d)[-1],\qquad R\Psi_{K/k}(H)\otimes^LR\Psi_{K/k}(H^\vee_K)\to\QQ_p/\ZZ_p(d)[-1]$$
    by \Lemma{duality restriction}.4. Hence both are also perfect in $D(S_k)$. The rest follows from \Theorem{duality for RPAU groups}.5.
\end{myproof}

\begin{remark}[for small et]
    By \Lemma{finite coeffs equivalence}.2 and \Proposition{cartier duality}.4, the morphism $u:S_K\to S_K'$ where $S_K'=\Spec K_\et$, induces an exact equivalence $u^*:D_p(S_K')\xrightarrow\sim D_p(S_K)$ which commutes with the Cartier duality and Yoneda functors from the category of $p$-primary finite étale $K$-groups. Similarly, by \Proposition{RPAU characterization} and \Proposition{definition of D0}, the obvious premorphism $v:S_k'\to S_k$ where $S_k'=\Spec k_\RP$ induces an exact equivalence $v^*:D_0(S_k)\xrightarrow\sim D_0(S_k')$ which commutes with the Serre duality and Yoneda functors from $D^b(\Alg_u^\RP(k))$. Thus all of \Theorem{mixed char}, as well as \Lemma{interpretation char mixed}.3 remain true if we replace $S_K$ by $S_K'$ and/or $S_k$ by $S_k'$, and accordingly replace $R\Psi_{K/k}$ by its composite on the left with $Lv^*$ and/or on the right with $u^*$.
\end{remark}


\subsection{Prime-to-characteristic steps}
\label{Prime-to-characteristic steps}

\begin{situation}[variety setup]
    Let $k$ be a field and $d$ any integer. Let $X/k$ be a proper, smooth, geometrically integral $k$-scheme of dimension $e\geq 0$. Let $U\subseteq X$ be any nonempty open. Set $S_U^v=U_\et$ and $S_k^v=\Spec k_\et$ and define a morphism of sites $\psi_{U/k}^v:S_U^v\to S_k^v$ by :
    $$(\psi_{U/k}^v)^{-1}:T/k\mapsto U_T=U\times_kT.$$
    We write $\Psi_{U/k}^v:\Sh(S_U^v)\to\Sh(S_k^v)$ the pushforward of $\psi_{U/k}^v$. We also define $R\Psi_{U/k,c}^v=R\Psi_{X/k}^v\circ j_!$ where $j_!:D(S_U^v)\to D(S_X^v)$ is the exceptional pushforward induced by the open immersion $j:U\to X$.
\end{situation}

We will as before often omit the superscript $(-)^v$ (which stands for "variety"). Note that $\psi_{U/k}^v$ is a morphism of sites is because $\Spec k_\et$ admits fiber products which are preserved by $(\psi_{U/k}^v)^{-1}$. As in section \ref{Finite coefficients}, over $S_U$ or $S_k$, for $\ell\neq\characteristic k$ a prime, by $\ZZ/\ell^n(r)$ and $\QQ_\ell/\ZZ_\ell(r)$ we mean the sheaves $\mu_{\ell^n}^{\otimes r}$ and $\varinjlim_n\mu_{\ell^n}^{\otimes r}$.

The following is a reformulation of Poincaré duality with the above notations.

\begin{theorem}[poincare duality]
    Consider \Situation{variety setup}. Let $\ell\neq\characteristic k$ be a fixed prime and $d$ a fixed integer.
    \begin{enumerate}
        \item The functors $R\Psi_{U/k}$ and $R\Psi_{U/k,c}$ send objects of $D_\ell(S_U)$ to objects of $D_\ell(S_k)$. If $G\in\Sh_\ell(S_U)$ then $R^q\Psi_{U/k}(G)=R^q\Psi_{U/k,c}(G)=0$ for $q\neq 0,\dots,2e$.
        \item There exists a canonical trace map $\tr:R\Psi_{U/k,c}(\QQ_\ell/\ZZ_\ell(d+e))\to\QQ_\ell/\ZZ_\ell(d)[-2e]$ in $D(S_k)$.
        \item For any $G\in D_\ell(S_U)$ with Cartier dual $G^\vee=R\HOM_{S_U}(G,\QQ_\ell/\ZZ_\ell(d+e))$, we have a perfect pairing :
        $$R\Psi_{U/k}(G)\otimes^LR\Psi_{U/k,c}(G^\vee)\to R\Psi_{U/k,c}(\QQ_\ell/\ZZ_\ell(d+e))\xrightarrow\tr\QQ_\ell/\ZZ_\ell(d)[-2e]$$
        which gives perfect Cartier dualities of finite étale $k$-groups for $q\in\ZZ$ :
        $$R^q\Psi_{U/k}(G)\otimes R^{2e-q}\Psi_{U/k,c}(G^\vee)\to \QQ_\ell/\ZZ_\ell(d).$$
        \item For $G\in D^+(S_U)$ we have $R\Psi_{U/k}(G)=R\Gamma(U_{\Adh k},G)$ and $R\Psi_{U/k,c}(G)=R\Gamma_c(U_{\Adh k},G)=R\Gamma(X_{\Adh k},j_{\Adh k,!}G)$, where we again write $G$ for its pullback to other étale sites and the morphism $j_{\Adh k}:U_{\Adh k}\to X_{\Adh k}$ is the base change of $j:U\to X$ by $\Spec \Adh k\to\Spec k$.
    \end{enumerate}
\end{theorem}

\begin{myproof}[of \Proposition{poincare duality}]
    \textbf{4.} As a complex of sheaves on the étale site of a field, $R\Psi_{U/k}(G)$ can be identified with the complex of Galois modules of stalks $R\Psi_{U/k}(G)_{\Adh k}$. By \Lemma{pushforward of sites}, $R^q\Psi_{U/k}(G)$ is the sheafification of the presheaf $T/k\mapsto H^q(U_T,G)$ for $q\in\ZZ$, so that : $$R^q\Psi_{U/k}(G)_{\Adh k}=\varinjlim_{l/k}H^q(U_l,G)=H^q\left(\varprojlim_{l/k}U_l,G\right)=H^q(U_{\Adh k},G)$$
    where $l$ ranges over finite separable extensions of $k$, by \cite[Lem. 1.16 of Ch. III]{Milne80}. Similarly $R^q\Psi_{U/k,c}(G)$ is the sheafification of the presheaf $T/k\mapsto H^q(X_T,(j_!G)|_{X_T})$, so $R^q\Psi_{U/k,c}(G)=H^q(X_{\Adh k},(j_!G)|_{X_{\Adh k}})$ and by base change \cite[Th. 0EYU]{stacks-project}, $R^q\Psi_{U/k,c}(G)=H^q(X_{\Adh k},j_{\Adh k,!}(G|_{X_{\Adh k}}))$ as desired.
    
    \textbf{1-3.} By \textbf{(4.)} these are reformulations of \cite[Eq. 3.2.6.1 of Exposé XVIII]{SGA4} in the case $d=0$, since $U_{\Adh k}$ is a smooth $\Adh k$-variety of dimension $e$ by the geometrically integral assumption. Case $d\in\ZZ$ follows from :
        \begin{align*}
                \QQ_\ell/\ZZ_\ell(d+e)
                &= \QQ_\ell/\ZZ_\ell(e)\otimes\QQ_\ell/\ZZ_\ell(d)\\
                H^q_c(U_{\Adh k},G\otimes\QQ_\ell/\ZZ_\ell(d))
                &= H^q_c(U_{\Adh k},G)\otimes\QQ_\ell/\ZZ_\ell(d)\\
                \HOM_{S_k}(G,\QQ_\ell/\ZZ_\ell(d+e))
                &= \HOM_{S_k}(G,\QQ_\ell/\ZZ_\ell(e))\otimes\QQ_\ell/\ZZ_\ell(d)
            \end{align*}
    for arbitrary $G\in\Sh(S_k)$.
\end{myproof}

Recall the definition of the étale lift $h^{K/k}(T)$, for $T/k$ étale, as a special case of \Definition{Kato lift definition}.

\begin{situation}[equal char 0 setup]
    Let $K$ be a complete discrete valuation field with residue field $k$ with equal characteristic $0$. Set $S_K^0=\Spec K_\et$ and $S_k^0=\Spec k_\et$ and define a premorphism of sites $\psi_{K/k}^0:S_K^0\to S_k^0$ by :
    $$(\psi_{K/k}^0)^{-1}:T/k\mapsto h^{K/k}(T).$$
    We write $\Psi_{K/k}^0:\Sh(S_K^0)\to\Sh(S_k^0)$ the pushforward of $\psi_{K/k}^0$.
\end{situation}

As before we will usually omit the superscript $(-)^0$ (which stands for "characteristic $0$"). 

Duality in \Situation{equal char 0 setup} is classical. For the purpose of duality for higher local fields, one could use the formalisms of \Situation{mixed char setup} and \Situation{positive equal char setup} exclusively for $p$-torsion coefficients and treat $\ell$-torsion coefficients more straightforwardly with just Galois cohomology and small étale sites, which can be treated identically to \Situation{equal char 0 setup}. In this section we will also show that $R\Psi_{K/k}$ satisfies prime-to-$p$ duality in \Situation{mixed char setup} and \Situation{positive equal char setup}. This is harder than when working directly with Galois cohomology, as can be seen in the proofs of \Lemma{RPsi computation} and \Theorem{equal char 0}. As noted above this is not strictly needed to obtain \Theorem{intro general} and \Theorem{intro finite}, but illustrates the kind of difficulties that appear in these formalisms, and shows that they work uniformly independantly of the characteristic. We need preliminary computations. The main tools are Kato and Suzuki's ind-smooth approximation of the Kato lift \cite[Prop. 4.1]{KatoSuzuki2019} and Fujiwara-Gabber's formal base change \cite[Cor. 1.18.(2)]{BhattMathew}, here in the form of \cite[Prop. 4.2]{KatoSuzuki2019}).

\begin{lemma}[general restriction]
    Let $K,k,S_K,S_k,R\Psi_{K/k}$ be as in \Situation{positive equal char setup}, \Situation{mixed char setup} or \Situation{equal char 0 setup}. Let $L/K$ be a finite separable extension with residue field $l$ finite (but not necessarily separable) over $k$, and $S_L,S_l,R\Psi_{L/l}$ the analogous objects. Then there are canonical isomorphisms for all $q\in\ZZ$ : $$R\Psi_{K/k}(\pi_{L/K,*}G)=\pi_{l/k,*}(R\Psi_{L/l}(G)),\qquad R^q\Psi_{K/k}(\pi_{L/K,*}G)=\pi_{l/k,*}(R^q\Psi_{L/l}(G))$$
    natural in $G\in D^+(S_L)$, where $\pi_{L/K}:S_L\to S_K$ and $\pi_{l/k}:S_l\to S_k$ are induced by base change.
\end{lemma}

\begin{myproof}[of \Lemma{general restriction}]
    The well-definedness of $\pi_{L/K}$ and $\pi_{l/k}$ and exactness of the pushforward is either by \Lemma{pushforward exactness} or properties of finite étale maps. For the mixed characteristic case \Situation{mixed char setup}, see \Lemma{mixed char restriction}. For the other cases, the proof is identical except we do not need \Lemma{interpretation char mixed}.1 : we directly use \Lemma{pushforward of sites} to show $R^q\Psi_{E/F}(G)$ is the sheafification of $T\mapsto H^q(h^{E/F}(T),G)$, for $(E,F)=(K,k),(L,l)$.
\end{myproof}

\begin{lemma}[indsmooth purity]
    Let $n\geq 1$ be an integer and $i_0:Z_0\to X_0$ a closed immersion of regular $\ZZ[1/n]$-schemes, of codimension $c\geq 0$. Let $\{X_\alpha\}_{\alpha\in A}$ be a filtered diagram of smooth $X_0$-schemes, with affine transition maps, and set $Z_\alpha=X_\alpha\times_{X_0}Z_0$, $X=\varprojlim_\alpha X_\alpha$, $Z=\varprojlim_\alpha Z_\alpha$. Label the maps as follows :
    $$\begin{tikzcd}
       Z \arrow[rr,bend left=30,"f_{Z,0}"]\arrow[r,"f_{Z,\alpha}"']\arrow[d,"i"] & Z_\alpha \arrow[r,"f_{Z,\alpha,0}"']\arrow[d,"i_\alpha"] & Z_0\arrow[d,"i_0"] \\
       X \arrow[rr,bend right=30,"f_{X,0}"]\arrow[r,"f_{X,\alpha}"] & X_\alpha \arrow[r,"f_{X,\alpha,0}"] & X_0
    \end{tikzcd}$$
    \begin{enumerate}
        \item For $F$ a sheaf of $\ZZ/n$-modules on $X_{0,\et}$ and $q\in\ZZ$, we have : $$R^qi_\alpha^!f_{X,\alpha,0}^*F = f_{Z,\alpha,0}^*R^qi_0^!F,\qquad
        R^qi^!f_{X,0}^*F = f_{Z,0}^*R^qi_0^!F=\varinjlim_\alpha f_{Z,\alpha}^*R^qi_\alpha^!f_{X,\alpha,0}^*F.$$
        \item For all $r\in\ZZ$ we have $R^qi^!(\ZZ/n(r)) = \ZZ/n(r-c)$ if $q=2c$ and $R^qi^!(\ZZ/n(r))=0$ if $q\neq 2c$.
    \end{enumerate}
\end{lemma}

\begin{myproof}[of \Lemma{indsmooth purity}]
    \textbf{1.} We start with $q=0$. Write $j:U\to X$, $j_\alpha:U_\alpha\to X_\alpha$, $j_0:U_0\to X_0$ the open complements of $Z$, $Z_\alpha$ and $Z_0$, and $f_{U,\alpha,0}:U_\alpha\to U_0$, $f_{U,0}:U\to U_0$, $f_{U,\alpha}:U\to U_\alpha$ the obvious maps, so that $U_\alpha$ and $U$ are base changes of $U_0$. We have :
    \begin{align*}
        i_\alpha^!f_{X,\alpha,0}^*F
        &= i_\alpha^*\ker(f_{X,\alpha,0}^*F\to j_{\alpha,*}j_\alpha^*f_{X,\alpha,0}^*F)\\
        &= \ker(i_\alpha^*f_{X,\alpha,0}^*F\to i_\alpha^*j_{\alpha,*}j_\alpha^*f_{X,\alpha,0}^*F)\\
        &= \ker(f_{Z,\alpha,0}^*i_0^*F\to i_\alpha^*j_{\alpha,*}f_{U,\alpha,0}^*j_0^*F)\\
        &= \ker(f_{Z,\alpha,0}^*i_0^*F\to i_\alpha^*f_{X,\alpha,0}^*j_{0,*}j_0^*F).
    \end{align*}
    The first equality is by \cite[Rem. 3.13 of Ch. II]{Milne80}. The second is by exactness of pullbacks. The third is by functoriality of pullbacks. For the fourth, by smoothness of $X_\alpha$ we have $f_{X,\alpha,0}^* = f_{X,\alpha,0}^!(-d)[-2d]$ on $\Sh(X_{0,\et},\ZZ/n)$, where $d$ is the (locally constant) dimension of $X_\alpha$ over $X_0$, and because $Rj_{\alpha,*}$ and $Rj_{0,*}$ commute with the Tate twist of $n$-torsion sheaves, when restricted to $\Sh(X_{0,\et},\ZZ/n)$ we have :
    $$Rj_{\alpha,*}f_{U,\alpha,0}^* =(Rj_{\alpha,*}f_{U,\alpha,0}^!)(-d)[-2d]=(f_{X,\alpha,0}^!Rj_{0,*})(-d)[-2d]=f_{X,\alpha,0}^*Rj_{0,*}$$
    and $j_{\alpha,*}f_{U,\alpha,0}^*=f_{X,\alpha,0}^*j_{0,*}$ as required. Then by functoriality of pullbacks and exactness of $f_{Z,\alpha,0}^*$ :
    $$i_\alpha^!f_{X,\alpha,0}^*F= \ker(f_{Z,\alpha,0}^*i_0^*F\to f_{Z,\alpha,0}^*i_0^*j_{0,*}j_0^*F)= f_{Z,\alpha,0}^*i_0^!F$$
    as desired for $q=0$. The claim for $q\geq 1$ follows because pullbacks are exact. For $i$, the proof is identical except in the fourth equality we use \cite[Lem. 1.16 of Ch. III]{Milne80} to have :
    $$Rj_*f_{U,0}^*=\varinjlim_\alpha f_{X,\alpha}^*Rj_{\alpha,*}f_{U,\alpha,0}^*=\varinjlim_\alpha f_{X,\alpha}^*f_{X,\alpha,0}^*Rj_{0,*}=f_{X,0}^*Rj_{0,*}.$$
    
    \textbf{2.} This follows from \textbf{(1.)} and the absolute purity theorem for the $i_\alpha$'s \cite[Th. 3.1.1]{Riou}.
\end{myproof}

\begin{lemma}[RPsi computation]
    Consider $K,k,S_K,S_k,R\Psi_{K/k}$ as in any of \Situation{positive equal char setup}, \Situation{mixed char setup} or \Situation{equal char 0 setup}. Let $n\geq 1$ be invertible in $k$ and $d\in\ZZ$ arbitrary. Then for $r,q\in\ZZ$ we have canonical isomorphisms :
    \begin{align*}
        R^q\Psi_{K/k}(\ZZ/\ell^n(r))
        &=\ZZ/\ell^n(r)\text{ if q=0,}\\
        &=\ZZ/\ell^n(r-1)\text{ if q=1,}\\
        &=0\text{ otherwise.}
    \end{align*}
\end{lemma}

\begin{myproof}[of \Lemma{RPsi computation}]
    The (restriction to affines of the) sheaf $R^q\Psi_{K/k}(\ZZ/n(r))\in\Sh(S_k)$ is the sheafification of the presheaf $T\in S_k\mapsto H^q(h^{K/k}(T),\ZZ/n(r))$ by \Lemma{pushforward of sites} in \Situation{positive equal char setup} and \Situation{equal char 0 setup}, and by \Lemma{interpretation char mixed} in \Situation{mixed char setup}. We have, for $r,s\in\ZZ$ :
    $$H^q(h^{K/k}(T),\ZZ/n(r+s))=\ZZ/n(s)\otimes H^q(h^{K/k}(T),\ZZ/n(r))$$
    so it suffices to prove the claim for one $r\in\ZZ$, for each $q\in\ZZ$. Here is a summary of the following proof.
    
    For $q=0$, we always easily have $\Psi_{K/k}(\mu_n)=\mu_n$, essentially by \Proposition{Kato lift structure}.2.
    
    For $q=1$, the valuation map gives an injection $\delta:\Psi_{K/k}(\GG_m)/n\to \ZZ/n$. In \Situation{equal char 0 setup} and \Situation{positive equal char setup} we get surjectivity of $\delta$ by Hilbert's theorem 90 for fields. In \Situation{mixed char setup} it is not clear that $R^1\Psi_{K/k}(\GG_m)$ is trivial so the same argument doesn't adapt. Instead we use a localization sequence and compute $H^q(h^{\O_K/k}(R),\ZZ/n(r))$ and $H^q_{\Spec R}(h^{\O_K/k}(R),\ZZ/n(r))$. The first group is computed using the fact that $(h^{\O_K/k}(R),\mm_Kh^{\O_K/k}(R))$ is Henselian, and the second by locally approximating $h^{\O_K/k}(R)$ with an ind-smooth $\O_K$-algebra and using absolute purity per \Lemma{indsmooth purity}.
    
    In \Situation{equal char 0 setup} and \Situation{positive equal char setup}, the cancellation $R^q\Psi_{K/k}(\ZZ/n(r))=0$ for $q>1$ follows from Lang's theorem for complete discrete valuation fields with separably closed residue field. In \Situation{mixed char setup} we use the same methods as for $q=1$ to conclude for $q>1$.
    
    \textbf{Consider \Situation{equal char 0}.} Then $R^q\Psi_{K/k}(\ZZ/n(r))$ is identified with the Galois module $H^q(K^\ur,\ZZ/n(r))$ (see also \Theorem{equal char 0}.4 below). By \cite[Cor. 1.18.(2)]{BhattMathew}, we have an equality $H^q(K^\ur,\ZZ/n(r))=H^q(\widehat{K^\ur},\ZZ/n(r))$. This reduces to computing the cohomology of a complete discrete valuation field with separably closed residue field. Cancellation for $q>1$ comes from $\cd_\ell(\widehat{K^\ur})\leq 1$ for $\ell\neq\characteristic k$ prime, by Lang's theorem. We have a split exact sequence given by the valuation map :
    $$0\to(\widehat{\O_K^\ur})^\times\to(\widehat{K^\ur})^\times\to\ZZ\to 0.$$
    By passing to $n$-torsion and quotient by $n$-powers, we get $\mu_n(\widehat{\O_K^\ur})=\mu_n(\widehat{K^\ur})$ and an exact sequence :
    $$0\to (\widehat{\O_K^\ur})^\times/n\to(\widehat{K^\ur})^\times\to\ZZ/n\to 0.$$
    By Hensel's lemma, $\mu_n(\widehat{\O_K^\ur})=\mu_n(\Adh k)$ which gives case $q=0$, $r=1$. Hensel's lemma also gives $(\widehat{\O_K^\ur})^\times/n=0$, so $(\widehat{K^\ur})^\times/n=\ZZ/n$. By Hilbert's theorem 90 we have $H^1(\widehat{K^\ur},\GG_m)=0$, so by Kummer theory $H^1(\widehat{K^\ur},\ZZ/n(1))=(\widehat{K^\ur})^\times/n$, which gives the case $q=1$, $r=1$. This concludes in \Situation{equal char 0}.
    
    \textbf{Consider \Situation{positive equal char setup} and \Situation{mixed char setup}.} We adapt the above proof. In \Situation{mixed char setup} we write $S_{\O_K}=\Spec \O_{K,\RPS}$ (\Definition{RPS site of O}) and $R\Psi_{\O_K/k}=i^*$. In \Situation{positive equal char setup} we write $S_{\O_K}=\Spec \O_{K,\RP}$ (\Definition{relatively perfect}) and $R\Psi_{\O_K/k}$ the pushforward of the premorphism $T\in S_k\mapsto h^{\O_K/k}(T)$. Thus we have a functor $R\Psi_{\O_K/k}:D(S_{\O_K})\to D(S_k)$ such that for $G\in\Sh(S_{\O_K})$ and $q\in\ZZ$, (the restriction to affines of) $R^q\Psi_{\O_K/k}(G)$ is the sheafification of $T\mapsto H^q(h^{\O_K/k}(T),G)$, either by \cite[Prop. 3.2]{KatoSuzuki2019} or \Lemma{pushforward of sites}.
    
    For any affine $\Spec R\in S_k$, $R$ is a finite product of integral domains, by definition of $\Spec k_\et^\perar$ in \Situation{positive equal char setup} or by \Lemma{local RP is good} in \Situation{mixed char setup}. Thus \Proposition{Kato lift structure}.2 yields an exact sequence in $\Sh(S_k)$ :
    $$0\to\Psi_{\O_K/k}(\GG_{m,\O_K})\to\Psi_{K/k}(\GG_{m,K})\to \ZZ\to 0$$
    Passing to $n$-torsion and quotient by $n$-powers, because $\Psi_{K/k}$ and $\Psi_{\O_K/k}$ are left exact, we have an isomorphism $\Psi_{K/k}(\mu_n)=\Psi_{\O_K/k}(\mu_n)$ and an exact sequence :
    $$\Psi_{\O_K/k}(\GG_{m,\O_K})/n\to\Psi_{K/k}(\GG_{m,K})/n\to\ZZ/n\to 0.$$
    For any affine $\Spec R\in S_k$, $(h^{\O_K/k}(R),\mm_Kh^{\O_K/k}(R))$ is a Henselian pair by \Proposition{Kato lift structure}.1, thus $\mu_n(h^{\O_K/k}(R))=\mu_n(R)$. After sheafifying we get $\Psi_{\O_K/k}(\mu_n)=\mu_n$, which proves case $q=0$, $r=1$.
    
    The map $\Psi_{\O_K/k}(\GG_m)\xrightarrow n\Psi_{\O_K/k}(\GG_m)$ is surjective in the étale topology : if $R\in S_k$ is affine and $x\in h^{\O_K/k}(R)^\times$, then $\Adh x\in R$ is a unit so $R\to R[x^{1/n}]$ is an étale faithfully flat $k$-algebra map and the restriction of $x$ to $h^{\O_K/k}(R[\Adh x^{1/n}])$ is an $n$-th power because $(h^{\O_K/k}(R),\mm_Kh^{\O_K/k}(R))$ is a Henselian pair. Thus $\Psi_{K/k}(\GG_m)/n=\ZZ/n$ and Kummer theory gives an injection $\delta:\Psi_{K/k}(\GG_{m,K})/n\to R^1\Psi_{K/k}(\mu_n)$. It remains to show $\delta$ is surjective and $R^q\Psi_{K/k}(\ZZ/n(r))=0$ for $q\geq 2$, which we check at stalks.
    
    \textbf{In \Situation{equal char p}}, the strict Henselization at a prime of an object of $S_k=\Spec k_\et^\perar$ is the algebraic closure of a perfect field extension of $k$. Let $F$ be such a field extension with algebraic closure $\Adh F$. Then :
    $$R^q\Psi_{K/k}(G)_{\Adh F}=\varinjlim_{E/F\text{ finite}}H^q(h^{K/k}(E),G)=H^q\left(\varinjlim_{E/F\text{ finite}}h^{K/k}(E),G\right)=H^q(h^{K/k}(\Adh F),G)$$
    where the last equality is by \cite[Cor. 1.18.(2)]{BhattMathew} : $h^{K/k}(F)$ is a discrete valuation field (\Proposition{Kato lift structure}.2), with maximal unramified extension $\varinjlim_{E/F}h^{K/k}(E)$ a Henselian discrete valuation field, whose completion is precisely $h^{K/k}(\Adh F)$. Since $h^{K/k}(\Adh F)$ is a complete discrete valuation field with algebraically closed residue field, we get $R^1\Psi_{K/k}(\GG_m)(\Adh F)=0$ by Hilbert's theorem 90, so $\delta$ is surjective, and $R^q\Psi_{K/k}(\ZZ/n(r))(\Adh F)=0$ for $q\geq 2$ by Lang's theorem. This proves the required properties in \Situation{equal char p}.
    
    \textbf{In \Situation{mixed char setup}}, consider $R$ a relatively perfectly smooth $k$-algebra and $R'$ the strict Henselization of $R$ at a prime identified with a point $x\in\Spec R$. By \cite[Prop. 4.1 and 4.2]{KatoSuzuki2019}, there exists an $\O_K$-algebra $B$ with the following properties.
    \begin{enumerate}
        \item[(i)] $B\otimes_{\O_K}k=R'$ ;
        \item[(ii)] $B$ is a filtered colimit $\varinjlim_\alpha B_\alpha$ of smooth $\O_K$-algebras $B_\alpha$ ;
        \item[(iii)] $(B,\mm_KB)$ is a Henselian pair ;
        \item[(iv)] We have isomorphisms $R^q\Psi_{K/k}(G)_{\Adh x}=H^q(B_K,G)=H^q(h^{K/k}(R),G)$ induced by restrictions, for $G\in\Sh_\fin(S_K)$ and $q\in\ZZ$, where $B_K=B\otimes_{\O_K}K$.
    \end{enumerate}
    
    By (i), (iii), \Proposition{Kato lift structure}.1 and \cite[Th. 09ZI]{stacks-project} we have for all $q\geq 0$ :
    $$H^q(h^{\O_K/k}(R'),\ZZ/n(1))=H^q(B,\ZZ/n(1))=H^q(R',\ZZ/n(1)).$$
    Write $Z=\Spec R'$, $X=\Spec h^{\O_K/k}(R')$, $U=\Spec h^{K/k}(R')$, $X'=\Spec B$ and $U'=\Spec B_K$. By (iv) and the above, restriction by $X\to X'$ gives an isomorphism of the exact sequences of \cite[Prop. 1.25 of Ch. III]{Milne80} :
    $$\begin{tikzcd}
       \cdots\arrow[r]& H^q_Z(X)\arrow[r]\arrow[d,"\sim"] & H^q(X)\arrow[r]\arrow[d,"\sim"] & H^q(U)\arrow[r]\arrow[d,"\sim"] & H^{q+1}_Z(X)\arrow[r]\arrow[d,"\sim"] & \cdots\\
       \cdots\arrow[r]& H^q_Z(X')\arrow[r] & H^q(X')\arrow[r] & H^q(U')\arrow[r] & H^{q+1}_Z(X')\arrow[r] & \cdots
    \end{tikzcd}$$
    where all cohomology groups are with coefficients $\ZZ/n(1)$.
    
    Writing $\Adh F$ the residue field of $R'$, by \cite[Th. 09ZI]{stacks-project} we have $H^q(X',\ZZ/n(1))=H^q(\Adh F,\ZZ/n(1))=0$ for $q\geq 1$ because $\Adh F$ is separably closed. Thus $H^q(U,\ZZ/n(1))=H^q(U',\ZZ/n(1))=H^{q+1}_Z(X',\ZZ/n(1))$ for $q\geq 1$. By \Lemma{indsmooth purity}.2 applied to the codimension $1$ closed immersion of regular schemes $\Spec k\to\Spec\O_K$ and the ind-smooth $\O_K$-scheme $X'$, we have $H^{q+1}_Z(X',\ZZ/n(1))=H^{q-1}(R',\ZZ/n)$ for $q\in\ZZ$. Thus $R^q\Psi_{K/k}(\ZZ/n(1))_{\Adh x}=H^{q-1}(R',\ZZ/n)= H^{q-1}(\Adh F,\ZZ/n)$ for $q\geq 1$.
    
    This shows $\delta:\ZZ/n\to R^1\Psi_{K/k}(\ZZ/n(1))$ is bijective, because at each stalk it is an injection of finite groups of equal orders $n$ (we do not claim that the isomorphism $R^q\Psi_{K/k}(\ZZ/n(1))_{\Adh x}=\ZZ/n$ of purity is induced by $\delta$). Also $R^q\Psi_{K/k}(\ZZ/n(1))=0$ for $q>1$, which concludes for \Situation{mixed char setup}.
\end{myproof}

\begin{theorem}[equal char 0]
    Consider $K,k,S_K,S_k,R\Psi_{K/k}$ as in any of \Situation{positive equal char setup}, \Situation{mixed char setup} or \Situation{equal char 0 setup} setup. Let $\ell\neq\characteristic k$ be a prime and $d\in\ZZ$ an arbitrary fixed integer.
    \begin{enumerate}
        \item The functor $R\Psi_{K/k}$ sends $D_\ell(S_K)$ to $D_\ell(S_k)$. If $G\in\Sh_\ell(S_K)$ then $R^q\Psi_{K/k}(G)=0$ for $q\neq 0,1$.
        \item There exists a canonical trace map $\tr:R\Psi_{K/k}(\QQ_\ell/\ZZ_\ell(d+1))\to\QQ_\ell/\ZZ_\ell(d)[-1]$ in $D(S_k)$.
        \item For any $G\in D_\ell(S_K)$ with Cartier dual $G^\vee=R\HOM_{S_K}(G,\QQ_\ell/\ZZ_\ell(d+1))$, we have a perfect pairing :
        $$R\Psi_{K/k}(G)\otimes^LR\Psi_{K/k}(G^\vee)\to R\Psi_{K/k}(\QQ_\ell/\ZZ_\ell(d+1))\xrightarrow\tr\QQ_\ell/\ZZ_\ell(d)[-1]$$
        which gives perfect Cartier dualities of finite étale $\ell$-groups for $q\in\ZZ$ :
        $$R^q\Psi_{K/k}(G)\otimes R^{1-q}\Psi_{K/k}(G^\vee)\to \QQ_\ell/\ZZ_\ell(d).$$
        \item For $G\in\Sh(S_K)$ and $q\in\ZZ$, $R^q\Psi_{K/k}(G)$ is the sheafification of $T\in S_k\mapsto H^q(h^{K/k}(T),G)$. In particular $R\Psi_{K/k}(G)_{\Adh k} = R\Gamma(K^\ur,G_{\Adh K})$. If in addition $G$ is torsion then $R\Psi_{K/k}(G)_{\Adh k} = R\Gamma(\widehat{K^\ur},G_{\Adh K})$.
    \end{enumerate}
\end{theorem}

\begin{myproof}[of \Theorem{equal char 0}]
    \textbf{4.} This is \Theorem{equal char p}.4 in \Situation{positive equal char setup}, \Lemma{interpretation char mixed}.3 in \Situation{mixed char setup}, and identical to \Theorem{equal char p}.4 in \Situation{equal char 0}. The comparison between $K^\ur$ and $\widehat{K^\ur}$ follows from \cite[Cor. 1.18.(2)]{BhattMathew}.
    
    \textbf{1.} Using \Lemma{general restriction} instead of \Lemma{mixed char restriction} and \Lemma{RPsi computation} instead of \cite{KatoSuzuki2019}, this is identical to \Theorem{mixed char}.1.
    
    \textbf{2.} This is a special case of \Lemma{RPsi computation}.
    
    \textbf{3.} By \Lemma{finite coeffs equivalence}, \Proposition{cartier duality} and \textbf{(1.)} and \textbf{(4.)} above, we can assume $S_K=\Spec K_\et$, $S_k=\Spec k_\et$ and $R\Psi_{K/k}=R\Gamma(K',-)$ where $K'=\widehat{K^\ur}$. By exactness of $(-)^\vee=\HOM_{k_\et}(-,\QQ_\ell/\ZZ_\ell(d))$, the cup-product pairing :
    $$R\Gamma(K',G)\otimes^L R\Gamma(K',G^\vee)\to R\Gamma(K',\QQ_\ell/\ZZ_\ell(d))\xrightarrow{\tr}\QQ_\ell/\ZZ_\ell(d-1)[-1]$$is perfect in $D(k_\et)$ if and only if for all $q\in\ZZ$, the pairing :
    $$H^q(K',G)\otimes H^{1-q}(K',G^\vee)\to H^1(K',\QQ_\ell/\ZZ_\ell(d+1))=\QQ_\ell/\ZZ_\ell(d)$$
    is perfect in $\Sh(k_\et)$. We can assume $G\in\Sh_\ell(k_\et)$ is concentrated in degree $0$, and we can assume $q=0$ by Lang's theorem and commutativity of the cup-product. We identify $G^\vee$ with $\HOM_{K_\et}(G,\QQ_\ell/\ZZ_\ell(d))$ concentrated in degree $0$.
    
    First consider the case $G=\Res_{L/K'}(\ZZ/\ell^n)=\ZZ/\ell^n[\Gal(L/K')]$ in $\Sh(K'_\et)$, for some $n\geq 1$ and finite Galois $L/K'$. Then $G^\vee=\ZZ/\ell^n[\Gal(L/K')]\otimes\ZZ/\ell^n(d+1)=\Res_{L/K'}(\ZZ/\ell^n(d+1))$, because $\ZZ/\ell^n(d+1)\cong\ZZ/\ell^n$ as a Galois module over $K'$. By Shapiro's lemma and \Lemma{RPsi computation} we have a diagram of pairings :
    $$\begin{tikzcd}
        H^0(K',G)&[-3em] \otimes&[-3em] H^1(K',G^\vee)\arrow[r]\arrow[d,"\sim"]& H^1(K',\QQ_\ell/\ZZ_\ell(d)) \arrow[d,"\sim"] \\
        H^0(L,\ZZ/\ell^n)\arrow[u,"\sim"]&[-3em] \otimes&[-3em]  H^1(L,\ZZ/\ell^n(d+1))\arrow[r]\arrow[d,"\sim"]& H^1(L,\QQ_\ell/\ZZ_\ell(d+1))\arrow[d,"\sim"]\\
        \ZZ/\ell^n \arrow[u,"\sim"]&[-3em]\otimes&[-3em] \ZZ/\ell^n(d)\arrow[r]& \QQ_\ell/\ZZ_\ell(d)
    \end{tikzcd}$$
    where the top and middle horizontal arrows are given by the cup-product, and the bottom arrow, given by multiplication, is a perfect pairing. The top half of the diagram commutes by the relation $x\cup\cores(y)=\cores(\res(x)\cup y)$ for cup products in group cohomology. The commutativity of the bottom half reduces to $d=1$ (because $\ZZ/\ell^n(d)$ is a trivial Galois module over $K'$ and $H^q(L,\ZZ/\ell^n(d))=H^q(L,\ZZ/\ell^n(1))\otimes\ZZ/\ell^n(d-1)$, since $L$ has all $\ell^n$-th roots of unity), which comes from the commutative diagrams given by Kummer theory and the valuation $v:L^\times/n\to\ZZ/n$ :
    $$\begin{tikzcd}
        H^0(L,\ZZ/\ell^n)&[-3em]\otimes&[-3em]  H^1(L,\ZZ/\ell^n(1))\arrow[r]\arrow[d,"\sim"]& H^1(L,\ZZ/\ell^n(1))\arrow[d,"\sim"]\\
        \ZZ/\ell^n \arrow[u,"\sim"]&[-3em]\otimes&[-3em] L^\times/\ell^n\arrow[r,"{(s,x)\mapsto x^s}"]& L^\times/\ell^n
    \end{tikzcd}
    \qquad
    \begin{tikzcd}
        \ZZ/\ell^n&[-3em]\otimes&[-3em] L^\times/\ell^n \arrow[r,"{(s,x)\mapsto x^s}"]\arrow[d,"v"]& L^\times/\ell^n\arrow[d,"v"]\\
        \ZZ/\ell^n\arrow[u,"\sim"]&[-3em]\otimes&[-3em] \ZZ/\ell^n \arrow[r,"{(s,t)\mapsto st}"]& \ZZ/\ell^n
    \end{tikzcd}$$
    Hence $H^q(K',G)\otimes H^{1-q}(K',G^\vee)\to H^1(K',\QQ_\ell/\ZZ_\ell(d+1))$ is isomorphic to a perfect pairing, thus perfect.
    
    Now consider $G\in\Sh_\ell(K'_\et)$ general. There is a finite Galois extension $L/K'$ over which $G$ is trivial. Write $P=\Res_{L/K'}(G|_L)$. Then we have a short exact sequence, for some $Q\in\Sh_\ell(K_\et)$ :
    $$0\to G\to P\to Q\to 0$$
    which gives diagrams with exact rows, commutative up to signs :
    $$\begin{tikzcd}
        0\arrow[r]& H^0(K',G) \arrow[r]\arrow[d,"\phi_G"] & H^0(K',P) \arrow[r]\arrow[d,"\phi_P"] & H^0(K',Q)\arrow[d,"\phi_Q"]\\
        0\arrow[r]& H^1(K',G^\vee_K)^\vee_k \arrow[r] & H^1(K',P^\vee_K)^\vee_k \arrow[r] & H^1(K',Q^\vee_K)^\vee_k
    \end{tikzcd}$$
    where $(-)^\vee_K$ and $(-)^\vee_k$ denote Cartier duality over $\Spec K_\et$ and $\Spec k_\et$ respectively. We have seen that $\phi_P$ is an isomorphism. The map $\phi_P$ is injective, so $\phi_G$ is injective. This applies to any $G\in\Sh_\ell(k_\et)$, so $\phi_Q$ is also injective. By the four lemma $\phi_G$ is surjective, thus an isomorphism as desired.
\end{myproof}


\subsection{Higher local fields}
\label{Higher local fields}

We assemble the previous steps to get results over higher local fields, by which we mean the following.

\begin{definition}[higher local field]
    For $d\geq 0$, a \emph{$d$-local field} is a field $K$ equipped with a sequence of fields $K=k_d,\dots,k_0$ where $k_0$ is perfect and for all $i>0$, $k_i$ is a Henselian discrete valuation field with residue field $k_{i-1}$, and $k_i$ is also complete whenever $\characteristic k_i>0$. The field $k_i$ is referred to as the \emph{$i$-local residue field of $K$}. We say $K$ has \emph{totally equal characteristic $p$} if $\characteristic k_d=\characteristic k_0=p$, or \emph{mixed characteristic $(0,p)$ at level $s$} for $1\leq s\leq d$ if $\characteristic k_d=\characteristic k_s=0$ and $\characteristic k_{s-1}=p>0$.
\end{definition}

If $K$ is a Henselian discrete valuation field and $L/K$ a finite extension, then the absolute value of $K$ extends uniquely to $L$. Thus if $L=l_d$ is a finite extension of a $d$-local field $K=k_d$, then by induction on $d$ one sees that $l_d$ is again $d$-local, and its $i$-local residue field $l_i$ is a finite extension of $k_i$ for all $0\leq i\leq d$. There exist integers $e_1,\dots,e_d\geq 1$ such that the valuation $v_{l_i}$ restricts to $e_iv_{k_i}$ on $k_i$, called the ramification indices of $L/K$, which satisfy\footnote{This is true more generally if $K=k_d,\dots,k_0$ is such that $k_i$ is a Henselian discrete valuation field with residue field $k_{i-1}$ and $\hat k_i/k_i$ is separable for each $i>0$, or equivalently $\O_{k_i}$ is excellent : see \cite[Prop. 10 of Sec. 1.4]{SerreLC}.} $e_d\cdots e_1[l_0:k_0]= [L:K]$. We say $L/K$ is \emph{totally unramified} if each $l_i/k_i$ is separable and $e_1=\dots=e_d=1$. The union of all totally unramified finite extensions of $K$ in a fixed separable closure is called the \emph{maximal totally unramified extension} of $K$, and written $K^\tur$.

Recall the étale or Kato lift $h^{k_{i+1}/k_i}(R_i)$ of \Definition{Kato lift definition}. If $R_i$ is an étale $k_i$-algebra, $h^{k_{i+1}/k_i}(R_i)$ is étale over $k_{i+1}$. If $k_{i+1}$ is complete of equicharacteristic $p>0$ and $R_i$ is relatively perfect over $k_i$, $h^{k_{i+1}/k_i}(R_i)$ is relatively perfect over $k_{i+1}$ by \Remark{Kato lift characterizations}. Hence we can iterate this construction as follows.

\begin{definition}[higher Kato lift]
    Let $K=k_d,\dots,k_0$ be a $d$-local field. Let $R_i$ be a $k_i$-algebra for some $0\leq i<d$.
    \begin{itemize}
        \item If $R_i$ is étale over $k_i$, define :
        $$h^{K/k_i}(R_i)=h^{k_d/k_{d_1}}\circ\cdots\circ h^{k_{i+1}/k_i}(R_i)$$
        the étale $K$-algebra obtained by successive étale lifts, called the \emph{iterated étale lift of $R_i$ over $K$}.
        \item If $\characteristic k_{d-1}>0$ and $R_i$ is relatively perfect over $k_i$, define :
        $$h^{\hat K/k_i}(R_i)=h^{\hat k_d/k_{d_1}}\circ\cdots\circ h^{k_{i+1}/k_i}(R_i)$$
        the complete algebra over $\hat K=\hat k_d$, the completion of $k_d$ with respect to its discrete valuation, obtained by successive Kato lifts, called the \emph{iterated Kato lift of $R_i$ over $\hat K$}.
    \end{itemize}
\end{definition}

Note that the iterated étale lift makes sense for $K$ a general higher local field (with Henselian discrete valuation fields of arbitrary characteristics at each step) but only étale $R_i$. In contrast the iterated Kato lift works for any relatively perfect $R_i$, but requires complete fields at each step and stops at the first level of characteristic $0$ ; in particular, we need the convention of \Definition{higher local field} that $k_i$ is complete whenever it has equal characteristic\footnote{If $K=k_d,\dots,k_0$ is a higher local field with mixed characteristic $(0,p)$ at level $0<s<d$, and each $k_i$ is complete, we can still lift a relatively perfect $k_0$-algebra $R_0$ to complete $k_i$-algebras $h^{k_i/k_0}(R_0)$ for each $1\leq i\leq d$. For $i\leq s$ this is Kato's construction, which is canonical. For $i>s$ the lift may not be unique, and depends on a choice of sections $k_{j-1}\to\O_{k_j}$ for $j>s$. Given such a choice we can \textit{ad hoc} define $h^{k_i/k_{i-1}}(-) = k_i\otimes_{\O_{k_i}} (-\otimes_{k_{i-1}}\O_{k_i})^\wedge$ to define the iterated lift.}. When $K=\hat K$, $\characteristic k_{d-1}>0$, and $R_i$ is étale over $k_i$, the two lifts coincide.

\begin{situation}[higher local setup]
    Let $K=k_d,\dots,k_0$ be a $d$-local field. Assume $d\geq 2$ and $K$ has mixed characteristic $(0,p)$ at level $2$. Consider $X$ a proper, smooth, geometrically integral $K$-variety of dimension $e\geq 0$ and $U\subseteq X$ a nonempty open. We define functors $R\Psi_{U/k_0},R\Psi_{U/k_0,c}:D(U_\et)\to D(k_{0,\et}^\perar)$ by :
    \begin{align*}
            R\Psi_{U/k_0} &=R\Psi_{k_1/k_0}^p\circ Lv^*\circ R\Psi_{k_2/k_1}^m\circ u^*\circ R\Psi_{k_3/k_2}^0\circ\cdots\circ R\Psi_{k_d/k_{d-1}}^0\circ R\Psi_{U/K}^v,\\
            R\Psi_{U/k_0,c}&=R\Psi_{k_1/k_0}^p\circ Lv^*\circ R\Psi_{k_2/k_1}^m\circ u^*\circ R\Psi_{k_3/k_2}^0\circ\cdots\circ R\Psi_{k_d/k_{d-1}}^0\circ R\Psi_{U/K,c}^v,
    \end{align*}
    with $u:\Spec k_{2,\Et}\to\Spec k_{2,\et}$ and $v:\Spec k_{1,\RP}\to\Spec k_{1,\RPS}$ the premorphisms defined by identity.
\end{situation}

\begin{theorem}[higher local duality]
    Consider \Situation{higher local setup}. For $G\in D(U_\et)$, define $G^\vee=R\HOM_{U_\et}(G,\QQ/\ZZ(d))$.
    \begin{enumerate}
        \item The functors $R\Psi_{U/k_0}$ and $R\Psi_{U/k_0,c}$ send $D_p(U_\et)$ to $D^b_{\W_{k_0}}(k_{0,\et}^\perar)$, and $D_\ell(U_\et)$ to $D_\ell(k_{0,\et}^\perar)$ for $\ell\neq p$ prime. If $G\in\Sh_\fin(U_\et)$ then $R^q\Psi_{U/k_0}(G)=R^q\Psi_{U/k_0,c}(G)=0$ for $q\neq 0,\dots,d+2e$.
        \item There exists a canonical trace morphism $\tr:R\Psi_{U/k_0,c}(\QQ/\ZZ(d+e))\to\QQ/\ZZ[-d-2e]$ in $D(k_{0,\et}^\perar)$.
        \item For $G\in D_\fin(U_\et)$, we have a perfect pairing in $D(k_{0,\et}^\perar)$ for $q\in\ZZ$ :
        $$R\Psi_{U/k_0}(G)\otimes^LR\Psi_{U/k_0,c}(G^\vee)\to R\Psi_{U/k_0,c}(\QQ/\ZZ(d+e))\xrightarrow\tr\QQ/\ZZ[-d-2e].$$
        \item If $G\in D_\ell(U_\et)$, we have perfect Cartier dualities of finite étale $\ell$-primary $k_0$-groups for $q\in\ZZ$ :
        $$R^q\Psi_{U/k_0}(G)\otimes R^{d+2e-q}\Psi_{U/k_0,c}(G^\vee)\to R^{d+2e}\Psi_{U/k_0,c}(\QQ_\ell/\ZZ_\ell(d+e))\xrightarrow\tr\QQ_\ell/\ZZ_\ell$$
        \item If $G\in D_p(U_\et)$, we have perfect Serre dualities of connected and finite étale parts :
        \begin{align*}
            R^q\Psi_{U/k_0}(G)^0&=\EXT^1_{k_{0,\et}^\perar}(R^{d+2e+1-q}\Psi_{U/k_0,c}(G^\vee)^0,\QQ_p/\ZZ_p),\\ \pi_0(R^q\Psi_{U/k_0}(G))&=\HOM_{k_{0,\et}^\perar}(\pi_0(R^{d+2e-q}\Psi_{U/k_0,c}(G^\vee)),\QQ_p/\ZZ_p).
        \end{align*}
    \end{enumerate}
\end{theorem}

\begin{myproof}[of \Theorem{higher local duality}]
    The factors of the composite functors $R\Psi_{U/k_0}$ and $R\Psi_{U/k_0,c}$ each satisfy dualities for $p$-torsion coefficients : $R\Psi_{U/K}^v$ and $R\Psi_{U/K,c}^v$ by \Theorem{poincare duality}, $R\Psi_{k_i/k_{i-1}}^0$ for $i>2$ by \Theorem{equal char 0}, $Lv^*\circ R\Psi_{k_2/k_0}^m\circ u^*$ by \Theorem{mixed char} and \Remark{for small et}, and $R\Psi_{k_1/k_0}^p$ by \Theorem{equal char 0}. They each satisfy dualities for $\ell$-torsion coefficients by \Theorem{poincare duality}, \Theorem{equal char 0}, \Lemma{finite coeffs equivalence} and \Theorem{cartier duality}.4. Perfect pairings of objects of $D_\ell(k_{0,\et}^\perar)$ and $D_{\W_k}(k_{0,\et}^\perar)$ decompose as in \textbf{(4.-5.)} by \Theorem{cartier duality}.1 and \Theorem{Serre duality}.5. 
    
    The theorem follows, except we only have $R^q\Psi_{U/k_0}(G)=R^q\Psi_{U/k_0,c}(G)=0$ for $q\neq 0,\dots,d+2e+1$ if $G\in \Sh_\fin(U_\et)$ has $p$-torsion. This is because the complex $R\Psi_{k_2/k_0}(G)=R\Psi_{k_1/k_0}^p\circ Lv^*\circ R\Psi_{k_2/k_1}^m\circ u^*(G)$ is \textit{a priori} concentrated in degrees $0,1,2,3$ and not $0,1,2$. We show $R^3\Psi_{k_2/k_0}(G)=0$ also for $G\in\Sh_p(k_{2,\et})$.
    
    Let $T\in\Spec k_{0,\et}^\perar$ and take $\Adh x\to T$ a geometric point. In \Proposition{char 0pp-0p interpretation} below we will see that :
    $$R^3\Psi_{k_2/k_0}(G)_{\Adh x}=\varinjlim_{T'/T} H^3(h^{\hat k_2/k_0}(R'),G)$$
    where $T'$ ranges over factorizations $\Adh x\to T'\to T$ with $T'\in k_{0,\et}^\perar$ étale over $T$. Since $h^{\hat k_2/k_0}$ commutes with finite products, we can assume $T=\Spec F$ with $F$ is a perfect field, $\Adh x$ is given by an algebraic closure $\Adh F$, and $T'$ runs over finite field $\Spec F'\to\Spec F$. Then :
    $$R^3\Psi_{k_2/k_0}(G)_{\Adh x}=H^3(h^{\hat k_2/k_0}(F)^\tur,G)$$
    where $M_2=h^{\hat k_2/k_0}(F)^\tur$ is a Henselian discrete valuation field with residue field $M_1=h^{k_1/k_0}(F)^\ur$, itself a Henselian discrete valuation field with algebraically closed residue field $M_0=\Adh F$. Clearly $M_0$ has cohomological dimension $0$ ; then inductively $M_i$ has cohomological dimension at most $\cd(\widehat{M_i})=\cd(M_{i-1})+1\leq i$ by \cite[Cor. 1.18.(2)]{BhattMathew} and \cite[Cor. of Th. 3]{Kato82}. Thus $R^3\Psi_{k_2/k_0}(G)_{\Adh x}=0$. This holds for any geometric point $\Adh x$ of $T\in\Spec k_{0,\et}^\perar$, thus $R^3\Psi_{k_2/k_0}(G)=0$ as desired.
\end{myproof}

\begin{remark}[same for indrat proet]
    Since the composite $\alpha_{k_0}R\Psi_{k_1/k_0}^p: D(k_{1,\RP})\to D(k_{0,\proet}^\indrat)$ also satisfies a duality by \Remark{for indrat proet}, where $\alpha_{k_0}:D(k_{0,\et}^\perar)\to D(k_{0,\proet}^\indrat)$ is the change of site functor, \Theorem{higher local duality} also holds if we replace $\Spec k_{0,\et}^\perar$ by $\Spec k_{0,\proet}^\indrat$ and $R\Psi_{U/k_0}$ and $R\Psi_{U/k_0,c}$ by their composites with $\alpha_{k_0}$.
\end{remark}

\begin{theorem}[higher local duality finite]
    Consider \Situation{higher local setup}. Assume $k_0$ is finite. Define functors $D(U_\et)\to D(*_\proet)$ :
    $$R\Psi_{U/*}=R\Psi_{k_0/*}\circ \alpha_{k_0}\circ R\Psi_{U/k_0},\qquad R\Psi_{U/*,c}=R\Psi_{k_0/*}\circ \alpha_{k_0}\circ R\Psi_{U/k_0,c}$$
    where $R\Psi_{k_0,*}:D(k_{0,\proet}^\indrat)\to D(*_\proet)$ is the functor of \Situation{finite field}.
    \begin{enumerate}
        \item The functors $R\Psi_{U/*}$ and $R\Psi_{U/*,c}$ send $D_p(U_\et)$ to $D^b_{\W_p}(*_\proet)$, and $D_\ell(U_\et)$ to $D_\ell(*_\proet)$ for $\ell\neq p$ prime. If $G\in\Sh_\fin(U_\et)$ then $R^q\Psi_{U/*}(G)=R^q\Psi_{U/*,c}(G)=0$ for $q\neq 0,\dots,d+2e+1$. In particular $R^q\Psi_{U/*}(G)$ and $R^q\Psi_{U/*,c}(G)$ are representable by locally compact Hausdorff groups for $q\in\ZZ$.
        \item There exists a canonical trace morphism $\tr:R\Psi_{U/*,c}(\QQ/\ZZ(d+e))\to\QQ/\ZZ[-d-2e-1]$ in $D(*_\proet)$.
        \item For $G\in D_\fin(U_\et)$ and $G^\vee=R\HOM_{U_\et}(G,\QQ/\ZZ(d))$ we have a perfect pairing in $D(*_\proet)$ :
        $$R\Psi_{U/*}(G)\otimes^LR\Psi_{U/*,c}(G^\vee)\to R\Psi_{U/*,c}(\QQ/\ZZ(d+e))\xrightarrow\tr\QQ/\ZZ[-d-2e-1]$$
        which induces perfect Pontryagin dualities of locally compact Hausdorff groups for $q\in\ZZ$ :
        $$R^q\Psi_{U/*}(G)\otimes R^{d+2e+1-q}R\Psi_{U/*,c}(G^\vee)\to R^{d+2e+1}R\Psi_{U/*,c}(\QQ/\ZZ(d+e))\xrightarrow\tr\QQ/\ZZ.$$
    \end{enumerate}
\end{theorem}

\begin{myproof}[of \Theorem{higher local duality finite}]
    This follows from \Theorem{higher local duality}, \Remark{same for indrat proet} and \Proposition{topological cohomology} for $p$-primary coefficients. For $G\in\Sh_\ell(k_\proet^\indrat)$ and $q\in\ZZ$, $R^q\Psi_{k_0/*}(G)$ is the sheaf representable by the finite group $H^q(k_\et,G)$. We have $H^q(k_\et,G)=0$ if $q\neq 0,1$, $H^1(k_\et,\QQ/\ZZ)=\QQ/\ZZ$, and perfect pairings :
    $$H^q(k_\et,G)\otimes H^{1-q}(k_\et,G^\vee)\to H^1(k_\et,\QQ/\ZZ)=\QQ/\ZZ$$
    for all $q\in\ZZ$ where $G^\vee=\HOM_{k_{0,\proet}^\indrat}(G,\QQ_\ell/\ZZ_\ell)$ is the Cartier dual of $G$. The corresponding properties for $R\Psi_{k_0,*}$ follow, and so does the theorem from composing with \Theorem{higher local duality}.
\end{myproof}

It remains to give arithmetic interpretations of \Theorem{higher local duality} and \Theorem{higher local duality finite}. The aim of the rest of the section is to relate the functors $R\Psi_{U/k_0}$, $R\Psi_{U/*}$, $R\Psi_{U/k_0,c}$, $R\Psi_{U/*,c}$ to the étale cohomology groups of $U$.

\begin{lemma}[monoidality premorphism]
    Let $v:S'\to S$ be a premorphism of sites.
    \begin{enumerate}
        \item For $F,F'\in D(S)$, there is a canonical morphism $Lv^*(F\otimes^LF')\to Lv^*F\otimes^LLv^*F'$ in $D(S')$. If $S'$ has finite products and $v^{-1}:S'\to S$ preserves those products, then this is an isomorphism.
        \item For $G,G'\in D(S')$ the cup-product $Rv_*G\otimes^LRv_*G'\to Rv_*(G\otimes^L G')$ corresponds under the adjunction of $Lv^*$ and $Rv_*$ to the composite :
        $$Lv^*(Rv_*G\otimes^LRv_*G')\xrightarrow{\text{\textbf{(1.)}}} (Lv^*Rv_*G)\otimes^L(Lv^*Rv_*G')\xrightarrow{\eps_G\otimes\eps_{G'}}G\otimes^LG'$$
        where $\eps:Lv^*Rv_*\to\id$ is the counit in $D(S')$.
    \end{enumerate}
\end{lemma}

\begin{myproof}[of \Lemma{monoidality premorphism}]
    \textbf{1.} The cup-product of $Rv_*$ and the unit $F'\to Rv_*Lv^*F'$ define a natural map $Rv_*R\HOM_{S'}(Lv^*F',T)\to R\HOM_S(F',Rv_*T)$ for $T\in D(S')$. By adjunctions, we get a morphism :
    \begin{align*}
        \Hom_{D(S')}((Lv^*F)\otimes^L (Lv^*F'),T)
        &= \Hom_{D(S)}(F,Rv_*R\HOM_{D(S')}(Lv^*F',T))\\
        &\to \Hom_{D(S)}(F,R\HOM_{D(S)}(F',Rv_*T))\\
        &=\Hom_{D(S')}(Lv^*(F\otimes^LF'),T).
    \end{align*}
    hence a canonical morphism $Lv^*(F\otimes^LF')\to (Lv^*F)\otimes^L (Lv^*F')$ by Yoneda's lemma. This is an isomorphism whenever the map $Rv_*R\HOM_{S'}(Lv^*F',T)\to R\HOM_S(F',Rv_*T)$ is an isomorphism, which is the case under the given additional assumption by \cite[Prop. 3.1.(1)]{SuzukiNeron}.
    
    \textbf{2.} Working out the definition of the morphism of \textbf{(1.)} for $F,F'\in D(S)$ reveals it corresponds, under the adjunction of $Lv^*$ and $Rv_*$, to the composite :
    $$F\otimes^L F'
    \xrightarrow{\eta_F\otimes^L\eta_{F'}}
    (Rv_*Lv^*F)\otimes^L(Rv_*Lv^*F'))
    \xrightarrow{\cup_{Lv^*F,Lv^*F'}}
    Rv_*((Lv^*F)\otimes^L(Lv^*F'))$$
    where $\cup_{G,G'}:Rv_*G\otimes^LRv_*G'\to Rv_*(G\otimes^LG')$ is the cup-product in $D(S)$ for $G,G'\in D(S')$, and $\eta_F:F\to Rv_*Lv^*F$ is the unit in $D(S)$. Thus we want a commutative square, for $G,G'\in D(S')$ :
    $$\begin{tikzcd}[column sep=60]
        Rv_*G\otimes^LRv_*G'\arrow[r,"\eta\otimes\eta"] & (Rv_*Lv^*Rv_*G)\otimes^L(Rv_*Lv^*Rv_*G') \arrow[r,"\cup"]\arrow[d,"Rv_*\eps\otimes^LRv_*\eps"]& Rv_*(Lv^*Rv_*G\otimes^LLv^*Rv_*G')\arrow[d,"Rv_*(\eps\otimes^L\eps)"]\\
        Rv_*G\otimes^LRv_*G'\arrow[u,equals]\arrow[r,equals]& Rv_*G\otimes^LRv_*G'\arrow[r,"\cup_{G,G'}"] & Rv_*(G\otimes^LG)
    \end{tikzcd}$$
    where $\eps:Lv^*Rv_*\to\id$ is the counit. The left-hand side commutes by unit-counit relations. The commutativity of the square :
    $$\begin{tikzcd}[column sep=60]
        (Rv_*Lv^*Rv_*G)\otimes^L(Rv_*Lv^*Rv_*G')\arrow[d,"Rv_*\eps_G\otimes^LRv_*\eps_{G'}"] \arrow[r,"\cup_{Lv^*Rv_*G,Lv^*Rv_*G'}"]& Rv_*(Lv^*Rv_*G\otimes^LLv^*Rv_*G')\arrow[d,"Rv_*(\eps_G\otimes^L\eps_{G'})"]\\
        Rv_*G\otimes^LRv_*G'\arrow[r,"\cup_{G,G'}"] & Rv_*(G\otimes^LG)
    \end{tikzcd}$$
    follows from the naturality of the cup-product $\cup_{G,G'}$ in $G$ and $G'$, clear by definition \cite[Prop. 2.4]{SuzukImp}.
\end{myproof}

\begin{proposition}[char 0pp-0p interpretation]
    Consider \Situation{higher local setup}.
    \begin{enumerate}
        \item Assume $d=2$ and $X=U=\Spec K$. Let $R\tilde\Psi_{K/k_0}:D(K_\Et)\to D(k_{0,\et}^\perar)$ be the derived pushforward of the premorphism $T/k_0\mapsto h^{\hat K/k_0}(T)$. Then for $G\in D_\fin(K_\et)$ we have :
        $$R\Psi_{K/k_0}(G) = R\tilde\Psi_{K/k_0}\circ u^*(G).$$
        In particular, $R^q\Psi_{K/k_0}(G)$ is the sheafification of $T\in \Spec k_{0,\et}^\perar\mapsto H^q(h^{\hat K/k_0}(T),G)$ for $q\in\ZZ$.
        \item In the situation of \textbf{(1.)}, we also have a commutative diagram of pairings :
        $$\begin{tikzcd}
            R\Psi_{K/k_0}(G)\arrow[d,"\sim"]&[-3em]\otimes^L&[-3em]R\Psi_{K/k_0}(G^\vee)\arrow[r]& R\Psi_{K/k_0}(\QQ/\ZZ(2))\arrow[d,"\sim"]\\
            u^*R\tilde\Psi_{K/k_0}(G)&[-3em]\otimes^L&[-3em]u^*R\tilde\Psi_{K/k_0}(G^\vee)\arrow[r]\arrow[u,"\sim"]& u^*R\tilde\Psi_{K/k_0}(\QQ/\ZZ(2))
        \end{tikzcd}$$
        where $G^\vee=R\HOM_{K_\et}(G,\QQ/\ZZ(2))$, the upper pairing is as in \Theorem{higher local duality} and the lower pairing is given by monoidality of $u^*$ and the cup-product for the derived pushforward $R\tilde\Psi_{K/k_0}$. 
        \item Consider the general case. Let $R\Psi^\et_{U/k_0}:D(U_\et)\to D(k_{0,\et})$ be the derived pushforward of the premorphism $T/k_0\mapsto U_{h^{K/k_0}(T)}$ and write $R\Psi^\et_{U/k_0,c}=R\Psi^\et_{X/k_0}\circ j_!$. Let $\beta_{k_0}:\Spec k_{0,\et}^\perar\to\Spec k_{0,\et}$ be the premorphism of sites defined by identity. Then for $G\in D_\fin(K_\et)$ we have :
        $$\beta_{k_0,*}\circ R\Psi_{U/k_0}(G) = R\Psi^\et_{U/k_0}(G),\qquad \beta_{k_0,*}\circ R\Psi_{U/k_0,c}(G) = R\Psi^\et_{U/k_0,c}(G).$$
        In particular we have identifications in $\Sh(k_{0,\et})$, for $q\in\ZZ$ :
        \begin{align*}
            R^q\Psi_{U/k_0}(G)_{\Adh k}&=H^q(U_{K^\tur},G|_{U_{K^\tur}})\\
            R^q\Psi_{U/k_0,c}(G)_{\Adh k}&=H^q_c(U_{K^\tur},G|_{U_{K^\tur}})=H^q(X_{K^\tur},j_{K^\tur,!}(G|_{U_{K^\tur}})).
        \end{align*}
        where $(-)|_{U_{K^\tur}}$ is the pullback to $D((U_{K^\tur})_\et)$ and $j_{K^\tur}:U_{K^\tur}\to X_{K^\tur}$ is the base change of $j$.
        \item In the situation of \textbf{(3.)}, we also have a commutative diagram of pairings for all $q,s\in\ZZ$ :
        $$\begin{tikzcd}
            R^q\Psi_{U/k_0}(G)_{\Adh k}\arrow[d,"\sim"]&[-3em]\otimes^L&[-3em]R^s\Psi_{U/k_0,c}(G^\vee)_{\Adh k}\arrow[r]& R^{q+s}\Psi_{U/k_0,c}(\QQ/\ZZ(d+e))_{\Adh k}\arrow[d,"\sim"]\\
            H^q(U_{K^\tur},G|_{U_{K^\tur}})&[-3em]\otimes^L&[-3em]H^s_c(U_{K^\tur},G^\vee|_{U_{K^\tur}})\arrow[u,"\sim"]\arrow[r]& H^{q+s}_c(U_{K^\tur},\QQ/\ZZ(d+e))
        \end{tikzcd}$$
        where $G^\vee=R\HOM_{U_\et}(G,\QQ/\ZZ(d+e))$, the upper pairing is of \Theorem{higher local duality} and the lower pairing is the cup product of étale cohomology.
    \end{enumerate}
\end{proposition}

\begin{myproof}[of \Proposition{char 0pp-0p interpretation}]
    \textbf{1.} Since $u^*$ defines an equivalence $D_\fin(K_\et)=D_\fin(K_\Et)$, we can identify $G$ with an object of $D_\fin(K_\Et)$ and it suffices to prove :
    $$R\Psi_{k_1/k_0}^p\circ Lv^*\circ R\Psi_{K/k_1}^m(G) = R\tilde \Psi_{K/k_0}(G).$$
    We reproduce an argument shared by Suzuki in a private communication. The idea is to eliminate the change of site step $Lv_{k_1}^*:D(k_{1,\RPS})\to D(k_{1,\RPS})$ by modifying $R\Psi_{K/k_1}^m$ and $R\Psi_{k_1/k_0}^p$ to pass through the same site $\Spec k_{1,\PRPS}$. Define $\Spec k_{1,\PRPS}$ as the category of $k_1$-schemes which are cofiltered limits of relatively perfectly smooth schemes with affine transition morphisms, with all $k_1$-scheme morphisms and the étale topology. Define premorphisms :
    $$\Spec k_{1,\RP}\xrightarrow{v_\P}\Spec k_{1,\PRPS}\xrightarrow w \Spec k_{1,\RPS}$$
    by identity on the underlying categories. Thus $v=w\circ v_\P$.
    
    For $F$ a sheaf of sets over $\Spec k_{1,\RPS}$ and $X=\varprojlim_iX_i\in \Spec k_{1,\PRPS}$, we have $w^*F(X)=\varinjlim_iF(X_i)$, from which it follows that the pullback of sets for $w$ is exact, and $w$ is a morphism of sites. Clearly $v_{\P,*}$ is exact and $v_{\P,*}v_\P^*=\id$ : in fact $v_{\P,*}F(X)=F(X)$ for $F\in\Sh(k_{1,\RP})$ and $X\in\Spec k_{1,\PRPS}$, and $v_\P^*F(X)=F(X)$ for $F\in\Sh(k_{1,\PRPS})$ and $X\in\Spec k_{1,\RP}$. Hence we also have $v_{\P,*}Lv_\P^*=\id$.
    
    By \cite[proof of Prop. 5.4]{SuzukiCFT}, for $T\in\Spec k_{0,\et}^\perar$, the $k_1$-scheme $(\psi_{k_1/k_0}^p)^{-1}(T)=h^{k_1/k_0}(T)$ belongs to $\Spec k_{1,\PRPS}$. Hence the premorphism $\psi_{k_1/k_0}:\Spec k_{1,\RP}\to \Spec k_{0,\et}^\perar$ factors through $v_P$ as $\psi_{k_1/k_0,\P}:\Spec k_{1,\PRPS}\to \Spec k_{0,\et}^\perar$ and we have :
    $$R\Psi_{k_1/k_0}^pLv^* = R\Psi_{k_1/k_0,\P}^pv_{\P,*}Lv_\P^*w^* = R\Psi_{k_1/k_0,\P}^pw^*$$
    where $R\Psi_{k_1/k_0,\P}^p:D(k_{1,\PRPS})\to D(k_{0,\et}^\perar)$ is the derived pushforward of $\psi_{k_1/k_0,\P}$.
    
    Define $\Spec \O_{K,\PRPS}$ as the category of flat $\O_K$-schemes with special fiber in $\Spec k_{1,\PRPS}$, with all $\O_K$-scheme morphisms and the étale topology. Define $w_\O:\Spec \O_{K,\PRPS}\to\Spec \O_{K,\RPS}$ the premorphism defined by identity. We have a commutative diagram of premorphisms of sites where all the horizontal maps defined by base change :
    $$\begin{tikzcd}
       \Spec K_\Et\arrow[r,"j_\P"]\arrow[d,equals] & \Spec \O_{K,\PRPS}\arrow[d,"w_\O"] & \Spec k_{1,\PRPS}\arrow[l,"i_\P"']\arrow[d,"w"]\\
       \Spec K_\Et\arrow[r,"j"] & \Spec \O_{K,\RPS}&\arrow[l,"i"']\Spec k_{1,\RPS}
    \end{tikzcd}$$
    which induces a natural map $Lw_\O^*Rj_*\to Rj_{\P,*}$. Writing $R\Psi_{K/k_1,\P}^m=Li_\P^*Rj_{\P,*}$, we get a natural map :
    $$w^*R\Psi_{K/k_1}^m = Li_\P^*Lw_\O^*Rj_*\to R\Psi_{K/k_1,\P}^m.$$
    
    By \cite[Cor. 3.1]{KatoSuzuki2019} the functor $A\mapsto A\otimes_{\O_K}k_1$ from the category of flat $\O_K$-algebras with relatively perfect special fibers is left adjoint to $R\mapsto h^{\O_K/k_1}(R)$. In particular the restriction of $i_\P^{-1}:\Spec k_{1,\PRPS}\to\Spec \O_{K,\PRPS}$ to affines has a right adjoint described by the Kato lift, so $i_\P^*$ coincides on affine objects with the pushforward of the premorphism $T/k_1\mapsto h^{\O_K/k_1}(T)$. Hence $i_\P^*$ is exact and $R\Psi_{K/k_1,\P}^m=i_\P^*Rj_{\P,*}$.
    
    This description of $i_\P^*$ also shows that $R^q\Psi_{K/k_1,\P}^m=H^qR\Psi_{K/k_1,\P}$ is the $q$-th left derived functor of $R^0\Psi_{K/k_1,\P}^m$ for $q\in\ZZ$, and by \Lemma{pushforward of sites}, $R^q\Psi_{K/k_1,\P}^m(F)$ coincides on affine objects with the sheafification of the presheaf $T/k_1\mapsto H^q(h^{K/k_1(T)},F)$ of $\Spec k_{1,\PRPS}$, for $F\in\Sh(K_\Et)$ and $q\in\ZZ$.
    
    We show $w^*R\Psi_{K/k_1}^m(G)\to R\Psi_{K/k_1,\P}^m(G)$ is an isomorphism for $G\in D_\fin(K_\Et)$. It suffices to consider $G\in\Sh_\fin(K_\Et)$. Let $T$ be an affine object of $\Spec k_{1,\PRPS}$ written as $T=\varprojlim_\lambda T_\lambda$ with $T_\lambda\in \Spec k_{1,\RPS}$. Let $x=\varprojlim_\lambda x_\lambda$ be a point of $T$, and $\Adh x=\varprojlim_\lambda\Adh x_\lambda$ a separable closure. Let $q\in\ZZ$. Then :
    \begin{align*}
        (w^*R^q\Psi_{K/k_1}^m(G))_{\Adh x}
        &= \varinjlim_\lambda R^q\Psi_{K/k_1}^m(G)_{\Adh x_\lambda} &\text{by description of }w^*,\\
        &= \varinjlim_\lambda \varinjlim_{S_\lambda/T_\lambda}H^q(h^{K/k_1}(S_\lambda),G)&\text{by \Lemma{interpretation char mixed}.2,}\\
        &= H^q(\varinjlim_\lambda\varinjlim_{S_\lambda/T_\lambda} h^{K/k_1}(S_\lambda),G)&\text{by \cite[Lem. 1.16 of Ch. III]{Milne80},}\\
        &=H^q(h^{K/k_1}(\varinjlim_\lambda\varinjlim_{S_\lambda/T_\lambda} S_\lambda),G)&\text{by \cite[Cor. 1.18.(2)]{BhattMathew} and \Lemma{Kato lift colimit},}\\
        &=H^q(h^{K/k_1}(T_x^{sh}),G)\\
        &=\varinjlim_{S/T} H^q(h^{K/k_1}(S),G) & \text{similarly,}\\
        &=(R^q\Psi_{K/k_1,\P}^m(G))_{\Adh x}&\text{by description of }R\Psi_{K/k_1,\P}^m
    \end{align*}
    where $S_\lambda$ ranges over factorizations $\Adh x_\lambda\to S_\lambda\to T_\lambda$ with $S_\lambda\in\Spec k_{0,\et}^\perar$ étale over $T_\lambda$ (and similarly for $S/T$). We repeatedly use \cite[Lem 04GW]{stacks-project} and the corollary that $T_x^{sh}=\varinjlim_\lambda T_{\lambda,x_\lambda}^{sh}$. Here it is not necessary that $G$ be upper bounded or with finite cohomology, but it is important that $G$ be the pullback of a lower bounded sheaf on $\Spec K_\et$ with torsion cohomology in order to use \cite[Cor. 1.18.(2)]{BhattMathew}.
    
    Combining the above with the previous $R\Psi_{k_1/k_0}^pLv^* = R\Psi_{k_1/k_0,\P}^pw^*$, we get :
    $$R\Psi_{K/k_0}(G) = R\Psi_{k_1/k_0}^pLv^*R\Psi_{K/k_1}^m(G)=R\Psi_{k_1/k_0,\P}^pR\Psi_{K/k_1,\P}^m(G).$$
    But from their descriptions on affine objects as pushforwards of $T/k_0\mapsto h^{k_1/k_0}(T_0)$ and $T/k_1\mapsto h^{K/k_1}(T)$, it follows that the composite of $R\Psi_{k_1/k_0,\P}^p$ and $R\Psi_{K/k_1,\P}^m$ coincides with the pushforward of the premorphism $T/k_0\mapsto h^{K/k_1}(h^{k_1/k_0}(T))=h^{K/k_0}(T)$. This is the functor $R\tilde\Psi_{K/k_0}$, finishing the proof.
    
    \textbf{2.} By \Lemma{monoidality premorphism}, $Lv^*$ and $Lv_\P^*$ are monoidal. As seen in \textbf{(1.)}, we have isomorphisms for $G\in D_\fin(K_\et)$ :
    \begin{align*}
        R\Psi_{K/k_0}(G)
        &= R\Psi_{k_1/k_0}^pLv^*R\Psi_{K/k_1}^mu^*(G),\\
        R\tilde\Psi_{K/k_0}
        &=R\Psi_{k_1/k_0,\P}^pR\Psi_{K/k_1,\P}^m,\\
        w^*R\Psi_{K/k_1}^m(G)
        &=R\Psi_{K/k_1,\P}^m(G),\\
        R\Psi_{k_1/k_0}^pLv^*
        &= R\Psi_{k_1/k_0,\P}^pw^*.
    \end{align*}
    
    In all cases the pairings are defined by using the monoidality of the (derived) pullbacks $u^*$, $Lv^*$, $i^*$, $i_\P^*$, and cup-products for pushforwards $R\Psi_{k_1/k_0}^p$, $Rj_*$, $R\Psi_{k_1/k_0,\P}$, $Rj_{\P,*}$ and $R\tilde\Psi_{K/k_0}$. In particular the natural multiplicative structures of the left and right sides of the first equality coincide by definition, and we need to show the other three isomorphisms are compatible with the natural mutiplicative structures.
    
    The formation of cup-products for pushforwards is compatible with the composition of premorphisms of sites, so multiplicative structures for $R\Psi_{k_1/k_0,\P}^pR\Psi_{K/k_1,\P}^m$ and for $R\tilde\Psi_{K/k_0}$ coincide.
    
    The equality $R\Psi_{k_1/k_0}^pLv^*= R\Psi_{k_1/k_0,\P}^pw^*$ comes from identifications $R\Psi_{k_1/k_0}^p=R\Psi_{k_1/k_0,\P}^pv_{\P,*}$ and $v_{\P,*} Lv^*=v_\P Lv_\P^*w^*=w^*$. The first is compatible with cup-products as before, by composition of pushforwards of premorphism of sites. The second comes from the unit $\id\to v_{\P,*}Lv_\P^*$, which is multiplicative in the sense that we have a commutative diagram for $F,F'\in\Sh(k_{1,\PRPS})$ :
    $$\begin{tikzcd}
       F\arrow[d] &[-3em] \otimes^L&[-3em] F'\arrow[d]\arrow[rr,equals]&& F \otimes^L F'\arrow[d]\\
       v_{\P,*}Lv_\P^*F&[-3em] \otimes^L&[-3em] v_{\P,*}Lv_\P^*F'\arrow[r]& v_{\P,*}(Lv_\P^*F \otimes^L Lv_\P^*F')\arrow[r,equals]&v_{\P,*}Lv_\P^*(F \otimes^L F')
    \end{tikzcd}$$
    where the lower pairing is given by the cup-product of $v_{\P,*}$ and monoidality of $Lv_\P^*$ ; this can be seen using \Lemma{monoidality premorphism}.2. Hence the multiplicative structures for $R\Psi_{k_1/k_0}^pLv^*$ and $R\Psi_{k_1/k_0,\P}^pw^*$ coincide.
    
    Let $G\in D_\fin(K_\Et)$. For $q\in\ZZ$, when restricted to affines the sheaf $w^*R^q\Psi_{K/k_1}^m(G)$ coincides with the sheafification of $T=\varprojlim_iT_i\mapsto \varinjlim_iH^q(h^{K/k_1}(T_i),G)$. On the other hand, $R^q\Psi_{K/k_1,\P}^m(G)$ when restricted to affines coincides with the sheafification of $T\mapsto H^q(h^{K/k_1}(T),G)$. By continuity of étale cohomology and \cite[Cor. 1.18.(2)]{BhattMathew}, we have that the restriction map :
    $$\varinjlim_iH^q(h^{K/k_1}(T_i),G)\to H^q(h^{K/k_1}(T),G)$$
    is an isomorphism, which describes the equality $w^*R\Psi_{K/k_1}^m(G)=R\Psi_{K/k_1,\P}^m(G)$ locally. On the other hand, the cup-products for the restrictions of $R\Psi_{K/k_1}^m$ and $R\Psi_{K/k_1,\P}^m$ to affines are both obtained by sheafifying the cup-products for these étale cohomology groups. In other words, for $F,F'\in\Sh(K_\Et)$ and $q,r\in\ZZ$ the following diagram :
    $$\begin{tikzcd}
        a(H^q(h^{K/k_1}(-),F) &[-3em] \otimes\arrow[d]&[-3em] H^r(h^{K/k_1}(-),F'))\arrow[r]& a(H^{q+r}(h^{K/k_1}(-),F \otimes F')\arrow[d]\\
        R^q\Psi_{K/k_1}^m(F) &[-3em] \otimes &[-3em] R^r\Psi_{K/k_1}(F')\arrow[r]& R^q\Psi_{K/k_1}^m(F \otimes F')
    \end{tikzcd}$$
    commutes when viewed as a diagram of sheaves over affine objects of $\Spec k_{1,\RPS}$ ; the analogous diagram for $R^q\Psi_{K/k_1,\P}^m$ also commutes over affine objects of $\Spec k_{1,\PRPS}$. Hence the compatibility of the isomorphism $w^*R\Psi_{K/k_1}^m(G)\xrightarrow\sim R\Psi_{K/k_1,\P}^m(G)$ with multiplicative structures comes from the commutativity of the diagram of cup-products :
    $$\begin{tikzcd}
        \varinjlim_iH^q(h^{K/k_1}(T_i),F) \arrow[d]&[-3em] \otimes&[-3em] \varinjlim_iH^r(h^{K/k_1}(T_i),F'))\arrow[r]\arrow[d]& \varinjlim_i H^{q+r}(h^{K/k_1}(T_i),F \otimes^L F')\arrow[d]\\
        H^q(h^{K/k_1}(T),F) &[-3em] \otimes&[-3em] H^r(h^{K/k_1}(T),F'))\arrow[r]& H^{q+r}(h^{K/k_1}(T),F \otimes F')
    \end{tikzcd}$$
    for $T=\varprojlim_iT_i$ an affine object of $\Spec k_{1,\PRPS}$, which follows from the compatibility of the tensor product with colimits and of the restrictions $H^q(h^{K/k_1}(T_i),-)\to H^q(h^{K/k_1}(T),-)$ with cup-products. This concludes.
    
    \textbf{3.} The exactness of $\beta_{k_0,*}$, which is simply the restriction to $\Spec k_{0,\et}$, is clear.  By \Lemma{pushforward of sites}, $R^q\Psi^\et_{U/k_0}(G)$ is the sheafification of $T/k_0\mapsto H^q(U_{h^{K/k_0}(T)},G)$ and $R^q\Psi^\et_{U/k_0,c}(G)$ is that of :
    $$T/k_0\mapsto H^q(X_{h^{K/k_0}(T)},(j_!G)|_{X_{h^{K/k_0}(T)}})=H^q_c(U_{h^{K/k_0}(T)},G|_{U_{h^{K/k_0}(T)}})$$
    where the latter equality is by base change \cite[Th. 0EYU]{stacks-project}. The identification of stalks follows and it remains to compare $R\Psi_{U/k_0}$ with $R\Psi_{U/k_0}^\et$ (and $R\Psi_{U/k_0,c}$ with $R\Psi_{U/k_0,c}^\et$). By \Lemma{finite coeffs equivalence} and \Theorem{equal char 0}.4, we have identifications in $D(k_{2,\et})$ :
    \begin{align*}
        \left[R\Psi_{k_{3/k_2}}^0\circ\cdots\circ R\Psi_{k_d/k_{d-1}}^0\circ R\Psi_{U/K}^v(G)\right]_{\Adh k_2}&
        = R\Gamma(k_3^\ur,-)\circ\cdots\circ R\Gamma(k_d^\ur,-)\circ R\Gamma(U_{\Adh K},-)(G),\\
        \left[R\Psi_{k_{3/k_2}}^0\circ\cdots\circ R\Psi_{k_d/k_{d-1}}^0\circ R\Psi_{U/K,c}^v(G)\right]_{\Adh k_2}&
        = R\Gamma(k_3^\ur,-)\circ\cdots\circ R\Gamma(k_d^\ur,-)\circ R\Gamma_c(U_{\Adh K},-)(G),
    \end{align*}
    while by \cite[Prop. 2.20 of Ch. III]{Milne80} we also have :
    \begin{align*}
        R\Gamma(U_{K^\tur},G|_{U^\tur})
        &= R\Gamma(k_2^\tur,-)\circ R\Gamma(k_3^\ur,-)\circ\cdots\circ R\Gamma(k_d^\ur,-)\circ R\Gamma(U_{\Adh K},-)(G),\\
        R\Gamma_c(U_{K^\tur},G|_{U^\tur})
        &= R\Gamma(k_2^\tur,-)\circ R\Gamma(k_3^\ur,-)\circ\cdots\circ R\Gamma(k_d^\ur,-)\circ R\Gamma_c(U_{\Adh K},-)(G).
    \end{align*}
    which reduces us to proving $\beta_{k_0,*}\circ R\Psi_{k_2/k_0}(G) = R\Psi^\et_{k_2/k_0}(G)$. This follows directly from \textbf{(1.)}.
    
    \textbf{4.} This is \textbf{(2.)} plus the fact that, in \textbf{(3.)}, all new identifications are compatible with cup-products.
\end{myproof}

\begin{proposition}[k0 finite interpretation]
    Consider \Situation{higher local setup}. Assume $k_0$ is finite. Recall $R\Psi_{U/*}$ and $R\Psi_{U/*,c}$ from \Theorem{higher local duality finite}.
    \begin{enumerate}
        \item For $G\in D_\fin(U_\et)$ and $q\in\ZZ$ we have canonical isomorphisms :
        $$R^q\Psi_{U/*}(G)(*)=H^q(U_\et,G),\qquad R^q\Psi_{U/*,c}(G)(*)=H^q_c(U_\et,G)=H^q(X_\et,j_!G).$$
        \item For $G\in D_\fin(U_\et)$ and $G^\vee=R\HOM_{U_\et}(G,\QQ/\ZZ(d+e))$, there is a commutative diagram of pairings :
        $$\begin{tikzcd}
            R^q\Psi_{U/*}(G)(*)\arrow[d,"\sim"] &[-3em] \otimes &[-3em] R^{d+2e+1-q}\Psi_{U/*,c}(G^\vee)(*)\arrow[d,"\sim"]\arrow[r]& R^{d+2e+1}\Psi_{U/*,c}(\QQ/\ZZ(d+e))(*)\arrow[d,"\sim"]\\
            H^q(U_\et,G) &[-3em] \otimes &[-3em] H^{d+2e+1-q}_c(U_\et,G^\vee)\arrow[r]& H^{d+2e+1}_c(U_\et,\QQ/\ZZ(d+e))
        \end{tikzcd}$$
    \end{enumerate}
\end{proposition}

\begin{myproof}[of \Proposition{k0 finite interpretation}]
    Consider the diagram :
    $$\begin{tikzcd}
        D(U_\et)\arrow[r,"R\Psi_{U/k_0}"]\arrow[d,equals]& D(k_{0,\et}^\perar)\arrow[r,"\alpha_{k_0}"]\arrow[d,"(-)_{\Adh k_0}"]& D(k_{0,\proet}^\indrat)\arrow[r,"R\Psi_{k_0/*}"]\arrow[d,"{(-)_{\Adh k_0}}"]& D(*_\proet)\arrow[d,"{R\Gamma(*,-)}"]\\
        D(U_\et)\arrow[r,"{R\Gamma(U_{K^\tur},-)}"']& D(k_{0,\et})\arrow[r,equals]& D(k_{0,\et})\arrow[r,"{R\Gamma(k_0,-)}"'] & D(\Ab)
    \end{tikzcd}$$
    where $\alpha_{k_0}:D(k_{0,\et}^\perar)\to D(k_{0,\proet}^\indrat)$ is the change of site functor. The leftmost square commutes by \Proposition{char 0pp-0p interpretation}.3. The second square commutes when restricted to $D^b_{\W_{k_0}}(k_{0,\et}^\perar)$ or $D_\fin(k_{0,\et}^\perar)$, by \Proposition{Yoneda for Wk}.4 or \Lemma{finite coeffs equivalence}.2. The third square clearly commutes by definition of $R\Psi_{k_0/*}$.
    
    The composite of the top arrows is exactly $R\Psi_{U/*}$. The composite of the bottom arrows is precisely $R\Gamma(U,-)$, by \cite[Prop. 2.20 of Ch. III]{Milne80}. By \Theorem{higher local duality}.1, \Proposition{Yoneda for Wk}.3 and \Lemma{finite coeffs equivalence}, we thus have :
    $$R\Gamma(*,R\Psi_{U/*}(G))=R\Gamma(U,G),\qquad G\in D_\fin(U_\et).$$
    By \Proposition{indproab over proet}.1, $\Gamma(*,-):\Sh(*_\proet)\to\Ab$ is exact hence $R^q\Psi_{U/*}(G)(*)=H^q(U,G)$ for $G\in D_\fin(U_\et)$. The identification $R^q\Psi_{U/*,c}(G)(*)=H^q_c(U,G)$ the same, with the addition of a base change argument \cite[Th. 0EYU]{stacks-project}, to show $R\Gamma(k_0,-)\circ R\Gamma_c(U_{K^\tur},-)=R\Gamma_c(U,-)$.
    
    This identification respects cup-products because the three following pairs of identifications do :
    \begin{align*}
        R\Psi_{U/k_0}(G)_{\Adh k_0} = R\Gamma(U_{K^\tur},G),\qquad & R\Psi_{U/k_0,c}(G^\vee)_{\Adh k_0} = R\Gamma_c(U_{K^\tur},G^\vee),\\
        \alpha_{k_0}(H)_{\Adh k_0} = H_{\Adh k_0},\qquad &\alpha_{k_0}(H^\vee)_{\Adh k_0} = (H^\vee)_{\Adh k_0},\\
        R\Gamma(*,R\Psi_{k_0/*}(W))=R\Gamma(k_0,W_{\Adh k_0}),\qquad &R\Gamma(*,R\Psi_{k_0/*}(W^\vee))=R\Gamma(k_0,(W^\vee)_{\Adh k_0}),\\
    \end{align*}
    for $G\in D_\fin(U_\et)$, $H\in D_\fin(k_{0,\et}^\perar)$ or $D^b_{\W_{k_0}}(k_{0,\et}^\perar)$, and $W\in D_\fin(k_{0,\proet}^\indrat)$ or $D^b_{\W_{k_0}}(k_{0,\proet}^\indrat)$, where $(-)^\vee=R\HOM(-,\QQ/\ZZ(r))$ with $r=d+e,0,0$ respectively for $G,H,W$. The compatibility for the first pair is by \Proposition{char 0pp-0p interpretation}.4, the second by \Theorem{Serre duality}.1, the third by composition of pushforwards.
\end{myproof}

The following corollary now gives a reformulation of \Theorem{higher local duality finite} with classical objects.

\begin{corollary}[higher local interpretation bis]
    Let $K$ be a $d$-local field of mixed characteristic at level $2$, with finite $0$-local residue field. Let $X$ be a smooth, proper, geometrically integral $K$-scheme of dimension $e$ and $U\subseteq X$ a nonempty open. Let $G\in D_\fin(U_\et)$ and $G^\vee=R\HOM_{U_\et}(G,\QQ/\ZZ(d+e))$. Then there exist canonical, locally compact Hausdorff group topologies on $H^q(U,G)$ and $H^q_c(U,G^\vee)$ for $q\in\ZZ$, functorial in $G$, such that the cup-product :
    $$H^q(U,G)\otimes H^{d+2e+1-q}_c(U,G^\vee)\to H^{d+2e+1}(U,\QQ/\ZZ(d+e))\to \QQ/\ZZ$$
    is a perfect Pontryagin duality. In particular, it is nondegenerate.
\end{corollary}

\begin{remark}[other characteristics]
    In \Situation{higher local field} we consider a $d$-local field $K=k_d,\dots,k_0$ with mixed characteristic $(0,p)$ at level $2$. With the same methods we can in fact consider the cases where $K$ has mixed characteristic at level $1$ or totally equal characteristic $0$.
    \begin{enumerate}
        \item If $K$ has mixed characteristic at level $1$, define :
        \begin{align*}
            R\Psi_{U/k_0} &=w_*\circ Lv^*\circ R\Psi_{k_1/k_0}^m\circ u^*\circ R\Psi_{k_2/k_1}^0\circ\cdots\circ R\Psi_{k_d/k_{d-1}}^0\circ R\Psi_{U/K}^v,\\
            R\Psi_{U/k_0,c} &=w_*\circ Lv^*\circ R\Psi_{k_1/k_0}^m\circ u^*\circ R\Psi_{k_2/k_1}^0\circ\cdots\circ R\Psi_{k_d/k_{d-1}}^0\circ R\Psi_{U/K,c}^v
        \end{align*}
        as functors $D(U_\et)\to D(k_{0,\et}^\perar)$, where $w:\Spec k_{0,\RP}\to\Spec k_{0,\et}^\perar$, $v:\Spec k_{0,\RP}\to\Spec k_{0,\RPS}$ and $u:\Spec k_{1,\Et}\to\Spec k_{1,\et}$ are induced by identity. Then \Theorem{higher local duality}, \Theorem{higher local duality finite}, \Proposition{char 0pp-0p interpretation}.3-4 and \Corollary{higher local interpretation bis} still hold verbatim. In fact in \Theorem{higher local duality}, $D_p(U_\et)$ maps to the smaller category $D^b(\Alg_u^\RP(k_0))\subseteq D^b_{\W_{k_0}}(k_{0,\et}^\perar)$, and in \Theorem{higher local duality finite}, $D_p(U_\et)$ maps to the smaller category $D_p(*_\proet)\subseteq D^b_{\W_p}(*_\proet)$. In particular, the pairing in \Corollary{higher local interpretation bis} is a perfect pairing of finite groups, recovering results of \cite{Izquierdo}. The only differences with the treated cases are that we need $w_*$ to satisfy a duality (this is \Proposition{Yoneda for Wk}.2 and \Theorem{Serre duality}.3.a), and that \Proposition{char 0pp-0p interpretation}.1-2 must be adapted for $d=1$.
        \item If $K$ has totally equal characteristic $0$, define :
        \begin{align*}
            R\Psi_{U/k_0} &=R\Psi_{k_1/k_0}^0\circ\cdots\circ R\Psi_{k_d/k_{d-1}}^0\circ R\Psi_{U/K}^v,\\
            R\Psi_{U/k_0,c} &=R\Psi_{k_1/k_0}^0\circ\cdots\circ R\Psi_{k_d/k_{d-1}}^0\circ R\Psi_{U/K}^v
        \end{align*}
        as functors $D(U_\et)\to D(k_{0,\et})$. Then \Theorem{higher local duality} adapts in the following form.
        \begin{itemize}
            \item If $G\in D_\fin(U_\et)$ then $R\Psi_{U/k_0}(G)$ and $R\Psi_{U/k_0,c}(G)$ belong to $D_\fin(k_{0,\et})$. If $G\in \Sh_\fin(U_\et)$, then $R\Psi_{U/k_0}(G)$ and $R\Psi_{U/k_0,c}(G)$ are concentrated in degrees $0,\dots,d+2e$.
            \item If $G\in D_\fin(U_\et)$ and $G^\vee=R\HOM_{U_\et}(G,\QQ/\ZZ(d+2))$ then we have a perfect Cartier duality :
            $$R\Psi_{U/k_0}(G)\otimes^L R\Psi_{U/k_0,c}(G^\vee)\to R\Psi_{U/k_0,c}(\QQ/\ZZ(d+e))\to\QQ/\ZZ[-d-2e].$$
            \item For $G\in D_\fin(U_\et)$, we have $R\Psi_{U/k_0}(G)=R\Gamma(U_{K^\tur},G)$ and $R\Psi_{U/k_0,c}(G)=R_c\Gamma(U_{K^\tur},G)$ in $\Sh(k_{0,\et})$.
        \end{itemize}
        If in addition $k_0$ is quasi-finite, meaning we have an isomorphism $\Gal(\Adh k_0/k_0)\cong\hat\ZZ$ (for instance, $k_0=\CC(\!(t)\!)$) then the functor $R\Gamma(k_0,-):D(k_{0,\et})\to D(\Ab)$ sends $D_\fin(k_{0,\et})$ to $D^b(\Ab_\fin)$, and Cartier dualities in $D_\fin(k_{0,\et})$ to perfect pairings of finite groups, via a trace map $H^1(k_0,\QQ/\ZZ)=\Hom_{\text{cont}}(\Gal(\Adh k_0/k_0),\QQ/\ZZ)\xrightarrow\sim\QQ/\ZZ$. Thus \Corollary{higher local interpretation bis} also holds in this case, where the pairing of étale cohomology groups is in fact a perfect pairing of (discrete) finite groups as in \cite{Izquierdo}.
    \end{enumerate}
    By contrast, the methods in this paper are insufficient to consider cases where $\characteristic k_3>0$, because we lack a positive equal characteristic step with imperfect residue field. Similarly, though we have a duality for $1$-local fields of equal characteristic $p>0$, we cannot consider varieties $U$ over them because we lack a $p$-torsion version of Poincaré duality.
\end{remark}


\section*{Appendix : Kato lifting}
\label{Appendix : Kato lifting}
\addcontentsline{toc}{section}{Appendix : Kato lifting}

In this appendix we review Kato's canonical lifting.

\begin{lemma}[Kato lift nilpotent uniqueness]
    Let $A$ be a ring and $I$ a nilpotent ideal. Write $k=A/I$ the quotient ring. Consider $A$-algebras $B$ and $B'$. Assume $B$ is formally étale, and assume $B'$ is either formally étale or flat, over $A$. Then any isomorphism $B/IB\xrightarrow\sim B'/IB'$ over $k$ lifts uniquely to an isomorphism $B\xrightarrow\sim B'$ over $A$.
\end{lemma}

\begin{myproof}[of \Lemma{Kato lift nilpotent uniqueness}]
    By Yoneda's lemma, a formally étale $A$-algebra $B$ is determined by the $k$-algebra $B/IB$ since :
    $$\Hom_A(B,C)=\Hom_A(B,C/IC)=\Hom_k(B/IB,C/IC)$$
    for any $A$-algebra $C$. This solves the case $B'$ is formally étale.
    
    Now assume $B$ is formally étale and $B'$ is flat. The formal étaleness of $B$ ensures there exists a unique lift $B\to B'$ of the isomorphism $B/IB\xrightarrow\sim B'/IB'$, as before. Once such a lift is given and $B'$ is flat, the morphism $B\to B'$ is an isomorphism which can be seen as follows. Let $K$ and $C$ be its kernel and cokernel as a morphism of $A$-modules. Then $C/IC$ is the cokernel of $B/IB\to B'/IB'$, by right-exactness of $-\otimes_A A/I$, so $C=IC$. It follows that $C=I^nC$ for all $n\geq 0$, and by nilpotence of $I$, $C=0$. Thus $B\to B'$ is surjective. By flatness of $B'$, we now have $K/IK=\ker(B/IB\to B'/IB')=0$ and as for $C$, $K=IK=I^2K=\dots=0$. Hence $B\to B'$ is bijective.\footnote{This argument, omitted from the original construction \cite{Kato82}, can be found in Jarod Alper's \emph{Notes on Deformation Theory} (\url{https://sites.math.washington.edu/\textasciitilde jarod/courses/math581J-fall21/def-theory-notes.pdf} at the time of writing.)}
\end{myproof}

\begin{proposition}[Kato lift nilpotent]
    Let $A$ be a ring, $I$ a nilpotent ideal, $k=A/I$ and $R$ a $k$-algebra. Assume that either :
    \begin{enumerate}
        \item $k$ has prime characteristic $p>0$ and $R$ is relatively perfect over $k$ ;
        \item the quotient $A\to k$ admits a ring section $k\to A$ and $R$ is formally étale over $k$.
    \end{enumerate}
    Then there exists a unique formally étale $A$-algebra $B$ equipped with a $k$-isomorphism $B\otimes_Ak\cong R$. If $R$ is $k$-flat, then $B$ is the unique $A$-flat algebra equipped with a $k$-isomorphism $B\otimes_Ak\cong R$.
\end{proposition}

\begin{myproof}[of \Proposition{Kato lift nilpotent}]
    In any case, uniqueness comes from \Lemma{Kato lift nilpotent uniqueness}. For the existence in case \textbf{(1.)}, see \cite[Lem. 1]{Kato82} : the algebra $B$ can be constructed by endowing $A$ with a certain $W_N(k)$-algebra structure and taking $B=A\otimes_{W_N(k)}W_N(R)$, for $N\geq 1$ greater than the nilpotence index of $I$, where $W_N$ denotes $p$-typical Witt vectors. For the existence in case \textbf{(2.)}, simply set $B=R\otimes_kA$ given any section $k\to A$ : it is formally étale (resp. flat) over $A$ by base change, if $R$ is.
\end{myproof}

\begin{proposition}[Kato lift complete]
    Let $A$ be a Noetherian, complete ring with respect to an ideal $I$ (meaning $A=\varprojlim_nA/I^n$), $k=A/I$ and $R$ a $k$-algebra. Consider one of the following cases :
    \begin{enumerate}
        \item $k$ has prime characteristic $p>0$ and $R$ is flat relatively perfect over $k$ ;
        \item the quotient $A\to k$ admits a ring section $k\to A$ and $R$ is flat formally étale over $k$.
    \end{enumerate}
    Then there exists a unique flat, $I$-adically complete $A$-algebra $B$ equipped with a $k$-isomorphism $B\otimes_Ak\cong R$.
\end{proposition}

\begin{myproof}[of \Proposition{Kato lift complete}]
    For each $n\geq 1$, there exists a unique flat $A/I^n$-algebra $B_n$ such that $B_n\otimes_Ak=R$, by \Proposition{Kato lift nilpotent}. By uniqueness, we have $B_m/I^nB_m=B_n$ for $m\geq n$, hence a map $B_m\to B_n$. Then $B=\varprojlim_n B_n$ is $A$-flat by \cite[Lem. 0912]{stacks-project}. For $n\geq 1$, by the same reference, $B/I^nB=\varprojlim_{m\geq n}B_m/I^nB_m$, so $B_n=B/I^nB$ and $B$ is $I$-complete. The uniqueness of $B$ follows from \Proposition{Kato lift complete} and from the fact that $B$ is determined by the $B/I^nB$'s.
\end{myproof}

\begin{remark}[case 2 of Kato lift]
    The assumption of \Proposition{Kato lift complete}.2 on $A$ is satisfied, for instance, if $A$ is a complete local ring of equal characteristic (for instance $k$ of characteristic $0$) since, by the Cohen structure theorem \cite[proof of Th. 032A]{stacks-project}, $A$ is a quotient of a power series ring over $k$.
\end{remark}

Recall that if $A$ is a Henselian local ring with residue field $k$, and $R$ is an étale $k$-algebra, there is a unique finite étale $A$-algebra $B$ such that $B\otimes_Ak=R$, called the \emph{étale lift of $R$ over $A$} \cite[Lem. 04GK]{stacks-project}.

\begin{definition}[Kato lift definition]
    Consider $A$ a ring, $I$ an ideal of $A$, $k=A/I$ and $R$ a $k$-algebra.
    \begin{itemize}
        \item If $A$ is Henselian local with maximal ideal $I$ and $R$ is $k$-étale, $h^{A/k}(R)$ is the étale lift of $R$ over $A$.
        \item If $A,I,k,R$ are as in \Proposition{Kato lift nilpotent} or \Proposition{Kato lift complete}, $h^{A/k}(R)$ is the $A$-algebra $B$ produced by that proposition, called the \emph{Kato lift of $R$ over $A$}.
    \end{itemize}
    If in addition $A$ is an integral domain with fraction field $K$, we write $h^{K/k}(R)=K\otimes_Ah^{A/k}(R)$.
\end{definition}

The two definitions of $h^{A/k}(R)$ above coincide when $A$ is a complete local ring and $R$ an étale $k$-algebra, by uniqueness of the Kato lift.

\begin{remark}[Kato lift generalization]
    The Kato lift can be partially generalized to non-affine $k$-schemes.
    \begin{enumerate}
        \item Bertapelle and Suzuki generalized the finite level Kato lift of \Proposition{Kato lift nilpotent} to non-affine schemes : if $X$ is a scheme, $I$ a sheaf of locally nilpotent ideals containing $p$, and $Z\to X$ the corresponding closed immersion, then any $T\in\RPSch/Z$ lifts uniquely to a formally étale $X$-scheme $h^{X/Z}(T)$ such that $h^{X/Z}(T)\times_XZ=T$, and $h^{X/Z}$ preserves fiber products and Zariski coverings \cite[Prop. 3.2]{BertapelleSuzuki}.
        \item By contrast, the infinite level Kato lift of \Proposition{Kato lift complete} does not extend reasonably to (even quasicompact quasiseparated) schemes. This is because the functor $h^{A/k}$ from affine relatively perfect $k$-schemes to affine $A$-schemes does not commute with fiber products : for relatively perfect $k$-algebras $R,R',S$ we have $h^{A/k}(R\otimes_SR')=[h^{A/k}(R)\otimes_{h^{A/k}(S)}h^{A/k}(R')]^\wedge$ and this differs from $h^{A/k}(R)\otimes_{h^{A/k}(S)}h^{A/k}(R')$ in general, for example if $S=k$ and $R$ and $R'$ are both infinite rank over $k$. Instead the infinite level Kato lift generalizes only as a ind-scheme (and a formal scheme, more precisely). Consider an ind-scheme $X=\varinjlim_nX_n$ whose transitions are infinitesimal thickenings, and such that the sheaf of ideals $I$ defined by $X_0$ contains $p$ and is Zariski-locally nilpotent over each $X_n$. Then any relatively perfect $X_0$-scheme $T$ lifts to an ind-scheme over $X$ by the formula $h^{X/X_0}(T)=\varinjlim_nh^{X_n/X_0}(T)$. It is uniquely determined as an $X$-ind-scheme $Y$ that is formally étale, complete in the sense that each $Y\times_XX_n$ is a scheme and $Y=\varinjlim_nY\times_XX_n$, and satisfies $Y\times_XX_0=T$.
    \end{enumerate}
\end{remark}

\begin{remark}[Kato lift characterizations]
    Let $A$ be a complete regular local ring of equicharacteristic $p>0$, with residue field $k$ such that $[k:k^p]<+\infty$. Then for $R$ relatively perfect over $k$, $h^{A/k}(R)$ is characterized as the unique complete, relatively perfect $A$-algebra $B$ such that $B\otimes_Ak=R$. To see this, first by \cite[Lem. 0C0S]{stacks-project} there exist isomorphisms $A\cong k[\![x_1,\dots,x_n]\!]$ and $K\cong k(\!(x_1,\dots,x_r)\!)$, and correspondingly :
    $$h^{A/k}(R)\cong R[\![x_1,\dots,x_r]\!]=\varprojlim_nR[x_1,\dots,x_r]/(x_1,\dots,x_r)^n$$
    so given $(a_1,\dots,a_s)$ any finite $p$-basis of $k$, the family $(a_1,\dots,a_s,x_1,\dots,x_r)$ is a strong $p$-basis of $A$ and $h^{A/k}(R)$ simultaneously. It follows that $h^{A/k}(R)$ is relatively perfect over $A$ by \Lemma{p-basis comparison}.2. Conversely, relatively perfect $A$-algebras are automatically flat by \cite[Prop. 5.2]{Kato86}, so any lift of $R$ relatively perfect and complete over $A$ must coincide with $h^{A/k}(R)$.
\end{remark}

\begin{remark}[usefulness of case 2]
    In this paper we use only case 1 of \Proposition{Kato lift complete}, which applies to both rings of positive equal characteristic and of mixed characteristic. Case 2 gives an analogous result for rings of equal characteristic $0$, but in practice it is less useful because the formally étale condition is too restrictive in characteristic $0$. For instance, if $A$ has equal characteristic $p>0$, the Kato lift of a relatively perfect $k$-algebra is relatively perfect over $A$ in many cases (\textit{e.g.} if $A$ is regular and $[k:k^p]<+\infty$, see \Remark{Kato lift characterizations}). By contrast if $A$ has characteristic $0$, the algebra $h^{A/k}(R)$ in either case 1 or case 2 is generally not formally étale.
\end{remark}

\begin{lemma}[local RP is good]
    Let $k$ be a field such that $[k:k^p]<\infty$ and $R$ is the relative perfection of a smooth $k$-algebra, then $R$ is a finite product of integral domains, which are relatively perfectly smooth over $k$.
\end{lemma}

\begin{myproof}[of \Lemma{local RP is good}]
    By \cite[Prop. 8.13]{BertapelleSuzuki}, $\Spec R$ is the (scheme theoretic) coproduct of finitely many irreducible components. In particular those are open in $\Spec R$, hence again relatively perfectly smooth over $k$. Thus $R=R_1\times\cdots\times R_n$ with $R_i$ relatively perfectly smooth over $k$ and irreducible. By \Lemma{p-basis comparison}.3, $R_i$ is also reduced, hence it is an integral domain by \cite[Lem. 01ON]{stacks-project}.
\end{myproof}

\begin{proposition}[Kato lift structure]
    Consider $A,I,k,R$ as in either case of \Proposition{Kato lift complete}. Let $J\subseteq R$ be an ideal and $I_{h^{A/k}(R),J}=\ker(h^{A/k}(R)\to R\to R/J)$.
    \begin{enumerate}
        \item The pair $(h^{A/k}(R),Ih^{A/k}(R))$ is Henselian. If $(R,J)$ is Henselian then so is $(h^{A/k}(R),I_{h^{A/k}(R),J})$. 
        \item Assume $A$ is a discrete valuation ring with fraction field $K$ and uniformizer $\pi$, and $I=\pi A$. Then :
        $$h^{A/k}(R)=\{0\}\sqcup\bigsqcup_{n\geq 0}\pi^n \left[h^{A/k}(R)\setminus \pi h^{A/k}(R)\right],\qquad h^{K/k}(R)=\{0\}\sqcup\bigsqcup_{n\in\ZZ}\pi^n \left[h^{A/k}(R)\setminus \pi h^{A/k}(R)\right].$$
        If in addition $R$ is a finite product of integral domains then there is a split exact sequence :
        $$0\to h^{A/k}(R)^\times \to h^{K/k}(R)^\times \to \Gamma(R,\ZZ)\to 0$$
        and if $R$ is a field then $h^{A/k}(R)$ is a discrete valuation ring with fraction field $h^{K/k}(R)$.
    \end{enumerate}
\end{proposition}

\begin{myproof}[of \Proposition{Kato lift structure}]
    \textbf{1.} The first statement comes from \cite[Lem. 0ALJ]{stacks-project} and the second from \cite[$(1)\Rightarrow(2)$ of Lem. 0DYD]{stacks-project}.
    
    \textbf{2.} For $x\in h^{A/k}(R)$, either $x=0$ or there exists some $n\geq 0$ such that $x\notin \pi^{n+1}h^{A/k}(R)$, since $h^{A/k}(R)$ is $I$-separated. Taking such an $n$ minimal, we have :
    $$h^{A/k}(R)=\{0\}\cup\bigcup_{n\geq 0}\left[\pi^nh^{A/k}(R)\setminus \pi^{n+1}h^{A/k}(R)\right]$$
    Since $h^{A/k}(R)$ is $A$-flat and $A$ is integral, the map $x\mapsto\pi x$ is injective on $h^{A/k}(R)$. Hence $\pi^nh^{A/k}(R)\setminus \pi^{n+1}h^{A/k}(R)=\pi^n\left[h^{A/k}(R)\setminus \pi h^{A/k}(R)\right]$. If $\pi^na=\pi^mb$ for some integers $n\geq m\geq 0$ and $a,b\in h^{A/k}(R)\setminus \pi h^{A/k}(R)$, then we have $a=\pi^{m-n}b$ so $m=n$ (because $a\notin\pi h^{A/k}(R)$) and $a=b$. Hence a disjoint union :
    $$h^{A/k}(R)=\{0\}\sqcup\bigsqcup_{n\geq 0}\pi^n\left[h^{A/k}(R)\setminus \pi h^{A/k}(R)\right].$$
    Since the tensor product of $A$-modules preserves colimits and $K=\bigcup_{n\geq 0}\pi^{-n}A$, we have $h^{K/k}(R)=\bigcup_{n\geq 0}\pi^{-n}h^{A/k}(R)$ and the previous decomposition extends as :
    $$h^{K/k}(R)=\{0\}\sqcup\bigsqcup_{n\in\ZZ}\pi^n\left[h^{A/k}(R)\setminus \pi h^{A/k}(R)\right].$$
    
    If $R$ is a finite product of integral domains, we reduce to $R$ integral by the observation that the factors of $R$ are also flat relatively perfect (resp. flat formally étale) $k$-algebras, and the functor $h^{A/k}$ commutes with finite products of such algebras. If $R$ is integral consider $x\in h^{K/k}(R)^\times$. By the above there exist unique $n,m\in\ZZ$ and elements $a,b\in h^{A/k}(R)\setminus\pi h^{A/k}(R)$ such that $x=\pi^na$ and $x^{-1}=\pi^{m}b$. Then $1=xx^{-1}=\pi^{n+m}ab$ and $ab\in h^{A/k}(R)\setminus\pi h^{A/k}(R)$ because $h^{A/k}(R)/\pi h^{A/k}(R)=R$ is integral. By uniqueness of the decomposition, $n=-m$ and $ab=1$. This shows $x$ decomposes uniquely as $x=\pi^na$ with $a\in h^{A/k}(R)^\times$. The required split exact sequence follows.
    
    If $R$ is a field then $h^{A/k}(R)$ is local with maximal ideal $\pi h^{A/k}(R)$ by \textbf{(1.)}, hence $h^{K/k}(R)=h^{A/k}(R)[\pi^{-1}]$ is its fraction field and the above shows any nonzero $x\in h^{A/k}(R)$ (which becomes invertible in $h^{K/k}(R)$) is of the form $\pi^na$ with $a\in h^{A/k}(R)^\times$ and some $n\geq 0$. This means $h^{A/k}(R)$ is a discrete valuation ring.
\end{myproof}

We finish with a couple useful lemmas about the Kato lift functor.

\begin{lemma}[Kato lift colimit]
    Consider $A,I,k$ as in \Proposition{Kato lift complete} and $\{R_\lambda\}$ a filtered diagram of $k$-algebras all satisfying the assumptions of one case of \Proposition{Kato lift complete}. Let $R=\varinjlim_\lambda R_\lambda$. Then $R$ satisfies the same assumptions and $h^{A/k}(R)$ is the $I$-completion of the $A$-algebra $\varinjlim_\lambda h^{A/k}(R_\lambda)$.
\end{lemma}

\begin{myproof}[of \Lemma{Kato lift colimit}]
    A filtered colimit of flat relatively perfect (resp. flat formally étale) $k$-algebras is again relatively perfect (resp. formally étale), and is flat by \cite[Lem. 05UU]{stacks-project}. Let $B_\lambda=h^{A/k}(R_\lambda)$ and $B=\varinjlim_\lambda B_\lambda$. Then $B$ is $A$-flat by \cite[Lem. 05UU]{stacks-project} and $B\otimes_Ak=\varinjlim_\lambda B_\lambda\otimes_Ak=R$. By \cite[Lem. 0912 and 031C]{stacks-project}, the $I$-completion $\hat B$ is $A$-flat, $I$-complete and $\hat B\otimes_Ak=B\otimes_Ak=R$. We conclude by uniqueness of the Kato lift.
\end{myproof}

\begin{lemma}[Kato lift restriction]
    Let $K$ be a complete discrete valuation field with residue field $k$. Consider $L/K$ a finite extension with residue field $l$ and $R/k$ an algebra. Assume $R$ is relatively perfect over $k$ if $\characteristic k=p>0$, or étale in general. Then there are canonical isomorphisms of $\O_L$- and $L$-algebras :
    $$h^{\O_L/l}(R\otimes_kl)=h^{\O_K/k}(R)\otimes_{\O_K}\O_L,\qquad h^{L/l}(R\otimes_kl)=h^{K/k}(R)\otimes_KL.$$
\end{lemma}

\begin{myproof}[of \Lemma{Kato lift restriction}]
    The second equality will follow from the first by associativity of the tensor product. Also by associativity of the tensor product, we get $(h^{\O_K/k}(R)\otimes_{\O_K}\O_L)\otimes_{\O_L}l=h^{\O_L/l}(R\otimes_kl)\otimes_{\O_L}l$
    and $h^{\O_K/k}(R)\otimes_{\O_K}\O_L$ is $\O_L$-flat by base change of the $\O_K$-flat algebra $h^{\O_K/k}(R)$ ; it is $\O_l$-étale if $R$ is étale. To see it is $\mm_L$-adically complete, remark that $\mm_L^e\O_L\subseteq\mm_K\O_L\subseteq\mm_L\O_L$ where $e$ is the ramification index of $L/K$, so by \cite[Lem. 0319]{stacks-project}, the notions of $\mm_K\O_L$- and $\mm_L\O_L$-completeness for $\O_L$-modules coincide. Then we have :
    \begin{align*}
        \varprojlim_n\left[(h^{\O_K/k}(R)\otimes_{\O_K}\O_L)/\mm_K^n\O_L\right]
        &= \varprojlim_n\left[(h^{\O_K/k}(R)/\mm_K^n\O_K)\otimes_{\O_K}\O_L\right] \\
        &=\left[\varprojlim_n(h^{\O_K/k}(R)/\mm_K^n\O_K)\right]\otimes_{\O_K}\O_L\\
        &=h^{K/k}(R)\otimes_{\O_K}\O_L
    \end{align*}
    where the second equality holds because $\O_L$ is a finite and flat, thus finite free by \cite[Lem. 02KB]{stacks-project}, $\O_K$-module. This proves $h^{\O_K/k}(R)\otimes_{\O_K}\O_L$ is an $\mm_L$-complete, flat (and étale if $R/k$ is étale) $\O_L$-algebra whose tensor to $l$ is $R\otimes_kl$, thus by \Proposition{Kato lift complete} we have $h^{\O_K/k}(R)\otimes_{\O_K}\O_L=h^{\O_L/l}(R\otimes_kl)$.
\end{myproof}


\renewcommand\refname{References}

\end{document}